\documentclass[reqno]{amsart} 

\author[A.~Appel]{Andrea Appel}
\address{Dipartimento di Scienze Matematiche, Fisiche e Informatiche, Universit\`a di Parma, INDAM GNSAGA, and INFN Gruppo Collegato di Parma, 43100 Parma PR, Italy}
\email{\href{mailto:andrea.appel@unipr.it}{andrea.appel@unipr.it}}
\author[J.~Li]{Jian-Rong Li}
\address{Faculty of Mathematics, University of Seville, Calle Tarfia s/n, 41012 Seville, Spain}
\email{\href{mailto:lijr07@gmail.com}{lijr07@gmail.com}}
\author[T.~Prze\'{z}dziecki]{Tomasz Prze\'{z}dziecki}
\address{Faculty of Mathematics, University of Vienna, Oskar Morgenstern Platz 1, 1090 Vienna, Austria}
\email{\href{mailto:tomasz.przezdziecki@univie.ac.at}{tomasz.przezdziecki@univie.ac.at}}
\author[B.~Vlaar]{Bart Vlaar}
\address{Beijing Institute of Mathematical Sciences and Applications, No. 544, Hefangkou Village, Huaibei Town, Huairou, Beijing, China}
\email{\href{mailto:b.vlaar@bimsa.cn}{b.vlaar@bimsa.cn}}

\title[Quantum loop symmetric pairs and abelian K-matrices]{Quantum loop symmetric pairs\\ and meromorphic abelian K-matrices}

\subjclass[2020]{
	Primary: 81R50. 
	Secondary: 16T25, 
	17B37, 
	81R12. 
}

\usepackage[english]{babel}
\usepackage[T1]{fontenc}
\usepackage{lmodern}
\usepackage[utf8x]{inputenc}
\usepackage{mathtools}
\mathtoolsset{showonlyrefs,showmanualtags}
\usepackage{empheq}
\usepackage{amsmath, amssymb}
\usepackage{stmaryrd}
\usepackage{enumerate,enumitem}
\usepackage[perpage,multiple,bottom]{footmisc}
\usepackage[usenames,dvipsnames]{color}
\usepackage{hyperref}
\hypersetup{
colorlinks=true,
citecolor=blue,
linkcolor=blue,
filecolor=blue, 
urlcolor=blue,
}
\usepackage[toc,page]{appendix}
\usepackage{xcolor}
\usepackage[breakable]{tcolorbox}
\usepackage{tikz}
\usetikzlibrary{shapes, backgrounds} 
\usepackage{tikz-cd}
\definecolor{orange}{rgb}{1,0.5,0}
\usepackage[all,cmtip,knot]{xy}

\usepackage{bm}

\newcommand{\nc}{\newcommand}
\newcommand{\rnc}{\renewcommand}

\nc{\tg}[1]{{\color{green!40!black}#1}}
\nc{\tb}[1]{{\color{blue!50!black}#1}}

\nc{\triend}{\parbox{2mm}{\hfill} \hfill\text{\hspace{0.2mm}}\hfill$\triangle$}
\nc{\ocend}{\parbox{2mm}{\hfill} \hfill\text{\hspace{0.2mm}}\hfill$\oslash$}
\nc{\black}{\parbox{2mm}{\hfill}\text{\hspace{0.2mm}}\hfill$\blacksquare$}
\nc{\rmkend}{$\triangle$}

\numberwithin{equation}{section}

\allowdisplaybreaks

\newtheorem*{theorem}{Theorem}
\newtheorem*{proposition}{Proposition}

\newtheorem*{lemma}{Lemma}
\newtheorem*{corollary}{Corollary}

\newenvironment{pf}{\medskip\paragraph{{\sc Proof}}}{\hfill$\square$\par\medskip}
\newenvironment{remark}{\medskip\paragraph{{\bf Remark.}}}{\par\medskip}
\newenvironment{remarks}{\medskip\paragraph{{\bf Remarks.} }}{\par\medskip}

\nc{\cA}{{\mathcal A}}
\nc{\cB}{{\mathcal B}}
\nc{\cC}{{\mathcal C}}
\nc{\cD}{{\mathcal D}}
\nc{\cE}{{\mathcal E}}
\nc{\cF}{{\mathcal F}}
\nc{\cG}{{\mathcal G}}
\nc{\cH}{{\mathcal H}}
\nc{\cI}{{\mathcal I}}
\nc{\cJ}{{\mathcal J}}
\nc{\cK}{{\mathcal K}}
\nc{\cL}{{\mathcal L}}
\nc{\cM}{\mathcal{M}}
\nc{\cN}{{\mathcal N}}
\nc{\cO}{{\mathcal O}}
\nc{\cP}{{\mathcal P}}
\nc{\cQ}{{\mathcal Q}}
\nc{\cR}{{\mathcal R}}
\nc{\cS}{\mathcal{S}}
\nc{\cT}{{\mathcal T}}
\nc{\cU}{\mathcal{U}}
\nc{\cV}{{\mathcal V}}
\nc{\cX}{{\mathcal X}}
\nc{\cY}{\mathcal{Y}}
\nc{\cW}{\mathcal{W}}
\nc{\cZ}{{\mathcal Z}}

\nc{\bbA}{{\mathbb{A}}}
\nc{\bbB}{{\mathbb{B}}}
\nc{\bbC}{{\mathbb{C}}}
\nc{\bbD}{{\mathbb{D}}}
\nc{\bbE}{{\mathbb{E}}}
\nc{\bbF}{{\mathbb{F}}}
\nc{\bbG}{{\mathbb{G}}}
\nc{\bbH}{{\mathbb{H}}}
\nc{\bbI}{{\mathbb{I}}}
\nc{\bbJ}{{\mathbb{J}}}
\nc{\bbK}{{\mathbb{K}}}
\nc{\bbL}{{\mathbb{L}}}
\nc{\bbM}{{\mathbb{M}}}
\nc{\bbN}{{\mathbb{N}}}
\nc{\bbO}{{\mathbb{O}}}
\nc{\bbP}{{\mathbb{P}}}
\nc{\bbQ}{{\mathbb{Q}}}
\nc{\bbR}{{\mathbb{R}}}
\nc{\bbS}{{\mathbb{S}}}
\nc{\bbT}{{\mathbb{T}}}
\nc{\bbU}{{\mathbb{U}}}
\nc{\bbV}{{\mathbb{V}}}
\nc{\bbX}{{\mathbb{X}}}
\nc{\bbY}{{\mathbb{Y}}}
\nc{\bbW}{{\mathbb{W}}}
\nc{\bbZ}{{\mathbb{Z}}}

\usepackage{mathrsfs}

\nc{\scrA}{{\mathscr A}}
\nc{\scrB}{{\mathscr B}}
\nc{\scrC}{{\mathscr C}}
\nc{\scrD}{{\mathscr D}}
\nc{\scrE}{{\mathscr E}}
\nc{\scrF}{{\mathscr F}}
\nc{\scrG}{{\mathscr G}}
\nc{\scrH}{{\mathscr H}}
\nc{\scrI}{{\mathscr I}}
\nc{\scrJ}{{\mathscr J}}
\nc{\scrK}{{\mathscr K}}
\nc{\scrL}{{\mathscr L}}
\nc{\scrM}{{\mathscr M}}
\nc{\scrN}{{\mathscr N}}
\nc{\scrO}{{\mathscr O}}
\nc{\scrP}{{\mathscr P}}
\nc{\scrQ}{{\mathscr Q}}
\nc{\scrR}{{\mathscr R}}
\nc{\scrS}{{\mathscr S}}
\nc{\scrT}{{\mathscr T}}
\nc{\scrU}{{\mathscr U}}
\nc{\scrV}{{\mathscr V}}
\nc{\scrX}{{\mathscr X}}
\nc{\scrY}{{\mathscr Y}}
\nc{\scrW}{{\mathscr W}}
\nc{\scrZ}{{\mathscr Z}}

\nc{\sfA}{{\mathsf A}}
\nc{\sfB}{{\mathsf B}}
\nc{\sfC}{{\mathsf C}}
\nc{\sfD}{{\mathsf D}}
\nc{\sfE}{{\mathsf E}}
\nc{\sfF}{{\mathsf F}}
\nc{\sfG}{{\mathsf G}}
\nc{\sfH}{{\mathsf H}}
\nc{\sfI}{{\mathsf I}}
\nc{\sfJ}{{\mathsf J}}
\nc{\sfK}{{\mathsf K}}
\nc{\sfL}{{\mathsf L}}
\nc{\sfM}{{\mathsf M}}
\nc{\sfN}{{\mathsf N}}
\nc{\sfO}{{\mathsf O}}
\nc{\sfP}{{\mathsf P}}
\nc{\sfQ}{{\mathsf Q}}
\nc{\sfR}{{\mathsf R}}
\nc{\sfS}{{\mathsf S}}
\nc{\sfU}{{\mathsf U}}
\nc{\sfV}{{\mathsf V}}
\nc{\sfX}{{\mathsf X}}
\nc{\sfY}{{\mathsf Y}}
\nc{\sfW}{{\mathsf W}}
\nc{\sfZ}{{\mathsf Z}}

\nc{\sfa}{{\mathsf a}}
\nc{\sfb}{{\mathsf b}}
\nc{\sfc}{{\mathsf c}}
\nc{\sfd}{{\mathsf d}}
\nc{\sfe}{{\mathsf e}}
\nc{\sff}{{\mathsf f}}
\nc{\sfg}{{\mathsf g}}
\nc{\sfh}{{\mathsf h}}
\nc{\sfi}{{\mathsf i}}
\nc{\sfj}{{\mathsf j}}
\nc{\sfk}{{\mathsf k}}
\nc{\sfl}{{\mathsf l}}
\nc{\sfm}{{\mathsf m}}
\nc{\sfn}{{\mathsf n}}
\nc{\sfo}{{\mathsf o}}
\nc{\sfp}{{\mathsf p}}
\nc{\sfq}{{\mathsf q}}
\nc{\sfr}{{\mathsf r}}
\nc{\sfs}{{\mathsf s}}
\nc{\sft}{{\mathsf t}}
\nc{\sfu}{{\mathsf u}}
\nc{\sfv}{{\mathsf v}}
\nc{\sfx}{{\mathsf x}}
\nc{\sfy}{{\mathsf y}}
\nc{\sfw}{{\mathsf w}}
\nc{\sfz}{{\mathsf z}}

\nc {\bfA}{{\mathbf A}}
\nc {\bfB}{{\mathbf B}}
\nc {\bfC}{{\mathbf C}}
\nc {\bfD}{{\mathbf D}}
\nc {\bfE}{{\mathbf E}}
\nc {\bfF}{{\mathbf F}}
\nc {\bfG}{{\mathbf G}}
\nc {\bfH}{{\mathbf H}}
\nc {\bfI}{{\mathbf I}}
\nc {\bfJ}{{\mathbf J}}
\nc {\bfK}{{\mathbf K}}
\nc {\bfL}{{\mathbf L}}
\nc {\bfM}{{\mathbf M}}
\nc {\bfN}{{\mathbf N}}
\nc{\bfO}{{\mathbf O}}
\nc {\bfP}{{\mathbf P}}
\nc {\bfQ}{{\mathbf Q}}
\nc {\bfR}{{\mathbf R}}
\nc {\bfS}{{\mathbf S}}
\nc {\bfT}{{\mathbf T}}
\nc {\bfU}{{\mathbf U}}
\nc {\bfV}{{\mathbf V}}
\nc {\bfX}{{\mathbf X}}
\nc {\bfY}{{\mathbf Y}}
\nc {\bfW}{{\mathbf W}}
\nc {\bfZ}{{\mathbf Z}}

\nc {\bfh}{{\mathbf h}}
\nc {\bfa}{{\mathbf a}}
\nc {\bfb}{{\mathbf b}}
\nc {\bfc}{{\mathbf c}}
\nc {\bfd}{{\mathbf d}}
\nc {\bfe}{{\mathbf e}}
\nc {\bff}{{\mathbf f}}
\nc {\bfg}{{\mathbf g}}
\nc {\bfz}{{\mathbf z}}

\nc {\fka}{{\mathfrak a}}
\nc {\fkb}{{\mathfrak b}}
\nc {\fkc}{{\mathfrak c}}
\nc {\fkd}{{\mathfrak d}}
\nc {\fke}{{\mathfrak e}}
\nc {\fkf}{{\mathfrak f}}
\nc {\fkg}{{\mathfrak g}}
\nc {\fkh}{{\mathfrak h}}
\nc {\fki}{{\mathfrak i}}
\nc {\fkj}{{\mathfrak j}}
\nc {\fkk}{{\mathfrak k}}
\nc {\fkl}{{\mathfrak l}}
\nc {\fkm}{{\mathfrak m}}
\nc {\fkn}{{\mathfrak n}}
\nc {\fko}{{\mathfrak o}}
\nc {\fkp}{{\mathfrak p}}
\nc {\fkq}{{\mathfrak q}}
\nc {\fkr}{{\mathfrak r}}
\nc {\fks}{{\mathfrak s}}
\nc {\fkt}{{\mathfrak t}}
\nc {\fku}{{\mathfrak u}}
\nc {\fkv}{{\mathfrak v}}
\nc {\fkx}{{\mathfrak x}}
\nc {\fky}{{\mathfrak y}}
\nc {\fkw}{{\mathfrak w}}
\nc {\fkz}{{\mathfrak z}}

\nc {\fkA}{{\mathfrak A}}
\nc {\fkB}{{\mathfrak B}}
\nc {\fkC}{{\mathfrak C}}
\nc {\fkD}{{\mathfrak D}}
\nc {\fkE}{{\mathfrak E}}
\nc {\fkF}{{\mathfrak F}}
\nc {\fkG}{{\mathfrak G}}
\nc {\fkH}{{\mathfrak H}}
\nc {\fkI}{{\mathfrak I}}
\nc {\fkJ}{{\mathfrak J}}
\nc {\fkK}{{\mathfrak K}}
\nc {\fkL}{{\mathfrak L}}
\nc {\fkM}{{\mathfrak M}}
\nc {\fkN}{{\mathfrak N}}
\nc {\fkO}{{\mathfrak O}}
\nc {\fkP}{{\mathfrak P}}
\nc {\fkQ}{{\mathfrak Q}}
\nc {\fkR}{{\mathfrak R}}
\nc {\fkS}{{\mathfrak S}}
\nc {\fkT}{{\mathfrak T}}
\nc {\fkU}{{\mathfrak U}}
\nc {\fkV}{{\mathfrak V}}
\nc {\fkX}{{\mathfrak X}}
\nc {\fkY}{{\mathfrak Y}}
\nc {\fkW}{{\mathfrak W}}
\nc {\fkZ}{{\mathfrak Z}}

\nc{\IN}{\bbN}
\nc{\IC}{\bbC}
\nc{\IR}{\bbR}
\nc{\IQ}{\bbQ}
\nc{\IZ}{\bbZ}
\nc{\IV}{\bbV}

\rnc{\a}{\fka}
\rnc{\b}{\fkb}
\nc{\fc}{\fkc}
\rnc{\d}{\fkd}
\nc{\e}{\fke}
\nc{\f}{\fkf}
\nc{\g}{\fkg}
\nc{\h}{\fkh}
\rnc{\i}{\fki}
\rnc{\j}{\fkj}
\rnc{\k}{\fkk}
\rnc{\l}{\fkl}
\nc{\m}{\fkm}
\nc{\n}{\fkn}
\rnc{\o}{\fko}
\nc{\p}{\fkp}
\nc{\q}{\fkq}
\rnc{\r}{\fkr}
\nc{\s}{\fks}
\rnc{\t}{\fkt}
\rnc{\u}{\fku}
\rnc{\v}{\fkv}
\nc{\x}{\fkx}
\nc{\y}{\fky}
\nc{\w}{\fkw}
\nc{\z}{\fkz}

\nc{\A}{\cA}
\nc{\B}{\cB}
\nc{\C}{\cC}
\nc{\D}{\cD}
\nc{\E}{\cE}
\nc{\F}{\cF}
\nc{\G}{\cG}
\rnc{\H}{\cH}
\nc{\I}{\cI}
\nc{\J}{\cJ}
\nc{\K}{\cK}
\rnc{\L}{\cL}
\nc{\M}{\cM}
\nc{\N}{\cN}
\rnc{\O}{\cO}
\rnc{\P}{\cP}
\nc{\Q}{\cQ}
\nc{\R}{\cR}
\rnc{\S}{\cS}
\nc{\T}{\cT}
\nc{\U}{\cU}
\nc{\V}{\cV}
\nc{\X}{\cX}
\nc{\Y}{\cY}
\nc{\W}{\cW}
\nc{\Z}{\cZ}

\nc{\ba}{{\bm\alpha}}
\nc{\bg}{{\bm\gamma}}
\nc{\bs}{{\bm\sigma}}
\nc{\bz}{{\bm\zeta}}
\nc{\bx}{{\bm\xi}}
\nc{\bpi}{{\bm\pi}}
\nc{\bphi}{{\bm\phi}}
\nc{\bpsi}{{\bm\psi}}
\nc{\be}{{\bm\eta}}
\nc{\beps}{{\bm\epsilon}}
\nc{\bom}{{\bm\omega}}
\nc{\brho}{{\bm\rho}}
\nc{\bvrho}{{\bm\varrho}}
\nc{\bxi}{{\bm\xi}}
\nc{\bk}{{\bm\kappa}}
\nc{\bTheta}{{\bm\Theta}}

\nc {\ul}{\underline}
\nc {\ol}{\overline}
\nc {\wt}{\widetilde}
\nc {\wh}{\widehat}
\nc{\scs}{\scriptscriptstyle}
\nc{\scsop}{\scriptscriptstyle\operatorname} 
\nc{\op}{{\operatorname{op}}}
\nc{\cop}{{\operatorname{cop}}}

\nc{\eps}{\epsilon}
\nc{\veps}{\varepsilon}

\nc{\ges}{\geqslant}
\nc{\les}{\leqslant}

\nc{\fext}[2]{{#1}[\negthinspace[#2]\negthinspace]}
\nc{\Lext}[2]{{#1}[#2]\negthinspace]}

\nc{\ten}{\otimes}

\nc {\ie}{{\emph{i.e.}}, }
\nc {\eg}{{\emph{e.g.}}, }
\nc {\qu}{\quad}
\nc {\qq}{\qquad}
\nc {\aand}{\qquad\mbox{and}\qquad}

\nc{\Hom}{\operatorname{Hom}}
\nc{\Ext}{\operatorname{Ext}}
\nc{\Aut}{\operatorname{Aut}}
\nc{\End}{\operatorname{End}}
\nc{\Tor}{\operatorname{Tor}}
\nc{\Pic}{\operatorname{Pic}}
\nc{\Ind}{\operatorname{Ind}}
\nc{\Fun}{\operatorname{Fun}}
\nc{\Res}{\operatorname{Res}}
\nc{\Obj}{\operatorname{Obj}}
\nc{\Rep}{\operatorname{Rep}}
\nc{\Ad}{\operatorname{Ad}}
\nc{\Ker}{\operatorname{Ker}}
\nc{\Mod}{\operatorname{Mod}}
\nc{\Sym}{\operatorname{Sym}}
\nc{\Tr}{\operatorname{Tr}}

\nc{\ev}{\operatorname{ev}}
\nc{\coev}{\operatorname{coev}}

\nc{\id}{{\operatorname{id}}}
\nc{\Spec}{\operatorname{Spec}}
\nc{\colim}{\operatorname{colim}}
\nc{\rank}{\operatorname{rk}}
\nc{\opspan}{\operatorname{span}}
\nc{\mult}{\operatorname{mult}} 
\nc{\norm}{\operatorname{norm}}

\nc {\Vect}{\operatorname{Vect}}

\rnc {\Re}{\operatorname{Re}}
\rnc {\Im}{\operatorname{Im}}

\nc{\ad}{\operatorname{ad}}

\nc{\Span}{\operatorname{span}}

\nc{\Phif}{\mathbf{\Phi}}
\nc{\Xf}{\mathbf{X}}
\nc{\Yf}{\mathbf{Y}}
\nc{\sW}{\scrW}
\nc{\sO}{\scrO}
\nc{\sM}{\scrM}
\nc{\sC}{\scrC}
\nc{\OUqg}{\O_{U_q(\g)}}
\nc{\pole}{\mathbf{p}}
\nc{\sK}{\scrK}
\nc{\DD}[1]{\underset{#1}{\Delta}}
\nc{\DT}[1]{\underset{#1}{\ten}}
\nc{\DDi}[1]{\underset{#1}{\Delta}^{\imath}}

\nc{\KMT}[1]{{\otimes}_{#1}^{\scsop{\mathrm KM}}}

\nc{\KMD}[1]{\underset{\scsop{\mathrm KM},{#1}}{\Delta}}

\nc{\FC}{{\bm\zeta}}

\nc{\GF}{\mathbf{\Gamma}}
\nc{\SF}{\mathbf{\Sigma}}
\nc{\Cc}{\mathbf{c}}

\nc{\La}{\Lambda}

\nc{\oD}{\operatorname{D}}

\nc{\pg}{\gamma}
\nc{\ps}{\sigma}

\nc{\SH}{\Sigma}

\nc{\RF}{\brho}
\nc{\ZF}{\bz}

\nc{\psa}[1]{{\mathbf s}_{#1,0}}
\nc{\psb}[1]{{\mathbf s}_{#1,-1}}

\nc{\pga}[1]{{\mathbf c}_{#1,0}}
\nc{\pgb}[1]{{\mathbf c}_{#1,-1}}

\nc{\fdel}{{\bm\delta}}

\nc{\Leta}{\varrho}
\nc{\deta}{\dot{\eta}}
\nc{\ddeta}{\ddot{\eta}}

\nc{\weta}{{\bm \eta}}
\nc{\wetad}{\dot{\bm \eta}}
\nc{\wetadd}{\ddot{\bm \eta}}
\nc{\etad}{\dot{\eta}}
\nc{\wy}{{\bf y}}
\nc{\ww}{{\bf w}}
\nc{\wx}{{\bf x}}
\nc{\wLetad}{\dot{\bvrho}}
\nc{\wLeta}{{\bvrho}}
\nc {\bfHd}{\dot{\mathbf H}}

\nc{\wgt}{\operatorname{wt}}

\nc{\bC}[2]{{\mathbf{C}^{#1}_{#2}}}
\nc{\bK}[2]{{{\bm\kappa}^{#1}_{#2}}}

\nc{\tk}{{\mathbf C}}
\nc{\tkt}{{\mathbf K}}
\nc{\tkb}{{}}

\nc{\RR}{{\mathscr R}}
\nc{\KK}{{\mathscr K}}
\nc{\bKK}{{\mathfrak K}}
\nc{\dKK}{\dot{\mathscr K}}
\nc{\LR}{\Lambda}
\nc{\qCB}{{\mathsf B}}
\nc{\qCC}{{\mathsf C}}
\nc{\qCBd}{\dot{\mathsf B}}
\nc{\qCG}{{\mathsf G}}
\nc{\qCD}{{\mathsf D}_q}

\nc{\Tb}{\T^{\theta}}

\nc{\oquad}{\operatorname{quad}}
\nc{\oprim}{\operatorname{prim}}

\nc{\odef}{\operatorname{def}}
\nc{\sgn}{\operatorname{sgn}}

\nc{\phalf}{\sfp^{\frac{1}{2}}}
\nc{\hc}{h_{\bfc}}

\nc{\ii}{\sqrt{-1}}

\begin{document}
\maketitle\thispagestyle{empty}

\begin{abstract}
Let $\g$ be a symmetrizable Kac--Moody algebra, let $q\in\bbC^\times$ not be a root 
of unity, and let $\bfU=U_q(L\g)$ be the corresponding quantum loop algebra, equipped 
with Hernandez' deformed Drinfeld coproduct.
For every diagram involution $\tau$ of $\g$, we introduce a coideal 
subalgebra $\bfB$ in a topological completion of $\bfU$, obtained by adapting the Letzter--Kolb construction of 
quantum symmetric pairs to the Drinfeld generators of $\bfU$. 
We derive necessary constraints on the parameters defining $\bfB$ from a
minimality condition on its intersection with the quantum loop Cartan subalgebra
$\bfU^0\subset\bfU$.
We show that $\bfB$ yields a topological realization of the quasi-split $\imath$quantum 
loop group $\bfU^\imath$ associated with $\tau$, in the presentation recently obtained 
by Lu--Pan, building on work of Chen--Lu--Pan--Ruan--Wang on $\imath$Hopf 
algebras.
By pullback, this gives rise to a rational restriction functor on category $\O$
representations. For modules with rational weights and generic values of the spectral
parameter, every $\bfU^\imath$-submodule of a restricted module is stable under the
derived quantum loop algebra $\bfU'$. In particular, irreducibility is generically
preserved under restriction for such modules.

For $|q|\neq1$ and under explicit assumptions on the quantum Cartan matrix, we 
construct meromorphic abelian K-matrices whose tensor counterparts yield boundary 
$q$-characters for both $\bfU$-modules and $\bfU^\imath$-modules.

Finally, we extend the coideal construction and Drinfeld coaction to 
shifted $\imath$quantum loop algebras with arbitrary coweight shifts, providing a 
new approach to $\imath$GKLO representations by restriction from ordinary shifted 
quantum loop algebras.

\end{abstract}

\setcounter{tocdepth}{1}
\tableofcontents


\section{Introduction}\label{s:introduction}

\subsection{}
Let $\g$ be a complex finite-dimensional simple Lie algebra, $q\in\bbC^\times$ not a 
root of unity, and $U_q(L\g)$ the corresponding quantum loop algebra.
In \cite{appel-vlaar-25a}, the first and fourth authors considered the spectral reflection 
equation on a tensor product of finite-dimensional $U_q(L\g)$-modules, 
with coefficients given by the action of the quantum R-matrix 
\cite{drinfeld-quantum-groups-87, frenkel-reshetikhin-92}, and constructed a family 
of universal solutions, depending on the choice of a quantum symmetric pair in 
$U_q(L\g)$.

Specifically, for each finite-dimensional $U_q(L\g)$-module $V$, there exists an 
operator $\bKK_V(z)$ such that
\begin{equation}\label{eq:intro-reflection}
	\begin{aligned}
		&\RR_{W^\psi,V^\psi}(w/z)_{21}(1\ten\bKK_W(w))
		\RR_{V^\psi,W}(zw)(\bKK_V(z)\ten1)\\
		&\qquad=(\bKK_V(z)\ten1)\RR_{W^\psi,V}(zw)_{21}
		(1\ten\bKK_W(w))\RR_{V,W}(w/z),
	\end{aligned}
\end{equation}
where $\RR_{V,W}(z)$ is the quantum R-matrix on $V\ten W$, $\psi$ is a prescribed twisting automorphism of $U_q(L\g)$, and $V^\psi$ denotes the pullback of $V$ through $\psi$.
The reflection equation was introduced by Cherednik and Sklyanin in the study of integrable systems with boundary conditions \cite{cherednik-84,sklyanin-88}, and expresses the compatibility of two-particle scattering with reflection at a boundary. 
The generalized form \eqref{eq:intro-reflection} appeared in \cite{cherednik-92}. Solutions of this equation are known as (spectral) K-matrices.

\subsection{}
In the representation-theoretic approach to the reflection equation a central role is played by quantum symmetric pairs.
These consist of a Drinfeld--Jimbo quantum group and a coideal subalgebra which quantizes the enveloping algebra of a symmetric Lie subalgebra \cite{letzter-02,kolb-14}.
Here, the choice of involution determines the K-matrix.
The construction of universal K-matrices for Drinfeld--Jimbo quantum groups of arbitrary Kac--Moody type was obtained in \cite{appel-vlaar-22}, extending the finite-type constructions of Bao--Wang \cite{bao-wang-18-kl} and Balagovi\'{c}--Kolb \cite{balagovic-kolb-19}. 

In affine type, this result led to a systematic construction of solutions of \eqref{eq:intro-reflection}.
For simplicity, we consider only the case where the underlying diagram involution $\tau$ of the quantum affine symmetric pair $U_q(\k)\subset U_q(L\g)$ fixes the affine node.
Proceeding as in \cite{drinfeld-quantum-groups-87}, a formal spectral K-matrix $\bKK_V(z)$
on	a finite-dimensional $U_q(L\g)$-module $V$ is readily obtained by applying the homogeneous shift $\SH_z$ to the universal K-matrix.
Specifically, $\bKK_V(z)$ satisfies the coproduct identity
	\begin{equation}\label{eq:intro-KM-coproduct}
	\bKK_{V{\KMT{w/z}}W}(z)
	=(1\ten\bKK_W(w))\RR_{V^\psi,W}(zw)(\bKK_V(z)\ten1),
\end{equation}
where ${\KMT{z}}$ denotes the tensor product
induced by the shifted Kac--Moody coproduct $\KMD{z}=(\id\ten\SH_z)\circ{\underset{\scsop{\mathrm{KM}}}{\Delta}}$
on $U_q(L\g)$,
and yields a $U_q(\k)$-intertwiner
\begin{equation}\label{eq:intro-KM-intertwiner}
	\bKK_V(z)\colon V(z)
	\longrightarrow V^\psi(z^{-1}),
\end{equation}
where $V(z)$ denotes the spectral shift of $V$, $\omega$ is the Chevalley involution, and $\psi=\omega\tau$.

\subsection{}
For irreducible finite-dimensional $V$, the formal K-matrix admits a normalized form
\begin{equation}\label{eq:intro-trigonometric-normalization}
	\bKK_V(z)=f_V(z)\,\bKK_V^{\mathrm{trig}}(z),
\end{equation}
where $f_V(z)\in\bbC(\negthinspace(z)\negthinspace)^\times$ and
$\bKK_V^{\mathrm{trig}}(z)\in\End(V)(z)$ is a trigonometric solution of
the reflection equation \eqref{eq:intro-reflection}.
The scalar function $f_V(z)$, and hence $\bKK_V(z)$, is meromorphic when $\g=\mathfrak{sl}_2$ and $|q|\neq1$.
For arbitrary finite-dimensional $V$, however, the meromorphicity of $\bKK_V(z)$ remains open even in this case.\\

In analogy with a similar factorization of the R-matrix \cite{khoroshkin-tolstoy-93, levendorskii-soibelman-stukopin-93, beck-kac-96, damiani-98}, it is natural to expect that, up to a gauge transformation, $\bKK(z)$ admits
a factorization
\begin{equation}\label{eq:intro-K-factorization}
	\bKK(z)=\bKK^+(z)\,\bKK^0(z)\,\bKK^-(z),
\end{equation}
where $\bKK^\pm(z)$ are rational operators supported on $U_q(L\n^\pm)$,
while $\bKK^0(z)$ is a meromorphic operator supported on 
$U_q(L\h)$.
In particular, such a decomposition should recover the boundary $q$-characters of the 
second and third authors \cite{Przez-23, li-przezdziecki-25, LP-eval, 
li-przezdziecki-26}, as a generalization of the Frenkel--Reshetikhin construction 
\cite{frenkel-reshetikhin-98}, with the modified tensor K-matrix of \cite{appel-vlaar-25c} 
in place of the R-matrix.\footnote{
More precisely, this amounts to taking the partial trace of the modified tensor K-matrix 
over a finite-dimensional $U_q(\k)$-module and then applying the loop Harish--Chandra 
projection.}
Evidence for the existence of the factorization \eqref{eq:intro-K-factorization}
appeared recently in the work of Baseilhac and Gainutdinov
\cite[Thm.~4.13]{baseilhac-gainutdinov-26}, for the semi-universal tensor K-matrix 
of the $q$-Onsager algebra evaluated on an irreducible two-dimensional representation of $U_q(L\mathfrak{sl}_2)$.

\subsection{}
The factorization of the R-matrix has the following interpretation in the category of 
finite-dimensional $U_q(L\g)$-modules.
The triangular factors $\RR^{\pm}(z)$ are twists relating the shifted Kac--Moody tensor 
product $\KMT{ s}$ to the Drinfeld tensor product
\cite{khoroshkin-tolstoy-94, enriquez-khoroshkin-pakuliak-07, damiani-24}.
These factors are rational and satisfy a unitarity condition. 
For $|q|\neq1$, the abelian factor $\RR^0(z)$ induces a meromorphic braiding on this category equipped with the Drinfeld tensor product \cite{gautam-toledano-17}.\\

Motivated by this interpretation, we expect conjugation by $\bKK^-(z)$ to carry 
$U_q(\k)$ to a {\em topological} right coideal subalgebra with respect to the Drinfeld 
coproduct, or equivalently to a subalgebra of natural operators 
on finite-dimensional $U_q(L\g)$-modules whose action on Drinfeld tensor products is 
right coideal. While $\bKK^{\pm}(z)$ would therefore intertwine such operators with the 
standard action of $U_q(\k)$, the operator $\bKK^0(z)$ would be an \emph{abelian 
K-matrix} with respect to the abelian R-matrix $\RR^0(z)$ and the Drinfeld tensor 
product. More precisely,
we expect a commutative diagram
\[
	\begin{tikzcd}[column sep=large,row sep=large]
		\underset{\scsop{KM}}{\Res}(V(z)) 
		\arrow[r,"\bKK^-_V(z)"]
		\arrow[d,"\bKK_V(z)"']
		&\underset{\scsop{D}}{\Res}(V(z))
		\arrow[d,"\bKK^0_V(z)"]\\
		\underset{\scsop{KM}}{\Res}(V^\psi(z^{-1}))
		&\underset{\scsop{D}}{\Res}(V^{\psi}(z^{-1}))
		\arrow[l,"\bKK^+_V(z)"]
	\end{tikzcd}
\]
where $\underset{\scsop{KM}}{\Res}$ is the standard restriction to $U_q(\k)$, $\underset{\scsop{D}}{\Res}$ denotes the expected topological restriction, and $\bKK^+_V(z)=\bKK^-_{V^\psi}(z^{-1})^{-1}$.\footnote{We expect a similar factorization for the tensor K-matrix, with its negative part acting as a comodule twist which twist-conjugates the standard coaction to a suitable {Drinfeld} coaction.}

\subsection{}
The goal of the present paper is to develop an independent approach to the abelian component. We construct a topological coideal subalgebra that yields meromorphic abelian $K$-matrices with respect to the Drinfeld tensor product and the abelian $R$-matrix.

Our construction uses only the Drinfeld presentation and therefore applies to the quantum affinisation of an arbitrary symmetrisable Kac--Moody algebra $\g$ \cite{hernandez-05,hernandez-07}. It thus extends the construction of abelian $K$-matrices beyond the setting in which we presently have a \emph{full} $K$-matrix as in \eqref{eq:intro-K-factorization}, encompassing, for instance, quantum toroidal algebras. More generally, we expect most of our results to extend almost verbatim to \emph{quiver} quantum loop algebras \cite{negut-arbitrary-quivers, negut-sala-schiffmann-25,  appel-gautam-26},\footnote{In that setting, however, we do not yet know the appropriate $\imath$Drinfeld--Serre relations, analogous to the zig-zag relations obtained by Negu\textcommabelow{t} in \cite{negut-symmetric-cartan}.}
see also \cite{ding-iohara-97}.
We now describe the construction and its consequences in detail.

\subsection{}
Henceforth, let $\g$ be an arbitrary symmetrizable Kac--Moody algebra,
$\bfU=U_q(L\g)$ the corresponding quantum affinization, $x_i^\pm(z)$ and $\phi_i^\pm(z)$ the series of its generators, and $\DD{s}$ the Drinfeld 
coproduct. Our first step is to adapt the Letzter--Kolb generators to the Drinfeld 
presentation of $\bfU$.
Let $\tau$ be a diagram involution of $\g$, and $\FC$ a formal central parameter.
We consider the fields
	\begin{align}
		\eta_i(z)&=\phi_i^+(z)\phi_{\tau(i)}^-(\FC^{-1}z^{-1}),\\ \label{eq:intro-coideal-generators}
		y_i(z)&=x_i^-(\FC^{-1}z^{-1})
		+\pg_i(z)x_{\tau(i)}^+(z)\phi_i^+(\FC^{-1}z^{-1})
		+\ps_i(z)\phi_i^+(\FC^{-1}z^{-1}),
	\end{align}
where $\pg_i(z),\ps_i(z)$ are scalar series subject initially only to conditions ensuring that the coefficients lie in $\Lext{\bfU}{\FC;\FC^{-1}}$.
Together with the elements $q^h$ for $h\in\h$ such that $\tau(h)=-h$, these coefficients generate a subalgebra $\bfB$ of $\Lext{\bfU}{\FC;\FC^{-1}}$.
The formulae for the Drinfeld coproduct immediately imply that $\bfB$ is a right coideal subalgebra (Proposition~\ref{ss:coideal}).

We further restrict the parameters by imposing necessary conditions for the intersection with the quantum loop Cartan subalgebra to be as small as possible. More precisely, let $\bfU^0$ denote the subalgebra of $\bfU$ generated by all $q^h$ and the coefficients of $\phi_i^\pm(z)$. We then determine necessary conditions on $\pg_i(z)$ and $\ps_i(z)$ for $\bfB$ to satisfy the minimality condition
\begin{equation}\label{eq:intro-minimality}
	\bfB\cap\Lext{\bfU^0}{\FC;\FC^{-1}}=\bfB^0,
\end{equation}
where $\bfB^0$ is generated by the elements $q^h$ satisfying $\tau(h)=-h$ and the coefficients of $\eta_i(z)$.
This is a loop version of the minimal-intersection requirement in the theory of 
quantum symmetric pairs \cite{letzter-02,kolb-14,regelskis-vlaar-20}.
Roughly, computing products of the generators produces additional terms in 
$\bfU^0$. Their vanishing forces $\pg_i(z),\ps_i(z)$ to arise from rational functions 
$\GF_i(z),\SF_i(z)$ depending on finitely many parameters.

The same computations show that, with these restrictions,
the coefficients of $y_i(z)$ and $\pg_i(\FC^{-1}z^{-1})\eta_i(z)$ satisfy the defining relations of the quasi-split $\imath$quantum loop algebra $\bfU^\imath$, yielding a morphism of algebras $\bfU^\imath\to\bfB$. Moreover, the right coideal property of $\bfB$ readily suggests the formulae for a $\bfU$-coaction on $\bfU^\imath$.
Thus, viewed \emph{ex post}, the construction of $\bfB$, guided by the minimal-intersection condition, offers an alternative perspective on the origin of the Drinfeld presentation of $\imath$-quantum groups  and their comodule structure
$\DDi{s}$ with respect to the Drinfeld coproduct \cite{lu-pan-26}.

\subsection{}
For the $q$-Onsager algebra realized as a coideal subalgebra of 
$U_q(L{\mathfrak{sl}}_2)$ with equal coideal parameters, Baseilhac and Kolb 
\cite{baseilhac-kolb-20} constructed real and imaginary root vectors, derived their 
commutation relations, and established a PBW basis, laying the groundwork for a Drinfeld-type presentation. Building on these results and 
motivated by the theory of $\imath$Hall algebras, Lu and Wang \cite{lu-wang-20} 
reformulated the rank-one relations in Drinfeld current form and established such a presentation in split affine $\mathsf{ADE}$ type. Their construction 
uses a universal $\imath$quantum group whose central specializations recover the 
coideal realizations with arbitrary nonzero parameters, including the two independent 
parameters of the $q$-Onsager embedding. Subsequent works \cite{zhang-22, 
	lu-wang-zhang-23, lu-wang-zhang-24, lu-pan-wang-zhang-25} extended this 
construction to the affinizations of arbitrary quasi-split Satake diagrams of finite type.

Lu and Pan \cite{lu-pan-26} subsequently extended this presentation to arbitrary Kac--Moody type and constructed a Drinfeld coaction on $\bfU^\imath$, see Theorem~\ref{thm:Lu-Pan-coaction}. Their use of the theory of $\imath$Hopf algebras developed by Chen--Lu--Pan--Ruan--Wang \cite{chen-pan-ruan-wang-25} yields a particularly conceptual and efficient proof of the coaction property.

\subsection{}
The coideal construction above recovers these results in a conceptually simpler form, 
with formulae obtained 
directly by restricting the Drinfeld coproduct to $\bfB$. We give an independent 
proof that these formulae define a coaction on $\bfU^\imath$ in 
Appendix~\ref{app:shifted-coaction}, verifying directly all the defining 
relations, including the $\imath$Drinfeld--Serre relations from \cite{lu-pan-26}. 
An advantage of this direct computational approach is that it extends readily 
to shifted quantum loop algebras, as we observe in \ref{ss:intro-shift}.

The relation between the coaction on $\bfU^\imath$ and the coideal subalgebra 
$\bfB$ becomes explicit through one-dimensional representations of 
$\bfU^\imath$. Indeed, the latter are parametrized by the same rational functions $\GF_i(z),\SF_i(z)$ that enter the construction of $\bfB$, see Proposition~\ref{prop:one-dimensional}.
This yields a surjective morphism of comodule algebras
\begin{equation}\label{eq:intro-character-realization}
	\xi_{\ps,\pg}
	=(\veps_{\ps,\pg}\ten\id)\circ\left.\DDi{s}\right|_{s=1}
	\colon\bfU^\imath\longrightarrow\bfB,
\end{equation}
where $\veps_{\ps,\pg}$ is the corresponding character.

\subsection{}
The morphisms $\xi_{\ps,\pg}$ give rise to a new family of restriction functors from 
suitable subcategories of $\bfU$-modules to $\bfU^\imath$-modules over $\bbC$.
The difficulty is that the coefficients of $y_i(z)$ are infinite series in $\FC^{-1}$, so 
$\FC$ cannot in general be specialized directly to a nonzero complex number.
We resolve this difficulty using the contour integral method of Gautam--Toledano 
Laredo \cite{gautam-toledano-16,gautam-toledano-17}.

Let $\sO$ denote the category of $\bfU$-modules whose restriction to $U_q(\g)$ lies 
in the category $\mathcal O \subset \Mod(U_q(\g))$, see Section \ref{sec:catO} for details.
On each weight space of a module in $\sO$, the series $\phi_i^\pm(z)$ and the two 
halves of $x_i^\pm(z)$ act as Laurent expansions of operator-valued rational 
functions, whose Laurent coefficients can therefore be recovered by contour integration.

For $\zeta\in\bbC^\times$, we define $\sO_\zeta\subset\sO$ by a genericity condition 
on these poles.
Specifically, for every node $i$, we require
\[
p_1p_2\notin\{\zeta^{-1},q_i^{a_{i\tau(i)}}\zeta^{-1}\}
\]
whenever $p_1,p_2$ are poles associated with $i$ or $\tau(i)$, taken over all weight 
spaces.
This separates the poles from their images under the involution 
$z\mapsto\zeta^{-1}z^{-1}$ and 
from the poles of the rational parameters $\GF_i(z),\SF_i(z)$.
The category $\sO_\zeta$ is dense in $\sO$ in the sense that every $V\in\sO$ belongs to 
$\sO_\zeta$ after a spectral shift outside an at most countable exceptional set.

Through contour integral formulae, we define the specialized action of the generators 
of $\bfB$ on every module in $\sO_\zeta$, yielding a restriction functor
\begin{equation}\label{eq:intro-restriction}
	\Res_{\zeta}^{\ps,\pg}\colon\sO_\zeta\longrightarrow\sM_\zeta.
\end{equation}
Here $\sM_\zeta$ denotes a subcategory of $\bfU^\imath$-modules, where the 
generating series act locally as Laurent expansions of rational functions.
Through the same method, we interpret the Drinfeld coaction as a rational, partially 
defined action of $\sO_\zeta$ on $\sM_\zeta$, analogous to the Drinfeld tensor product.
The restriction functor $\Res_{\zeta}^{\ps,\pg}$ is then recovered by acting on the one-dimensional module corresponding to $\GF_i(z), \SF_i(z)$.

\subsection{}
Under additional genericity assumptions, the action of the Drinfeld generators of 
$\bfU$ can be recovered directly from the restricted action of $\bfU^\imath$ through 
$\bfB$. This observation yields a surprisingly strong result on the structure of 
$\bfU^\imath$-modules obtained by restriction. Specifically, in 
Theorem~\ref{thm:generic-submodules-cyclicity-restriction}, we prove that, for a
module $V\in\sO$ with rational weights and a generic spectral shift
$V(s)\in\sO_\zeta$, every $\bfU^\imath$-submodule of
$\Res_{\zeta}^{\ps,\pg}(V(s))$ is stable under the derived quantum loop algebra 
$\bfU'$. When $\bfU'=\bfU$, \eg when $\g$ is of finite type, restriction 
gives a bijection on the entire lattice of submodules.
We also show that, in arbitrary type and for the same generic shifts, 
irreducible modules remain irreducible under restriction, while highest weight modules 
remain cyclic on their highest weight vectors.
The former result gives an analogue in category $\sO$ of the finite-dimensional generic 
irreducibility theorem of \cite{appel-vlaar-25a}. 
For infinite-dimensional modules, it is entirely new, even when $\g$ is finite-dimensional.

\subsection{}
We now turn to the abelian K-matrix, which we construct by adapting the methods 
developed in \cite{appel-gautam-26}.
For this construction, we assume $|q|\neq1$ and require that the determinant of the 
modified quantum Cartan matrix does not vanish at any nonzero integral power of 
$q$, and that its two extreme Laurent coefficients belong to $\{\pm1\}$.
The conditions on the Cartan matrix hold in finite type and in affine type other than 
$\sfA_1^{(1)}$.

In Theorem~\ref{thm:basic-K}, we construct two basic abelian K-matrices 
$\bKK^{0,\pm}_V(s)$ on modules in $\sO_\zeta$.
Their formal expressions are exponentials of series in $s^{\pm1}$ with coefficients in 
$\bfU^0$, multiplied by an $s$-independent operator diagonal on weight spaces.
The compatibility with the Drinfeld tensor product takes the form
\begin{equation}\label{eq:intro-basic-coproduct}
	\bKK^{0,\pm}_{V_1\DT{s_2/s_1}V_2}(s_1)
	=\bigl(\bKK^{0,\pm}_{V_1}(s_1)\ten\bKK^{0,\pm}_{V_2}(s_2)\bigr)
	\RR^{0,\pm}_{V_1^\bpsi,V_2}(s_1s_2),
\end{equation}
where $\bpsi$ is the involution obtained by composing the loop Chevalley involution 
with $\tau$ and the spectral shift by $\zeta^{-1}$.
After fixing the diagonal normalization, equation~\eqref{eq:intro-basic-coproduct}, 
together with compatibility with spectral shifts, determines the exponential factor up 
to a group-like series.
We determine this remaining factor by requiring that
\begin{equation}\label{eq:intro-basic-intertwiner}
	\bKK^{0,\pm}_V(s)\colon
	\Res_{\zeta}^{\ps,\pg}(V(s))
	\longrightarrow
	\Res_{\zeta}^{\ps,\pg}(V^\bpsi(s^{-1}))
\end{equation}
be an intertwiner.
We establish this relation on representations by contour integration.
It holds after rational continuation of the currents, while it fails as a coefficientwise 
identity between their formal expansions.

On each finite-dimensional weight space, the resulting series converge in 
neighbourhoods of zero and infinity, respectively, and extend meromorphically to 
$\bbC^\times$.
Moreover, they satisfy twisted unitarity:
\[
\bKK^{0,+}_V(s)^{-1}=\bKK^{0,-}_{V^{\bpsi}}(s^{-1}).
\]
Finally, if $V$ is highest weight with rational weights, dividing each operator by its 
eigenvalue on a highest weight vector yields a common normalised abelian 
$K$-matrix, rational on every weight space, and therefore unique.
Indeed, its uniqueness with this normalization follows directly from the generic cyclicity 
established above.

\subsection{}
By \cite{kolb-20,appel-vlaar-25b}, the basic K-matrix $\bKK^{0,\pm}_{V}(s)$ 
readily yields a tensor K-matrix $\KK^{0,\pm}_{M,V}(s)$ acting on the tensor product 
of a $\bfU^\imath$-module $M$ and $\bfU$-module $V$. 
Specifically, for $V_1,V_2\in\sO_\zeta$,
\begin{equation}\label{eq:intro-tensor-K}
	\KK^{0,\pm}_{V_1,V_2}(s)
	=\RR^{0,\mp}_{V_1,V_2^\bpsi}(s^{-1})^{-1}
	\bigl(1\ten\bKK^{0,\pm}_{V_2}(s)\bigr)
	\RR^{0,\pm}_{V_1,V_2}(s).
\end{equation}
Note that, since 
$\bKK^{0,\pm}_{V}(s)$ depends on a choice of parameters $\GF_i(z), \SF_i(z)$, so 
does $\KK^{0,\pm}_{M,V}(s)$.
Let $\sM_\zeta^{\mathrm{res}}$ denote the full subcategory of $\sM_\zeta$ 
consisting of subquotients of representations in the image of the corresponding 
restriction functor.
In Theorem~\ref{thm:tensor-K}, we prove that these operators descend to $\sM_\zeta^{\mathrm{res}}$ as meromorphic $\bfU^\imath$-intertwiners from $M\DT{s}V$ to $M\DT{s^{-1}}V^{\bpsi}$.

The construction of the abelian tensor K-matrix $\KK^{0,\pm}_{M,V}(s)$ gives a 
natural definition of boundary $q$-characters for both $\bfU$-modules and
$\bfU^\imath$-modules. In analogy with the construction of the boundary transfer 
matrices carried out in \cite[Sec.~6]{appel-vlaar-25c}, we set
\begin{equation}\label{eq:intro-reduced-tensor-K}
	\dKK^{0,\pm}(s)
	=\bigl(1\ten\bKK^{0,\pm}(s)\bigr)^{-1}\KK^{0,\pm}(s).
\end{equation}
Tracing over the second factor gives the boundary 
$q$-character of a $\bfU$-module, which is readily multiplicative with respect to the 
Drinfeld tensor product.
Tracing over the first factor gives the boundary $q$-character of a
$\bfU^\imath$-module, which satisfies a module identity for the Drinfeld action 
involving a folding of the usual $q$-character on the $\bfU$-module.
Both identities readily follow from the coproduct identities of $\dKK^{0,\pm}(s)$.

These constructions give, respectively, the abelian counterpart of the boundary 
transfer matrices of the first and fourth authors \cite{appel-vlaar-25c} and a trace 
realization of the boundary $q$-characters 	for $\bfU^\imath$-modules of the second 
and third authors \cite{Przez-23, li-przezdziecki-25, LP-eval, li-przezdziecki-26}.

\subsection{}\label{ss:intro-shift}
Finally, we extend our construction to \emph{shifted} $\imath$quantum loop
algebras $\bfU^\imath_\mu$. These algebras were introduced in simply-laced
split type by the second and third authors \cite{LP-gklo}, and in quasi-split
$\mathsf{ADE}$ type by Lu--Wang--Weekes \cite{lu-wang-weekes-26}.
In \cite{LP-Lev}, the second and third authors obtained minimalistic
presentations of Levendorskii type for $\bfU^\imath_\mu$ with anti-dominant
shift $\mu$, and used them to construct embeddings into shifted quantum
affine algebras, yielding a graded coideal structure.
Here we obtain analogous results in the topological setting by extending
our direct proof of the Drinfeld coaction to the shifted setting.
In Theorem~\ref{thm:shifted-coaction}, we construct a coassociative family
of algebra homomorphisms
\[
\DDi{s}{}_{\mu_1,\mu_2}\colon
\bfU^\imath_{\mu_1+\mu_2}\longrightarrow
\Lext{(\bfU^\imath_{\mu_1}\ten\bfU_{\mu_2})[\FC^{\pm1}]}{s^{-1};s},
\]
where $\mu_1,\mu_2$ are arbitrary coweights, inducing algebra homomorphisms
\[
\Psi_\mu\colon
\bfU^\imath_\mu\longrightarrow
\Lext{\bfU_\mu}{\FC;\FC^{-1}}.
\]
In both cases, the formulae on the generating series agree with those in the
unshifted case.

As an application, the shifted coaction provides a conceptual explanation
for the existence and structure of $\imath$GKLO representations.
The homomorphisms constructed in \cite{LP-gklo,lu-wang-weekes-26} give
infinite-dimensional representations by difference operators and are
expected to define truncations quantizing trigonometric analogues of
fixed point loci and symmetric quotients of the affine Grassmannian
\cite{BPT,lu-wang-weekes-25}.

In Corollary~\ref{ss: iGKLO}, we obtain a remarkably simple construction of $\imath$GKLO-type representations
 by pulling back the ordinary GKLO representations of $\bfU_\mu$
from \cite{GKLOq,FT19} along $\Psi_\mu$. The three types of 
summands in the $\imath$GKLO formulae (positive, negative and no difference operator) naturally correspond to the three summands in 
\eqref{eq:intro-coideal-generators}. Moreover, the admissibility constraint on 
the parameters $\ps_i(z)$ exactly matches the condition on the appearance of 
the no-difference-operator summand observed earlier in \cite{lu-wang-weekes-25,lu-wang-weekes-26}.

\subsection{Outline of the paper}
In Section~\ref{s:QL}, we review quantum loop algebras, the rational action of their 
generating series, and the Drinfeld coproduct.
In Sections~\ref{s:coideal-QLSP} and~\ref{s:comodule-realization}, we construct the 
topological coideal subalgebra, derive the parameter constraints, and relate it to 
characters of $\bfU^\imath$.
We then construct restriction and the Drinfeld action on representations, and prove 
generic irreducibility and cyclicity in Section~\ref{s:representations-rationality}.
After recalling the abelian R-matrix in Section~\ref{s:abelian-R-matrix}, we construct the 
basic abelian K-matrix in Section~\ref{s:basic-abelian-K-matrix} and its tensor version 
and boundary $q$-characters in Section~\ref{s:tensor-abelian-K-matrix}.
In Section~\ref{s:shifted-generalizations}, we extend the constructions to shifted 
quantum affine algebras and obtain $\imath$GKLO representations.
We give the direct proof of the Drinfeld coaction for shifted $\imath$quantum loop groups in Appendix~\ref{app:shifted-coaction}.


\subsection*{Acknowledgements}
This work began during two Research in Groups stays: one at ICMS in Edinburgh in January 2025, and one at the Erwin Schr\"{o}dinger Institute in Vienna in June 2025 and January 2026. The authors are grateful to both institutions for their support and wonderful working conditions.

AA thanks Sachin Gautam for many valuable discussions about this project during his long-term visit to the University of Parma in 2025–2026. AA also thanks Simone Castellan and Valerio Toledano Laredo for helpful discussions.
\\

AA was partially supported by the INdAM Research Project E53C23001740001 and the PRIN Grant 2022HMBTTL.
JL was supported by the Beatriz Galindo Senior grant BG24/00114 from the Spanish Ministry of Science, Innovation and Universities. 
TP was supported by the EPSRC grant No. EP/W022834/1, as well as the Austrian Science Fund (FWF)  grant PAT9039323, Grant-DOI 10.55776 /PAT9039323.
BV was supported by the International Scientists Program of the Beijing Natural Sciences Foundation (grant number IS24003).\\

During the preparation of this manuscript, the authors used OpenAI’s GPT-6 Astra to 
assist with proofreading, layout editing, and literature searches.
All suggestions were critically evaluated by the authors, and all arguments, computations, and references in the final manuscript were independently verified. The authors take full responsibility for the content and final form of the paper.
	

\section{Quantum loop algebras}\label{s:QL}

\subsection{Symmetrizable Kac--Moody algebras \cite{kac-90}}\label{ss:sym-km-summary}

\nc{\crt}[1]{h_{#1}}
\nc{\drc}[1]{\delta_{#1}}
\nc{\symd}[1]{d_{#1}}
\nc{\symdi}[1]{d_{#1}^{-1}}
\nc{\iip}[2]{\langle#1,#2\rangle}
\nc{\Diag}{\operatorname{diag}}
\nc{\rt}[1]{\alpha_{#1}}
\nc{\tinv}[1]{t_{#1}}

Let $\bfI$ be a finite set, and let $\sfA=(a_{ij})_{i,j\in\bfI}$ be a generalized Cartan matrix, that is,
$a_{ii}=2$, $a_{ij}\in\IZ_{\leqslant 0}$ for $i\neq j$, and $a_{ij}=0$ implies $a_{ji}=0$.
Let $(\h,\Pi,\Pi^\vee)$ be a realization of $\sfA$, namely:
\begin{itemize}[leftmargin=1.5em]\itemsep0.25cm
	\item $\h$ is a finite--dimensional complex vector space\footnote{In contrast to \cite{kac-90}, we do not impose the minimality condition $\dim \h = 2|\bfI|-\operatorname{rank}(\sfA)$, see \cite[Sec.~11]{appel-toledano-19b}.}
	\item $\Pi=\{\alpha_i\}_{i\in\bfI}$ is a linearly independent subset of $\h^*$
	\item $\Pi^\vee=\{\crt{i}\}_{i\in\bfI}$ is a linearly independent subset of $\h$
	\item $\alpha_i(\crt{j})=a_{ji}$ for all $i,j\in\bfI$
\end{itemize}
The Kac--Moody algebra attached to $\sfA$ and the realization $(\h,\Pi,\Pi^\vee)$ is the Lie algebra $\g$
generated by $\h$ together with elements $\{e_i,f_i\}_{i\in\bfI}$, subject to the relations $[\h,\h]=0$ and
\begin{gather}
	[h,e_i]=\alpha_i(h)e_i
	\qquad\qquad
	[h,f_i]=-\alpha_i(h)f_i
	\qquad\qquad
	[e_i,f_j]=\drc{ij}\crt{i}
\end{gather}
as well as
\[
\mathsf{ad}(e_i)^{1-a_{ij}}(e_j)=0=\mathsf{ad}(f_i)^{1-a_{ij}}(f_j),
\qquad i\neq j.
\]
We denote by $\n^\pm\subset\g$ the Lie subalgebras generated by $\{e_i\}_{i\in\bfI}$ and
$\{f_i\}_{i\in\bfI}$, respectively.\\

From now on, assume that $\sfA$ is indecomposable and symmetrizable. Fix an invertible diagonal matrix
$\sfD=\Diag(\symd{i})_{i\in\bfI}$, with coprime integers $\symd{i}\in\IZ_{>0}$, such that
$\sfD\sfA$ is symmetric. Then there exists a symmetric non--degenerate bilinear form
$\iip{\cdot}{\cdot}$ on $\h$ satisfying $\iip{\crt{i}}{-}=\symdi{i}\alpha_i$.
The associated identification $\nu:\h\to\h^*$ is equivariant with respect to the action of the Weyl group $W$ and satisfies
$\nu(\crt{i})=\symdi{i}\alpha_i$. In particular, it restricts to an isomorphism
$\h'\simeq\sfQ\ten_{\IZ}\IC$, where $\h'$ is the span of
$\{\crt{i}\}_{i\in\bfI}$ and ${\mathsf{Q}}=\bigoplus_{i\in\bfI}\IZ{\alpha}_i\subseteq\h^*$ is the root lattice.
Observe that $\iip{\crt{i}}{\crt{i}}=2\symdi{i}$, and that the induced form on $\h^*$ satisfies
$\iip{\rt{i}}{\rt{i}}=2\symd{i}\in 2\IZ_{>0}$. Given the choice of a subspace $\h''$ such that $\h=\h'\oplus\h''$, the form $\iip{\cdot}{\cdot}$ extends uniquely to a non--degenerate invariant
symmetric bilinear form on $\g$, still denoted by $\iip{\cdot}{\cdot}$. It is characterized by the property
$\iip{e_i}{f_j}=\drc{ij}\symdi{i}$, and one has
$[x,y]=\iip{x}{y}\cdot\tinv{\alpha}$ for all $x\in\g_{\alpha}$ and $y\in\g_{-\alpha}$, where
$\tinv{\alpha}=\nu^{-1}(\alpha)$. Finally, we fix $\rho\in\h^*$ such that $\rho(\crt{i})=1$ for every $i\in\bfI$.

\subsection{Quantum loop algebras \cite{drinfeld-88, hernandez-05}}\label{ss:QL}	
Fix $\hbar\in\bbC\setminus\bbQ$, set $q=e^{\pi\ii\hbar}$ and $q_i=q^{d_i}$.\footnote{
For representation--theoretic purposes, it is convenient to work
with the so-called \emph{unrolled} quantum group, whose definition naturally
involves the logarithmic deformation parameter $\hbar$. The condition $\hbar\notin\bbQ$ is equivalent to requiring that
$q=e^{\pi\ii\hbar}$ is not a root of unity. The stronger analytic condition $|q|\neq1$ will be
imposed in Section \ref{ss:quantum-Cartan} for the construction of the
meromorphic abelian R-matrix.}
The quantum loop algebra $\bfU=U_q(L\g)$ is the unital associative
$\bbC$-algebra generated by the elements
\[
\{q^h\}_{h\in\h},
\qquad
\{\phi^\pm_{i,\pm k}\}_{i\in\bfI,\ k\in\bbZ_{\geqslant0}},
\qquad
\{x^\pm_{i,k}\}_{i\in\bfI,\ k\in\bbZ},
\]
subject to the following relations. 
Then, for all $i,j\in\bfI$, $h\in\h$, $k,l\in\bbZ_{\geqslant0}$, $r,s\in\bbZ$,
and $\veps,\veps'\in\{+,-\}$,\footnote{
We use the conventions $\phi^+_{i,n}=0$ ($n<0$) and $\phi^-_{i,n}=0$ ($n>0$).}
\begin{align}
	\label{eq:QL-full-rels-cartan}
	q^0=1,
	\qquad
	q^hq^{h'}&=q^{h+h'},
	\qquad
	\phi^\pm_{i,0}=q_i^{\pm h_i},
	\\
	\label{eq:QL-full-rels-phi-h}
	q^h\phi^\pm_{i,\pm k}
	=
	\phi^\pm_{i,\pm k}q^h,
	\quad&\qquad
	[\phi^\veps_{i,\veps k},\phi^{\veps'}_{j,\veps' l}]=0,
	\\[1em]
	\label{eq:QL-full-rels-weight}
	q^h x^\pm_{i,r}q^{-h}
	&=
	q^{\pm\alpha_i(h)}x^\pm_{i,r},
	\\
	\label{eq:QL-full-rels-phi-x-plus}
	\phi^\veps_{i,\veps k}x^+_{j,r-1}
	-
	q_i^{a_{ij}}{\phi^\veps_{i,\veps k-1}}x^+_{j,r}
	&=
	q_i^{a_{ij}}x^+_{j,r-1}\phi^\veps_{i,\veps k}
	-
	x^+_{j,r}{\phi^\veps_{i,\veps k-1}},
	\\
	\label{eq:QL-full-rels-phi-x-minus}
	{\phi^\veps_{i,\veps k-1}}x^-_{j,r}
	-
	q_i^{a_{ij}}\phi^\veps_{i,\veps k}x^-_{j,r-1}
	&=
	q_i^{a_{ij}}x^-_{j,r}{\phi^\veps_{i,\veps k-1}}
	-
	x^-_{j,r-1}\phi^\veps_{i,\veps k},
	\\[1em]
	\label{eq:QL-full-rels-x-x-plus}
	q_i^{a_{ij}}x^+_{i,r-1}x^+_{j,s}
	-
	x^+_{i,r}x^+_{j,s-1}
	&=
	x^+_{j,s}x^+_{i,r-1}
	-
	q_i^{a_{ij}}x^+_{j,s-1}x^+_{i,r},
	\\
	\label{eq:QL-full-rels-x-x-minus}
	q_i^{a_{ij}}x^-_{i,r}x^-_{j,s-1}
	-
	x^-_{i,r-1}x^-_{j,s}
	&=
	x^-_{j,s-1}x^-_{i,r}
	-
	q_i^{a_{ij}}x^-_{j,s}x^-_{i,r-1},
	\\[1em]
	\label{eq:QL-full-rels-double}
	[x^+_{i,r},x^-_{j,s}]
	&=
	\frac{\delta_{ij}}{q_i-q_i^{-1}}
	\left(
	\phi^+_{i,r+s}
	-
	\phi^-_{i,r+s}
	\right),
\end{align}
and, for $i\neq j$, $N=1-a_{ij}$, and $r_1,\dots, r_N\in\bbZ$,
\begin{align}
	\label{eq:QL-full-rels-serre}
	\Sym_{\ul{r}}
	\sum_{k=0}^{N}
	(-1)^k
	\begin{bmatrix}N\\ k \end{bmatrix}_{q_i}
	x^\pm_{i,r_{1}}\cdots x^\pm_{i,r_k}
	x^\pm_{j,s}
	x^\pm_{i,r_{k+1}}\cdots x^\pm_{i,r_{N}}
	=0
\end{align}
where $\Sym_{\ul{r}}$ denotes symmetrization in  the indices $r_1,\dots,r_N$.\footnote{For an indeterminate $p$, the quantum integer and quantum factorial are
\[
[n]_p=\sum_{j=0}^{n-1}p^{-n+1+2j}\aand
[n]_p!=\prod_{j=1}^{n}[j]_p,
\]
with $[0]_p=0$ and $[0]_p!=1$.
The quantum binomial coefficient is
\[
\begin{bmatrix}n\\k\end{bmatrix}_p
=\frac{[n]_p!}{[k]_p!\,[n-k]_p!}
\qquad (0\les k\les n).
\]
These expressions are Laurent polynomials in $p$, hence admit
specialization at every $p\in\bbC^\times$. At $p=1$, they recover
the ordinary integers, factorials, and binomial coefficients, respectively.}

\begin{remark}
The subalgebra generated by $\{q^h\}_{h\in\h}$, and $\{x^{\pm}_{i,0}\}_{i\in\bfI}$
is isomorphic to the (unrolled) Drinfeld--Jimbo quantum group $U_q(\g)$ associated to $\g$. 
\hfill\rmkend
\end{remark}

\subsection{Fields}\label{ss:QL-fields}

For every $i\in\bfI$, define $\phi^{+}_i(z)\in\fext{\bfU}{z}$, $\phi^{-}_i(z)\in\fext{\bfU}{z^{-1}}$, and $x^{\pm}_i(z)\in\fext{\bfU}{z,z^{-1}}$ by
\[
\phi_{i}^{\pm}(z) = \sum_{k \ges 0} \phi_{i, \pm k}^{\pm} z^{\pm k}
\qq\qq 
x_{i}^{\pm}(z) = \sum_{k \in \bbZ} x_{i,k}^{\pm} z^{k}. 
\]
For every $i,j\in\bfI$, set
\[
\ZF_{ij}(z,w)=q_i^{a_{ij}}z-w
\aand 
\RF_{ij}(z,w)
=
-\frac{\ZF_{ij}(z,w)}{\ZF_{ji}(w,z)}
=
\frac{q_i^{a_{ij}}z-w}{z-q_i^{a_{ij}}w}.
\]
Note that $\RF_{ij}(z,w)$ depends only on $z/w$, hence we shall often write $\RF_{ij}(z/w)=\RF_{ij}(z,w)$. Since $d_ia_{ij}=d_ja_{ji}$, one has $q_i^{a_{ij}}=q_j^{a_{ji}}$, hence
\[
\ZF_{ij}(z,w)=\ZF_{ji}(z,w),
\qquad
\RF_{ij}(z,w)=\RF_{ji}(z,w),
\qquad
\RF_{ij}(z,w)^{-1}=\RF_{ij}(w,z).
\]
Let $\fdel(z)=\sum_{k\in\bbZ}z^k$
be the formal $\delta$-distribution. Then the defining relations imply the following current relations. For all $i,j\in\bfI$,
\begin{align}
	\label{eq:commutative-rel-field}
	[\phi^{\pm}_i(z),\phi^{\pm}_j(w)]=0=[\phi^{\pm}_i(z),\phi^{\mp}_j(w)],
\end{align}
\begin{align}
	\label{eq:weight-rel-fields-plus}
	\phi^\pm_i(z)x^{+}_j(w)=\RF_{ij}(w,z)x^{+}_j(w)\phi^{\pm}_i(z),\\
	\label{eq:weight-rel-fields-minus}
	\phi^\pm_i(z)x^{-}_j(w)=\RF_{ij}(z,w)x^{-}_j(w)\phi^{\pm}_i(z),
\end{align}
\begin{align}
	\label{eq:quadratic-rel-fields-plus}
	\ZF_{ji}(z,w)x^{+}_i(z)x^{+}_j(w)=-\ZF_{ij}(w,z)x^{+}_j(w)x^{+}_i(z),\\
	\label{eq:quadratic-rel-fields-minus}
	\ZF_{ji}(w,z)x^{-}_i(z)x^{-}_j(w)=-\ZF_{ij}(z,w)x^{-}_j(w)x^{-}_i(z),
\end{align}
\begin{align}
	\label{eq:double}
	[x^+_i(z),x^-_j(w)]
	=
	\frac{\delta_{ij}}{q_i-q_i^{-1}}
	\left(
	\fdel(w/z)\phi^+_i(w)
	-
	\fdel(z/w)\phi^-_i(z)
	\right),
\end{align}
and, for $i\neq j$, with $N=1-a_{ij}$,
\begin{align}\label{eq:serre-relations}
	\operatorname{Sym}_{\ul{z}}
	\sum_{k=0}^{N}
	(-1)^k
	\begin{bmatrix}N\\ k\end{bmatrix}_{q_i}
	x^{\pm}_i(z_1)\cdots x^{\pm}_i(z_k)x^{\pm}_j(w)
	x^{\pm}_i(z_{k+1})\cdots x^{\pm}_i(z_N)
	=0,
\end{align}
where $\operatorname{Sym}_{\ul{z}}$ denotes symmetrization in the variables $z_1,\dots,z_N$.\\ 

Finally, define $\{h_{i,\pm k}\}_{k>0}$ by
\begin{equation}\label{eq:log-phi}
\phi^{\pm}_i(z)
=
\phi^{\pm}_{i,0}
\exp\left(
\pm(q_i-q_i^{-1})
\sum_{k>0}h_{i,\pm k}z^{\pm k}
\right),
\end{equation}
and set $h^\pm_{i}(z)=\pm(q_i-q_i^{-1})
\sum_{k>0}h_{i,\pm k}z^{\pm k}$.
In particular, we have $[h_{i, k}, x^{\pm}_{j,r}]=\pm\frac{1}{k}[ka_{ij}]_{q_i}x^{\pm}_{j,r+k}$, hence
\begin{equation}\label{eq:adjoint-h}
[h_{i, k}, x^{\pm}_{j}(w)]=\pm\frac{1}{k}[ka_{ij}]_{q_i}w^{-k}x^{\pm}_j(w).
\end{equation}

\subsection{Representations and rationality \cite{beck-kac-96, hernandez-07}}  \label{sec:catO}
Define $\phi^{+}_i(z), x^{\pm}_i(z)^+\in\fext{\bfU}{z}$ and  $\phi^{-}_i(z), x^{\pm}_i(z)^-\in\fext{\bfU}{z^{-1}}$ by
\begin{align}\label{eq:QL-fields}
	\phi_{i}^{+}(z) &= \sum_{k \ges 0} \phi_{i,k}^{+} z^{k},
	\qq\qq \hspace{4.5mm}
	x_{i}^{\pm}(z)^+ = \sum_{k \ges 0} x_{i,k}^{\pm} z^{k},\\
	\phi_{i}^{-}(z) &= \sum_{k \ges 0} \phi_{i,-k}^{-} z^{-k},
	\qq\qq
	x_{i}^{\pm}(z)^- = -\sum_{k < 0} x_{i,k}^{\pm} z^{k}.
\end{align}
Let $\sO$ be the full subcategory of $\bfU$-modules $V$ whose restriction to $U_q(\g)$
is in category $\OUqg$, that is
\begin{enumerate}[label=(O\arabic*)]\itemsep0.25cm
\item the action of the operators $\{q^h\}_{h\in\h}$ is semisimple, \ie
\begin{equation}
	V=\bigoplus_{\mu\in\h^*}V_\mu	
\end{equation}
where $V_\mu=\{v\in V\,\vert\, \forall h\in\h,\,q^hv=q^{\mu(h)}v\}$;
\item for every $\mu\in\h^*$, $\dim(V_\mu)<\infty$;
\item the set of weights $P(V)=\{\mu\in\h^*\,\vert\,V_\mu\neq0\}$ is contained in a finite union of cones of the form $\lambda-\sfQ_+$.\\
\end{enumerate}
For every $V\in\sO$, $\mu\in P(V)$, $i\in\bfI$, and $\veps\in\{\pm\}$, the formal operators
\begin{equation}
\phi^\pm_i(z)\in\fext{\End(V_\mu)}{z^{\pm1}}
\aand
x^{\veps}_i(z)^{\pm}\in\fext{\Hom(V_\mu,V_{\mu+\veps\alpha_i})}{z^{\pm1}}
\end{equation}
are the expansions of rational operators $\Phif_i(z), \Xf_i^{\veps}(z)$ at $z=0$ and $z=\infty$. Specifically, set 
\begin{equation}\label{eq:H1}
\wt{h}_{i,\pm 1}
=[2]_{q_i}^{-1}h_{i,\pm1}=\pm\frac{\phi_{i,0}^{\mp}\phi_{i,\pm1}^{\pm}}{q^2_i-q^{-2}_i}.
\end{equation} Then,
\begin{align}\label{eq:QL-rational-gen-X}
\Xf_i^{\veps}(z)&=(1-\veps z\ad(\wt{h}_{i,1}))^{-1}(x_{i,0}^{\veps}),\\
\Phif_i(z)&=\phi^{-}_{i,0}+(q_i-q_i^{-1})[\Xf_i^+(z), x_{i,0}^-].
\end{align}

\begin{remarks}
	There is a slight abuse of notation in the definition of the rational functions $\Phif_i(z)$ and $\Xf_i^{\pm}(z)$. Since we do not assume that representations are finite-dimensional, $\Phif_i(z)$ and $\Xf_i^{\pm}(z)$ can only be defined at the level of individual weight spaces. For simplicity, we omit the dependence on the weight $\mu \in P(V)$.\hfill\rmkend
\end{remarks}

\subsection{Contour integral formulae \cite{gautam-toledano-17}}\label{ss:QL-contour}
Let $V\in\sO$. For every $\mu\in P(V)$ and $i\in\bfI$, let $\pole_i(V_\mu)\subset\bbC^\times$
be the set of poles of the rational functions $\Phif_i(z)$ and $\Xf^\veps_i(z)$ on $V_\mu$, and $\cC_{i,\mu}\subset\bbC^\times$ a Jordan curve enclosing $\pole_i(V_\mu)$ and not $0$.\footnote{By a Jordan curve we shall mean a finite disjoint union of simple closed curves, each oriented clockwise, whose bounded interior domains are pairwise disjoint. Our orientation convention is opposite to that of \cite{gautam-toledano-17}. This is due to our use of the inverse spectral parameter, which interchanges the roles of zero and infinity and therefore reverses the residue convention.} Then, for any $k\in\bbZ$ and $r\in\bbZ_{>0}$,
\begin{equation}\label{eq:contour-X-Phi}
x_{i,k}^{\pm}\vert_{V_\mu}=\oint_{\cC_{i,\mu}}\Xf_i^{\pm}(z)z^{-k-1}dz,
\qq\qq
\phi_{i,\pm r}^{\pm}\vert_{V_\mu}=\pm\oint_{\cC_{i,\mu}}\Phif_i(z)z^{\mp r-1}dz,
\end{equation}
and
\begin{equation}\label{eq:contour-Phi-0}
(\phi_{i,0}^+-\phi_{i,0}^-)\vert_{V_{\mu}}=	\oint_{\cC_{i,\mu}}\Phif_i(z)\frac{dz}{z},
\end{equation}
where $\oint_{\cC} f=\frac{1}{2\pi\ii}\int_{\cC}$.\footnote{
The same result holds for any weight module with finite-dimensional weight spaces. In particular, the condition (O3) is not necessary.}

\subsection{Drinfeld coproduct \cite{hernandez-05,hernandez-07, damiani-24}}
\label{ss:drinfeld-coproduct}
The algebra $\bfU$ is equipped with a morphism of algebras
\begin{equation}
\DD{s}\colon \bfU\to\Lext{\bfU\ten\bfU}{s^{-1};s}
\end{equation}
defined on the generators as follows:
\begin{equation}\label{eq:Drinfeld-generators}
	\begin{split}
	\DD{s}(\phi_{i,\pm k}^{\pm}) &= \sum_{a=0}^{k} \phi_{i,\pm a}^{\pm} \otimes \phi_{i,\pm(k-a)}^{\pm} s^{\pm(k-a)},\\
	\DD{s}(x_{i,k}^{+}) &= x_{i,k}^{+} \otimes 1 + \sum_{a \ges 0} \phi_{i,-a}^{-} \otimes x_{i,k+a}^{+} s^{k+a}, \\
	\DD{s}(x_{i,k}^{-}) &= 1\otimes x_{i,k}^{-} s^{k} + \sum_{a \ges 0} x_{i,k-a}^{-} \otimes \phi_{i,a}^{+} s^{a}.
	\end{split}
\end{equation}
Hence, in terms of fields, one has
\begin{equation}\label{eq:Drinfeld-fields}
\begin{split}
	\DD{s}(x_{i}^{+}(z)) &= x_{i}^{+}(z) \otimes 1 + \phi_{i}^{-}(z) \otimes x_{i}^{+}(zs), \\
	\DD{s}(x_{i}^{-}(z)) &= 1\otimes x_{i}^{-}(zs) + x_i^{-}(z) \otimes \phi_{i}^{+}(zs), \\
	\DD{s}(\phi_{i}^{\pm}(z)) &= \phi_{i}^{\pm}(z) \otimes \phi_{i}^{\pm}(zs).
\end{split}
\end{equation}
The map $\DD{s}$ is coassociative in the sense that
\begin{equation}\label{eq:Drinfeld-coassociative}
(\DD{s_1}\ten\id)\circ\DD{s_1s_2} = (\id\ten\DD{s_2})\circ\DD{s_1}.
\end{equation} 

\subsection{Rationality of Drinfeld tensor product \cite{hernandez-07, gautam-toledano-17}}
For every $V,W\in\sO$, the Drinfeld coproduct yields an action of $\bfU$ on 
$\Lext{V\ten W}{s^{-1};s}$. This is the Laurent expansion at $s=0$ of a family 
of actions of $\bfU$ on $V\ten W$, denoted $V\DT{s}W$, whose matrix coefficients are rational functions of $s$ \cite[Lemma~3.10]{hernandez-07}.
While for the commuting generators $\{\phi_i^\pm(z)\}_{i\in\bfI}$ this is clear, 
the rationality of the action of $\{x_{i}^{\pm}(z)\}_{i\in\bfI}$ on $V\DT{s}W$ 
becomes evident when expressed in terms of contour integrals \cite[Sec.~4.2]{gautam-toledano-17}.

Let $s\in\bbC^\times$ be such that, for every $\lambda\in P(V), \mu\in P(W)$ and $i\in\bfI$, 
\begin{equation}
	\pole_i(V_{\lambda})\cap s^{-1}\pole_i(W_\mu)=\emptyset
\end{equation} 
where $\pole_i(V_\lambda), \pole_i(W_\mu)\subset\bbC^\times$ are the sets of poles, see \ref{ss:QL-contour}.\footnote{
Since the weight spaces are countable and finite-dimensional, this condition excludes only a countable subset of Lebesgue measure zero in $\bbC^\times$.}

For every $\lambda\in P(V), \mu\in P(W)$,
let $\cC_{i,\lambda},\cC_{i,\mu}\subset\bbC^\times$ be two Jordan curves which do not enclose $0$ and are such that
\begin{itemize}\itemsep0.25cm
	\item $\cC_{i,\lambda}$ encloses $\pole_i(V_\lambda)$ and none of the points in $s^{-1}\pole_i(W_\mu)$;
	\item $\cC_{i,\mu}$ encloses $s^{-1}\pole_i(W_\mu)$ and none of the points in $\pole_i(V_\lambda)$.
\end{itemize}
Then, the action of $x_{i,k}^\pm$ on $V\DT{s}W$ is given by
	\begin{equation}\label{eq:Drinfeld-gen-contour}
	\begin{split}
		\DD{s}(x_{i,k}^{+}) 
		&= x_{i,k}^{+} \otimes 1 + \oint_{\cC_{i,\mu}}
		\left(\sum_{a \ges 0} \phi_{i,-a}^{-}w^{-a} \right)\otimes \Xf_i^+(sw)w^{-k-1}dw \\
		&=x_{i,k}^{+} \otimes 1 + \oint_{\cC_{i,\mu}}\Phif_i(w)\otimes \Xf_i^+(sw)w^{-k-1}dw,\\ 
		\DD{s}(x_{i,k}^{-})
		&= 1\otimes x_{i,k}^{-} s^{k} + \oint_{\cC_{i,\lambda}}\Xf_i^-(w)\ten
		\left(\sum_{a \ges 0} \phi_{i,a}^{+} s^{a}w^a\right)
		w^{-k-1}dw\\
		&=1\otimes x_{i,k}^{-} s^{k} + \oint_{\cC_{i,\lambda}}\Xf_i^-(w)\ten\Phif_i(sw)w^{-k-1}dw.
	\end{split}
\end{equation}
In terms of currents, the action on $V\DT{s}W$ reads
\begin{equation}\label{eq:Drinfeld-current-contour}
	\begin{split}
		\DD{s}(x_{i}^{+}(z)^{\pm}) 
		&=\Xf_i^+(z)\ten 1 + \oint_{\cC_{i,\mu}}\Phif_i(w)\otimes \Xf_i^+(sw)
		\frac{dw}{w-z},\\ 
		\DD{s}(x_{i}^{-}(z)^{\pm})
		&=1\otimes \Xf_{i}^{-}(sz) + \oint_{\cC_{i,\lambda}}\Xf_i^-(w)\ten\Phif_i(sw)\frac{dw}{w-z}.\\
	\end{split}
\end{equation}

\section{Coideal realization of quantum loop symmetric pairs}\label{s:coideal-QLSP}

We adapt the definition of quasi--split quantum symmetric pairs for Drinfeld--Jimbo quantum groups due to Letzter and Kolb \cite{letzter-02,kolb-14} to the setting of quantum loop algebras. We derive necessary constraints on the defining parameters from a minimality condition on the intersection with the quantum loop Cartan subalgebra, in analogy with the Letzter-Kolb theory of quantum symmetric pairs, and its generalization in \cite{regelskis-vlaar-20, regelskis-vlaar-22}.

In Sections~\ref{ss:B-def}--\ref{ss:psi-invariance}, we define the coideal subalgebra 
$\bfB$ depending on the parameters $\ps_i(z),\pg_i(z)$ ($i\in\bfI$). 
In Sections~\ref{ss:SF+GF}--\ref{ss:eta-aux}, we introduce a first restriction on the 
parameters, inducing weight and quadratic relations.
Finally, in Section~\ref{ss:DS-simply-laced}, we introduce further parameter 
restrictions, arising from the Drinfeld--Serre relations in simply-laced types.

	\subsection{The algebra $\bfB$}\label{ss:B-def}
	Let $\tau\colon\bfI\to\bfI$ be an involutive map such that $a_{\tau(i)\tau(j)}=a_{ij}$. In particular, one has $d_{\tau(i)}=d_i$ for every $i\in\bfI$.\footnote{
	Namely, the tuple $(d'_i)_{i\in\bfI}$ with $d'_i=d_{\tau(i)}$ is formed by coprime positive integers such that $d'_ia_{ij}=d_{\tau(i)}a_{\tau(i)\tau(j)}=d_{\tau(j)}a_{\tau(j)\tau(i)}=d'_ja_{ji}$. Hence, by uniqueness, $d'_i=d_i$.}
	The involution $\tau\colon\bfI\to\bfI$ naturally extends to an automorphism on $\bfU$. Specifically, we set
	$\tau(q^h)=q^{\tau(h)}$, relying on the lift of $\tau$ to $\h$ constructed by \cite{kac-wang-92}, and
	\begin{equation}
		\tau(x^{\pm}_i(z))=x^{\pm}_{\tau(i)}(z)\aand \tau(\phi^{\pm}_i(z))=\phi^{\pm}_{\tau(i)}(z).
	\end{equation}
	In particular, $\tau(h^{\pm}_i(z))=h^{\pm}_{\tau(i)}(z)$. Note also that $\tau$ preserves
	$U_q(\g)\subset\bfU$.\\
	
	For every $i\in\bfI$, we consider
	\begin{equation}
		\ps_{i}(z)=\sum_{k\in\bbZ}\ps_{i,k}z^k\in\fext{\bbC[\FC^{\pm1}]}{z^{\pm1}}
		\aand
		\pg_{i}(z)=\sum_{k\ges0}\pg_{i,-k}z^{-k}\in\fext{\bbC[\FC^{\pm1}]}{z^{-1}}
	\end{equation}
such that 
\begin{equation}\label{eq:cond-sigma-gamma}
\lim_{k\to\infty}\operatorname{deg}_{\FC}(\ps_{i,k})-k=-\infty
\quad\mbox{and}\quad
\lim_{k\to\infty}\operatorname{ord}_{\FC^{-1}}(\pg_{i,-k})=\infty,
\end{equation}
where $\operatorname{ord}_{\FC^{-1}}$ denotes the lowest exponent with respect to $\FC^{-1}$.
\\

We denote by $\bfB$ the subalgebra of $\Lext{\bfU}{\FC;\FC^{-1}}$ generated 
over $\bbC[\FC^{\pm1}]$ by the elements $\{q^{h}\}_{h\in\h^{-\tau}}$, $\{\eta_{i,r}\}_{i\in\bfI,r\ges0}$ and $\{y_{i,k}\}_{i\in\bfI, k\in\bbZ}$ defined as follows:
\begin{equation}\label{eq:B-generators}
\begin{aligned}
	\eta_{i,k} &= \sum_{a=0}^{k} \phi_{i,a}^{+} \phi_{\tau(i),-k+a}^{-} \FC^{k-a},\\
	y_{i,k} &= x_{i, -k}^{-} \FC^k + \sum_{a,b \ges 0} \pg_{i,-a} x_{\tau(i),k+a+b}^{+} \phi_{i,b}^{+} \FC^{-b} + \sum_{a \ges 0} \ps_{i,k+a} \phi_{i,a}^{+} \FC^{-a}.
\end{aligned}
\end{equation}
In particular, $\eta_{i,0} = q_i^{h_i-h_{\tau(i)}}$.
Note that $\eta_{i,k}\in\bfU[\FC]$, while $y_{i,k}\in\Lext{\bfU}{\FC;\FC^{-1}}$ thanks to the condition \eqref{eq:cond-sigma-gamma}. Specifically, the conditions on $\pg_i(z)$ and $\ps_i(z)$ guarantee, respectively, that the second and third terms in $y_{i,k}$ lie in $\Lext{\bfU}{\FC;\FC^{-1}}$.\\

Define the generating series $\eta_{i}(z)\in\fext{\bfB}{z}$ and $y_i(z)\in\fext{\bfB}{z,z^{-1}}$ by
\begin{equation}
	\eta_{i}(z)=\sum_{r\ges0}\eta_{i,r}z^r\aand y_i(z)=\sum_{k\in\bbZ}y_{i,k}z^k.
\end{equation}
Then, 
\begin{align} \label{eq: eta defi}
	\eta_{i}(z) &= \phi_{i}^{+}(z) \phi_{\tau(i)}^{-}(\FC^{-1}{z}^{-1}),\\ \label{eq: y defi}
	y_{i}(z) &= x_{i}^{-}(\FC^{-1}{z}^{-1}) + \pg_{i}(z) x_{\tau(i)}^{+}(z) \phi_{i}^{+}(\FC^{-1} z^{-1}) + \ps_{i}(z) \phi_{i}^{+}(\FC^{-1}z^{-1}).
\end{align}

\begin{lemma}
	Let $i\in\bfI$. Then,
\begin{equation}\label{eq:weight-B}
q^{h}y_i(z)=q^{-\alpha_i(h)}y_i(z)q^{h}
\end{equation}
for every $h\in\h^{\theta}$ if and only if $\sigma_i(z)=0$ whenever $\tau(i)\neq i$.
\end{lemma}
\begin{pf} It is enough to observe that 
$\alpha_i(h)=\theta(\alpha_i)(h)=-\alpha_{\tau(i)}(h)$ if $h\in\h^\theta$.
When $\tau(i)=i$, one has $\alpha_i(h)=0$ and \eqref{eq:weight-B}
reduces to $q^{h}y_i(z)=y_i(z)q^{h}$. Instead, when $\tau(i)\neq i$, the condition 
$\sigma_i(z)=0$ is equivalent to \eqref{eq:weight-B}.
\end{pf}

\subsection{Coideal property}\label{ss:coideal}
By imposing linearity in $\FC$, we extend the Drinfeld coproduct
to a morphism of algebras
\begin{align}
	\DD{s}\colon\Lext{\bfU}{\FC;\FC^{-1}}\to\Lext{\Lext{(\bfU\ten\bfU)}{s^{-1};s}}{\FC;\FC^{-1}},
\end{align}
which, by a slight abuse of notation, we denote by the same symbol.
The subalgebra $\bfB\subset\Lext{\bfU}{\FC;\FC^{-1}}$ is a right coideal with respect to $\DD{s}$. More precisely, we have the following result.

\begin{proposition}
	The restriction of the Drinfeld coproduct to $\bfB$ yields a morphism of algebras
	\begin{equation}
		\DD{s}\vert_{\bfB}\colon\bfB\to\Lext{(\bfB\ten\bfU)[\FC^{\pm1}]}{s^{-1};s}.
	\end{equation}	
\end{proposition}

\begin{pf}
	It is enough to consider the explicit formulae of the Drinfeld coproduct on the generators. In terms of fields, we have:
	\begin{align}\label{eq:DD-eta-field}
		\DD{s}(\eta_{i}(z))
		&=\phi_i^+(z)\phi_{\tau(i)}^-(\FC^{-1}z^{-1})\ten\phi_i^+(zs)\phi_{\tau(i)}^-(\FC^{-1}z^{-1}s)\\		
		&=\eta_{i}(z)\ten\phi_i^+(zs)\phi_{\tau(i)}^-(\FC^{-1}z^{-1}s)	
	\end{align}
	and
	\begin{align}\label{eq:DD-y-field}
		\DD{s}(y_{i}(z))&=
		1\ten x^{-}_i(\FC^{-1}z^{-1}s)
		+x^{-}_i(\FC^{-1}z^{-1})\ten\phi_i^+(\FC^{-1}z^{-1}s)\\
		&\quad+\pg_i(z)x_{\tau(i)}^+(z)\phi_i^+(\FC^{-1}z^{-1})\ten\phi^+_i(\FC^{-1}z^{-1}s)\\
		&\quad+\pg_i(z)\phi^-_{\tau(i)}(z)\phi_i^+(\FC^{-1}z^{-1})\ten x_{\tau(i)}^+(zs)\phi^+_i(\FC^{-1}z^{-1}s)\\
		&\quad+\ps_i(z)\phi^+_i(\FC^{-1}z^{-1})\ten\phi^+_i(\FC^{-1}z^{-1}s)\\
		&=1\ten x^{-}_i(\FC^{-1}z^{-1}s)+y_i(z)\ten\phi_i^+(\FC^{-1}z^{-1}s)\\
		&\qquad+\pg_i(z)\eta_i(\FC^{-1}z^{-1})\ten x_{\tau(i)}^+(zs)\phi^+_i(\FC^{-1}z^{-1}s)
	\end{align}
	Hence, on the generators $\eta_{i,k}$ and $y_{i,k}$, it reads: 
	\begin{align}
		\label{eq:DD-eta-gen}
		\DD{s}(\eta_{i,k})&=	\sum_{\substack{a,b,c\ges0\\a+b+c=k}}\eta_{i,a}\ten\phi^+_{i,b}\phi^-_{\tau(i),-c}\FC^{c}s^{b-c}\\
		\label{eq:DD-y-gen}
		\DD{s}(y_{i,k})&=
		1\otimes x^-_{i,-k}\FC^k s^{-k}
		+\sum_{a\ges 0}y_{i,k+a}\otimes \phi^+_{i,a}\FC^{-a}s^a\\
		&\quad+
		\sum_{a,b,c\ges 0}
		\pg_{i,-c}\eta_{i,a}\otimes x^+_{\tau(i),\,k+a+b+c}\phi^+_{i,b}
		\,\FC^{-a-b}s^{k+a+c+2b}.
	\end{align}
	Note that only finitely many negative powers of $s$ appear, and, for every $\ell\in\bbZ$, the coefficient of $s^{\ell}$ lives in $(\bfB\ten\bfU)[\FC^{\pm1}]$.
\end{pf}

\subsection{$\psi$-invariance}\label{ss:psi-invariance}
The (loop) Chevalley involution $\bom\colon\bfU\to\bfU$ is the algebra map given by
$\bom(q^h)=q^{-h}$ ($h\in\h$), and
\begin{equation}
	\bom(x^{\pm}_i(z))=-x^{\mp}_i(z^{-1})
	\aand
	\bom(\phi^{\pm}_i(z))=\phi^{\mp}_i(z^{-1}).
\end{equation}
Specifically, 
	\begin{equation}
	\bom(x^{\pm}_{i,k})=-x^{\mp}_{i,-k}
	\aand
	\bom(\phi^{\pm}_{i,\pm r})=\phi^{\mp}_{i,\mp r}
\end{equation}
for every $k\in\bbZ$ and $r\in\bbZ_{\ges0}$.
Note that $\bom(h^{\pm}_{i}(z))=h^{\mp}_{i}(z^{-1})$ and $\bom(h_{i,\pm r})=-h_{i,\mp r}$ for every $r\in\bbZ_{>0}$.
 Set $\theta=\tau\circ\bom$. 
\\

We extend both $\tau$ and $\bom$ to $\Lext{\bfU}{\FC;\FC^{-1}}$ by imposing linearity in $\FC$, and we set 
\begin{equation}
\bpsi=\SH_\FC^{-1}\circ\tau\circ\bom
\end{equation}
where $\SH_\FC\colon\bfU\to\bfU[\FC^{\pm1}]$ is the formal shift operator in $\FC$ given by
\begin{equation}
	\SH_\FC(x^{\pm}_i(z))=x^{\pm}_i(\FC z)
	\aand
	\SH_\FC(\phi^{\pm}_i(z))=\phi^{\pm}_i(\FC z).
\end{equation}
Note that $\tau\circ\bom=\bom\circ\tau$, $\tau\circ\SH_\FC=\SH_\FC\circ\tau$, and
$\SH_\FC^{-1}\circ\bom=\bom\circ\SH_\FC$. Hence, $\bpsi^2=\id$.

Define $\{\Leta_{i,r}\}_{r>0}$ by
\begin{equation}\label{eq:def-rho}
	\eta_i(z)=\phi^+_{i,0}\phi^-_{\tau(i),0}\exp((q_i-q_i^{-1})\sum_{r>0}\Leta_{i,r}z^r).
\end{equation}
Hence, for every $r\in\bbZ_{>0}$, $\Leta_{i,r}=h_{i,r}-h_{\tau(i),-r}\FC^{r}$.
Set $\Leta_i(z)=\sum_{r>0}\Leta_{i,r}z^r\in\fext{\bfB}{z}$.\\

Both $\bfU$ and $\bfB$ have large abelian subalgebras.
Let $\bfU^0\subset\bfU$ be the subalgebra generated by $\{q^h\}_{h\in\h}$ and $\{\phi^{\pm}_{i,\pm k}\}_{i\in\bfI,k>0}$ (or, equivalently, $\{h_{i,\pm k}\}_{i\in\bfI, k>0}$).
Similarly, let $\bfB^0\subset\bfB$ be the subalgebra generated by $\{q^{h}\}_{h\in\h^{\theta}}$ and $\{\eta_{i,k}\}_{i\in\bfI, k> 0}$ (or, equivalently, $\{\Leta_{i,k}\}_{i\in\bfI, k>0}$).

\begin{lemma}
$\bfB^0\subset\bfU^0[\FC^{\pm1}]^{\bpsi}$.
\end{lemma}

\begin{pf}
It is enough to observe that every element in $\h^{\theta}$ has the form $h-\tau(h)=h+\bpsi(h)$ for some $h\in\h$ and, for every $i\in\bfI$ and $r>0$, $\Leta_{i,r}=h_{i,r}-h_{\tau(i),-r}\FC^r=h_{i,r}+\bpsi(h_{i,r})$.
\end{pf}	

\subsection{Rational parameters}\label{ss:SF+GF}
In analogy with the approach to quantum symmetric pairs developed in \cite{letzter-02,kolb-14}, it is natural to ask that the inclusion $\bfB^0\subseteq\bfB\cap\Lext{\bfU^0}{\FC;\FC^{-1}}$ be an equality, or at least that 
the intersection be as close as possible to $\bfB^0$. 
In Sections~\ref{ss:weight-quadratic-relations}--\ref{ss:DS-simply-laced}, 
we shall exhibit explicit computations, showing that certain products of 
elements of $\bfB$ produce elements of 
$\Lext{\bfU^0}{\FC;\FC^{-1}}$ not lying in $\bfB^0$. The vanishing of 
these unwanted terms yields necessary constraints on the fields $\ps_i(z)$ 
and the series $\pg_i(z)$ ($i\in\bfI$), which we address in the proposition below. Specifically, they force
\[
\ps_i(z)=\SF_i(z)^+-\SF_i(z)^-,
\qquad
\pg_i(z)=\GF_i(z)^-,
\]
where $\SF_i(z)$ and $\GF_i(z)$ are rational functions depending on finitely many parameters, and the superscripts $+$ and $-$ denote expansions at $z=0$ and $z=\infty$, respectively. In particular, $\ps_i(z)$ and $\pg_i(z)$
automatically satisfy \eqref{eq:cond-sigma-gamma}. We conjecture that, 
with this choice of rational parameters, these necessary conditions are 
also sufficient, that is
\begin{equation}\label{eq:minimality}
	\bfB\cap\Lext{\bfU^0}{\FC;\FC^{-1}}=\bfB^0.
\end{equation}
We refer to \eqref{eq:minimality} as the \emph{minimality condition}.

\begin{proposition}\label{prop:sigma-gamma}
\hfill
\begin{enumerate}[leftmargin=2em]\itemsep0.25cm
\item 	
For every $i\in\bfI$, let $\ps_{i}(z)=\sum_{k\in\bbZ}\ps_{i,k}z^k\in\fext{\bbC[\FC^{\pm1}]}{z^{\pm1}}$
be such that
\begin{align}\label{cond:sigma}
	\ps_{i,k+1}=\ps_{i,k-1}\FC
\end{align}
for every $k\in\bbZ$. Then, there exist $\psa{i}, \psb{i}\in\bbC[\FC^{\pm1}]$ ($i\in\bfI$) such that
\begin{equation}
\ps_{i,2k}=\psa{i}\FC^k\quad\mbox{and}\quad\ps_{i,2k+1}=\psb{i}\FC^{k+1}
\end{equation}
and therefore
\begin{equation}\label{eq:sigma-delta}
\ps_{i}(z)=(\psa{i}+\psb{i}z\FC)\fdel(\FC z^2).
\end{equation}
In particular, $\ps_i(z)$ satisfies 
\eqref{eq:cond-sigma-gamma}. Moreover, set
\begin{equation}\label{eq:Sigma-expansions}
	\ps_i(z)^+=\sum_{k\geqslant 0}\ps_{i,k}z^k
	\qquad\text{and}\qquad
	\ps_i(z)^-=-\sum_{k<0}\ps_{i,k}z^k.
\end{equation} 
Then, $\ps_i(z)^{\pm}$ are, respectively, the expansion at $z=0$ and $z=\infty$ of the
rational function over $\bbC(\FC)$
\begin{equation}\label{eq:SF}
	\SF_i(z)=\frac{\psa{i}+\psb{i}\FC z}{1-\FC z^2}.
\end{equation}
\item 
Consider the system of equations
\begin{equation}\label{cond:gamma}
	\begin{aligned}
		\delta_{j\tau(i)}\,\fdel(\FC zw)
		((q_i^{a_{i\tau(i)}}z-w)&\gamma_{\tau(i)}(w)
		+(q_i^{a_{i\tau(i)}}w-z)\gamma_i(z))=\\
		=&(q_i^{a_{ij}}-1)(q_i-q_i^{-1})(z+w)\sigma_i(z)\sigma_j(w)
	\end{aligned}
\end{equation}
where $i,j\in\bfI$ and $\pg_{i}(z)=\sum_{k\ges0}\pg_{i,-k}z^{-k}\in\fext{\bbC[\FC^{\pm1}]}{z^{-1}}$.
Then, \eqref{cond:gamma} admits a solution if and only if 
\begin{equation}\label{cond:sigma-vanishing}
	\S=\{i\in\bfI\mid\SF_i(z)\neq0\}
	\,\mbox{ is a completely disconnected subset of }\,
	\bfI^\tau,
\end{equation}
\ie for every $i,j\in\S$, $\tau(i)=i$, $\tau(j)=j$, and $a_{ij}=0$ if $i\neq j$. 

\item 
Let $\pg_i(z)\in\fext{\bbC[\FC^{\pm1}]}{z^{-1}}$ ($i\in\bfI$) be a solution of \eqref{cond:gamma}.
Then, there exist
$\pga{i}, \pgb{i}\in\bbC[\FC^{\pm1}]$ ($i\in\bfI$) with
\begin{equation}
\pgb{i}=-\pgb{\tau(i)}
\end{equation}
such that the following holds.\\
\begin{enumerate}[leftmargin=2em]\itemsep0.25cm 
\item If $i\in\bfI^\tau$, 
$\pg_i(z)$ is given by $\pg_{i,0}=\pga{i}$ and
\begin{equation}
\begin{aligned}
	\gamma_{i,-2(r+1)}&=(
(q_i^2-1)q_i^{2r}\pga{i}
-(q_i-q_i^{-1})(q_i^{2r+2}-1)(\psa{i}^2+\FC \psb{i}^2))\FC^{-r-1},\\
\gamma_{i,-(2r+1)}&=-(q_i-q_i^{-1})(q_i^{2r+2}+q_i^{2r}-2)\psa{i}\psb{i}\,\FC^{-r}
\end{aligned}
\end{equation}
for $r\geqslant 0$. In particular, $\pg_i(z)$ satisfies \eqref{eq:cond-sigma-gamma}.
Moreover, $\pg_i(z)$ is the expansion at $z=\infty$ of the rational function over $\bbC(\FC)$
\begin{equation}\label{eq:GF-fixed}
\GF_i(z)=\frac{1-\FC z^2}{q_i^2-\FC z^2}
\left(\pga{i}+q_i(q_i-q_i^{-1})^2\,\SF_i(z)\SF_i(\FC^{-1}z^{-1})\right).
\end{equation}
\item If $i\in\bfI\setminus\bfI^\tau$, $\pg_i(z)$ is given by $\pg_{i,0}=\pga{i}$,
$\pg_{i,-1}=\pgb{i}$, and
\begin{equation}
	\begin{aligned}
		\gamma_{i,-2(r+1)}&=(q_i^{a_{i\tau(i)}}\pga{i}-\pga{\tau(i)})q_i^{ra_{i\tau(i)}}\FC^{-r-1},\\
		\gamma_{i,-(2r+1)}&=\pgb{i}q_i^{ra_{i\tau(i)}}\FC^{-r}
	\end{aligned}
\end{equation}
for $r\geqslant 0$.
In particular, $\pg_i(z)$ satisfies \eqref{eq:cond-sigma-gamma}. Moreover, $\pg_i(z)$
is  the expansion at $z=\infty$ of the rational function over $\bbC(\FC)$
\begin{equation}\label{eq:GF-not-fixed}
\GF_i(z)=\frac{\pga{\tau(i)}-\pgb{i}\FC z-\pga{i}\FC z^2}{q_i^{a_{i\tau(i)}}-\FC z^2}
\end{equation}
with $\GF_i(0)=q_i^{-a_{i\tau(i)}}\pga{\tau(i)}$ and $\GF_i(\infty)=\pga{i}$.
\end{enumerate}

\end{enumerate}
\end{proposition}

Part (1) is clear. The proofs of (2) and (3) are carried out in \ref{sss:pf-parameter-2} and \ref{sss:pf-parameter-3a}--\ref{sss:pf-parameter-3b}, respectively.

\begin{corollary}
	Let $\SF_i(z), \GF_i(z)$ ($i\in\bfI$) be defined as in Proposition~\ref{ss:SF+GF}.
	Then,
	\begin{equation}\label{eq:gamma-symmetry}
		\GF_i(\FC^{-1}z^{-1})
		=
		\frac{q_i^{a_{i\tau(i)}}-\FC z^2}{1-q_i^{a_{i\tau(i)}}\FC z^2}
		\GF_{\tau(i)}(z)
	\end{equation}
	for every $i\in\bfI$.
\end{corollary}

\begin{pf}
If $\tau(i)=i$, the relation \eqref{eq:gamma-symmetry} is clear. If $\tau(i)\neq i$,
\begin{equation}
	\GF_{\tau(i)}(z)
	=\frac{\pga{i}+\pgb{i}\FC z-\pga{\tau(i)}\FC z^2}{q_i^{a_{i\tau(i)}}-\FC z^2}
	=\frac{\pga{\tau(i)}-\pgb{i}z^{-1}-\pga{i}\FC^{-1}z^{-2}}{1-q_i^{a_{i\tau(i)}}\FC^{-1}z^{-2}}
\end{equation}
and the result follows.
\end{pf}

\subsubsection{Proof of Proposition~\ref{ss:SF+GF} (2)}\label{sss:pf-parameter-2}
First, we observe that
\begin{equation}
(z+w)\ps_i(z)\ps_j(w)=0\,\mbox{ if and only if }\, \ps_i(z)=0\,\mbox{ or }\,\ps_j(w)=0.	
\end{equation}
Indeed, the coefficient of $z^mw^n$ in $(z+w)\ps_i(z)\ps_j(w)$ is given by
\begin{equation}
\begin{array}{lcr}
	(\psa{i} \psb{j}+\psb{i} \psa{j})\,\FC^{(m+n)/2} & \mbox{when} & m+n \text{ even},\\[2mm]
	(\psa{i} \psa{j}+\FC \psb{i} \psb{j})\,\FC^{(m+n-1)/2} & \mbox{when} &  m+n \text{ odd}.
\end{array}
\end{equation}
Assume $(z+w)\ps_i(z)\ps_j(w)=0$ and $\ps_i(z)\neq0$. Then, one has
\begin{equation}
\begin{bmatrix}
\psa{i} & \psb{i}\\ \FC \psb{i} & \psa{i}
\end{bmatrix}	
\begin{bmatrix}
\psb{j} \\ \psa{j}	
\end{bmatrix}
=
\begin{bmatrix}
0\\0	
\end{bmatrix}
\end{equation}
whose only solution is given by $\psa{j}=0=\psb{j}$.\footnote{
More precisely, the determinant of the system is 
$\psa{i}^2-\FC \psb{i}^2$. If $\psb{i}=0$, then $\psa{i}^2\neq0$. If $\psb{i}\neq 0$ and $\psa{i}^2-\FC \psb{i}^2=0$, then $\FC$ should admit a square root in $\bbC(\FC)$.} Hence, $\ps_j(w)=0$.

Let $i\in\bfI$ with $\tau(i)\neq i$. For $j=i$, \eqref{cond:gamma} reduces to 
\[
(q_i^{2}-1)(q_i-q_i^{-1})(z+w)\sigma_i(z)\sigma_i(w)=0,
\]
hence $\ps_i(z)=0$ and $i\not\in\S$. Moreover, let $i,j\in\bfI^\tau$ with $j\neq i=\tau(i)$. Then,
\eqref{cond:gamma} reduces to 
\[
(q_i^{a_{ij}}-1)(q_i-q_i^{-1})(z+w)\sigma_i(z)\sigma_j(w)=0.
\]
Thus, if $a_{ij}\neq 0$, then $\sigma_i(z)=0$ or $\sigma_j(w)=0$. Hence, if $i,j\in\S$ with $i\neq j$, then $a_{ij}=0$.

\subsubsection{Proof of Proposition~\ref{ss:SF+GF} (3a)}\label{sss:pf-parameter-3a}
Let $i\in\bfI^\tau$. Then, \eqref{cond:gamma} with $j=i$ reads
\begin{equation}
	\begin{aligned}
		\fdel(\FC zw)
		((q_i^{2}z-w)\pg_{i}(w)
		+(q_i^{2}w-z)\pg_i(z))=q_i(q_i-q_i^{-1})^2(z+w)\sigma_i(z)\sigma_i(w).
	\end{aligned}
\end{equation}
First, we compute the solution of the associated homogeneous system:
\begin{equation}
	\begin{aligned}
		\L(z,w)=\fdel(\FC zw)((q_i^{2}z-w)\pg'_{i}(w)+(q_i^{2}w-z)\pg'_i(z))=0.
	\end{aligned}
\end{equation}
Set $\L(z,w)=\sum_{m,n\in\bbZ}\L_{m,n}z^mw^n$. Then, a direct expansion gives
\begin{equation}
\L_{m,n}=
\begin{cases}
-2\pg'_{i,-1}\FC^n & m=n\\
((q_i^2-1)\pg'_{i,0}-\FC \pg'_{i,-2})\FC^n & m=n+1\\
(q_i^2\pg'_{i,-(2r+1)}-\FC \pg'_{i,-(2r+3)})\FC^{n+2r+1} & m=n+2r+2\\
(q_i^2\pg'_{i,-(2r+2)}-\FC \pg'_{i,-(2r+4)})\FC^{n+2r+2} & m=n+2r+3
\end{cases}	
\end{equation}
for every $n\in\bbZ$ and $r\geqslant 0$. Hence, one has $\pg'_{i,-(2r+1)}=0$ ($r\geqslant 0$) and 
\begin{equation}
\begin{aligned}
\pg'_{i}(z)
=\pga{i}\left(1+\sum_{r\geqslant 1}(q_i^2-1)q_i^{2r-2}\FC^{-r}z^{-2r}\right)
\end{aligned}
\end{equation}
where $\pga{i}=\pg'_{i,0}$ is a free parameter in $\bbC[\FC^{\pm1}]$. One readily checks
that $\pg'_{i}(z)$ is the expansion at $z=\infty$ of the rational function
	\begin{equation}
	\begin{aligned}
		\pga{i}\frac{1-\FC^{-1}z^{-2}}{1-q_i^2\FC^{-1}z^{-2}}
		=\pga{i}\frac{1-\FC z^{2}}{q_i^2-\FC z^{2}}.
	\end{aligned}
\end{equation}
Let $\pg_i'''(z)$ be the expansion at $z=\infty$ of
\begin{equation}
	\begin{aligned}
		\frac{1-\FC z^{2}}{q_i^2-\FC z^{2}}\SF_i(z)\SF_i(\FC^{-1}z^{-1})
		&=-z^{-1}
		\frac{(\psa{i}+\psb{i}z^{-1})(\psb{i}+\FC^{-1}\psa{i}z^{-1})}
		{(1-\FC^{-1}z^{-2})(1-q_i^2\FC^{-1}z^{-2})}\\
		&=-
		\psa{i}\psb{i}z^{-1}\frac{1+\FC^{-1}z^{-2}}{(1-\FC^{-1}z^{-2})(1-q_i^2\FC^{-1}z^{-2})}
		\\
		&\hspace{1.5cm}
		-
		(\psa{i}^2+\FC \psb{i}^2)\frac{\FC^{-1}z^{-2}}{(1-\FC^{-1}z^{-2})(1-q_i^2\FC^{-1}z^{-2})}.
	\end{aligned}
\end{equation}
It remains to prove that $\pg_i''(z)=q_i(q_i-q_i^{-1})^2\pg'''_i(z)$ is a particular solution of \eqref{cond:gamma}. Set $\pg_i''(z)=\sum_{r\geqslant 0}\pg''_{i,-r}z^{-r}$.
We have $\pg''_{i,0}=0$, $\pg''_{i,-1}=-q_i(q_i-q_i^{-1})^2\psa{i}\psb{i}$, and
\begin{equation}
\begin{aligned}
	\pg''_{i,-(2r+1)}&=-(q_i-q_i^{-1})((q_i^{2}+1)q_i^{2r}-2)\psa{i}\psb{i}\FC^{-r},\\
	\pg''_{i,-2r}&=-(q_i-q_i^{-1})(q_i^{2r}-1)(\psa{i}^2+\psb{i}^2\FC)\FC^{-r}
\end{aligned}
\end{equation}
for $r\geqslant 1$.\footnote{
Note that
\[
\frac{1+t}{(1-t)(1-q^2t)}
=
\sum_{r\ges 0}\frac{(q^2+1)q^{2r}-2}{q^2-1}\,t^r,
\qquad
\frac{t}{(1-t)(1-q^2t)}
=
\sum_{r\ges 1}\frac{q^{2r}-1}{q^2-1}\,t^r.
\]
}
Hence,
\begin{equation}
\begin{aligned}
q_i^2\FC^r\pg''_{i,-(2r+1)}-\FC^{r+1}\pg''_{i,-(2r+3)}=q_i(q_i-q_i^{-1})^2\cdot 2\psa{i}\psb{i},\\
q_i^2\FC^{r+1}\pg''_{i,-(2r+2)}-\FC^{r+2}\pg''_{i,-(2r+4)}=q_i(q_i-q_i^{-1})^2\cdot (\psa{i}^2+\psb{i}^2\FC)
\end{aligned}
\end{equation}
which coincides with \eqref{cond:gamma}.

\subsubsection{Proof of Proposition~\ref{ss:SF+GF} (3b)}\label{sss:pf-parameter-3b}
Let $i\in\bfI$ with $\tau(i)\neq i$. Set $j=\tau(i)$. Note that $a_{ij}=a_{ji}$ and $q_i^{a_{ij}}=q_j^{a_{ij}}$. Then, by \ref{sss:pf-parameter-2}, $\ps_i(z)=0=\ps_j(z)$,
and \eqref{cond:gamma} reads
\begin{equation}
	\fdel(\FC zw)
	((q_i^{a_{ij}}z-w)\gamma_{j}(w)
	+(q_i^{a_{ij}}w-z)\gamma_i(z))=0.
\end{equation}
In terms of the coefficients of $\pg_i(z),\pg_j(w)$, this is equivalent to the conditions
\begin{equation}
	\pg_{i,-1}+\pg_{j,-1}=0,\quad
	\FC\pg_{i,-2}=q_i^{a_{ij}}\pg_{i,0}-\pg_{j,0},\,\quad
	\FC\pg_{j,-2}=q_i^{a_{ij}}\pg_{j,0}-\pg_{i,0},
\end{equation}
and, for all $r\geqslant0$,	
	\begin{equation}
	\begin{aligned}
		q_i^{a_{ij}}& \pg_{k, -(2r+1)}-\FC \pg_{k, -(2r+3)}=0
		\quad\mbox{and}\quad
		q_i^{a_{ij}}\pg_{k,-2(r+1)}-\FC \pg_{k,-2(r+2)}=0
	\end{aligned}	
\end{equation}
for $k=i,j$.
Set $\pga{i}=\pg_{i,0}$, $\pga{j}=\pg_{j,0}$, and $\pgb{i}=\pg_{i,-1}=-\pg_{j,-1}=-\pgb{j}$.
Then, for $r\geqslant0$,
\begin{equation}
\pg_{i,-(2r+1)}=q_i^{ra_{ij}}\FC^{-r}\pgb{i}=-\pg_{j,-(2r+1)}	
\end{equation}
and
\begin{equation}
\begin{aligned}
\pg_{i,-2(r+1)}=q_i^{ra_{ij}}(q_i^{a_{ij}}\pga{i}-\pga{j})\FC^{-r-1},\\
\pg_{j,-2(r+1)}=q_i^{ra_{ij}}(q_i^{a_{ij}}\pga{j}-\pga{i})\FC^{-r-1}.
\end{aligned}
\end{equation}

\subsection{Weight and quadratic relations}\label{ss:weight-quadratic-relations}
Henceforth, for every $i\in\bfI$, we fix the rational functions over $\bbC(\FC)$
\begin{equation}\label{eq:SF+GF}
\begin{aligned}
	\SF_i(z)&=\delta_{i\tau(i)}\frac{\psa{i}+\psb{i}\FC z}{1-\FC z^2},\\
	\GF_i(z)&=
		\frac{1}{q_i^{a_{i\tau(i)}}-\FC z^2}\left(
		\pga{\tau(i)}-\pgb{i}\FC z-\pga{i}\FC z^2\right.\\
		&\hspace{3cm}\left.+
		(1-\FC z^2)q_i(q_i-q_i^{-1})^2\,\SF_i(z)\SF_{i}(\FC^{-1}z^{-1})\right),
\end{aligned}
\end{equation}
satisfying the conditions of Proposition~\ref{ss:SF+GF}. That is,
\begin{itemize}\itemsep0.25cm
	\item $\S=\{i\in\bfI\mid\SF_i(z)\neq0\}$ is a completely disconnected subset of $\bfI^\tau$;
	\item for every $i\in\bfI$, $\pgb{i}=-\pgb{\tau(i)}$.
\end{itemize}
Finally, we set $\ps_i(z)=\SF_i(z)^+-\SF_i(z)^-$ and $\pg_i(z)=\GF_i(z)^-$.

\begin{remark}
Eventually, we shall specialize $\FC$ to a complex number
$\zeta\in\bbC^\times$. For simplicity, we therefore assume that
$\pga{i}\in\bbC^\times$ for every $i\in\bfI$.\footnote{
	The condition $\pga{i}\neq0$ will be needed for the
	construction of the abelian K-matrix in
	Section~\ref{s:basic-abelian-K-matrix}.
	Our argument establishing the rationality of the generators
	$\eta_i(z)$ on representations, see Proposition~\ref{ss:eta-aux}, only requires the weaker condition $\GF_i(z)\neq0$ for every $i\in\bfI$. 
	This requirement can be avoided by modifying the definition of $\bfB$, 
	replacing the Cartan generators $\eta_i(z)$ with
	$\deta_i(z)=\pg_i(\FC^{-1}z^{-1})\eta_i(z)$.
	If $\GF_i(z)=0$ for all $i\in\bfI$, then also $\SF_i(z)=0$
	for all $i\in\bfI$, and the modified algebra is generated over
	$\bbC[\FC^{\pm1}]$ by $\{x^-_{i,k}\}_{i\in\bfI, k\in\bbZ}$
	and $\{q^h\}_{h\in\h^{\theta}}$.
}
\hfill\rmkend
\end{remark}

This choice of parameters is instrumental to the following result.
Recall that $\RF_{ij}(z,w)=\frac{q_i^{a_{ij}}z-w}{z-q_i^{a_{ij}}w}=\RF_{ij}(z/w)$, see \ref{ss:QL}.
\begin{proposition}
For any $i,j\in\bfI$,
\begin{equation}\label{eq:weight-rel-eta-y}
\eta_i(z)y_j(w)=\RF_{ij}(\FC zw)\RF_{\tau(i)j}(w/z)y_j(w)\eta_i(z)	
\end{equation}
and
\begin{equation}\label{eq:quadratic-y-y}
\small
\begin{aligned}
	(q_i^{a_{ij}}z-w)&y_i(z)y_{j}(w)+(q_i^{a_{ij}}w-z)y_{j}(w)y_i(z)=\\
	&\hspace{-1.15cm}=\delta_{i\tau(j)}\frac{\fdel(\FC zw)}{q_i-q_i^{-1}}\left((q_i^{a_{i\tau(i)}}z-w)\pg_{\tau(i)}(w)\eta_{\tau(i)}(z)
	+(q_i^{a_{i\tau(i)}}w-z)\pg_i(z)\eta_{i}(w)\right).
\end{aligned}
\end{equation}
\end{proposition}
The proof is carried out in \ref{sss:QSP-weight-rel}--\ref{sss:QSP-quadratic-rel}.
This result is equivalently expressed in terms of the generators $\eta_{i,k}$ ($k\in\bbZ_{\geqslant0}$) and $y_{i,r}$ ($r\in\bbZ$) as follows.
\begin{corollary}\label{cor:explicit-rels}
For any $i,j\in\bfI$, $k\in\bbZ_{\geqslant0}$ and $r,s\in\bbZ$,
	\begin{equation}\label{eq:weight-y-modes}
	\begin{aligned}
		&
		\FC\,\eta_{i,k-1}y_{j,r-2}
		-
		\FC q_i^{a_{\tau(i)j}}\eta_{i,k-2}y_{j,r-1}
		-
		q_i^{a_{ij}}\eta_{i,k}y_{j,r-1}
		+
		q_i^{a_{ij}+a_{\tau(i)j}}\eta_{i,k-1}y_{j,r}
		\nonumber\\
		&\qquad =
		q_i^{a_{ij}+a_{\tau(i)j}}\FC\,y_{j,r-2}\eta_{i,k-1}
		-
		q_i^{a_{ij}}\FC\,y_{j,r-1}\eta_{i,k-2}
		-
		q_i^{a_{\tau(i)j}}y_{j,r-1}\eta_{i,k}
		+
		y_{j,r}\eta_{i,k-1}.
	\end{aligned}
	\end{equation}
	and
	\begin{equation}\label{eq:quadratic-y-modes}
	\begin{aligned}
		&q_i^{a_{ij}}
		\left(
		y_{i,r-1}y_{j,s}
		+
		y_{j,s-1}y_{i,r}
		\right)
		-
		\left(
		y_{i,r}y_{j,s-1}
		+
		y_{j,s}y_{i,r-1}
		\right)
		\\
		&\quad =
		\frac{\delta_{i\tau(j)}}{q_i-q_i^{-1}}
		\sum_{n\in\bbZ}
		\FC^n
		\Big(
		q_i^{a_{i\tau(i)}}
		\eta_{\tau(i),r-n-1}
		\pg_{\tau(i),s-n}
		-
		\eta_{\tau(i),r-n}
		\pg_{\tau(i),s-n-1}
		\\
		&\hspace{4.5cm}
		+
		q_i^{a_{i\tau(i)}}
		\pg_{i,r-n}
		\eta_{i,s-n-1}
		-
		\pg_{i,r-n-1}
		\eta_{i,s-n}
		\Big).
	\end{aligned}
	\end{equation}
\end{corollary}


\subsubsection{Proof of relation \eqref{eq:weight-rel-eta-y}}\label{sss:QSP-weight-rel}
By \eqref{eq:weight-rel-fields-minus}, one has
\begin{align}
\phi_i^+(z)\phi_{\tau(i)}^-&(\FC^{-1}z^{-1})x^-_j(\FC^{-1}w^{-1})=\\
&=\RF_{ij}(\FC zw)\RF_{\tau(i)j}(w/z)
x^-_j(\FC^{-1}w^{-1})\phi_i^+(z)\phi_{\tau(i)}^-(\FC^{-1}z^{-1})
\end{align}
and, by \eqref{eq:weight-rel-fields-plus}, 
\begin{align}
	\phi_i^+(z)\phi_{\tau(i)}^-&(\FC^{-1}z^{-1})x^+_{\tau(j)}(w)\phi_{j}^+(\FC^{-1}w^{-1})=\\
	&=\RF_{\tau(i)\tau(j)}(\FC zw)\RF_{i\tau(j)}(w/z)
	x^+_{\tau(j)}(w)\phi_{j}^+(\FC^{-1}w^{-1})\phi_i^+(z)\phi_{\tau(i)}^-(\FC^{-1}z^{-1})\\
	&=\RF_{ij}(\FC zw)\RF_{\tau(i)j}(w/z)
	x^+_{\tau(j)}(w)\phi_{j}^+(\FC^{-1}w^{-1})\phi_i^+(z)\phi_{\tau(i)}^-(\FC^{-1}z^{-1})
\end{align}
where the second equality follows from the identities
\begin{equation}\label{eq:rho-symmetry}
\RF_{ij}(x)=\RF_{\tau(i)\tau(j)}(x)
\quad\mbox{and}\quad
\RF_{i\tau(j)}(x)=\RF_{\tau(i)j}(x).
\end{equation}
It remains to show that
\begin{equation}\label{eq:cond-weight-rel}
	\left(1-\RF_{ij}(\FC zw)\RF_{\tau(i)j}(w/z)\right)\ps_j(w)
	\phi_i^+(z)\phi_{\tau(i)}^-(\FC^{-1}z^{-1})\phi_j^+(\FC^{-1}w^{-1})
	=0.
\end{equation}
Recall that, by definition \eqref{eq:sigma-delta},
\begin{equation}
	\ps_j(w)=\delta_{j\tau(j)}(\ps_{j,0}+\ps_{j,-1}\FC w)\,\fdel(\FC w^2).
\end{equation}
Hence, \eqref{eq:cond-weight-rel} is obvious when $\tau(j)\neq j$.
On the other hand, 
	\begin{align}\label{eq:sigma-delta-2}
	\RF_{ij}(\FC zw)\RF_{ij}(w/z)\fdel(\FC w^2)=
	\RF_{ij}(z/w)\RF_{ij}(w/z)\fdel(\FC w^2)=\fdel(\FC w^2),
\end{align}
where the second identity follows since $\RF_{ij}(z,w)=\RF_{ij}(w,z)^{-1}$.
Thus, \eqref{eq:cond-weight-rel} follows when $\tau(j)=j$.

\subsubsection{Proof of relation~\eqref{eq:quadratic-y-y}}\label{sss:QSP-quadratic-rel}
We shall expand \eqref{eq:quadratic-y-y} and observe that many components in the LHS vanish.\\
First, one has
\begin{align}\label{eq:quadratic-A}
(q_i^{a_{ij}}z-w)&x^-_i(\FC^{-1}z^{-1})x^-_j(\FC^{-1}w^{-1})\\
&+(q_i^{a_{ij}}w-z)x^-_j(\FC^{-1}w^{-1})x^-_i(\FC^{-1}z^{-1})=0	
\end{align}
since it reduces to the quadratic relation \eqref{eq:quadratic-rel-fields-minus}.\\
Similarly, one has
\begin{equation}\label{eq:quadratic-B}
\tiny
\begin{split}
	(q_i^{a_{ij}}&z-w)\pg_i(z)x^+_{\tau(i)}(z)\phi^+_i(\FC^{-1}z^{-1})
	\pg_j(w)x^+_{\tau(j)}(w)\phi^+_j(\FC^{-1}w^{-1})\\
	&\qquad+(q_i^{a_{ij}}w-z)\pg_j(w)x^+_{\tau(j)}(w)\phi^+_j(\FC^{-1}w^{-1})\pg_i(z)x^+_{\tau(i)}(z)\phi^+_i(\FC^{-1}z^{-1})=\\
	&=(q_i^{a_{ij}}z-w)x^+_{\tau(i)}(z)	x^+_{\tau(j)}(w)
	\RF_{i\tau(j)}(\FC zw)\pg_i(z)\pg_j(w)\phi^+_i(\FC^{-1}z^{-1})\phi^+_j(\FC^{-1}w^{-1})\\
	&\qquad\qquad+(q_i^{a_{ij}}w-z)x^+_\tau(j)(w)x^+_\tau(i)(z)
	\RF_{\tau(i)j}(\FC zw)\pg_i(z)\pg_j(w)\phi^+_j(\FC^{-1}w^{-1})\phi^+_i(\FC^{-1}z^{-1})=\\
	&=((q_i^{a_{\tau(i)\tau(j)}}z-w)x^+_{\tau(i)}(z)	x^+_{\tau(j)}(w)
	+(q_i^{a_{\tau(i)\tau(j)}}w-z)x^+_{\tau(i)}(w)x^+_{\tau(j)}(z))\cdot\\
	&\hspace{4.45cm}\cdot
	\RF_{i\tau(j)}(\FC zw)\pg_i(z)\pg_j(w)\phi^+_i(\FC^{-1}z^{-1})\phi^+_j(\FC^{-1}w^{-1})=0
\end{split}
\end{equation}
where the second equality follows from $a_{ij}=a_{\tau(i)\tau(j)}$ and \eqref{eq:rho-symmetry}, and the third equality follows from the quadratic relation \eqref{eq:quadratic-rel-fields-plus}.\\
For the mixed terms, one has
\begin{equation}\label{eq:quadratic-AC}
\tiny
\begin{split}
(q_i&^{a_{ij}}z-w)
(x_i^-(\FC^{-1}z^{-1})\ps_j(w)\phi^+_j(\FC^{-1}w^{-1})+\ps_i(z)\phi^+_i(\FC^{-1}z^{-1})x_j^-(\FC^{-1}w^{-1}))\\
&+(q_i^{a_{ij}}w-z)
(x_j^-(\FC^{-1}w^{-1})\ps_i(z)\phi^+_i(\FC^{-1}z^{-1})+\ps_j(w)\phi^+_j(\FC^{-1}w^{-1})x_i^-(\FC^{-1}z^{-1}))=\\
=&\ps_j(w)((q_i^{a_{ij}}z-w)
x_i^-(\FC^{-1}z^{-1})\phi^+_j(\FC^{-1}w^{-1})
+(q_i^{a_{ij}}w-z)\phi^+_j(\FC^{-1}w^{-1})x_i^-(\FC^{-1}z^{-1}))\\
&+\ps_i(z)\left((q_i^{a_{ij}}z-w) \phi^+_i(\FC^{-1}z^{-1})x_j^-(\FC^{-1}w^{-1}) +(q_i^{a_{ij}}w-z)
x_j^-(\FC^{-1}w^{-1})\phi^+_i(\FC^{-1}z^{-1})\right)\\
=&0
\end{split}
\end{equation}
since both summands vanish by the weight relation \eqref{eq:weight-rel-fields-minus}.\\
Finally, one has
\begin{equation}\label{eq:quadratic-BC}
\tiny
\begin{split}
(q_i^{a_{ij}}&z-w)
(
\pg_i(z)x_{\tau(i)}^+(z)\phi_i^+(\FC^{-1}z^{-1})
\ps_j(w)\phi^+_j(\FC^{-1}w^{-1})\\
&\hspace{4.1cm}+\ps_i(z)\phi^+_i(\FC^{-1}z^{-1})\pg_j(w)x_{\tau(j)}^+(w)\phi^+_j(\FC^{-1}w^{-1}))\\
+&(q_i^{a_{ij}}w-z)
(\pg_j(w)x_{\tau(j)}^+(w)\phi^+_j(\FC^{-1}w^{-1})\ps_i(z)\phi^+_i(\FC^{-1}z^{-1})\\
&\hspace{5.1cm}+\ps_j(w)\phi^+_j(\FC^{-1}w^{-1})\pg_{i}(z)x_{\tau(i)}^+(z)\phi^+_i(\FC^{-1}z^{-1}))=\\
=&((q_i^{a_{ij}}z-w)
+(q_i^{a_{ij}}w-z)\RF_{j\tau(i)}(\FC zw))
\pg_i(z)\ps_j(w)x_{\tau(i)}^+(z)\phi_i^+(\FC^{-1}z^{-1})\phi^+_j(\FC^{-1}w^{-1})\\
&+((q_i^{a_{ij}}z-w)\RF_{i\tau(j)}(\FC zw)
+(q_i^{a_{ij}}w-z))
\ps_i(z)\pg_j(w)x_{\tau(j)}^+(w)\phi_i^+(\FC^{-1}z^{-1})\phi^+_j(\FC^{-1}w^{-1})\\
=&0
\end{split}
\end{equation}
since
\begin{itemize}\itemsep0.25cm
\item if $\tau(i)\neq i$ and $\tau(j)\neq j$, then $\ps_i(z)=0=\ps_j(w)$;
\item if $\tau(i)=i$ and $\tau(j)=j$, then $\RF_{ij}(\FC zw)\ps_i(z)=\RF_{ij}(w/z)\ps_i(z)$;
\item if $\tau(i)=i$ and $\tau(j)\neq j$, then $\ps_j(w)=0$,
\begin{align}
	\RF_{i\tau(j)}(\FC zw)=\RF_{\tau(i)j}(\FC zw)=\RF_{ij}(\FC zw),
\end{align}
and therefore $\RF_{i\tau(j)}(\FC zw)\ps_i(z)=\RF_{ij}(w/z)\ps_i(z)$;
\item if $\tau(j)=j$ and $\tau(i)\neq i$, then $\ps_i(z)=0$,
\begin{align}
	\RF_{j\tau(i)}(\FC zw)=\RF_{\tau(j)i}(\FC zw)=\RF_{ji}(\FC zw)=\RF_{ij}(\FC zw),
\end{align}
and therefore $\RF_{j\tau(i)}(\FC zw)\ps_j(w)=\RF_{ij}(z/w)\ps_j(w)$.
\end{itemize}
Therefore, it remains to study the two mixed products with a non-trivial component along $\Lext{\bfU^0}{\FC;\FC^{-1}}$.
One has
\begin{equation}\label{eq:quadratic-A-B}
\tiny
\begin{split}
(q_i^{a_{ij}}z-w)&(x_i^-(\FC^{-1}z^{-1})\pg_j(w)x_{\tau(j)}^+(w)\phi^+_j(\FC^{-1}w^{-1})\\
&\hspace{4cm}+\pg_i(z)x_{\tau(i)}^+(z)\phi^+_i(\FC^{-1}z^{-1})x_j^-(\FC^{-1}w^{-1}))\\
&\hspace{-1.25cm}+(q_i^{a_{ij}}w-z)(x_j^-(\FC^{-1}w^{-1})\pg_i(z)x_{\tau(i)}^+(z)\phi^+_i(\FC^{-1}z^{-1})\\
&\hspace{4cm}+\pg_j(w)x_{\tau(j)}^+(w)\phi^+_j(\FC^{-1}w^{-1})x_i^-(\FC^{-1}z^{-1}))=\\
&\hspace{-1.5cm}=(q_i^{a_{ij}}z-w)(\pg_j(w)x_i^-(\FC^{-1}z^{-1})x_{\tau(j)}^+(w)\phi^+_j(\FC^{-1}w^{-1})\\
&\hspace{2cm}+\pg_i(z)\RF_{ij}(w/z)x_{\tau(i)}^+(z)x_j^-(\FC^{-1}w^{-1})\phi^+_i(\FC^{-1}z^{-1}))\\
&\hspace{-1.25cm}+(q_i^{a_{ij}}w-z)(\pg_i(z)x_j^-(\FC^{-1}w^{-1})x_{\tau(i)}^+(z)\phi^+_i(\FC^{-1}z^{-1})\\
&\hspace{2cm}+\pg_j(w)\RF_{ij}(z/w)x_{\tau(j)}^+(w)x_i^-(\FC^{-1}z^{-1})\phi^+_j(\FC^{-1}w^{-1}))\\
&\hspace{-1.5cm}=(q_i^{a_{ij}}z-w)\pg_j(w)\left[x_i^-(\FC^{-1}z^{-1}),x_{\tau(j)}^+(w)\right]\phi^+_j(\FC^{-1}w^{-1})\\
&\hspace{0cm}+(q_i^{a_{ij}}w-z)(\pg_i(z)\left[x_j^-(\FC^{-1}w^{-1}),x_{\tau(i)}^+(z)\right]\phi^+_i(\FC^{-1}z^{-1}))\\
&\hspace{-1.5cm}=\frac{\delta_{i\tau(j)}}{q_i-q_i^{-1}}(q_i^{a_{ij}}z-w)\pg_j(w)\left(\fdel(\FC zw)\phi^-_i(w)
-\fdel((\FC zw)^{-1})\phi^+_i(\FC^{-1}z^{-1})\right)\phi^+_j(\FC^{-1}w^{-1})\\
&\hspace{-1.25cm}+\frac{\delta_{j\tau(i)}}{q_i-q_i^{-1}}(q_i^{a_{ij}}w-z)\pg_i(z)\left(\fdel(\FC zw)\phi^-_{j}(z)
-\fdel((\FC zw)^{-1})\phi^+_{j}(\FC^{-1}w^{-1})\right)\phi^+_i(\FC^{-1}z^{-1})\\
&\hspace{-1.5cm}=\delta_{i\tau(j)}\frac{\fdel(\FC zw)}{q_i-q_i^{-1}}\left((q_i^{a_{i\tau(i)}}z-w)\pg_{\tau(i)}(w)\eta_{\tau(i)}(\FC^{-1}w^{-1})
+(q_i^{a_{i\tau(i)}}w-z)\pg_i(z)\eta_{i}(\FC^{-1}z^{-1})\right)\\
&\hspace{-1.25cm}-\delta_{i\tau(j)}\frac{\fdel(\FC zw)}{q_i-q_i^{-1}}
\left(
(q_i^{a_{i\tau(i)}}z-w)\pg_{\tau(i)}(w)+(q_i^{a_{i\tau(i)}}w-z)\pg_i(z)
\right)
\phi^+_{\tau(i)}(\FC^{-1}w^{-1})\phi^+_i(\FC^{-1}z^{-1})
\end{split}	
\end{equation}
where the second equality follows from \eqref{eq:weight-rel-fields-plus}--\eqref{eq:weight-rel-fields-minus} and the identity
$\RF_{ij}(x)=\RF_{ji}(x)$, the third equality follows from the identities
\begin{equation}
(q_i^{a_{ij}}z-w)\RF_{ij}(w/z)=-(q_i^{a_{ij}}w-z)\quad\mbox{and}\quad	(q_i^{a_{ij}}w-z)\RF_{ij}(z/w)=-(q_i^{a_{ij}}z-w),
\end{equation}
and the fourth equality follows from \eqref{eq:double} and the identity $\fdel(x)=\fdel(x^{-1})$.\\

Finally, the last contribution is given by
	\begin{equation}
	\small
	\begin{split}
		(q_i^{a_{ij}}-1)(z+w)\ps_i(z)\ps_j(w)\phi_i^+(\FC^{-1}z^{-1})\phi_j^+(\FC^{-1}w^{-1}).
	\end{split}
	\end{equation}
Hence, the result follows from Proposition~\ref{prop:sigma-gamma}.


\subsection{Recursion formula for $y_i(z)$}
For every $i\in\bfI$, set $\oD_i=\ad(\Leta_{i,1})$, see \eqref{eq:def-rho}. 

\begin{proposition}\label{prop:y-rationality}
For every $i\in\bfI$ and $k\in\bbZ$,
\begin{equation}\label{eq:yik-recursion}
	\oD_i(y_{i,k})
	=[a_{i\tau(i)}]_{q_i}y_{i,k+1}-[2]_{q_i}y_{i,k-1}\FC
\end{equation}
In particular,
	\begin{equation}\label{eq:yi-recursion}
	(\oD_i-[a_{i\tau(i)}]_{q_i} z^{-1}
	+[2]_{q_i}\FC z)(y_i(z)^\pm)
	=
	-
	[a_{i\tau(i)}]_{q_i} z^{-1}y_{i,0}
	-
	[2]_{q_i}\FC\, y_{i,-1}
\end{equation}
where $y_i(z)^+=\sum_{k\geqslant 0}y_{i,k}z^k$ and 
$y_i(z)^-=-\sum_{k<0}y_{i,k}z^k$.
\end{proposition}

\begin{pf}
By \eqref{eq:def-rho}, $\Leta_{i,1}=(q_i-q_i^{-1})^{-1}\eta_{i,0}^{-1}\eta_{i,1}$. Hence, by Corollary~\ref{cor:explicit-rels}, one has
\begin{equation}
	[\Leta_{i,1},y_i(z)]
	=
	\left(
	[a_{i\tau(i)}]_{q_i} z^{-1}
	-
	[2]_{q_i}\FC z
	\right)y_i(z)
\end{equation}
and \eqref{eq:yik-recursion} follows. 
The formulae for $y_i(z)^{\pm}$ are obtained by summation over non-negative and negative powers, respectively.
\end{pf}

\begin{remark}
In $\bfU$, one has
\begin{equation}
	[h_{i,m},x^{\pm}_{j,k}]=\pm\frac{1}{m}[ma_{ij}]_{q_i}x^{\pm}_{j,k+m}	
\end{equation}
for every $i,j\in\bfI$, $m\in\bbZ_{\neq0}$, and $k\in\bbZ$. Hence,
\begin{align}
	\oD_i(x^-_{i,k})&=-[2]_{q_i}x^-_{i,k+1}+[a_{i\tau(i)}]_{q_i}\FC x^-_{i,k-1},\\
	\oD_i(x^+_{\tau(i),k})&=[a_{i\tau(i)}]_{q_i}x^+_{\tau(i),k+1}-[2]_{q_i}\FC x^+_{\tau(i),k-1}	
\end{align}	
for every $i\in\bfI$ and $k\in\bbZ$. Therefore,
	\begin{equation}
	\begin{aligned}
	\oD_i(y_{i,k})
	&=[a_{i\tau(i)}]_{q_i}y_{i,k+1}-[2]_{q_i}y_{i,k-1}\FC
	-[2]_{q_i}\sum_{a\ges0}\left(\ps_{i,k+1+a}-\ps_{i,k-1+a}\FC\right)\phi^+_{i,a}\FC^{-a}.
	\end{aligned}
\end{equation}	
The choice of $\SF_i(z)$ as in \eqref{eq:SF+GF}, see also Proposition~\ref{ss:SF+GF} (1), 
is therefore determined by the vanishing of the last term. As before, this is required by the minimality condition \eqref{eq:minimality}.
\hfill\rmkend
\end{remark}

\subsection{Auxiliary relation for $\eta_i(z)$}\label{ss:eta-aux}
We shall use the following observation.
\begin{proposition}\label{prop:eta-rationality}
For every $i\in\bfI$, one has
\begin{equation}\label{eq:etai-rationality}
	\small
	\begin{split}
	w(1-q_i^{a_{i\tau(i)}}\FC z^2)&
	\pg_{\tau(i)}(\FC^{-1}z^{-1})\eta_{\tau(i)}(z)
	\\
	&+z(1-q_i^{a_{i\tau(i)}}\FC w^2)
	\pg_i(\FC^{-1}w^{-1})\eta_i(w)
	=F_i(z,w),
	\end{split}
\end{equation}
where
\begin{equation}
\small
\begin{aligned}
	F_i(z,w)
	={}&
	-(q_i-q_i^{-1})\FC zw(1-\FC zw)
	\Big[
	(q_i^{a_{i\tau(i)}}z-w)y_i(z)^+y_{\tau(i)}(w)^+
	\\
	&\hspace{3cm}
	+
	(q_i^{a_{i\tau(i)}}w-z)y_{\tau(i)}(w)^+y_i(z)^+
	\\
	&\hspace{3.25cm}
	+
	q_i^{a_{i\tau(i)}}y_{i,-1}y_{\tau(i)}(w)^+
	+
	q_i^{a_{i\tau(i)}}y_{\tau(i),-1}y_i(z)^+
	\\
	&\hspace{3.5cm}
	-
	y_i(z)^+y_{\tau(i),-1}
	-
	y_{\tau(i)}(w)^+y_{i,-1}
	\Big]
	\\
	&\hspace{3.75cm}
	+
	(1-\FC zw)
	\left(
	w\pg_{\tau(i),0}\eta_{\tau(i),0}
	+
	z\pg_{i,0}\eta_{i,0}
	\right).
\end{aligned}
\end{equation}
In particular,
\begin{equation}
(1-q_i^{a_{i\tau(i)}}\FC w^2)
\pg_i(\FC^{-1}w^{-1})\eta_i(w)
=
\frac{\partial F_i}{\partial z}(0,w)
-w\bigl(\pg_{\tau(i),0}\eta_{\tau(i),1}
+\FC\pg_{\tau(i),-1}\eta_{\tau(i),0}\bigr).
\end{equation}
\end{proposition}

\begin{pf}
The result follows by multiplying \eqref{eq:quadratic-y-modes} in
Corollary~\ref{cor:explicit-rels} by $z^rw^s$ and summing over $r,s\ges0$.
\end{pf}

\subsection{Drinfeld--Serre relations}
\label{ss:DS-simply-laced}
We shall impose the following additional restrictions on the parameters defined in \ref{ss:weight-quadratic-relations}.
\begin{enumerate}[leftmargin=2em]\itemsep0.25cm
\item For every $i\in\bfI\setminus\bfI^\tau$ with $a_{i\tau(i)}=0$,
\begin{equation}\label{eq:serre-orthogonal-parameters}
\pga{i}=\pga{\tau(i)}\aand \pgb{i}=0,
\end{equation}
that is, $\GF_i(z)=\pga{i}=\pga{\tau(i)}=\GF_{\tau(i)}(z)$.
\item 
For every $i\in\bfI^\tau$ with $a_{ji}\not\in 2\bbZ$ for some $j\in\bfI^\tau$,
\begin{equation}\label{eq:serre-sigma-parity}
\psa{i}=0=\psb{i},
\end{equation}
that is, $\SF_i(z)=0$.
\end{enumerate}

\begin{remark}
In simply-laced type, condition \eqref{eq:serre-sigma-parity} implies that a fixed vertex
carrying a nonzero $\SF_i(z)$ has no other fixed vertex as a neighbour.\hfill\rmkend
\end{remark}

These restrictions on the parameters are reminiscent of those occurring in standard quantum symmetric pairs \cite{regelskis-vlaar-20, appel-vlaar-22} and, as before, they naturally appear by imposing the minimality condition \eqref{eq:minimality}. 
Namely, for 
$i,j\in\bfI$ with $i\neq j$ and $N=1-a_{ij}$, set
\begin{equation}\label{eq:serre-expression}
\bfS_{ij}(\ul z,w)=\operatorname{Sym}_{\ul z}
\sum_{p=0}^{N}(-1)^p\begin{bmatrix}N\\p\end{bmatrix}_{q_i}
y_i(z_1)\cdots y_i(z_p)y_j(w)y_i(z_{p+1})\cdots y_i(z_N),
\end{equation}
where symmetrization is the sum over all permutations, without
division by $N!$. In contrast with the Drinfeld--Serre relations of $\bfU$, 
the formal expression $\bfS_{ij}(\ul z,w)$ does not always vanish
in $\fext{\bfU}{\ul{z},w}$, but it can be expressed as a linear combination of \emph{lower degree} terms, which either are in $\fext{\bfB}{\ul{z},w}$ or have a non-trivial component along $\bfU^0$. Conditions \eqref{eq:serre-orthogonal-parameters}--\eqref{eq:serre-sigma-parity} guarantee that such terms vanish.

\begin{proposition}\label{prop:serre-relations}
Set $\deta_i(z)=\pg_i(\FC^{-1}z^{-1})\eta_i(z)$ and
\begin{equation}\label{eq:serre-realized-cartan}
\ddeta_i(z)=
\frac{q_i^{a_{i\tau(i)}}\FC z^2-1}
  {q_i^{a_{i\tau(i)}}-\FC z^2}\,
\deta_{\tau(i)}(z),
\end{equation}
where the scalar factor is expanded at $z=0$.
Under \eqref{eq:serre-orthogonal-parameters} and
\eqref{eq:serre-sigma-parity}, the following identities hold.
\begin{enumerate}
\item If $a_{ij}=0$ and $i\neq\tau(j)$, then $\bfS_{ij}(z,w)=0$.
\item If $a_{ij}=0$ and $j=\tau(i)$, then\footnote{Note that in this case $\ddeta_i(z)=-\deta_{\tau(i)}(z)$, hence $\ddeta_i(z)-\ddeta_j(z)=\deta_i(z)-\deta_j(z)$.}
\begin{equation}\label{eq:serre-orthogonal-correction}
\bfS_{ij}(z,w)=\frac{\fdel(\FC zw)}{q_i-q_i^{-1}}
\bigl(\deta_i(w)-\deta_j(z)\bigr).
\end{equation}
\item If $a_{ij}=-1$ and $i\neq\tau(i),\tau(j)$, then
$\bfS_{ij}(\ul{z},w)=0$.
\item If $a_{ij}=-1$ and $i=\tau(i)$, then
\begin{equation}\label{eq:serre-fixed-minus-one-full}
\begin{aligned}
\bfS_{ij}(\ul z,w)
={}&
\Sym_{\ul{z}}\frac{\fdel(\FC z_1z_2)}{q_i-q_i^{-1}}\Bigl(
y_j(w)\ddeta_i(z_2)\\
&\qquad+
\bigl(1-[2]_{q_i}\RF_{ij}(z_2,w)\bigr)
\ddeta_i(z_2)y_j(w)\Bigr).
\end{aligned}
\end{equation}
\item If $a_{ij}=-1$ and $j=\tau(i)$, then
\begin{equation}\label{eq:serre-exchanged-minus-one-full}
\begin{aligned}
\bfS_{ij}(\ul z,w)
={}&\Sym_{\ul{z}}
\frac{\fdel(\FC z_1w)}{q_i-q_i^{-1}}\Bigl(
\bigl(1+\RF_{ii}(z_1,z_2)\bigr)
\ddeta_i(z_1)y_i(z_2)\\
&\hspace{-1cm}+
\bigl(\RF_{ji}(w,z_2)-[2]_{q_i}\bigr)
\ddeta_j(w)y_i(z_2)
	+y_i(z_2)\bigl(\ddeta_j(w)-[2]_{q_i}\ddeta_i(z_1)\bigr)
\Bigr).
\end{aligned}
\end{equation}
\end{enumerate}
Here the rational functions on the right-hand sides are expanded at
$z_2/w=0$ in (4), and at $z_1/z_2=0$ and $w/z_2=0$, respectively,
in (5), before symmetrization.
\end{proposition}

\begin{remark}
A weaker version of the Drinfeld--Serre relations for $\bfB$ can be obtained by 	
multiplying (2) by $(z-w)$ and (3), (4), (5) by $(z_1-q_iw)(z_2-q_iw)$. 
One easily shows that these follow directly from the
quadratic relation \eqref{eq:quadratic-y-y} and require no restriction on the parameters.
\hfill\rmkend
\end{remark}

As in the case of the weight relations \eqref{eq:weight-rel-eta-y} and the quadratic 
relations \eqref{eq:quadratic-y-y}, the proof of Proposition~\ref{prop:serre-relations} is 
by a direct computation, which we present in a generalized form in Appendix 
\ref{app:shifted-coaction}.\\ 

However, it follows readily from the result above that, in simply-laced types, after a 
non--scalar renormalization and a relabeling of the $\imath$Cartan generators, the 
algebra $\bfB$ is a quotient of the affine quasi--split $\imath$quantum group of the same 
type, see, \eg \cite[Def. 2.1]{lu-wang-weekes-26}. Building on this observation, in the 
next section, we shall provide an alternative proof of the Drinfeld--Serre relations in 
$\bfB$ for arbitrary types. 


\section{Comodule realization of quantum loop symmetric pairs}
\label{s:comodule-realization}

We recall the definition of the quasi--split $\imath$quantum group $\bfU^\imath$ 
associated to $\tau$
and the topological $\bfU$-comodule structure  recently constructed in \cite{lu-pan-26} building on the theory of $\imath$Hopf algebras introduced in \cite{chen-pan-ruan-wang-25}.
We then show that  every choice of admissible parameters $\SF_i(z), \GF_i(z)$ ($i\in\bfI$) determines a character of $\bfU^\imath$, and yields a surjective morphism from $\bfU^\imath$ to $\bfB$.

\subsection{$\imath$Quantum loop algebras \cite{lu-wang-20, zhang-22,lu-pan-26}}\label{ss:i-quantum-group}
Let $\bfU^\imath$ be the algebra generated over $\bbC$ by $\FC^{\pm1}$, $\{q^h\}_{h\in\h^{\theta}}$, $\{\wetad_{i,k}\}_{i\in\bfI,k\geqslant0}$, $\{\wetad_{i,0}^{-1}\}_{i\in\bfI}$, and $\{\wy_{i,r}\}_{i\in\bfI,r\in\bbZ}$, subject to the following relations. 
Set 
\begin{equation}
	\wetad_{i}(z)=\sum_{k\ges0}\wetad_{i,k}z^k \aand \wy_i(z)=\sum_{r\in\bbZ}\wy_{i,r}z^r.
\end{equation}
Then, $\FC$ is central with inverse $\FC^{-1}$, and, for every $i,j\in\bfI$ and $h\in\h^{\theta}$,
\begin{gather}
\label{eq:cartan-abs}
q^0=1,
\qquad
q^hq^{h'}=q^{h+h'},\\
\label{eq:eta-h-abs}
\wetad_{i,0}\wetad_{i,0}^{-1}=\wetad_{i,0}^{-1}\wetad_{i,0}=1, \qquad
q^h\wetad_{i,k}=\wetad_{i,k}q^h,
\quad
[\wetad_{i,k},\wetad_{j,l}]=0,
\end{gather}
\begin{align}
\label{eq:weight-abs}
q^h \wy_{i}(z)q^{-h}&=q^{-\alpha_i(h)}\wy_{i}(z),\\
\label{eq:weight-rel-eta-y-abs}
\wetad_i(z)\wy_j(w)&=\RF_{ij}(\FC zw)\RF_{\tau(i)j}(w/z)\wy_j(w)\wetad_i(z)\\[1em]
\label{eq:quadratic-y-y-abs}
\begin{split}
(q_i^{a_{ij}}z-w)\wy_i(z)\wy_{j}(w)&+(q_i^{a_{ij}}w-z)\wy_{j}(w)\wy_i(z)=\\
&\hspace{-1.75cm}=\delta_{i\tau(j)}\frac{\fdel(\FC zw)}{q_i-q_i^{-1}}\left((q_i^{a_{i\tau(i)}}z-w)\wetad_{\tau(i)}(z)
+(q_i^{a_{i\tau(i)}}w-z)\wetad_{i}(w)\right).
\end{split}
\end{align}
In \eqref{eq:weight-rel-eta-y-abs}, the product of rational functions is expanded at $z=0$.
Moreover, for every $i\neq j$,
\begin{equation}\label{eq:drinfeld-serre-relations-abs}
\operatorname{Sym}_{\ul{z}}
\sum_{p=0}^N(-1)^p
\begin{bmatrix}N\\p\end{bmatrix}_{q_i}\ww(\ul{v}_p)=0,
\end{equation}
where $N=1-a_{ij}$, $\ul{v}_p=((i,z_1),\ldots,(i,z_p),(j,w),
(i,z_{p+1}),\ldots,(i,z_N))$, and $\ww(\ul{v}_p)$ is recursively defined as follows.

For an $n$-tuple $\ul{v}=((i_1,z_1),\ldots,(i_n,z_n))$, we denote by
$\ul{v}^k$ the tuple obtained by removing its $k$th entry. Multiple
superscripts indicate the simultaneous deletion of the corresponding entries, preserving the order of the other ones. As in \ref{ss:DS-simply-laced}, we introduce the auxiliary Cartan currents
\begin{equation}\label{eq:abstract-serre-cartan}
\wetadd_i(z)=
\frac{q_i^{a_{i\tau(i)}}\FC z^2-1}
  {q_i^{a_{i\tau(i)}}-\FC z^2}\,
\wetad_{\tau(i)}(z),
\end{equation}
with the scalar factor expanded at $z=0$.
Set $\ww(\varnothing)=1$ and define recursively
\begin{equation}\label{eq:abstract-serre-W-recursion}
\begin{aligned}
\ww(\ul{v})={}&\wy_{i_1}(z_1)\ww(\ul{v}^1)\\
&-\sum_{k=2}^n
\delta_{i_1,\tau(i_k)}\frac{\fdel(\FC z_1z_k)}
 {q_{i_1}-q_{i_1}^{-1}}
\left(\prod_{s=2}^{k-1}\RF_{i_k i_s}(z_k/z_s)\right)
\wetadd_{i_k}(z_k)\ww(\ul{v}^{1,k}),
\end{aligned}
\end{equation}
where $\RF_{i_k i_s}(z_k/z_s)$ is expanded at $z_k/z_s=0$.
In particular,
\[
\ww((i,z),(j,w))=
\wy_i(z)\wy_j(w)-
\delta_{i,\tau(j)}\frac{\fdel(\FC zw)}{q_i-q_i^{-1}}
\wetadd_j(w).
\]
\begin{remarks}
\hfill
\begin{enumerate}[leftmargin=2em]\itemsep0.25cm
\item	
If $i\neq\tau(i),\tau(j)$, $\ww(\ul{v}_p)=\wy_i(z_1)\cdots \wy_i(z_p)\wy_j(w)
\wy_i(z_{p+1})\cdots \wy_i(z_N)$ and \eqref{eq:drinfeld-serre-relations-abs} is the usual Drinfeld--Serre relation.
\item 
If $a_{ij}=0$ or $a_{ij}=-1$, expanding the recursion gives precisely the
formulae in Proposition~\ref{prop:serre-relations}, with
$y_i,\deta_i,\ddeta_i$ replaced by
$\wy_i,\wetad_i,\wetadd_i$, respectively.
\hfill\rmkend
\end{enumerate}
\end{remarks}

We denote by $\bfU^{\imath,0}\subset\bfU^\imath$ the {\em loop $\imath$Cartan subalgebra}, generated by $\FC^{\pm1}$, the elements $\{q^h\}_{h\in\h^{\theta}}$, 
$\{\wetad_{i,k}\}_{i\in\bfI, k\ges 0}$, and $\{\wetad_{i,0}^{-1}\}_{i\in\bfI}$.
Set
\begin{equation}\label{eq:def-rho-bis}
\wetad_i(z)=\wetad_{i,0}\exp\left((q_i-q_i^{-1})\sum_{k>0}\wLetad_{i,k}z^k\right)
\end{equation}
and $\wLetad_i(z)=\sum_{k>0}\wLetad_{i,k}z^k\in\fext{\bfU^{\imath,0}}{z}$.
We shall refer to  \eqref{eq:drinfeld-serre-relations-abs} as the {\em $\imath$Drinfeld--Serre relations}.

\subsection{Comparison with the original presentation}
\label{ss:Lu-Pan-presentation}
The presentation of $\bfU^\imath$ introduced above differs slightly
from those in \cite{lu-wang-weekes-26, lu-pan-26}.
We enlarge the loop $\imath$Cartan subalgebra by adjoining the elements
$\{q^h\}_{h\in\h^{\theta}}$, and relabel and renormalize the other generators.
Specifically, in the notation of
\cite[Def.~2.11, Sec.~4.3]{lu-pan-26}, the change of generators is
given by $C\longmapsto\FC$, $B_i(z)\longmapsto\wy_i(z)$, 
\begin{equation}\label{eq:Lu-Pan-dictionary}
\bbK_i\longmapsto-q_i^{-\frac{a_{i\tau(i)}}{2}}\wetad_{i,0},
\aand
\Theta_i(z)\longmapsto
\frac{1-q_i^{a_{i\tau(i)}}\FC z^2}
{1-q_i^{-a_{i\tau(i)}}\FC z^2}\,
\wetad_{\tau(i),0}^{-1}\wetad_{\tau(i)}(z).
\end{equation}
The scalar rational factor in \eqref{eq:Lu-Pan-dictionary}
is expanded at $z=0$ and has constant term $1$.
Since the Cartan zero modes are invertible, this change of
generators is invertible.\\

The compatibility with the relations is straightforward.
The Cartan relations \eqref{eq:cartan-abs}, \eqref{eq:eta-h-abs}, \eqref{eq:weight-abs} are immediate, while the weight and quadratic relations
\eqref{eq:weight-rel-eta-y-abs}--\eqref{eq:quadratic-y-y-abs} are readily 
obtained by substitution in \cite[Eqs.~(2.44)--(2.46)]{lu-pan-26}.
For the quadratic relations, observe that
\begin{equation}\label{eq:Lu-Pan-Cartan-combination}
q_i^{-\frac{a_{i\tau(i)}}{2}}\bbK_{\tau(i)}\Theta_i(z)
\longmapsto
\frac{q_i^{a_{i\tau(i)}}\FC z^2-1}
{q_i^{a_{i\tau(i)}}-\FC z^2}\,
\wetad_{\tau(i)}(z),
\end{equation}
and, on the support of $\fdel(\FC zw)$, the scalar factor in
\eqref{eq:Lu-Pan-Cartan-combination} becomes
\[
\frac{q_i^{a_{i\tau(i)}}z-w}
{q_i^{a_{i\tau(i)}}w-z}.
\]
Finally, the recursive form of the Drinfeld--Serre relations
\eqref{eq:drinfeld-serre-relations-abs} follows from \cite[Thm.~4.7]{lu-pan-26}, which identifies the $\imath$quantum loop group with the centrally extended loop Borel subalgebra $\bfU^{\ges0}$ endowed with the $\tau$-twisted product defined in \cite[Eq.~(4.36)]{lu-pan-26}. This approach is based on the theory of $\imath$Hopf algebras developed in \cite{chen-pan-ruan-wang-25}. Under this identification, the relations in $\bfU^\imath$ are obtained by expressing the ordinary Drinfeld--Serre relations of $\bfU^{\ges0}$ in terms of the twisted product, and applying the change of generators above.

Namely, set
$\wx(\ul{v})=x_{i_1}^+(z_1)\cdots x_{i_n}^+(z_n)$ for every $\ul{v}=((i_1,z_1),\ldots,(i_n,z_n))$. The Drinfeld--Serre relations in $\bfU^{\ges0}$ read
\begin{align}
\operatorname{Sym}_{\ul{z}}
\sum_{p=0}^{N}
(-1)^p
\begin{bmatrix}N\\ p\end{bmatrix}_{q_i}
\wx(\ul{v}_p)
=0,
\end{align}
where $N=1-a_{ij}$ and $\ul{v}_p=((i,z_1),\ldots,(i,z_p),(j,w),
(i,z_{p+1}),\ldots,(i,z_N))$. By definition of the twisted product, one has
\begin{equation}
\begin{aligned}
	x_{i_1}^+(z_1)*\wx(\ul{v}^1)
	=\wx(\ul{v})
	+\sum_{k=2}^n \delta_{i_1,\tau(i_k)}
	\frac{\fdel(Cz_1z_k)}
	{q_{i_1}-q_{i_1}^{-1}}
	\left(\prod_{s=2}^{k-1}\RF_{i_k i_s}(z_k/z_s)\right)
	\phi^+_{i_k}(z_k)\wx(\ul{v}^{1,k}),
\end{aligned}
\end{equation}
where $\RF_{i_k i_s}(z_k/z_s)$ is expanded at $z_k/z_s=0$.
Under the change of generators \eqref{eq:Lu-Pan-dictionary}, 
the isomorphism with $\bfU^\imath$ maps $x_i^+(z)$, $\phi^+_i(z)$ respectively to 
$\wy_i(z)$, $\wetadd_i(z)$.
Hence, the recursion \eqref{eq:abstract-serre-W-recursion} follows, since $\phi^+_{i_k}(z_k)\wx(\ul{v}^{1,k})=\phi^+_{i_k}(z_k)*\wx(\ul{v}^{1,k})$.

\subsection{Comodule structure \cite{lu-pan-26}}\label{ss:comodule-structure}
In \cite[Prop.~4.15]{lu-pan-26}, the algebra $\bfU^\imath$ is endowed with a 
comodule structure over $\bfU$. The existence of the coaction follows as before 
by adapting the $\imath$Hopf algebra approach from \cite{chen-pan-ruan-wang-25}
to the case of quantum loop algebras with their Drinfeld coproduct.
\begin{theorem}\cite{lu-pan-26}\label{thm:Lu-Pan-coaction}
There is a morphism of algebras
\[
\DDi{s}\colon\bfU^\imath\longrightarrow
\Lext{(\bfU^\imath\ten\bfU)[\FC^{\pm1}]}{s^{-1};s}
\]
given by $\DDi{s}(\FC)=\FC\ten1$,
$\DDi{s}(q^h)=q^h\ten q^h$ ($h\in\h^{\theta}$), and
\begin{align}
	\label{eq:DDi-weta-field}
	\DDi{s}(\wetad_i(z))
	&=\wetad_i(z)\ten\phi_i^+(zs)
	\phi_{\tau(i)}^-(\FC^{-1}z^{-1}s),\\
	\label{eq:DDi-wy-field}
	\DDi{s}(\wy_i(z))
	&=1\ten x_i^-(\FC^{-1}z^{-1}s)
	+\wy_i(z)\ten\phi_i^+(\FC^{-1}z^{-1}s)\notag\\
	&\quad+\wetad_i(\FC^{-1}z^{-1})\ten
	x_{\tau(i)}^+(zs)\phi_i^+(\FC^{-1}z^{-1}s).
\end{align}
Moreover, $(\id\ten\veps)\circ\DDi{s}=\id$ and
\[
(\DDi{s_1}\ten\id)\circ\DDi{s_1s_2}
=(\id\ten\DD{s_2})\circ\DDi{s_1}.
\]
\end{theorem}
The coaction $\Delta^\imath$ is not the one introduced in \cite{lu-pan-26}. However, as we explain below, it can be recovered from the latter by twisting with an anti-isomorphism on $\bfU^\imath$ and taking into account the different conventions on the Drinfeld coproduct on $\bfU$. This choice makes the identification with the coideal structure on $\bfB$ immediate, see \ref{ss:coideal}.\\

We denote by $\bfU_{\operatorname{LP}}$ the quantum loop algebra
and the Drinfeld coproduct $\Delta_{\operatorname{LP}}$ considered in \cite{lu-pan-26}, specialized by setting the two affine central elements $C,C'$ equal to $1$, $K_iK'_i=1$, and $v=q$, and extended with the elements $\{q^h\}_{h\in\h}$.
There is an algebra anti-isomorphism $\bff$ from $\bfU_{\operatorname{LP}}$ to $\bfU$ given by $\bff(q^h)=q^{-h}$,
\begin{equation}\label{eq:LP-coproduct-dictionary}
\bff(X_i^\pm(u))=\mp x_i^\pm(u^{-1})
\aand
\bff(K_i^{\pm1}\psi_i^\pm(u))=\phi_i^\mp(u^{-1}).
\end{equation}
In particular, $\bff(K_i^{\pm1})=\phi^\mp_{i,0}$. 
This accounts for the weight relations and
the denominator $v_i^{-1}-v_i$ in \cite[Eq.~(2.14)]{lu-pan-26}.
Then, by \cite[Eqs.~(4.4)--(4.7)]{lu-pan-26}, it follows that
\begin{equation}\label{eq:LP-coproduct-intertwiner}
\DD{s}\circ\bff
=(\id\ten\SH_s)\circ(\bff\ten\bff)\circ\Delta_{\mathrm{LP}},
\end{equation}
where $\SH_s$ is the formal shift operator, hence $\SH_s(x_i^\pm(u))=x_i^\pm(su)$ and
$\SH_s(\phi_i^\pm(u))=\phi_i^\pm(su)$, see \ref{ss:psi-invariance}.\\

We denote by $\bfU^\imath_{\operatorname{LP}}$ the $\imath$quantum  loop algebra in the presentation of \cite{lu-pan-26}, extended with the elements $\{q^h\}_{h\in\h^{\theta}}$, and by $\Delta^\imath_{\operatorname{LP}}$ the corresponding comodule structure. Let $\bfg'\colon\bfU^\imath_{\operatorname{LP}}\to\bfU^\imath$ be the isomorphism from \ref{ss:Lu-Pan-presentation} and $\bfg''$ the anti-involution\footnote{
The weight and quadratic relations are easily checked. For the Serre relations it suffices to check the relations at mode zero, by \cite[Sec.~4.3]{lu-pan-26}, which follow since, in this case, $\bfg''$ reduces to the standard anti-involution of the underlying quantum symmetric pair in $U_q(\g)$, see \cite[Lem.~2.6(1)]{lu-pan-26}.
} of $\bfU^\imath$ given by
$ \bfg''(\FC)=\FC$, $\bfg''(q^h)=q^{-h}$, 
\begin{equation}\label{eq:iota-spectral-reflection}
\bfg''(\wy_i(z))=\wy_i(\FC^{-1}z^{-1})
\aand
\bfg''(\wetad_i(z))=\wetad_{\tau(i)}(z).
\end{equation}
Set $\bfg=\bfg''\circ\bfg'$. Then, $\bfg(C)=\FC$, 
$\bfg(\bbK_i)=-q_i^{-a_{i\tau(i)}/2}\wetad_{\tau(i),0}$,
\begin{equation}\label{eq:LP-coaction-dictionary}
\bfg(B_i(u))=\wy_i(\FC^{-1}u^{-1}),
\aand
\bfg(\Theta_i(u))=
\frac{1-q_i^{a_{i\tau(i)}}\FC u^2}{1-q_i^{-a_{i\tau(i)}}\FC u^2}
\wetad_{i,0}^{-1}\wetad_i(u),
\end{equation}
where the rational function is expanded at $u=0$. Then, one readily checks that
\begin{equation}\label{eq:LP-coaction-transport}
\DDi{s}=(\id\ten\SH_s)\circ (\bfg\ten\bff)\circ
\Delta^\imath_{\mathrm{LP}}\circ\bfg^{-1},
\end{equation}
hence it is a morphism of algebras, satisfying the counit identity and the shifted coassociativity identity. Finally, note that, proceeding as in Proposition~\ref{ss:coideal},
one checks that the image of $\Delta^\imath$ lies in $\Lext{(\bfU^\imath\ten\bfU)[\FC^{\pm1}]}{s^{-1};s}$.


\begin{remark}
We give a direct proof of Theorem~\ref{thm:Lu-Pan-coaction} in
Appendix~\ref{app:shifted-coaction}, in the more general setting of
shifted quantum loop algebras.\hfill\rmkend
\end{remark}

\subsection{Characters and parameters}
\label{ss:characters-admissible-parameters}
The definition of the coideal subalgebra $\bfB\subset\Lext{\bfU}{\FC;\FC^{-1}}$
in \ref{ss:B-def} relies on a choice of formal series $\pg_{i}(z)\in\fext{\bbC[\FC^{\pm1}]}{z^{-1}}$, 
$\ps_i(z)\in\fext{\bbC[\FC^{\pm1}]}{z^{\pm1}}$ ($i\in\bfI$).
By imposing the minimality condition \eqref{eq:minimality}, we constrained
these parameters: they are encoded by the rational functions $\SF_i(z), \GF_i(z)$ ($i\in\bfI$) described in \eqref{eq:SF+GF} and subject to the additional constraints \eqref{eq:serre-orthogonal-parameters}--\eqref{eq:serre-sigma-parity}.\\

We now observe that the same rational functions classify the characters of $\bfU^\imath$.
In this subsection and the next, we temporarily allow the parameters $\psa{i}, \psb{i}$, $\pga{i}, \pgb{i}$ determining $\SF_i(z), \GF_i(z)$ to vary over $\bbC[\FC^{\pm1}]$, see Remark~\ref{ss:weight-quadratic-relations}. More precisely, we regard the
parameters $\pga{i}$ as formal invertible elements, with
the identifications required by \eqref{eq:serre-orthogonal-parameters},
and set\footnote{We shall prove below that the construction of non-trivial characters requires imposing in $\bbC[\pga{i}^{\pm1}\mid i\in\bfI]$ the relation $\pga{i}=\pga{\tau(i)}$ whenever $i\neq\tau(i)$ and $a_{i\tau(i)}=0$.
}
\begin{equation}\label{eq:parameter-ring}
R_0=\bbC[\pga{i}^{\pm1}\mid i\in\bfI]
\big/\bigl(\pga{i}-\pga{\tau(i)}\mid
i\neq\tau(i),\ a_{i\tau(i)}=0\bigr),
\quad
R=R_0[\FC^{\pm1}],
\end{equation}
and  $\bfU^\imath_R=R\otimes_{\bbC[\FC^{\pm1}]}\bfU^\imath$.
The remaining parameters $\pgb{i},\psa{i},\psb{i}$ are allowed to be arbitrary elements
of $R$, subject to the conditions specified below.\footnote{This convention
replaces the condition $\psa{i}, \psb{i}, \pga{i}, \pgb{i}\in\bbC$ of
Remark~\ref{ss:weight-quadratic-relations} only in these two subsections.}\\

By a character of $\bfU^\imath_R$ we mean an $R$-algebra homomorphism
$\veps\colon\bfU^\imath_R\to R$, or equivalently a representation
on the free rank-one $R$-module.

\begin{proposition}\label{prop:one-dimensional}
\hfill
\begin{enumerate}[leftmargin=2em]\itemsep0.25cm
\item 
Let $\GF_i(z),\SF_i(z)$ ($i\in\bfI$) be as in \eqref{eq:SF+GF}, and set
$\pg_i(z)=\GF_i(z)^-$ and $\ps_i(z)=\SF_i(z)^+-\SF_i(z)^-$.
Then, the assignments
$\veps_{\ps,\pg}(q^h)=1$ ($h\in\h^{\theta}$),
\begin{equation}\label{eq:admissible-character}
\veps_{\ps,\pg}(\wetad_i(z))=\pg_i(\FC^{-1}z^{-1}),
\aand
\veps_{\ps,\pg}(\wy_i(z))=\ps_i(z)
\end{equation}
define a character
$\veps_{\ps,\pg}\colon\bfU^\imath_R\to R$ if and only if
\eqref{eq:serre-orthogonal-parameters} and
\eqref{eq:serre-sigma-parity} hold.
\item 
Every character $\veps\colon\bfU^\imath_R\to R$ such that
\begin{equation}\label{eq:character-normalization}
\veps(q^h)=1\quad (h\in\h^{\theta}) \aand \veps(\wetad_{i,0})=\pga{i}
\end{equation}
has the form \eqref{eq:admissible-character} for unique $\SF_i(z), \GF_i(z)$ ($i\in\bfI$).
\end{enumerate}
\end{proposition}

\begin{pf}
(1) The relations \eqref{eq:cartan-abs}--\eqref{eq:eta-h-abs} are immediate.
By \eqref{eq:SF+GF}, $\ps_i(z)=\delta_{i\tau(i)}(\psa{i}+\psb{i}\FC z)\fdel(\FC z^2)$. 
In particular, when $i\neq\tau(i)$, $\ps_i(z)=0$ and \eqref{eq:weight-abs} is trivial.
When $i=\tau(i)$, then $\alpha_i(h)=0$ for $h\in\h^{\theta}$ and \eqref{eq:weight-abs}
reduces to $\ps_i(z)=\ps_i(z)$.
Similarly, \eqref{eq:weight-rel-eta-y-abs} is trivial when $j\neq\tau(j)$, while, 
when $j=\tau(j)$, it is enough to observe that $(1-\RF_{ij}(\FC zw)\RF_{ij}(w/z))\fdel(\FC w^2)=0$, see \eqref{eq:sigma-delta-2}. 
The quadratic relations \eqref{eq:quadratic-y-y-abs} reduce to the defining system 
of series $\pg_i(z)$ ($i\in\bfI$), see \eqref{cond:gamma}, hence they automatically hold.\\

It remains to prove the $\imath$Drinfeld--Serre relations
\eqref{eq:drinfeld-serre-relations-abs}. Through $\veps$, the recursion 
\eqref{eq:abstract-serre-W-recursion} reduces to 
\begin{equation}\label{eq:serre-scalar-recursion}
\begin{aligned}
	\lambda(\ul{v})={}&\ps_{i_1}(z_1)\lambda(\ul{v}^1)
	-\sum_{k=2}^n
	\RF_{1k}(\ul{v})\ddot{\pg}_{i_k}(z_k)
	\lambda(\ul{v}^{1,k}),
\end{aligned}
\end{equation}
where $\lambda(\ul{v})=\veps(\ww(\ul{v}))$,
\begin{align}
\label{eq:serre-contraction-kernel}
\RF_{1k}(\ul{v})&=
\delta_{i_1,\tau(i_k)}\frac{\fdel(\FC z_1z_k)}
{q_{i_1}-q_{i_1}^{-1}}
\prod_{s=2}^{k-1}
\RF_{i_k i_s}(z_k/z_s)^+,\\
\label{eq:serre-scalar-cartan}
\ddot{\pg}_i(z)&=\veps(\wetadd_i(z))=
\frac{q_i^{a_{i\tau(i)}}\FC z^2-1}
{q_i^{a_{i\tau(i)}}-\FC z^2}\,
\pg_{\tau(i)}(\FC^{-1}z^{-1})
\stackrel{\eqref{eq:gamma-symmetry}}{=}-\GF_i(z)^+.
\end{align}
In particular, $\ddot{\pg}_i(0)=-q_i^{-a_{i\tau(i)}}\pga{\tau(i)}$. Hence, \eqref{eq:drinfeld-serre-relations-abs} becomes
\begin{equation}\label{eq:scalar-serre}
	\S_{ij}(\ul{z}, w)=\operatorname{Sym}_{\ul{z}}
	\sum_{p=0}^{N}
	(-1)^p
	\begin{bmatrix}N\\ p\end{bmatrix}_{q_i}
	\lambda(\ul{v}_p)
	=0.
\end{equation}
For $i\neq\tau(i)$ and $j\neq\tau(i)$, one has 
$\lambda(\ul{v}_p)=\ps_j(w)\prod_{k=1}^N\ps_{i}(z_k)=0$, and \eqref{eq:drinfeld-serre-relations-abs} is trivial.
For $i\neq\tau(i)$ and $j=\tau(i)$ with $a_{ij}=0$, one has
\begin{equation}
\S_{i\tau(i)}(z, w)=\frac{\fdel(\FC zw)}{q_i-q_i^{-1}}(\ddot{\pg}_{\tau(i)}(z)-\ddot{\pg}_i(w)),
\end{equation}
and \eqref{eq:drinfeld-serre-relations-abs} is readily seen to be equivalent to \eqref{eq:serre-orthogonal-parameters}. 
For $i\neq\tau(i)$ and $j=\tau(i)$ with $a_{ij}<0$, every term in the recursive expansion of $\lambda(\ul{v}_p)$ contains at least $N-1\geqslant1$ factors of $\ps_i$, since at most one $i$-entry can be paired with the unique $j$-entry. Thus $\S_{i\tau(i)}(\ul{z}, w)$ vanishes automatically,
see, \eg \eqref{eq:serre-exchanged-minus-one-full}.

It remains to consider the case $i=\tau(i)$. In this case, every term in the recursive expansion of $\lambda(\ul{v}_p)$ contains the factor $\ps_j(w)$. Hence, if $\ps_j(w)=0$, there is nothing to prove. Suppose $\ps_j(w)\neq0$, hence $j=\tau(j)$. If $a_{ij}=0$,
$\S_{ij}(z, w)=[\ps_i(z),\ps_j(w)]=0$. If $a_{ij}<0$, the condition \eqref{eq:SF+GF}
forces $\ps_i(z)=0$. Moreover, each $\lambda(\ul{v}_p)$ has the form
\begin{equation}\label{eq:vanishing-factor}
\lambda(\ul{v}_p)=\ps_j(w)\sum_{M\in\M}(-1)^{|M|}\prod_{k\not\in M}\ps_i(z_k)
\cdot L^p_{M}(\ul{z}, w),
\end{equation}
where $\M$ is the set of matchings of $\{1,\dots, N\}$. 
Hence, if $a_{ij}$ is even, then $N$ is odd, every matching leaves at least one index out,
and \eqref{eq:vanishing-factor} vanishes automatically.
On the other hand, if $a_{ij}$ is odd, by computing explicitly the coefficient of $w^s$ in $\S_{ij}(\ul{z},w)$, one observes that $\S_{ij}(\ul{z},w)=0$ implies $\ps_j(z)=0$, see, \eg \eqref{eq:serre-fixed-minus-one-full}. Hence, \eqref{eq:drinfeld-serre-relations-abs} is equivalent to \eqref{eq:serre-sigma-parity}.\\

(2) We proceed as in (1). Let $\veps$ be a character with the prescribed values, and set
\[
\pg_i(z)=\veps\bigl(\wetad_i(\FC^{-1}z^{-1})\bigr),
\qquad \ps_i(z)=\veps(\wy_i(z)).
\]
The weight relation \eqref{eq:weight-rel-eta-y-abs} with $i=j$ forces
$\ps_j(z)=0$ when $j\neq\tau(j)$. On the other hand, when $j=\tau(j)$, it yields
$(1-\FC z^2)\ps_j(z)=0$. Hence
\[
\ps_j(z)=(\psa{j}+\psb{j}\FC z)\fdel(\FC z^2).
\]
The quadratic relation \eqref{eq:quadratic-y-y-abs} reduces to \eqref{cond:gamma}, so
$\ps_i(z)$, $\pg_i(z)$ are forced to be of the form \eqref{eq:SF+GF}. Finally, by (1), the $\imath$Drinfeld--Serre relations are equivalent to \eqref{eq:serre-orthogonal-parameters} and \eqref{eq:serre-sigma-parity}.
\end{pf}

\begin{remark}
The result above extends to arbitrary characters corresponding to non-trivial
group homomorphisms $\h^{\theta}\to\bbC$.\hfill\rmkend
\end{remark}

\subsection{Central quotients and the realization in $\bfB$}\label{ss:central-specializations}
For every $i\in\bfI$, set
\begin{equation}\label{eq:central-zero-modes}
\bfz_i=\wetad_{i,0}q^{-d_ih_i+d_{\tau(i)}h_{\tau(i)}}.
\end{equation}
One readily shows that the elements $\bfz_i$ are central and invertible, and
satisfy
\begin{equation}\label{eq:central-coinvariance}
	\DDi{s}(\bfz_i)=\bfz_i\ten1.
\end{equation}	
Hence, the coaction $\Delta^\imath$ descends to the quotient
\begin{equation}\label{eq:central-quotient}
\bfU^\imath_{R,\mathbf c}
=\bfU^\imath_R/(\bfz_i-\pga{i})_{i\in\bfI},
\end{equation}
where we are imposing the relation $\wetad_{i,0}=\pga{i}q^{d_ih_i-d_{\tau(i)}h_{\tau(i)}}$.
Moreover, the characters in Proposition~\ref{prop:one-dimensional} also descend
to $\bfU^\imath_{R,\mathbf c}$, since, by \eqref{eq:admissible-character},
$\veps_{\ps,\pg}(\bfz_i)=\pga{i}=\GF_i(\infty)$.\\

For a pair $(\pg,\ps)$ satisfying Proposition~\ref{prop:one-dimensional},
define $\bfB_R$ to be the $R$-subalgebra of
$\Lext{(\bfU\ten_{\bbC} R_0)}{\FC;\FC^{-1}}$ generated by the elements
in \eqref{eq:B-generators} and $q^h$ for $h\in\h^{\theta}$.

\begin{proposition}\label{prop:central-realization}
Let $\veps_{\ps,\pg}$ be the character \eqref{eq:admissible-character}.
Then, there is a surjective morphism of comodule algebras $\xi_{\ps,\pg}\colon\bfU^\imath_R\to\bfB_R$
given by
\begin{equation}\label{eq:xi}
\xi_{\ps,\pg}(q^h)=q^h,\qquad
\xi_{\ps,\pg}(\wetad_i(z))=\pg_i(\FC^{-1}z^{-1})\eta_i(z),
\qquad
\xi_{\ps,\pg}(\wy_i(z))=y_i(z).
\end{equation}
\end{proposition}

\begin{pf}
It is enough to observe that
\begin{equation}\label{eq:morphism-xi}
	\xi_{\ps,\pg}=(\veps_{\ps,\pg}\ten\id)\circ\DDi{s}\vert_{s=1}\colon \bfU^\imath_R\to\Lext{(\bfU\ten_{\bbC} R_0)}{\FC;\FC^{-1}}.
\end{equation} 
Indeed, the formulae \eqref{eq:xi} follow readily by applying $\veps_{\ps,\pg}$ to the first tensor
factor of \eqref{eq:DDi-weta-field}--\eqref{eq:DDi-wy-field}, and then setting
$s=1$.
Since the constant coefficient of $\pg_i(\FC^{-1}z^{-1})$ is the invertible element $\pga{i}$, its inverse power series recovers every coefficient of
$\eta_i(z)$ from the image of $\wetad_i(z)$. Thus $\xi_{\ps,\pg}$ is surjective.
Moreover, from the coassociativity of the coaction, one has
$\DD{s}\circ\xi_{\ps,\pg}=(\xi_{\ps,\pg}\ten\id)\circ\DDi{s}$, that is, $\xi_{\ps,\pg}$
is a morphism of comodule algebras. 
\end{pf}

\begin{corollary}
The currents $y_i(z)$ and
$\deta_i(z)=\pg_i(\FC^{-1}z^{-1})\eta_i(z)$ in $\bfB_R$
satisfy the $\imath$Drinfeld--Serre relations \eqref{eq:drinfeld-serre-relations-abs}, with $\wy_i,\wetad_i, \wetadd_i$ replaced by $y_i,\deta_i, \ddeta_i$.
\end{corollary}

\begin{remarks}
By \eqref{eq:central-coinvariance}, one has
$\xi_{\ps,\pg}(\bfz_i)=\veps_{\ps,\pg}(\bfz_i)=\pga{i}$, so the
morphism factors through \eqref{eq:central-quotient}, inducing a surjective
morphism of comodule algebras 
$\xi_{\ps,\pg}\colon\bfU^\imath_{R,\mathbf c}\to\bfB_R$.	
Moreover, every evaluation morphism $\ev\colon R\to\bbC[\FC^{\pm1}]$ over $\bbC[\FC^{\pm1}]$ induces a surjective map 
\begin{equation}\label{eq:embedding}
\xi_{\ps,\pg}\colon\bbC[\FC^{\pm1}]\otimes_{R,\ev}\bfU^{\imath}_{R,\bfc}
\longrightarrow\bfB.
\end{equation}
Hence, Corollary~\ref{ss:central-specializations} extends to $\bfB\subset\Lext{\bfU}{\FC;\FC^{-1}}$.\hfill\rmkend
\end{remarks}


\section{Representations and rationality}\label{s:representations-rationality}
In this section, we prove that the morphism of comodule algebras $\xi_{\ps,\pg}$
constructed in Section~\ref{ss:central-specializations} yields a restriction functor on a one-parameter family of dense subcategories in $\sO$.
We then show that, for modules with rational weights and generic values of the spectral parameter, every $\bfU^\imath$-submodule of a restricted module is stable under the derived quantum loop algebra $\bfU'$. In particular, irreducibility is generically preserved under restriction for such modules. Finally, we describe the corresponding module structure on a category of $\bfU^\imath$-modules over $\bbC$.

\subsection{The category $\sM_\zeta$}\label{ss:category-M}
Fix $\zeta\in\bbC^\times$. Let $\sM_\zeta$ be the full subcategory of $\bfU^\imath$-modules over $\bbC$ such that
\begin{enumerate}[label=(M\arabic*)]\itemsep0.25cm
	\item\label{cond:M1} $M$ is at most countably dimensional;
	\item\label{cond:M2} $\FC$ acts on $M$ as the scalar $\zeta\in\bbC^\times$;
	\item\label{cond:M3} the action of $\{q^h\}_{h\in\h^{\theta}}$ on $M$ is semisimple, \ie
	\[M=\bigoplus_{\lambda\in(\h^{\theta})^*}M_\lambda\]
	where $M_\lambda=\{m\in M\mid \forall h\in\h^{\theta}, q^hm=q^{\lambda(h)}m\}$;
	\item\label{cond:M4} the action of $\bfU^{\imath,0}$ on $M$ is locally finite;
	\item\label{cond:M5} for every $m\in M$, $\dim(\sum_{i\in\bfI,k\in\bbZ}{\bfU^{\imath,0}_\zeta}\wy_{i,k}m)<\infty$.
\end{enumerate}
By (M4), $M$ decomposes into the sum of generalized eigenspaces with respect to $\bfU^{\imath,0}_\zeta$. Set $P(M)=\{\lambda\in(\h^\theta)^*\mid M_\lambda\neq0\}$.

\begin{proposition}
Let $M\in\sM_\zeta$.
\begin{enumerate}[leftmargin=2em]\itemsep0.25cm
\item The action of
\begin{equation}
	\wy_i(z)^+=\sum_{r\ges0}\wy_{i,r}z^r \aand \wy_i(z)^-=-\sum_{r<0}\wy_{i,r}z^r
\end{equation}
is locally given by the expansion at $z=0$ and at $z=\infty$, respectively, of the rational operator:
\begin{equation}
	\Yf_i(z)=([a_{i\tau(i)}]_{q_i}-\ad(\wLetad_{i,1})z-[2]_{q_i}\zeta z^2)^{-1}([a_{i\tau(i)}]_{q_i}\wy_{i,0}
	+[2]_{q_i}\zeta\, \wy_{i,-1}z).
\end{equation}
\item
The action of $\wetad_i(z)$ is locally given by the expansion at $z=0$ of the rational operator:
\begin{equation}
	\begin{aligned}
	\bfHd_i(z)
	&=
	\frac{\wetad_{i,0}+z\wetad_{i,1}-\zeta z^2\wetad_{\tau(i),0}}
	{1-q_i^{a_{i\tau(i)}}\zeta z^2}
	\\
	&\quad-
	\frac{q_i^{a_{i\tau(i)}}\zeta z^2(q_i-q_i^{-1})}
	{1-q_i^{a_{i\tau(i)}}\zeta z^2}
	\Big[
	\wy_{i,-1}
	z^{-1}\bigl(\Yf_{\tau(i)}(z)-\wy_{\tau(i),0}\bigr)
	\\
	&\hspace{2.25cm}
	-
	q_i^{-a_{i\tau(i)}}z^{-1}\bigl(\Yf_{\tau(i)}(z)-\wy_{\tau(i),0}\bigr)\wy_{i,-1}
	\\
	&\hspace{2.75cm}
	+
	\bigl(
	\Yf_{\tau(i)}(z)\wy_{i,0}
	-
	q_i^{-a_{i\tau(i)}}\wy_{i,0}\Yf_{\tau(i)}(z)
	\bigr)
	\Big].
\end{aligned}
\end{equation}
Moreover, $\bfHd_i(z)$ satisfies the identity
\begin{equation}\label{eq:eta-symmetry}
	\bfHd_i(\zeta^{-1}z^{-1})=\frac{1-q_i^{a_{i\tau(i)}}\zeta z^2}{q_i^{a_{i\tau(i)}}-\zeta z^2}\bfHd_{\tau(i)}(z).
\end{equation}
\end{enumerate}
\end{proposition}

\begin{pf}
Part (1) follows from Proposition~\ref{prop:y-rationality} and (M4)--(M5).
For part (2), apply Proposition~\ref{prop:eta-rationality} to the abstract
currents, replacing $\pg_i(\zeta^{-1}z^{-1})\eta_i(z)$ by $\wetad_i(z)$.
Take the coefficient of $z$ in \eqref{eq:etai-rationality}, and use the
$z^0w^0$ coefficient of \eqref{eq:quadratic-y-y-abs} to eliminate
$\wetad_{\tau(i),1}$. This gives the displayed rational formula.
Finally, evaluating \eqref{eq:etai-rationality} at
$w=\zeta^{-1}z^{-1}$ gives
\begin{align}
	\zeta^{-1}z^{-1}(1-q_i^{a_{i\tau(i)}}\zeta z^2)
	\bfHd_{\tau(i)}(z)
	+
	z(1-q_i^{a_{i\tau(i)}}\zeta^{-1}z^{-2})
	\bfHd_i(\zeta^{-1}z^{-1})=0
\end{align}
and \eqref{eq:eta-symmetry} follows.
\end{pf}

\begin{remark}\label{rmk:eta-normalization}
In view of \eqref{eq:eta-symmetry}, it may be useful to correct the generating series $\{\wetad_i(z)\}_{i\in\bfI}$ by any family of rational functions $\{g_i(z)\}_{i\in\bfI}$ such that 
\[g_i(\zeta^{-1}z^{-1})=
\frac{q_i^{a_{i\tau(i)}}-\zeta z^2}{1-q_i^{a_{i\tau(i)}}\zeta z^2}
g_{\tau(i)}(z).\]
Then, the action of $g_i(z)\wetad_i(z)$ on any $M\in\sM_{\zeta}$ is given by
a rational operator $\bfH_i(z)$ satisfying the identity $\bfH_i(\zeta^{-1}z^{-1})=\bfH_{\tau(i)}(z)$.\\

We shall consider the following normalization in $\bfU^\imath$.
Let $\ps_i(z),\pg_i(z)$ be defined as in \ref{ss:SF+GF}. Then, by \eqref{eq:gamma-symmetry}, one can choose $g_i(z)=\pg_i(\FC^{-1}z^{-1})^{-1}$ and set $\weta_i(z)=g_i(z)\wetad_i(z)$. Then, $\xi_{\ps,\pg}(\weta_i(z))=\eta_i(z)$.	
\hfill\rmkend
\end{remark}

\subsection{The category \texorpdfstring{$\sO_\zeta$}{O-zeta}}\label{ss:cat-O-zeta}
The elements of the subalgebra $\bfB\subset\Lext{\bfU}{\FC;\FC^{-1}}$
do not automatically act on any complex $\bfU$-module.
While the generators
$\{\eta_{i,r}\}_{i\in\bfI,r\ges0}$ are Laurent polynomials in $\FC$, and
therefore can be evaluated at any $\zeta\in\bbC^\times$, the generators
$\{y_{i,k}\}_{i\in\bfI,k\in\bbZ}$ are Laurent series in $\FC^{-1}$.\\

Fix $\zeta\in\bbC^\times$. We denote by $\sO_\zeta$ the full subcategory of
$\bfU$-modules $V\in\sO$ whose poles satisfy the following pairwise
non-congruence condition:
\begin{equation}\label{eq:cat-W-zeta}
\forall p,q\in\pole_i(V_\mu)\cup\pole_{\tau(i)}(V_\nu)\colon
pq\notin
\{\zeta^{-1},q_i^{a_{i\tau(i)}}\zeta^{-1}\}
\end{equation}
for every $i\in\bfI$ and every $\mu,\nu\in P(V)$, see \ref{ss:QL-contour}.\\

Let $\ps_i(z),\pg_i(z)$ ($i\in\bfI$) be 
defined as in \eqref{eq:SF+GF}, together with
\eqref{eq:serre-orthogonal-parameters} and \eqref{eq:serre-sigma-parity}. Let
$\pole_{i,\SF},\pole_{i,\GF}\subset\bbC^\times$ 
be the sets of poles of the corresponding rational functions $\SF_i(z)$ and $\GF_i(z)$ at
$\FC=\zeta$, see \eqref{eq:SF+GF}. Thus
\begin{equation}\label{eq:SF+GF-poles}
\pole_{i,\SF}\subseteq\{\pm\zeta^{-1/2}\}
\qquad\mbox{and}\qquad
\pole_{i,\GF}\subseteq
\{\pm q_i^{a_{i\tau(i)}/2}\zeta^{-1/2},\pm\zeta^{-1/2}\}.
\end{equation}
For $\mu\in P(V)$, set
\begin{equation}\label{eq:inverted-pole-set}
\pole_i^\zeta(V_\mu)
=
\{\zeta^{-1}p^{-1}\mid p\in\pole_i(V_\mu)\}.
\end{equation}
The map $\iota_\zeta(w)=\zeta^{-1}w^{-1}$ is an involution, and
$\iota_\zeta(w)=p$ if and only if $w=\zeta^{-1}p^{-1}$. Thus
$\pole_i^\zeta(V_\mu)$ contains the poles of
$\Phif_i(\zeta^{-1}w^{-1})$ and
$\Xf_i^\pm(\zeta^{-1}w^{-1})$ on $V_\mu$.
As a consequence of \eqref{eq:cat-W-zeta}, we obtain the following conditions of separation of poles, which are needed in \ref{ss:QSP-contour-restriction} to define the action of $\bfU^\imath$ on $V\in\sO_\zeta$ by contour integral formulae.

\begin{lemma}\label{eq:pole-symmetry}
For every $V\in\sO_\zeta$, every $i\in\bfI$, and every
$\mu,\nu\in P(V)$, the following conditions hold.
\begin{enumerate}[label=(P\arabic*)]\itemsep0.25cm
\item $\pole_i^\zeta(V_\mu)\cap\pole_{\tau(i)}(V_\nu)=\emptyset$;

\item $\pole_i^\zeta(V_\mu)\cap\pole_{i,\SF}=\emptyset$;

\item $\pole_{\tau(i)}(V_\nu)\cap\pole_{i,\GF}=\emptyset$.

\item $\pole_{\tau(i)}(V_\nu)\cap\pole_{i,\SF}=\emptyset$.
\end{enumerate}
\end{lemma}

\begin{pf}
(P1) follows from the condition $ab\neq\zeta^{-1}$ with $a\in\pole_i(V_\mu)$,  $b\in\pole_{\tau(i)}(V_\nu)$.
(P2) follows from \eqref{eq:SF+GF-poles} and the condition $a^2\neq\zeta^{-1}$ with $a\in\pole_i(V_\mu)$.
(P3) follows from \eqref{eq:SF+GF-poles} and the condition $a^2\not\in\{\zeta^{-1}, q_i^{a_{i\tau(i)}}\zeta^{-1}\}$ with $a\in\pole_{\tau(i)}(V_\nu)$.
(P4) follows from \eqref{eq:SF+GF-poles} and the condition
$a^2\neq\zeta^{-1}$ with $a\in\pole_{\tau(i)}(V_\nu)$.
\end{pf}

\begin{remark}
Since $\pole_{i,\SF}=\pole_{\tau(i),\SF}$, we also have
$\pole_{\tau(i)}^\zeta(V_\mu)\cap\pole_{i,\SF}=\emptyset$.
\hfill\rmkend	
\end{remark}

The subcategories $\sO_\zeta\subset\sO$ are dense in the following sense. Let
$V\in\sO$. Since $P(V)$ is countable and, for each weight $\mu$, $\pole_i(V_\mu)$ is finite, the set of values of $\zeta\in\bbC^\times$
for which the condition \eqref{eq:cat-W-zeta} fails is countable. Equivalently, for a
fixed $\zeta\in\bbC^\times$, all spectral shifts $V(s)$ belong to
$\sO_\zeta$ except for at most countably many values of
$s\in\bbC^\times$.

\subsection{Restriction of $\bfU$-modules}\label{ss:QSP-contour-restriction}
Let $\xi_{\ps,\pg}\colon\bfU^\imath\to\bfB$ be the surjective morphism of comodule algebras constructed in \ref{ss:central-specializations}.\footnote{
For simplicity, we disregard the issue arising from the possible non-invertibility of $\pga{i}$ in $\bbC[\FC^{\pm1}]$. Alternatively, one may assume directly that $\pga{i}\in\bbC^\times$.}
We now prove that $\xi_{\ps,\pg}$ yields an action of $\bfU^\imath$ on every module 
in $\sO_\zeta$.

Let $\sM_{\zeta,\bfc}$ be the full subcategory consisting of modules
on which the central elements $\bfz_i$ act as $\pga{i}=\pg_{i,0}|_{\FC=\zeta}\in\bbC^\times$ for every $i\in\bfI$. Equivalently,
	these are the objects of $\sM_\zeta$ that factor through the central
	quotient with $\FC=\zeta$ and $\bfz_i=\pga{i}$. 

\begin{proposition}\label{prop:restriction}
For every $\zeta\in\bbC^\times$, pullback through $\xi_{\ps,\pg}$ gives a
restriction functor from $\sO_\zeta$ to $\sM_{\zeta,\mathbf c}$. In particular, for every $V\in\sO_\zeta$
and $\mu\in P(V)$, the operator
\[\wy_{i,k}\in{\Hom(V_{\mu}, V_{\mu+\alpha_{\tau(i)}}\oplus V_\mu\oplus V_{\mu-\alpha_i})},\] 
where $i\in\bfI$ and $k\in\bbZ$,
is given by
	\begin{align}
		\wy_{i,k}\vert_{V_{\mu}}
		&=-\oint_{\La_{i,\mu}^\zeta}\Xf^-_{i}(\zeta^{-1}w^{-1})\frac{dw}{w^{k+1}}
		+\oint_{\La_{i,\mu}^{\SF}}\SF_i(w)\Phif_i(\zeta^{-1}w^{-1})\frac{dw}{w^{k+1}}\\
	&\qquad\qquad+\oint_{\La_{i,\mu}^\tau}\GF_i(w)\Xf_{\tau(i)}^+(w)\Phif_i(\zeta^{-1}w^{-1})\frac{dw}{w^{k+1}}
\end{align}
where
$\La^\SF_{i,\mu}, \La^\tau_{i,\mu},\La^\zeta_{i,\mu}\subset\bbC^\times$ are three Jordan curves which do not enclose $0$ and are such that
\begin{itemize}\itemsep0.25cm
	\item $\La^\SF_{i,\mu}$ encloses $\pole_{i,\SF}$ and none of the points in $\pole_i^\zeta(V_\mu)$;
	\item $\La^\tau_{i,\mu}$ encloses $\pole_{\tau(i)}(V_\mu)$ and none of the points in $\pole_i^\zeta(V_\mu)\cup\pole_{i,\GF}$;
		\item $\La^\zeta_{i,\mu}$ encloses $\pole_i^\zeta(V_\mu)$.
\end{itemize}
Moreover, the operators
\begin{equation}
	\wetad_i(z)\in\fext{\End(V_\mu)}{z}
	\aand
	\wy_i(z)^\pm\in\fext{\Hom(V_\mu, V_{\mu-\alpha_i}\oplus V_\mu\oplus V_{\mu+\alpha_{\tau(i)}})}{z^{\pm1}}
\end{equation}
are the expansions at either $z=0$ for $\wetad_i(z)$ and $\wy_i(z)^+$ or $z=\infty$ for $\wy_i(z)^-$ of the rational operators
\begin{align}
	\bfHd_i(z)&=\GF_i(\zeta^{-1}z^{-1})\Phif_i(z)\Phif_{\tau(i)}(\zeta^{-1}z^{-1}),\\[0.5em]
		\Yf_i(z)&=-\oint_{\La_{i,\mu}^\zeta}\Xf^-_{i}(\zeta^{-1}w^{-1})\frac{dw}{w-z}
	+\oint_{\La_{i,\mu}^\SF}\SF_i(w)\Phif_i(\zeta^{-1}w^{-1})\frac{dw}{w-z}\\
	&\qquad\qquad+\oint_{\La_{i,\mu}^{\tau}}\Xf_{\tau(i)}^+(w)\GF_i(w)\Phif_i(\zeta^{-1}w^{-1})\frac{dw}{w-z}.
\end{align}	
\end{proposition}

\begin{pf}
	The existence of the Jordan curves $\La^\SF_{i,\mu}, \La^\tau_{i,\mu}$ is guaranteed by the conditions (P1)--(P4) above.
	We first record the effect of the inversion of the spectral
	parameter. By \eqref{eq:contour-X-Phi}, followed by the change of
	variable $u=\zeta^{-1}w^{-1}$, one has
	\begin{equation}\label{eq:inverted-X-contour}
		\zeta^k x^-_{i,-k}\big|_{V_\mu}
		=
		-\oint_{\La^\zeta_{i,\mu}}
		\Xf_i^-(\zeta^{-1}w^{-1})\frac{dw}{w^{k+1}}.
	\end{equation}
	Indeed, $du=-\zeta^{-1}w^{-2}dw$. The map
	$u\mapsto\zeta^{-1}u^{-1}$ is holomorphic away from $0$ and
	preserves the clockwise orientation of the contours under
	consideration, so the minus sign from the differential remains.
	Relying on \eqref{eq:inverted-X-contour} and on the residue
	representations of $x^+_{\tau(i),k+a+b}$ and $\ps_{i,k+a}$, see
	\ref{ss:QL-contour} and \ref{ss:SF+GF}, we have
	\begin{align}
	y_{i,k}&=x_{i,-k}^-\zeta^k
+\sum_{a,b\ges0}\pg_{i,-a}x^+_{\tau(i),k+a+b}\phi^+_{i,b}\zeta^{-b}
+\sum_{a\ges0}\ps_{i,k+a}\phi^+_{i,a}\zeta^{-a}\\
	&=-\oint_{\La_{i,\mu}^\zeta}\Xf^-_{i}(\zeta^{-1}w^{-1})\frac{dw}{w^{k+1}}
+\oint_{\La_{i,\mu}^\SF}\SF_i(w)\left(\sum_{a\ges0}\phi_{i,a}^+\zeta^{-a}w^{-a}\right)\frac{dw}{w^{k+1}}\\
&\qquad\qquad+\oint_{\La_{i,\mu}^{\tau}}\Xf_{\tau(i)}^+(w)\left(\sum_{a\ges0}\pg_{i,-a}w^{-a}\right)\left(\sum_{b\ges0}\phi^+_{i,b}\zeta^{-b}w^{-b}\right)\frac{dw}{w^{k+1}}\\
	&=-\oint_{\La_{i,\mu}^\zeta}\Xf^-_{i}(\zeta^{-1}w^{-1})\frac{dw}{w^{k+1}}
+\oint_{\La_{i,\mu}^\SF}\SF_i(w)\Phif_i(\zeta^{-1}w^{-1})\frac{dw}{w^{k+1}}\\
&\qquad\qquad+\oint_{\La_{i,\mu}^{\tau}}\Xf_{\tau(i)}^+(w)\GF_i(w)\Phif_i(\zeta^{-1}w^{-1})\frac{dw}{w^{k+1}}	
\end{align}
where the second and third equalities follow since $\GF_i(w)$ is holomorphic inside $\La_{i,\mu}^{\tau}$ and $\Phif_i(\zeta^{-1}w^{-1})$ is holomorphic inside both $\La_{i,\mu}^\SF$ and $\La_{i,\mu}^{\tau}$. The conditions (M1)--(M5) are clearly satisfied, hence the restriction of $V$ to $\bfU^\imath$ lies in $\sM_\zeta$.
Since $\xi_{\ps,\pg}(\bfz_i)=\pga{i}$, it belongs to
$\sM_{\zeta,\mathbf c}$. Finally, the rational expression of the operators $\wetad_i(z)$ and $\wy_i(z)^\pm$ are readily obtained.
\end{pf}


\nc{\FP}{\mathscr{F}_{\pole}}
\nc{\Sub}{\operatorname{Sub}}

\subsection{Principal parts and contour integrals}\label{ss:principal-parts-moments}

We record two elementary facts about rational functions which will be used in the next section, see also \cite{appel-gautam-26}.\\

Let $E$ be a finite-dimensional vector space. Let $\pole\subset\bbC^\times$ be finite and $\FP(E)$ be the set of $E$-valued rational functions that vanish at infinity and whose poles are in $\pole$. 
Recall that any such function can be expressed as the sum of its principal parts, hence
\begin{equation}
\FP(E)
=
\bigoplus_{\alpha\in\pole}\bigoplus_{r\ges1}
E\,(w-\alpha)^{-r}.
\end{equation}

\begin{lemma}\label{lem:principal-parts-moments}
Let 
$F(w)=
\sum_{\alpha\in\pole}F_{\alpha}(w)
\in\FP(E)$ with $F_{\alpha}(w)=\sum_{r=1}^{N_\alpha}
F_{\alpha,r}(w-\alpha)^{-r}$.
\begin{enumerate}
\item 
Set
\[
F_k=\oint_{\C} F(w)\frac{dw}{w^{k+1}}
\]
where $\C\subset\bbC^\times$ is a Jordan curve enclosing $\pole$ and not $0$.
Then,
\begin{equation}\label{eq:modes-principal-parts-equivalence}
	\Span\{F_k\mid k\in\bbZ\}=
	\Span\{F_{\alpha,r}\mid \alpha\in\pole,r=1,\dots, N_{\alpha}\}.
\end{equation}
\item Let $\C_\alpha$ be a Jordan curve around $\alpha$, and $H(w)$ holomorphic within $\C_{\alpha}$ with $H(\alpha)\in\End(E)$ invertible. Then, for $z$ outside $\C_\alpha$,
\begin{align}
\oint_{\C_\alpha}H(w)F_\alpha(w)\frac{dw}{w-z}=\sum_{r=1}^{N_\alpha}\frac{G_r}{(z-\alpha)^r}	
\end{align}
where $(G_1,\dots, G_{N_{\alpha}})$
is obtained from $(F_{\alpha,1},\dots, F_{\alpha, N_{\alpha}})$ through an invertible upper triangular matrix with coefficients in $\End(E)$.
\end{enumerate}
\end{lemma}

\begin{pf}
(1) One has
\begin{equation}
	\oint_\C\frac{1}{(w-\alpha)^r}\frac{dw}{w^{k+1}}=
(-1)^{r}\binom{k+r-1}{r-1}\alpha^{-k-r},
\end{equation} 
hence
\begin{align}\label{eq:modes-principal parts}
F_k&
=\sum_{\alpha\in \pole}\sum_{r=1}^{N_\alpha}
(-1)^{r}\binom{k+r-1}{r-1}\alpha^{-k-r}F_{\alpha,r}
\end{align}
and one inclusion is obvious. For the opposite one, it is enough to observe that taking sufficiently many values of $k$, \eqref{eq:modes-principal parts} produces a linear system whose coefficient matrix is an invertible generalized Vandermonde matrix.\\

(2) Set $t=w-\alpha$ and $u=z-\alpha$. Consider the expansion $H(w)=\sum_{k\ges 0}H_kt^k$ near $\alpha$. Then,
\begin{align}
\oint_{\C_\alpha}H(w)F_\alpha(w)\frac{dw}{w-z}&=
\oint_{\C_\alpha}\Big(\sum_{k\ges0}H_kt^k\Big)\Big(\sum_{r=1}^{N_{\alpha}} F_{\alpha,r}t^{-r}\Big)\Big(-\sum_{n\ges0}t^nu^{-n-1}\Big)dw\\
&=\sum_{r=1}^{N_{\alpha}}\Big(\sum_{\ell=0}^{N_{\alpha}-r}H_\ell F_{\alpha, r+\ell}\Big)(z-\alpha)^{-r},	
\end{align}
and we obtain 
\begin{equation}
G_r=\sum_{\ell=0}^{N_{\alpha}-r}H_\ell F_{\alpha, r+\ell}.	
\end{equation}
Since $H_0=H(\alpha)$ is invertible, the corresponding triangular map is invertible.
\end{pf}

\subsection{Category $\sO$ modules with rational weights}\label{ss:cat-O-rational}
Let $\sO^\bbQ\subset\sO$ be the full subcategory of modules $V\in\sO$ such that $\lambda(h_i)\in\bbQ$ for every $\lambda\in P(V)$ and $i\in\bfI$. We shall prove in Theorem~\ref{thm:generic-submodules-cyclicity-restriction} that irreducible modules in $\sO^\bbQ_\zeta=\sO^\bbQ\cap\sO_\zeta$ are generically irreducible over $\bfU^\imath$.

We start with the following straightforward observation.

\begin{lemma}\label{lem:semisimple}
Let $X$ be a semisimple operator on a vector space $V$, $\Spec(X)\subset\bbC$ the set of its eigenvalues, $f\colon\Spec(X)\to\bbC$ a function, and define $f(X)\in\End(V)$ by $f(X)|_{V_\lambda}=f(\lambda)\id_{V_\lambda}$. If $f$ is injective, then every subspace preserved by $f(X)$ is preserved by $X$.
\end{lemma}

\begin{pf}
Let $v\in W$, and write $v=\sum_{\lambda\in S}v_\lambda$, where
$S\subseteq\Spec(X)$ is finite and $v_\lambda\in\Ker(X-\lambda)$.
It is enough to prove that $v_\lambda\in W$ for each $\lambda\in S$.

We set
\[
p_\lambda(x)=
\prod_{\substack{\mu\in S\\ \mu\neq\lambda}}
\frac{x-f(\mu)}{f(\lambda)-f(\mu)}\in\bbC[x],
\]
which is well-defined since $f$ is injective on $\Spec(X)$.
Then, for every $\mu\in S$, $p_\lambda(f(\mu))=\delta_{\lambda\mu}$.
Therefore,
\[
v_\lambda=\sum_{\mu\in S} p_\lambda(f(\mu))v_\mu=p_\lambda(f(X))v\in W,
\]
since $W$ is preserved by $f(X)$.
\end{pf}

\begin{proposition}\label{prop:cat-O-rational}
Let $V\in\sO^\bbQ$, and $W\subseteq V$ a subspace preserved by the operators $\phi^+_{i,0}-\phi^-_{i,0}$ ($i\in\bfI$). Then, $W$ is preserved by $\phi^\pm_{i,0}$ ($i\in\bfI$).
\end{proposition}

\begin{pf}
For every $i\in\bfI$, let $\Spec_i(V)\subset\bbC^\times$ be the set of eigenvalues of $\phi^+_{i,0}$ on $V$. Specifically, $\Spec_i(V)=\{q_i^{\lambda(h_i)}\mid \lambda\in P(V)\}\subset \bbC^\times$.
It is enough to observe that, since $\hbar\not\in\bbQ$ and $V\in\sO^\bbQ$, the function $f(z)=z-z^{-1}$ is injective on $\Spec_i(V)$. Namely, $f(q^a)=f(q^b)$ if and only if either $e^{\pi\ii\hbar(a-b)}=1$ or $e^{\pi\ii\hbar(a+b)}=-1$. Since $\hbar\not\in\bbQ$, one has that, whenever $a,b\in\bbQ$, then $f(q^a)=f(q^b)$ if and only if $a=b$. Hence the result follows from Lemma~\ref{lem:semisimple}. 
\end{pf}

\subsection{Generic irreducibility under restriction}
\label{ss:generic-irr-restriction}
Let $\bfU'\subset\bfU$ be the subalgebra generated by $\{q^h\}_{h\in\h'}$, $x_{i,r}^\pm$, and $\phi^\pm_{i,\pm k}$ ($i\in\bfI$, $r\in\bbZ$, $k\in\bbZ_{\ges 0}$).
\begin{theorem}\label{thm:generic-submodules-cyclicity-restriction}
Let $V\in\sO^\bbQ_\zeta$ be such that
\begin{itemize}[leftmargin=3em]\itemsep0.25cm
\item for every $i\in\bfI$, $\mu\in P(V)$, and
$\alpha\in\pole_{\tau(i)}(V_\mu)$, the operator
\begin{equation}\label{eq:generic-invertibility-GPhi}
	\GF_i(\alpha)\Phif_i(\zeta^{-1}\alpha^{-1})\big|_{V_\mu}
\end{equation}
is invertible;
\item 
for every $i\in\bfI$ and $\mu,\nu\in P(V)$,
\begin{equation}\label{eq:generic-Cartan-separation}
	\widehat{\pole}_i(V_\mu)\cap
	\{\zeta^{-1}p^{-1}\mid
	p\in\widehat{\pole}_{\tau(i)}(V_\nu)\}=\emptyset,
\end{equation}
where $\widehat{\pole}_i(V_\mu)$ is $\pole_i(V_\mu)$ enlarged by
the poles of $\Phif_i(z)^{-1}|_{V_\mu}$.
\end{itemize}
Consider the action of $\bfU^\imath$ on $V$ given by restriction through $\xi_{\ps,\pg}$.
\begin{enumerate}[leftmargin=3em]\itemsep0.25cm
\item Every $\bfU^\imath$-submodule in $V$ is a
$\bfU'$-submodule. 
\item If $V$ is irreducible as a $\bfU$-module, then its
restriction to $\bfU^\imath$ is irreducible.
\item If $V$ is highest weight over $\bfU$ with highest weight vector $v\in V$, then
its restriction to $\bfU^\imath$ is cyclic on $v$.
\end{enumerate}
\end{theorem}

\begin{pf}
(1)
Let $W\subset V$ be a $\bfU^\imath$-submodule. We first prove that $W$
is stable under the Drinfeld generators of $\bfU$.
First let $E$ be a finite sum of ordinary weight spaces. The rational operator
\[
\bfH_i(z)=\GF_i(\zeta^{-1}z^{-1})^{-1}\bfHd_i(z)
=\Phif_i(z)\Phif_{\tau(i)}(\zeta^{-1}z^{-1})
\]
preserves $W\cap E$, as does the logarithmic derivative $\bfH_i'(z)\bfH_i(z)^{-1}$. Since \eqref{eq:generic-Cartan-separation} separates the poles of the
logarithmic derivatives of the two factors in $\bfH_i(z)$, taking principal parts
at $\bigcup_{V_\mu\subseteq E}\widehat{\pole}_i(V_\mu)$ therefore
recovers $\Phif_i'(z)\Phif_i(z)^{-1}$, which vanishes at
infinity. Thus its coefficients preserve $W\cap E$.
Set $F_i(z)=\Phif_i(z)\phi_{i,0}^-$. Then, $F_i(0)=1$ and recursion on the coefficients in
\[
F_i'(z)=(\Phif_i'(z)\Phif_i(z)^{-1})F_i(z)
\]
shows that $F_i$  preserves $W\cap E$.
In particular, $F_i(\infty)=(\phi_{i,0}^{-})^{2}$ preserves $W\cap E$.
Since $V$ has rational weights and $\hbar\notin\bbQ$, the map
$t\mapsto t^{2}$ is injective on $\Spec(\phi^-_{i,0})$, hence
Lemma~\ref{lem:semisimple} gives stability under $\phi^-_{i,0}$.
Consequently all Laurent coefficients of $\Phif_i(z)$ and its
inverse preserve $W\cap E$. Since $E$ was arbitrary, the ordinary
Cartan currents preserve $W$.\\

Let $v\in W$ and write
\[
v=\sum_\mu v_\mu\in\bigoplus_\mu V_\mu.
\]
	For each weight component, by Proposition~\ref{prop:restriction}, the action of the current $y_i(z)^+$ is given by
\begin{equation}\label{eq:Y-action-three-pieces-generic}
	\begin{aligned}
	\Yf_i(z)\big|_{V_\mu}
	&=
	-\oint_{\Lambda^\zeta_{i,\mu}}
	\Xf_i^-(\zeta^{-1}w^{-1})
	\frac{dw}{w-z}
	\\
	&\quad\quad+
	\oint_{\Lambda^\Sigma_{i,\mu}}
	\SF_i(w)\Phif_i(\zeta^{-1}w^{-1})
	\frac{dw}{w-z}
	\nonumber\\
	&\quad\quad\quad+
	\oint_{\Lambda^\tau_{i,\mu}}
	\GF_i(w)\Xf_{\tau(i)}^+(w)\Phif_i(\zeta^{-1}w^{-1})
	\frac{dw}{w-z}.
	\end{aligned}
\end{equation}
Clearly, $W$ is stable under the action of $y_{i,k}$ ($i\in\bfI$, $k\in\bbZ$).
Hence, by Lemma~\ref{lem:principal-parts-moments} (1), all principal part
coefficients of the rational function $\Yf_i(z)v$ belong to $W$.\\

Now, fix $\mu$ such that $v_{\mu}\neq0$ and $\alpha\in\bbC^\times$ is a pole of
$\Xf_i^-(\zeta^{-1}w^{-1})v_{\mu}$.  Hence, $\alpha\in\pole_i^\zeta(V_\mu)$ and therefore, 
by \eqref{eq:cat-W-zeta}, $\Xf_{\tau(i)}^+(w)v_\nu$ must be regular at $w=\alpha$
	for any weight component $v_\nu$ of $v$. By (P1) and (P2), at
	$z=\alpha$ the second and third summands in
	\eqref{eq:Y-action-three-pieces-generic}, whose contours enclose the poles of $\SF_i(w)$
	and $\Xf_{\tau(i)}^+(w)$, respectively, are regular.

Hence the principal part of $\Yf_i(z)v$ at $z=\alpha$ comes only from the first summand.
By Lemma~\ref{lem:principal-parts-moments} (2), applied with $H(w)=1$, the
	principal part coefficients of $-\Xf_i^-(\zeta^{-1}w^{-1})v$,
	and therefore also those of $\Xf_i^-(\zeta^{-1}w^{-1})v$,
at $w=\alpha$ belong to $W$. Varying $\alpha$, and then applying
Lemma~\ref{lem:principal-parts-moments} (1), we conclude that 
$x_{i,r}^-v\in W$ for all $r\in\bbZ$. Thus, since $v\in W$ is arbitrary,
we have $x_{i,r}^-W\subset W$ for all $i\in\bfI$ and $r\in\bbZ$.\\

	Let now $\alpha$ be a pole of $\Xf_{\tau(i)}^+(w)v_\mu$ 
	for some weight component $v_\mu$, thus $\alpha\in\pole_{\tau(i)}(V_\mu)$. 
By \eqref{eq:cat-W-zeta}, $\alpha$ is not a pole of
$\Xf_i^-(\zeta^{-1}w^{-1})v_\nu$
	for any weight component $v_\nu$ of $v$. By (P1) and (P4), in this case
	the first two summands in
\eqref{eq:Y-action-three-pieces-generic} are regular at $w=\alpha$, and
the principal part of $\Yf_i(z)v$ at $z=\alpha$ comes only from the third summand.
Set
\[
\Phif^{\GF}_i(w)=\GF_i(w)\Phif_i(\zeta^{-1}w^{-1}).
\]
Choose a finite sum $E$ of ordinary weight spaces containing $v$. 
The preceding argument, applied to every
vector in $W\cap E$, shows that the principal part coefficients of
$\Xf_{\tau(i)}^+(w)\Phif^{\GF}_i(w)$ map $W\cap E$ into $W$.
By \eqref{eq:generic-Cartan-separation} and
\eqref{eq:generic-invertibility-GPhi}, $\Phif^{\GF}_i(w)$ is invertible on $E$,
and the Taylor coefficients of $\Phif^{\GF}_i(w)^{-1}$ at $\alpha$ preserve $W\cap E$.
Multiplying the principal part on the right by this inverse and
taking its principal part again recovers that of
$\Xf_{\tau(i)}^+(w)|_E$.
Hence these coefficients belong to $W$ and, by
Lemma~\ref{lem:principal-parts-moments} (1), we get
$x_{\tau(i),r}^+v\in W$ for all $r\in\bbZ$. This gives $x_{i,r}^+W\subset W$ for all $i\in\bfI$ and $r\in\bbZ$.\\

From \eqref{eq:QL-full-rels-double}, it follows that the action of $\phi^{\pm}_{i,k}$ ($k\neq0$) and $\phi^+_{i,0}-\phi^-_{i,0}$ preserve $W$ for every $i\in\bfI$. 
Hence, the result follows from Proposition~\ref{prop:cat-O-rational}.\\

(2) It is enough to observe that, if $V$ is irreducible over $\bfU$,
then it is irreducible over $\bfU'$.\footnote{
Indeed, every irreducible
	module in $\sO$ is the irreducible quotient of a Verma module associated
	with a loop weight \cite{hernandez-05}. The Verma modules for $\bfU$ and
	$\bfU'$ are isomorphic as $Q^+$-graded vector spaces and the operators $q^h$ ($h\in\h''$) of $\bfU$
	act diagonally on the homogeneous components. 
	Consequently, the maximal proper $\bfU'$-submodule is
	automatically preserved, and hence it coincides with the maximal proper $\bfU$-submodule.}
Hence, the result follows from~(1).\\

(3) As in (2), it is enough to observe that, if $V$ is highest weight over $\bfU$, then 
it is highest weight over $\bfU'$. Let $v\in V$ be a highest weight vector such that $V=\bfU\cdot v$. Then, 
\[
\bfU^\imath\cdot v=\bfU'\cdot v=\bfU\cdot v=V,
\]
where the first equality follows from (1).
\end{pf}

\begin{remark}
In \cite[Sec.~3.2]{negut-26}, Negu\textcommabelow{t} proved that, in the case of a strongly symmetrizable Kac--Moody algebra, the Drinfeld pairing on $\bfU$ is non-degenerate.
By a result of \cite{appel-gautam-26}, it then follows that irreducible modules in $\sO^{\bbQ}$ separate the points of $\bfU$.
A PBW argument then shows that the morphism $\xi_{\ps,\pg}$ from the central quotient in \eqref{eq:embedding} is injective.
We expect this to be true in arbitrary type, see also the special case described in \cite[Thm.~4.13]{lu-pan-26}.
\hfill\rmkend
\end{remark}


\subsection{Action of $\bfU$-modules}
For any $V\in\sO_\zeta$, we denote by $\bfHd_i(z)$ and $\Yf_i(z)$ the rational operators
describing the action of $\etad_{i}(z)$ and $y_i(z)^{\pm}$ on $V$.

\begin{proposition}\label{prop:action}
	Let $V,W\in\sO_\zeta$, and $s\in\bbC^\times$ be such that, for every $\lambda\in P(V), \mu\in P(W)$ and $i\in\bfI$, 
	\begin{equation}\label{eq:action-pole-separation}
		\begin{gathered}
		\pole^\bfB_i(V_{\lambda})
		\cap s\pole^\zeta_i(W_\mu)=\emptyset,\\
		\pole^{\bfB,\zeta}_i(V_{\lambda})
		\cap s^{-1}\pole_{\tau(i)}(W_\mu)=\emptyset,\\
		s\pole^\zeta_i(W_\mu)
		\cap s^{-1}\pole_{\tau(i)}(W_\mu)=\emptyset
		\end{gathered}
	\end{equation} 
	where $\pole^\bfB_i(V_\lambda)\subset\bbC^\times$ is the set of poles relative to the action of $y_i(z),\deta_i(z)=\pg_i(\FC^{-1}z^{-1})\eta_i(z)\in\bfB$ on $V_\lambda$, $\pole_i(W_\mu)\subset\bbC^\times$ is as above, and
	$\pole_i^{\bfB,\zeta}(V_\lambda)
	=(\pole_i^\bfB(V_\lambda))^\zeta$.\footnote{
	As before, for any $S\subset\bbC^\times$, $S^\zeta=\{\zeta^{-1}z^{-1}\mid z\in S\}$.}
The action of $y_{i,k}$ on $V\DT{s}W$ is given by
\begin{equation}\label{eq:rational-DD-y}
	\begin{split}
		\DD{s}(y_{i,k})
		&=1\otimes x^-_{i,-k}\zeta^k s^{-k}
		+\oint_{\cC^\bfB_{i,\lambda}}\Yf_i(w)\ten\Phif_i(\zeta^{-1}w^{-1}s)w^{-k-1}dw\\
		&\qquad+\oint_{\cC_{i,\mu}^\bfB}
			\bfHd_i(\zeta^{-1}w^{-1})\ten\Xf^+_{\tau(i)}(sw)\Phif_i(\zeta^{-1}w^{-1}s)
		w^{-k-1}dw\\
	\end{split}
\end{equation}
for every $\lambda\in P(V), \mu\in P(W)$, where $\cC^\bfB_{i,\lambda},\cC^\bfB_{i,\mu}\subset\bbC^\times$ are two Jordan curves which do not enclose $0$ and are such that
		\begin{itemize}\itemsep0.25cm
			\item $\cC^\bfB_{i,\lambda}$ encloses $\pole^\bfB_i(V_\lambda)$ and none of the points in $s\pole^{\zeta}_i(W_\mu)$;
			\item $\cC^\bfB_{i,\mu}$ encloses $s^{-1}\pole_{\tau(i)}(W_\mu)$ and none of the points in $\pole^{\bfB,\zeta}_i(V_\lambda)$ or $s\pole^{\zeta}_i(W_\mu)$.
		\end{itemize}
In terms of currents, the action on $V\DT{s}W$ reads
\begin{equation}\label{eq:rational-DD-y-current}
		\begin{split}
			\DD{s}(\deta_i(z))&=\bfHd_i(z)\ten\Phif_i(zs)\Phif_{\tau(i)}(\zeta^{-1}z^{-1}s),\\
			\DD{s}(y_{i}(z)^{\pm})
			&=1\otimes
			\left(x^-_{i,0}-\Xf_i^-(\zeta^{-1}z^{-1}s)\right)
			+\oint_{\cC^\bfB_{i,\lambda}}\Yf_i(w)\ten\Phif_i(\zeta^{-1}w^{-1}s)\frac{dw}{w-z}\\
			&\hspace{2cm}+\oint_{\cC_{i,\mu}^\bfB}
			\bfHd_i(\zeta^{-1}w^{-1})\ten\Xf^+_{\tau(i)}(ws)\Phif_i(\zeta^{-1}w^{-1}s)
		\frac{dw}{w-z}.
	\end{split}
\end{equation}
\end{proposition}

\begin{pf}
The existence of the Jordan curves $\cC^\bfB_{i,\lambda},
\cC^\bfB_{i,\mu}$ is guaranteed by
\eqref{eq:action-pole-separation}. Then, we proceed as in
Proposition~\ref{prop:restriction}:
\begin{equation}
\begin{split}
	\DD{s}(y_{i,k})&=
	1\otimes x^-_{i,-k}\zeta^k s^{-k}
	+\sum_{a\ges 0}y_{i,k+a}\otimes \phi^+_{i,a}\zeta^{-a}s^a\\
	&\qquad+
	\sum_{a,b,c\ges 0}
	\pg_{i,-c}\eta_{i,a}\otimes x^+_{\tau(i),\,k+a+b+c}\phi^+_{i,b}
	\,\zeta^{-a-b}s^{k+a+c+2b}\\
	&=1\otimes x^-_{i,-k}\zeta^k s^{-k}
	+\oint_{\cC^\bfB_{i,\lambda}}\Yf_i(w)\ten
	\left(\sum_{a\ges0}\phi^+_{i,a}\zeta^{-a}s^aw^{-a}\right)w^{-k-1}dw\\
	&\hspace{-1.5cm}+\oint_{\cC_{i,\mu}^\bfB}
	\left(\sum_{a,c\ges 0}\pg_{i,-c}\eta_{i,a}\zeta^{-a}w^{-a-c}\right)
	\ten
	\Xf^+_{\tau(i)}(sw)\left(\sum_{b\ges 0}\phi^+_{i,b}\zeta^{-b}s^bw^{-b}\right)
	w^{-k-1}dw\\
	&=1\otimes x^-_{i,-k}\zeta^k s^{-k}
	+\oint_{\cC^\bfB_{i,\lambda}}\Yf_i(w)\ten\Phif_i(\zeta^{-1}w^{-1}s)w^{-k-1}dw\\
	&\qquad+\oint_{\cC_{i,\mu}^\bfB}
	\bfHd_i(\zeta^{-1}w^{-1})\ten\Xf^+_{\tau(i)}(sw)\Phif_i(\zeta^{-1}w^{-1}s)
	w^{-k-1}dw.
\end{split}
\end{equation}
The result follows.
\end{pf}

\begin{remarks}
\hfill
\begin{enumerate}[leftmargin=2em]\itemsep0.25cm
\item
Relying on the same formulae, one obtains an action of $\sO_\zeta$ on $\sM_\zeta$.
Specifically, the morphism $\DDi{s}$
induces an action of $\bfU^\imath$ on $M\DT{s}V$ where $M\in\sM_\zeta$ and 
$V\in\sO_\zeta$. By \ref{cond:M1}--\ref{cond:M5}, this is described by
\eqref{eq:rational-DD-y-current}, where the Jordan curves are defined locally.
By \eqref{eq:central-coinvariance}, this action preserves
$\sM_{\zeta,\mathbf c}$ for each choice of the central scalars.
\item 
By the same pole-separation argument as in Theorem~\ref{thm:generic-submodules-cyclicity-restriction},
for $V,W\in\sO_\zeta^\bbQ$, the restriction of $V\DT{s}W$ to
$\bfU^\imath$ is generically cyclic if $V$ and $W$ are highest weight,
and generically irreducible if both remain irreducible under restriction to $\bfU^\imath$.

Similarly, for $M\in\sM_\zeta$ and $V\in\sO_\zeta^\bbQ$, the module
$M\DT{s}V$ is generically irreducible if $M$ is irreducible over $\bfU^\imath$
and $V$ remains irreducible under restriction to $\bfU^\imath$.
\hfill\rmkend
\end{enumerate}
\end{remarks}


\section{The abelian R-matrix}
\label{s:abelian-R-matrix}

In this section, we briefly review the representation-theoretic construction of 
meromorphic abelian R-matrices on category $\O$ modules over $\bfU$ following 
\cite{appel-gautam-26}.

\subsection{Quantum Cartan matrices}\label{ss:quantum-Cartan}
Henceforth, we assume $|q|\neq1$, equivalently $\hbar\notin\bbR$.
For $\sfp \in \bbC^\times$, let $\qCB(\sfp)=(\qCB_{ij}(\sfp))_{i,j\in\bfI}$ be the
{\em symmetrized $\sfp$-Cartan matrix} defined by
\[\qCB_{ij}(\sfp)=[d_i a_{ij}]_\sfp=[d_j a_{ji}]_\sfp=\qCB_{ji}(\sfp)\]
for every $i,j\in\bfI$. Then,
\[
[h_{i,k},x_j^\pm(w)]=\pm\frac{q-q^{-1}}{q_i-q_i^{-1}}\frac{[k]_q}{k}\qCB_{ij}(q^k) w^{-k} x_j^\pm(w).
\]
In the following, we shall consider the \emph{modified} symmetrized $\sfp$-Cartan matrix
$\qCBd(\sfp)=(\qCBd_{ij}(\sfp))_{i,j\in\bfI}$ defined by
\begin{equation}
\qCBd_{ij}(\sfp)=\sfp^{d_ia_{ij}}-\sfp^{-d_ia_{ij}}=(\sfp-\sfp^{-1})\qCB_{ij}(\sfp).
\end{equation}
We denote by $\qCBd^*(\sfp)$ the adjugate matrix. Note that
\[
\qCBd^*(\sfp)
=
(\sfp-\sfp^{-1})^{|\bfI|-1}\qCB^*(\sfp).
\]
Set
\begin{equation}\label{eq:def-D(p)}
\scrD(\sfp)=\det(\qCBd(\sfp))
=\sum_{r=m}^{M}c_r\sfp^r
\end{equation}
with $c_mc_M\neq0$.
Since
$\scrD(\sfp^{-1})=(-1)^{|\bfI|}\scrD(\sfp)$, one has
$m=-M$ and $c_m=(-1)^{|\bfI|}c_M$.\\

Henceforth, we assume that the symmetrizable Cartan matrix $\sfA$
satisfies the following conditions:
\begin{enumerate}[label=(C\arabic*)]\itemsep0.25cm
\item\label{cond:qCartan-non-vanishing} for every $k\in\bbZ_{\neq0}$, $\scrD(q^k)\neq0$;
\item\label{cond:qCartan-dominant} $c_m,c_M\in\{\pm 1\}$.
\end{enumerate}

\begin{remark}
Conditions~\ref{cond:qCartan-non-vanishing}
and~\ref{cond:qCartan-dominant} are satisfied when $\sfA$ is of finite
type, or of affine type other than $A_1^{(1)}$, see
\cite{khoroshkin-tolstoy-93} and \cite[App.~A]{gautam-toledano-17}
for finite type, and \cite[App.~C]{appel-gautam-wendlandt-26} for
affine type.
\hfill\rmkend
\end{remark}

For every $k\in\bbZ_{\neq0}$, 
\[
\det(\qCBd(q^k))
=
(q^k-q^{-k})^{|\bfI|}\det(\qCB(q^k))
\qquad(k\in\bbZ_{\neq0}).
\] 
Thus, since $q$ is not a root of unity, the condition~\ref{cond:qCartan-non-vanishing} is equivalent to the
invertibility of $\qCB(q^k)$ for every $k\in\bbZ_{\neq0}$. We write
\[
\qCC(q^k)
=
\bigl(\qCC_{ij}(q^k)\bigr)_{i,j\in\bfI}
=
\qCB(q^k)^{-1}.
\]
For every $i\in\bfI$ and $k\in\bbZ_{\neq0}$, set
\[
\wh{h}_{i,k}
=
\sum_{a\in\bfI}
\frac{q_a-q_a^{-1}}{q-q^{-1}}
\frac{k}{[k]_q}
\qCC_{ia}(q^k)h_{a,k}.
\]
Then, for every $j\in\bfI$, one has
\[
\begin{aligned}
[\widehat h_{i,k},x_j^\pm(w)]
&=
\sum_{a\in\bfI}
\frac{q_a-q_a^{-1}}{q-q^{-1}}
\frac{k}{[k]_q}
\qCC_{ia}(q^k)[h_{a,k},x_j^\pm(w)]\\
&=
\pm\delta_{ij}w^{-k}x_j^\pm(w).
\end{aligned}
\]

\subsection{The abelian tensor Casimir operator}\label{ss:operator-C}
Let $f\in\End(\h)$ be a linear endomorphism which is self-adjoint with respect to the symmetric bilinear form on $\h$. With a slight abuse of notation, we shall use the same letter to denote its dual $f\in\End(\h^*)$. Set
\begin{equation}\label{eq:operator-C}
\bC{f}{}=q^{(\id\otimes f)(\Omega_0)}=q^{(f\otimes\id)(\Omega_0)}	
\end{equation}
where $\Omega_0=\sum_iu_i\ten u^i\in\h\ten\h$ is the canonical tensor. Hence, $\bC{f}{}$ naturally acts on any tensor product in $\sO$. Specifically, for every $V,W\in\sO$, $\lambda\in P(V)$ and $\mu\in P(W)$, one has
\begin{align}\label{eq:operator-C-weight}
\bC{f}{V_\lambda, W_\mu}=\bC{f}{}\mid_{V_\lambda\otimes W_\mu}=q^{(\lambda,f(\mu))}=q^{(f(\lambda),\mu)}.
\end{align}
Moreover, one has
\begin{equation}
\label{eq:coproduct-operator-C}
\begin{aligned}
	\bC{f}{V_1\DT{s}V_2,V_3}&=\bC{f}{V_1,V_3}\bC{f}{V_2,V_3},\\
	\bC{f}{V_1,V_2\DT{s}V_3}&=\bC{f}{V_1,V_3}\bC{f}{V_1,V_2},
\end{aligned}
\end{equation}
and
\begin{align}\label{eq:twisted-operator-C}
\bC{f}{V_1,V_2}=(1\ 2)\circ\bC{f}{V_2,V_1}\circ(1\ 2)	
\end{align}
for every $V_1,V_2,V_3\in\sO$.
Note that, for every homogeneous $x\in\bfU$ with $\wgt(x)=\beta$, 
one has 
\begin{equation}
\label{eq:adjoint-C}
\Ad(\bC{f}{})(x\ten 1)=x\ten q^{\tinv{f(\beta)}}
\aand
\Ad(\bC{f}{})(1\ten x)=q^{\tinv{f(\beta)}}\ten x,
\end{equation}
where $\tinv{\alpha}=\nu^{-1}(\alpha)$, see \ref{ss:sym-km-summary}.
Note that $\tinv{f(\alpha)}=\nu^{-1}(f(\alpha))=f(\nu^{-1}(\alpha))=f(\tinv{\alpha})$.
We set $\bC{\pm}{}=\bC{\pm\id}{}$.

\subsection{The abelian R-matrix}\label{ss:abelian-R-matrix}
In \ref{ss:QL-fields}, we set
$h^{\pm}_i(z)=\pm(q_i-q_i^{-1})\sum_{k>0}h_{i,\pm k}z^{\pm k}$.
Let $V,W\in\sO$, $\lambda\in P(V)$, $\mu\in P(W)$. 
For every $i,j\in\bfI$, set
\begin{equation}\label{eq:Tij}
\T_{ij}^{\pm}(s)
=\sum_{k>0}(q_i-q_i^{-1})(q_j-q_j^{-1})kh_{i,\mp k}\ten h_{j, \pm k}s^{\pm k}.
\end{equation}
Let $\sfp$ be the q-shift operator given by $\sfp\cdot f(s)=f(qs)$. We consider the following multiplicative difference equation
\begin{equation}\label{eq:R-diff-eq}
\scrD(\sfp)\cdot\Lambda^{\pm}(s)=\sum_{i,j\in\bfI}\qCBd^*_{ij}(\sfp)\cdot\T_{ij}^{\pm}(s)
\end{equation}
where $\Lambda^+(s)\in\fext{\bfU^0\ten\bfU^0}{s}$ and $\Lambda^-(s)\in\fext{\bfU^0\ten\bfU^0}{s^{-1}}$.

\begin{theorem}[\cite{appel-gautam-26}]\label{thm:appel-gautam-26}
Assume that Conditions~\ref{cond:qCartan-non-vanishing}
and~\ref{cond:qCartan-dominant} hold.
Then, the system \eqref{eq:R-diff-eq} admits two unique formal solutions $\Lambda^{\pm}(s)\in s^{\pm1}\fext{\bfU^0\ten\bfU^0}{s^{\pm1}}$. Moreover, the operators
\begin{equation}
\RR^{0,\pm}(s)=\bC{\pm}{}\exp(\Lambda^{\pm}(s))
\end{equation}
induce two meromorphic braidings on $\sO$ related by unitarity.
Specifically, for every $V_1,V_2\in\sO$, $\RR_{V_1V_2}^{0,\pm}(s)$
define two meromorphic $\End(V_1\otimes V_2)$--valued functions
with the following properties:
\begin{enumerate}[font=\upshape]\itemsep0.25cm
	\item\label{R0:int} For every $V_1,V_2$,
	\[
	(1\ 2) \circ \RR^{0,\pm}_{V_1,V_2}(s) : V_1\DT{s} V_2 \to
	(V_2\DT{s^{-1}}V_1)(s)
	\]
	is a $\bfU$--intertwiner.
	\item\label{R0:shift} For every $V_1,V_2\in\sO$
	and $a,b\in\bbC^\times$, one has
	\[
	\RR^{0,\pm}_{V_1(a),V_2(b)}(s) = \RR^{0,\pm}_{V_1,V_2}(sa^{-1}b).
	\]
	\item\label{R0:cabling} For every $V_1,V_2,V_3\in\sO$, one has
	\[
	\begin{aligned}
		\RR^{0,\pm}_{V_1\DT{s_1}V_2,V_3}(s_1s_2) &= \RR^{0,\pm}_{V_1,V_3}(s_1s_2)
		\RR^{0,\pm}_{V_2,V_3}(s_2), \\
		\RR^{0,\pm}_{V_1,V_2\DT{s_2}V_3}(s_1) &= \RR^{0,\pm}_{V_1,V_3}(s_1s_2)
		\RR^{0,\pm}_{V_1,V_2}(s_1).
	\end{aligned}
	\]
	\item\label{R0:unitary} For every $V_1,V_2\in\sO$, one has
	\[
	\RR^{0,+}_{V_1,V_2}(s)^{-1} = (1\ 2)\circ
	\RR^{0,-}_{V_2,V_1}(s^{-1})\circ (1\ 2).
	\]
	\item\label{R0:rat} For all highest weight representations $V_1,V_2\in\sO$, there are meromorphic functions $f^{\pm}_{V_1V_2}(s)$ such that the normalized operator
	\[
	\RR^{0,\norm}_{V_1V_2}(s)=f^{\pm}_{V_1V_2}(s)^{-1}\RR^{0,\pm}_{V_1V_2}(s)
	\]
	is independent of $\pm$, rational and unitary in $s$.
\end{enumerate}
\end{theorem}

The two solutions are easily computed. Namely, it is enough to observe that, if
$f(s)=\sum_{k\ges0}f_ks^{\pm k}$, then
\[\D(\sfp)\cdot f(s)=\sum_{k\ges 0}\D(q^{\pm k})f_ks^{\pm k}\]
for every $\D(\sfp)\in\bbC[\sfp^{\pm1}]$. Since 
\begin{equation}
\frac{1}{q^k-q^{-k}}\qCC_{ij}(q^k)
=\frac{1}{q^k-q^{-k}}\frac{\qCB^*_{ij}(q^k)}{\det(\qCB(q^k))}
=\frac{\qCBd^*_{ij}(q^k)}{\det(\qCBd(q^k))},
\end{equation} 
one has
\begin{equation}\label{eq:log-R}
\begin{aligned}
	\LR^{\pm}(s)&=\pm\sum_{k>0}s^{\pm k}\sum_{i,j\in\bfI}
	\bfa_{ij,k}\,h_{i,\mp k}\otimes h_{j,\pm k},
\end{aligned}
\end{equation}
where
\begin{equation}\label{eq:def-coeff-a}
\begin{aligned}
	\bfa_{ij, k}=(q_i-q_i^{-1})(q_j-q_j^{-1})\frac{k}{q^k-q^{-k}}\qCC_{ij}(q^{k}).
\end{aligned}
\end{equation}
Note that $\bfa_{ij,k}=\bfa_{ji,k}$ and $\bfa_{ij,-k}=\bfa_{ij,k}$.
This result shows that the tensor provided in \cite{hernandez-05} 
gives rise to a meromorphic abelian R-matrix on the category $\sO$ for quantum affinizations and recovers, in particular, the action of the abelian R-matrix for quantum affine algebras obtained in \cite{khoroshkin-tolstoy-93,levendorskii-soibelman-stukopin-93,damiani-98, gautam-toledano-17}.

\begin{remark}
While the cabling identities and the unitarity property are clearly valid at the universal level, the intertwining identity of $\RR^{0,\pm}(s)$ is valid only at the level of representations.
\hfill\rmkend 
\end{remark}

\section{The basic abelian K-matrix}
\label{s:basic-abelian-K-matrix}

In this section, we adapt the strategy from \cite{appel-gautam-26} to construct 
a meromorphic abelian K-matrix on restricted $\bfU$-modules in $\sO_\zeta$, intertwining the action of $\bfU^\imath$. In finite type, this K-matrix is expected to coincide with the weight-preserving component of the trigonometric K-matrices constructed in \cite{appel-vlaar-25a} for finite-dimensional representations of the corresponding quasi-split quantum symmetric pairs of affine type.

\subsection{The abelian Casimir operator}\label{ss:Casimir}
Let $\Omega_0=\sum_iu_i\ten u^i\in\h\ten\h$ be defined as in \ref{ss:operator-C}.
For every self-adjoint endomorphism $f\in\End(\h)$,
we set
\begin{equation}\label{eq:operator-k}
\bK{f}{}=q^{C^f_0/2},
\end{equation}
where $C^f_0=\sum_if(u_i)u^i$. Hence, for every $V\in\sO$ and $\lambda\in P(V)$, one has
\begin{equation}
\bK{f}{V_\lambda}=q^{(\lambda,f(\lambda))/2}.	
\end{equation}
In particular, one has
\begin{align}\label{eq:operator-k-twisted}
\bK{f}{V_{\lambda}}=q^{(\lambda, f(\lambda))/2}=q^{(f(\lambda), f^2(\lambda))/2}=\bK{f}{V_{f(\lambda)}}
\end{align}
whenever $f^2=\id$.
Moreover, one has
\begin{equation}
\label{eq:coproduct-operator-k}
	\bK{f}{V_1\DT{s}V_2}=\bC{f}{V_1,V_2}\bK{f}{V_1}\bK{f}{V_2}
\end{equation}
for every $V_1,V_2\in\sO$.
Note that, for every homogeneous $x\in\bfU$ with $\wgt(x)=\beta$, 
one has 
\begin{align}
\label{eq:adjoint-operator-k}
\Ad(\bK{f}{})(x)&=q^{(\beta,f(\beta))/2}xq^{\tinv{f(\beta)}},
\end{align}
where $\tinv{\alpha}=\nu^{-1}(\alpha)$, see \ref{ss:sym-km-summary}.

\subsection{The logarithmic parameters}\label{ss:bK0-coeff-b}
Fix $\zeta\in\bbC^\times$. Let $\ps_i(z),\pg_i(z)$ ($i\in\bfI$) be 
defined as in \eqref{eq:SF+GF}, together with
\eqref{eq:serre-orthogonal-parameters} and \eqref{eq:serre-sigma-parity}. 
We consider the resulting complex parameters and rational functions as in Section~\ref{ss:cat-O-zeta}.
In particular, $\pga{i}\neq0$ for every $i\in\bfI$. Thus,
$\GF_i(0)\GF_i(\infty)\neq0$. We shall need the following functions.
\begin{enumerate}[leftmargin=2em]\itemsep0.25cm
\item We define
$\bfG_{i}^\pm(z)=\sum_{k>0}g_{i,k}^{\pm}z^{\pm k}\in\fext{\bbC}{z^{\pm1}}$ by
\begin{equation}\label{eq:log-GF-pm}
\exp(\bfG_i^\pm(z))
=\GF^\pm_i(z),
\end{equation}
where $\GF^+_i(z)=\frac{\GF_i(z)}{\GF_i(0)}$, $\GF^-_i(z)=\frac{\GF_i(z)}{\GF_i(\infty)}$, and  the two expressions are expanded at $z=0$ and $z=\infty$,
respectively. 
\item We define $\bfD_i^\pm(z)=\sum_{k>0}d_{i,k}^\pm z^{\pm k}\in\fext{\bbC}{z^{\pm1}}$ by
\begin{equation}\label{eq:def-d-exp}
\exp(2\bfD_i^\pm(z))=\GF_i^\pm(z)^2\left(
\frac{1-q_i^{\mp a_{i\tau(i)}}\zeta^{\pm 1} z^{\pm 2}}
{1-q_i^{\pm a_{i\tau(i)}}\zeta^{\pm 1} z^{\pm 2}}
\right),
\end{equation}
where as before the two expressions are expanded at $z=0$ and $z=\infty$,
respectively. More explicitly, we set
\begin{equation}\label{eq:def-d-bK}
d_{i,k}^{\pm}
=
\begin{cases}
\displaystyle
g_{i,k}^{\pm}\pm\frac{1}{k}
\left(
q_i^{\frac{k}{2}a_{i\tau(i)}}
-q_i^{-\frac{k}{2}a_{i\tau(i)}}
\right)
\zeta^{\pm k/2}& \mbox{if $k$ is even},\\[1.1em]
\hfill g_{i,k}^{\pm}\hfill & \mbox{if $k$ is odd}.
\end{cases}
\end{equation}
\end{enumerate}
Finally, we set
\begin{equation}\label{eq:bK0-coefficients}
\bfb_{i,k}^{\pm}
=
-\zeta^{\mp k}
k(q_i-q_i^{-1})
\frac{1}{q^k-q^{-k}}
\sum_{j\in\bfI}
\qCC_{ij}(q^k)d_{j,k}^{\pm}.
\end{equation}

\begin{lemma}\label{lem:bK0-coeff-b}
Let $\bfb_{i,k}^\pm\in\bbC$ be defined as in \eqref{eq:bK0-coefficients}.
Then, $\bfb_{i,k}^{+}=\zeta^{-k}\bfb_{\tau(i),k}^{-}$.
\end{lemma}
\begin{pf}
Relation \eqref{eq:gamma-symmetry}, together with
$\GF_\ell(0)=q_\ell^{-a_{\ell\tau(\ell)}}
\GF_{\tau(\ell)}(\infty)$, gives
\begin{equation}
	\log\left(
	\frac{\GF_\ell(\zeta^{-1}z^{-1})}{\GF_\ell(0)}
	\right)
	={}
	\log\left(
	\frac{\GF_{\tau(\ell)}(z)}
	{\GF_{\tau(\ell)}(\infty)}
	\right)
	+
	\log\left(
	q_\ell^{a_{\ell\tau(\ell)}}
	\frac{q_\ell^{a_{\ell\tau(\ell)}}-\zeta z^2}
	{1-q_\ell^{a_{\ell\tau(\ell)}}\zeta z^2}
	\right).
\end{equation}
The last logarithm has expansion at infinity
\[
-\sum_{k>0}
\frac{q_\ell^{ka_{\ell\tau(\ell)}}
	-q_\ell^{-ka_{\ell\tau(\ell)}}}{k}
\zeta^{-k}z^{-2k}.
\]
Hence, by comparing coefficients and using \eqref{eq:def-d-bK}, one has
\begin{equation}\label{eq:d-bK-unitarity}
	d_{\tau(\ell),k}^{-}=\zeta^{-k}d_{\ell,k}^{+}.
\end{equation}
By \eqref{eq:bK0-coefficients},
\begin{align}
	\bfb_{\tau(i),k}^{-}
	&=
	-\zeta^{k}
	k(q_i-q_i^{-1})
	\frac{1}{q^k-q^{-k}}
	\sum_{j\in\bfI}
	\qCC_{\tau(i)j}(q^k)d_{j,k}^{-}\\
	&=
	-k(q_i-q_i^{-1})
	\frac{1}{q^k-q^{-k}}
	\sum_{j\in\bfI}
	\qCC_{i\tau(j)}(q^k)d_{\tau(j),k}^{+}\\
	&=
	\zeta^{k}\bfb^+_{i,k},
\end{align}
where the second equality follows from \eqref{eq:d-bK-unitarity} and 
$\qCC_{ij}(q^k)=\qCC_{\tau(i)\tau(j)}(q^k)$. The result follows.
\end{pf}

\begin{remark}
The definitions \eqref{eq:log-GF-pm}, \eqref{eq:def-d-exp}, yielding the definition \eqref{eq:bK0-coefficients} of the coefficients $\bfb_{i,k}^\pm\in\bbC$, are motivated by Lemma~\ref{lem:b-gamma}.
\hfill\rmkend
\end{remark}

\subsection{The series $\Upsilon_{\bfa}^\pm$ and $\Upsilon_{\bfb}^\pm$}\label{ss:Upsilon-a-b}
Set
\begin{align}
\Upsilon_{\bfa}^\pm(s)
&=
\mp\frac{1}{2}
\sum_{k>0}s^{\pm2k}
\sum_{i,j\in\bfI}
\bfa_{ij,k}\zeta^{\mp k}
h_{\tau(i),\pm k}h_{j,\pm k},
\label{eq:bK0-Upsilon-a}\\
\Upsilon_{\bfb}^\pm(s)
&=
\sum_{k>0}s^{\pm k}
\sum_{i\in\bfI}\bfb_{i,k}^{\pm}h_{i,\pm k}
\label{eq:bK0-Upsilon-b},
\end{align}
where $\bfa_{ij,k}$ is defined as in \eqref{eq:def-coeff-a}.
These series admit a difference-equation description analogous to
\eqref{eq:R-diff-eq}. Let $\phalf$ denote the
half--shift $\phalf\cdot f(s)=f(q^{1/2}s)$, and set
\begin{align}
\mathcal Q_{ij}^{\pm}(s)
&=
-\frac{1}{2}\sum_{k>0}
k(q_i-q_i^{-1})(q_j-q_j^{-1})
\zeta^{\mp k}
h_{\tau(i),\pm k}h_{j,\pm k}s^{\pm2 k},
\label{eq:bK0-source-Q}\\
\mathcal P_{ij}^{\pm}(s)
&=
\mp\sum_{k>0}
k(q_i-q_i^{-1})d_{j,k}^{\pm}
\zeta^{\mp k}h_{i,\pm k}s^{\pm k}.
\label{eq:bK0-source-P}
\end{align}

Let $\bfU^{0,\pm} \subset \bfU^0$ be the subalgebras generated by $\{h_{i,\pm k}\}_{i\in\bfI,k>0}$, respectively.

\begin{lemma}\label{prop:bK0-difference-equation}
Assume that $\det(\qCB(q^k))\neq0$ for every
$k\in\bbZ_{\neq0}$.
\begin{enumerate}[leftmargin=2em]\itemsep0.25cm
	\item 
	The series $\Upsilon_{\bfa}^\pm(s)$ are the unique formal solutions in ${s^{\pm2}}\fext{\bfU^{0,\pm}}{s^{\pm2}}$ of the difference equations
	\begin{equation}\label{eq:bK0-form-diff-eq-a}
		\begin{aligned}
			\scrD(\phalf)\cdot\Upsilon_{\bfa}^{\pm}(s)
			&=
			\sum_{i,j\in\bfI}
			\qCBd^*_{ij}(\phalf)\cdot
			\mathcal Q_{ij}^{\pm}(s).
		\end{aligned}
	\end{equation}
	\item 
	The series $\Upsilon_{\bfb}^\pm(s)$ are the unique formal solutions in ${s^{\pm1}}\fext{\bfU^{0,\pm}}{s^{\pm1}}$ of the difference equations
	\begin{equation}\label{eq:bK0-form-diff-eq-b}
		\begin{aligned}
			\scrD(\sfp)\cdot\Upsilon_{\bfb}^{\pm}(s)
			&=
			\sum_{i,j\in\bfI}
			\qCBd^*_{ij}(\sfp)\cdot
			\mathcal P_{ij}^{\pm}(s).
		\end{aligned}
	\end{equation}
\end{enumerate}
\end{lemma}

\begin{pf}
The uniqueness of the formal solutions follows from the invertibility of
$\qCBd(q^{\pm k})$ for $k>0$. Then, it is enough to observe that,
since $\qCB(q^{-k})=\qCB(q^k)$, one has
\begin{equation}\label{eq:bK0-adjugate-ratio}
	\frac{\qCBd^*_{ij}(q^{\pm k})}
	{\det(\qCBd(q^{\pm k}))}
	=
	\pm\frac{1}{q^k-q^{-k}}\qCC_{ij}(q^k).
\end{equation}
Thus, the formal solutions of \eqref{eq:bK0-form-diff-eq-a} and  \eqref{eq:bK0-form-diff-eq-b} are readily given by \eqref{eq:bK0-Upsilon-a} and \eqref{eq:bK0-Upsilon-b}, respectively.
\end{pf}

\subsection{The abelian basic K-matrix}\label{ss:bK0-form-thm}
Fix an element $\hc^-\in\h$ such that
\begin{equation}\label{eq:bK0-H-minus}
q^{-\alpha_i(\hc^-)}=-q_i^{-a_{i\tau(i)}/2}\pga{i}=-q_i^{-a_{i\tau(i)}/2}\GF_i(\infty)
\end{equation}
for every $i\in\bfI$, and set
$\hc^+=-\bpsi(\hc^-)=\tau(\hc^-)$. Note that
\begin{equation}\label{eq:bK0-H-plus}
q^{-\alpha_i(\hc^+)}
=-q_i^{-a_{i\tau(i)}/2}\pga{\tau(i)}
=-q_i^{a_{i\tau(i)}/2}\GF_i(0)
\end{equation}
for every $i\in\bfI$. This choice is not canonical since it requires choosing a
logarithm of $\pga{i}\in\bbC^\times$ ($i\in\bfI$).

\begin{theorem}\label{thm:basic-K}
Set $\Upsilon^\pm(s)=\Upsilon^\pm_{\bfa}(s)+\Upsilon_{\bfb}^{\pm}(s)$ where $\Upsilon_{\bfa}^\pm(s), \Upsilon_{\bfb}^\pm(s)$ are the solutions of \eqref{eq:bK0-form-diff-eq-a} and \eqref{eq:bK0-form-diff-eq-b}, respectively. 
Then the operators
\begin{equation}\label{eq:bK0-form}
	\bKK^{0,\pm}(s)
	=
	q^{\hc^\pm}\bK{\pm\theta}{}\exp(\Upsilon^{\pm}(s))
\end{equation}
induce two meromorphic basic K-matrices on the category $\sO_\zeta$
related by twisted unitarity.\footnote{
By definition, $q^{\hc^\pm}\bK{\pm\theta}{}=q^{\hc^\pm\pm C_0^\theta/2}$.}
Specifically, for every $V\in\sO_\zeta$, the action of $\bKK^{0,\pm}(s)$ on $V$
defines two meromorphic $\End(V)$-valued functions
$\bKK_{V}^{0,\pm}(s)$ with the following properties:
\begin{enumerate}[font=\upshape]\itemsep0.25cm
	\item\label{bK0:int} For every $V\in\sO_\zeta$,
	\begin{equation}\label{eq:bK0-intertwining}
		\bKK^{0,\pm}_{V}(s) : V(s) \to
		V^{\bpsi}(s^{-1})
	\end{equation}
	is a $\bfB$-intertwiner.
	\item\label{bK0:shift} For every $V\in\sO_\zeta$
	and $a\in\bbC^\times$, one has
	\begin{equation}\label{eq:bK0-shift}
		\bKK^{0,\pm}_{V(a)}(s) = \bKK^{0,\pm}_{V}(sa).
	\end{equation}
	\item\label{bK0:coproduct} For every $V_1,V_2\in\sO_\zeta$, one has
	\begin{align}
		\label{eq:bK0-coproduct}
		\bKK^{0,\pm}_{V_1\DT{s_2/s_1}V_2}(s_1) &=
		\bKK^{0,\pm}_{V_1}(s_1)\ten\bKK^{0,\pm}_{V_2}(s_2)
		\RR^{0,\pm}_{V_1^\bpsi,V_2}(s_1s_2).
	\end{align}
	\item\label{bK0:unitary} For every $V\in\sO_\zeta$, one has
	\begin{equation}\label{eq:bK0-unitary}
		\bKK^{0,+}_{V}(s)^{-1} =
		\bKK^{0,-}_{V^{\bpsi}}(s^{-1}).
	\end{equation}
	\item\label{bK0:rat} For every highest weight representation $V\in\sO^\bbQ_{\zeta}$, there is a meromorphic function $g^{\pm}_{V}(s)$ such that the normalized operator
	\begin{equation}\label{eq:bK0-norm}
		\bKK^{0,\norm}_{V}(s)=g^{\pm}_{V}(s)^{-1}\bKK^{0,\pm}_{V}(s)
	\end{equation}
	is independent of $\pm$, rational and unitary in $s$.
\end{enumerate}
\end{theorem}

\begin{remarks}
\hfill
\begin{enumerate}[leftmargin=2em]\itemsep0.25cm
	\item 
	The property \eqref{eq:bK0-shift} follows readily from \eqref{eq:bK0-form}.
	\item 
	The existence and the properties of the normalization \eqref{eq:bK0-norm} follow 
	from Theorem~\ref{thm:generic-submodules-cyclicity-restriction} (3) and \eqref{eq:bK0-intertwining}. Namely, since $V\in\sO_\zeta^\bbQ$, by Theorem~\ref{thm:generic-submodules-cyclicity-restriction} (3) the
	restriction of $V$ to $\bfB$ remains generically cyclic over the highest weight vector $v_0\in V$. Hence, by \eqref{eq:bK0-intertwining}, 
	the morphisms $\bKK^{0,\pm}_V(s)$ are determined by their value on $v_0$.
	By \eqref{eq:bK0-form}, $\bKK^{0,\pm}_V(s)$ must preserve $v_0$, \ie there exist
	two meromorphic functions $g_V^{\pm}(s)$ such that $\bKK^{0,\pm}_V(s)v_0=g_V^{\pm}(s)v_0$. Hence, $g_V^{\pm}(s)^{-1}\bKK^{0,\pm}_V(s)$ yields the unique rational unitary $\bfB$-intertwiner $\bKK^{0,\norm}_{V}(s)$ determined by the condition $\bKK^{0,\norm}_{V}(s)v_0=v_0$.
	\hfill\rmkend
\end{enumerate}
\end{remarks}

The proof is carried out in \ref{ss:bK0-meromorphic}--\ref{ss:bK0-intertwining-main}.
Specifically, the meromorphicity of the operators $\bKK^{0,\pm}_V(s)$ is proved in 
\ref{ss:bK0-meromorphic}. Then, we observe in \ref{ss:bK0-ansatz} that the properties \eqref{eq:bK0-shift} and \eqref{eq:bK0-coproduct} force $\Upsilon^\pm(s)=\Upsilon_{\bfa}^\pm(s)+\Upsilon_{\bfb}^\pm(s)$ for some $\bfb^\pm_{i,k}\in\bbC$. Finally, we prove in \ref{ss:bK0-intertwining}--\ref{ss:bK0-intertwining-main} that the coefficients $\bfb^\pm_{i,k}$
defined in \eqref{eq:bK0-coefficients} are in fact determined by the intertwining property \eqref{eq:bK0-intertwining}. The unitarity property \eqref{eq:bK0-unitary} is proved in \ref{ss:bK0-ansatz}.

\subsection{Meromorphicity}\label{ss:bK0-meromorphic}

We first prove the analytic part of Theorem~\ref{thm:basic-K}.

\begin{proposition}\label{prop:bK0-meromorphic}
For every $V\in\sO$ and $\lambda\in P(V)$, the action of the formal series $\exp(\Upsilon^\pm(s))$ on $V_\lambda$ is given by the 
expansions at $s=0$ and $s=\infty$, respectively, of
$\End(V_\lambda)$-valued meromorphic functions.
\end{proposition}

The proof is carried out in \ref{sss:complex-analysis-lemma}--\ref{sss:pf-mero-2}.

\subsubsection{}\label{sss:complex-analysis-lemma} We shall need the following general results proved in \cite{gautam-toledano-17,appel-gautam-26}.

\begin{lemma}\label{lem:bK0-complex-analysis}
Let $W$ be a finite-dimensional vector space and $\scrA\subseteq\End(W)$ a commutative subalgebra. 
\begin{enumerate}[leftmargin=2em]\itemsep0.25cm
	\item \label{lem:bK0-meromorphic-summation}
	Let $\mathcal S^\pm(t)\in t^{\pm1}\fext{\scrA}{t^\pm1}$ be
	such that $\exp(\mathcal S^\pm(t))$ is the expansion at $t=0$,
	respectively at $t=\infty$, of a meromorphic function which is
	holomorphic there and takes the value $1$.
	
	Fix $q\in\bbC^\times$ with $|q|\neq1$ and let $\sfp$ be the corresponding multiplicative difference operator. Let $a\in\bbZ_{>0}$ and write
	$\D(\sfp)=\sum_{r=m}^{M}c_r\sfp^r\in{\bbZ}[\sfp,\sfp^{-1}]$, where
	$c_m,c_M\in\{\pm1\}$. Assume that
	$\D(q^{\pm ak})\neq0$ for every $k>0$.
	Let
	$\F^\pm(t)\in{t^{\pm1}}\fext{\scrA}{t^{\pm1}}$ be the unique formal solution of
	\begin{equation}\label{eq:bK0-general-difference-equation}
		\D(\sfp^a)\cdot \F^\pm(t)=\mathcal S^\pm(t).
	\end{equation}
	Then $\F^\pm(t)$ is convergent at $t=0$,
	respectively at $t=\infty$, and
	$\exp(\F^\pm(t))$ extends to a meromorphic function on
	$\bbC^\times$.
	
	\item \label{lem:bK0-logarithmic-pairing}
	Let $F, G$ be $\mathscr{A}^\times$-valued rational functions, with $G$ holomorphic and equal to $1$ at zero or infinity. Let $d$ be the greatest common divisor of the nonzero orders of the zeroes and poles of the eigenvalues of $F(w)$ on $W$, with the convention $d=1$ if all these eigenvalues are constant, and let $\C$ be a clockwise Jordan curve enclosing all poles of $F^{\pm 1}$ in $\bbC^\times$. Then, for any $m\in\bbZ$ and $s$ sufficiently close to $0$ or $\infty$, respectively, the operator-valued integral
	\[
	\exp\left(
	-\frac{m}{d}\oint_{\C}
	F(w)^{-1}F'(w)\log G(sw^{-1})\,dw
	\right)
	\]
	defines an analytic function of $s$ which extends to an
	$\End(W)$-valued rational function of $s$.
\end{enumerate}
\end{lemma}

\begin{pf}
Part (1) is the one-parameter specialization of the result proved in \cite{appel-gautam-26}.
Part (2) is proved in \cite[Sec.~8.6]{gautam-toledano-17} for $m=d$.
The same proof applies for arbitrary $m$: the orders of the scalar
divisors are divisible by $d$, while the unipotent contribution is the
exponential of a rational nilpotent function.
\end{pf}

\subsubsection{}\label{sss:pf-mero-1}
We shall further refine the series $\Upsilon^\pm(s)=\Upsilon_{\bfa}^\pm(s)+\Upsilon_{\bfb}^\pm(s)$, relying on Lemma~\ref{sss:complex-analysis-lemma} to prove the meromorphicity of its exponential.
Set $\bfE_i^\pm(z)=\bfD_i^\pm(z)-\bfG_i^\pm(z)$.
By \eqref{eq:log-GF-pm} and \eqref{eq:def-d-exp}, one has
\begin{equation}\label{eq:bK0-E-rational}
\exp(2\bfE_i^\pm(z))
=
\frac{1-q_i^{\mp a_{i\tau(i)}}\zeta^{\pm1}z^{\pm2}}
{1-q_i^{\pm a_{i\tau(i)}}\zeta^{\pm1}z^{\pm2}}.
\end{equation}
In particular, $\bfE_i^\pm(z)$ is even. Accordingly, we consider the linear decomposition 
\[
\P_{ij}^\pm
=
\P_{ij,\bfG}^\pm+\P_{ij,\bfE}^\pm
\qquad\text{and}\qquad
\Upsilon_{\bfb}^\pm
=
\Upsilon_{\bfb,\bfG}^\pm+\Upsilon_{\bfb,\bfE}^\pm.
\]
Then, \eqref{eq:bK0-form-diff-eq-b} gives
\begin{equation}\label{eq:bK0-source-G}
\scrD(\sfp)\Upsilon_{\bfb,\bfG}^\pm(s)
=
\mathcal S_{\bfG}^\pm(s),
\end{equation}
where $\mathcal S_{\bfG}^\pm(s)=
\sum_{i,j}\qCBd^*_{ij}(\sfp)\P_{ij,\bfG}^\pm(s)$.\\

The series $\Upsilon_{\bfa}^\pm$ and
$\Upsilon_{\bfb,\bfE}^\pm$ are even in $s$. 
We use a hat to denote the corresponding
series in $t^{\pm1}=s^{\pm2}$. Hence,
the equations \eqref{eq:bK0-form-diff-eq-a} and
\eqref{eq:bK0-form-diff-eq-b} read
\begin{align}
\scrD(\sfp)\widehat\Upsilon_{\bfa}^\pm(t)
&=
\sum_{i,j\in\bfI}\qCBd^*_{ij}(\sfp)\widehat\Q_{ij}^\pm(t),
\label{eq:bK0-even-diff-a}\\
\scrD(\sfp^2)\widehat\Upsilon_{\bfb,\bfE}^\pm(t)
&=
\sum_{i,j\in\bfI}\qCBd^*_{ij}(\sfp^2)
\widehat\P_{ij,\bfE}^\pm(t).
\label{eq:bK0-even-diff-b}
\end{align}
Consequently, their sum satisfies
\begin{equation}\label{eq:bK0-even-diff}
\scrD(\sfp)\scrD(\sfp^2)
\left(
\widehat\Upsilon_{\bfa}^\pm(t)
+\widehat\Upsilon_{\bfb,\bfE}^\pm(t)
\right)
= \mathcal S_{\bfQ,\bfE}^\pm(t),
\end{equation}
where
\begin{equation}
\mathcal S_{\bfQ,\bfE}^\pm(t)
=
\scrD(\sfp^2)\sum_{i,j\in\bfI}\qCBd^*_{ij}(\sfp)
\widehat\Q_{ij}^\pm(t)
+
\scrD(\sfp)\sum_{i,j\in\bfI}\qCBd^*_{ij}(\sfp^2)
\widehat\P_{ij,\bfE}^\pm(t).
\end{equation}
{Set $\widetilde{\scrD}(\sfp)=\scrD(\sfp)\scrD(\sfp^2)$.
This Laurent polynomial has integer coefficients and extreme
coefficients equal to $1$ by Condition~\ref{cond:qCartan-dominant}.
Moreover, $\widetilde{\scrD}(q^{\pm k})\neq0$ for every $k>0$
by Condition~\ref{cond:qCartan-non-vanishing}.
It is therefore enough to prove that}
\begin{equation}\label{eq:bK0-meromorphic-sources}
\exp(\mathcal S_{\bfG}^\pm(s))
\qquad\text{and}\qquad
\exp(\mathcal S_{\bfQ,\bfE}^\pm(t))
\end{equation}
are meromorphic on $V_\lambda$, holomorphic at zero and infinity, and equal to $1$ there. Indeed,
Lemma~\ref{lem:bK0-complex-analysis}\eqref{lem:bK0-meromorphic-summation},
applied directly to \eqref{eq:bK0-even-diff} with
$\D=\widetilde{\scrD}$ and $a=1$, gives the meromorphicity of
$\exp(\wh\Upsilon_{\bfa}^\pm(t)
+\wh\Upsilon_{\bfb,\bfE}^\pm(t))$, while its
application to \eqref{eq:bK0-source-G} with $a=1$ gives the
meromorphicity of $\exp(\Upsilon_{\bfb,\bfG}^\pm(s))$.

\subsubsection{}\label{sss:pf-mero-2}
Following \cite[Sec.~8.6]{gautam-toledano-17}, see also
\cite{appel-gautam-wendlandt-26,appel-gautam-26}, we use logarithmic
pairings of rational functions. Let $\C_{i,\lambda}$ be a clockwise
Jordan curve which encloses the poles of
$\Phif_i(w)^{\pm1}$ on $V_\lambda$ and does not enclose $0$.
For $s$ sufficiently close to $0$, respectively to $\infty$, the
second factor in each of the following integrals is holomorphic on and
inside the contour. Then,
on $V_\lambda$,
\begin{align}
\Q_{ij}^\pm(s)
&=
-\frac12\oint_{\C_{\tau(i),\lambda}}
\Phif_{\tau(i)}(w)^{-1}\Phif_{\tau(i)}'(w)
\bfH_j^\pm(s^2\zeta^{-1}w^{-1})\,dw,
\label{eq:bK0-contour-Q}\\
\P_{ij}^\pm(s)
&=
-\oint_{\C_{i,\lambda}}
\Phif_i(w)^{-1}\Phif_i'(w)
\bfD_j^\pm(s\zeta^{-1}w^{-1})\,dw,
\label{eq:bK0-contour-P}
\end{align}
where $\bfH^\pm_i(z)$ denotes the action on $V_\lambda$ of the series $h^\pm_i(z)$ defined in \eqref{eq:log-phi}.
Indeed, for the expansion at zero the integrals extract the residue at
$w=0$. For the expansion at infinity, the clockwise contour extracts
the residue at infinity. Using \eqref{eq:log-phi}, the two computations
give, respectively, \eqref{eq:bK0-source-Q} and
\eqref{eq:bK0-source-P}.

Since $\exp(\bfG_i^\pm(z))=\GF_i^\pm(z)$ is rational,
Lemma~\ref{lem:bK0-complex-analysis}\eqref{lem:bK0-logarithmic-pairing},
applied with $m=d$, shows that
$\exp(\P_{ij,\bfG}^\pm(s))$ is rational. More generally, an integral
multiple of a shifted pairing is covered by taking the corresponding
integral multiple of $d$ for $m$. Since the coefficients of
$\qCBd^*_{ij}(\sfp)$ are integers, it follows that
$\exp(\mathcal S_{\bfG}^\pm(s))$ is rational, holomorphic at the
corresponding expansion point, and equal to $1$ there.

Both $\exp(\bfH_i^\pm(z))$ and
$\exp(2\bfE_i^\pm(z))$ are rational, by the defining identity for
$\bfH_i^\pm$ and \eqref{eq:bK0-E-rational}. The same lemma, again with
$m=d$, therefore shows that $2\Q_{ij}^\pm$ and
$2\P_{ij,\bfE}^\pm$ have rational exponentials. Since all the
coefficients occurring in \eqref{eq:bK0-even-diff} are integers, the
function
\begin{equation}\label{eq:bK0-doubled-source}
R^\pm(t)=\exp\left(2\mathcal S_{\bfQ,\bfE}^\pm(t)\right)
\end{equation}
is rational, holomorphic at the corresponding expansion point, and
equal to $1$ there.\\

We claim that the normalized local square root
$\exp(\mathcal S_{\bfQ,\bfE}^\pm(t))$ of $R^\pm(t)$ is rational. We
use the following elementary criterion. On a joint generalized
eigenspace for a commutative algebra, write a rational operator as
$r(t)(1+N(t))$, where $r(t)$ is its scalar rational eigenvalue and $N(t)$ is
nilpotent, see \cite[Lem. 4.12]{gautam-toledano-16}. Its normalized local square root is rational provided every
zero and pole of $r(t)$ has even order. Indeed, the scalar square root
is then rational, and
\[
(1+N(t))^{1/2}
=
\sum_{k\ges0}\binom{1/2}{k}N(t)^k
\]
is a finite sum. Hence, it remains to verify such evenness condition for \eqref{eq:bK0-doubled-source}.\\

Set
\begin{equation}\label{eq:bK0-parity-data}
\begin{aligned}
	c_i(\sfp)&=\qCBd_{i,\tau(i)}(\sfp),\\
	M(\sfp)&=\scrD(\sfp^2)P_\tau\qCBd^*(\sfp),\\
	\ell_i(\sfp)&=\scrD(\sfp)
	\sum_j\qCBd^*_{ij}(\sfp^2)c_j(\sfp),
\end{aligned}
\end{equation}
where $P_\tau=(\delta_{i,\tau(j)})_{i,j\in\bfI}$ is the permutation
matrix of the diagram involution $\tau$.
The Laurent polynomial matrix $M(\sfp)$ is symmetric and integral.
Fix a joint generalized eigenspace, let $\varphi_i(z)$ be the scalar
eigenvalue of $\Phif_i(z)$, and let $r^\pm(t)$ be that of $R^\pm(t)$.
Write the coefficients in
\eqref{eq:bK0-parity-data} as
\[
M_{ij}(\sfp)=\sum_r M_{ij,r}\sfp^r,
\qquad
\ell_i(\sfp)=\sum_r\ell_{i,r}\sfp^r,
\]
and set $n_{i,a}=\operatorname{ord}_{z=a}\varphi_i(z)$, $a\in\bbC^\times$.
Only finitely many of these integers are nonzero. By the
logarithmic pairing formula,\footnote{For nonzero scalar rational
		functions $f,g$, put $n_a=\operatorname{ord}_{w=a}f(w)$. If $\C$ is
		a Jordan curve, enclosing the zeroes and poles of $f$ in $\bbC^\times$,
		excluding $0$, and $\log g(s/w)$ is holomorphic on and inside $\C$,
		the residue theorem gives
		\[
		\exp\left(-\oint_{\C}\frac{f'(w)}{f(w)}\log g(s/w)\,dw\right)
		=\prod_a g(s/a)^{n_a},
		\]
		where the product runs over the enclosed zeroes and poles.} 
for any $\alpha\in\bbC^\times$,
the orders at $t=\alpha$ of the scalar eigenvalues of
$\exp(2\widehat\Q_{ij}^\pm(t))$ and
$\exp(2\widehat\P_{ij,\bfE}^\pm(t))$ are, respectively,
\[
\sum_{\zeta ab=\alpha}n_{\tau(i),a}n_{j,b}
\qquad\text{and}\qquad
\sum_{\zeta a^2q^{\pm d_ja_{j\tau(j)}}=\alpha}\pm n_{i,a}.
\]
The shift $\sfp^r$ moves a divisor point $\alpha$ to $q^{-r}\alpha$.
Consequently, for every $t_0\in\bbC^\times$,
\begin{equation}\label{eq:bK0-source-order-parity}
	\begin{aligned}
		\operatorname{ord}_{t=t_0}r^\pm(t)
		&=\left(\sum_{i,j,r}M_{ij,r}
		\sum_{\substack{a,b\in\bbC^\times\\\zeta ab=q^rt_0}}
		n_{i,a}n_{j,b}\right)-\left(\sum_{i,r}\ell_{i,r}
		\sum_{\substack{a\in\bbC^\times\\\zeta a^2=q^rt_0}}n_{i,a}\right).
	\end{aligned}
\end{equation}
In characteristic two, the diagonal of $M(\sfp)$ and the linear term are related by
\begin{equation}\label{eq:bK0-characteristic-congruence}
\ell_i(\sfp)
\equiv
M_{ii}(\sfp)
=
\scrD(\sfp^2)\qCBd^*_{\tau(i),i}(\sfp)
\pmod 2.
\end{equation}
Indeed, $\qCBd(\sfp)P_\tau$ is symmetric, where $P_\tau$ is the
permutation matrix of $\tau$, and for every invertible symmetric matrix
$B$ over a field of characteristic $2$ one has
\[
(B^{-1})_{ii}
=
\sum_jB_{jj}(B^{-1}_{ij})^2.
\]
By Condition~\ref{cond:qCartan-dominant}, the reduction of
$\scrD(\sfp)$ modulo $2$ is nonzero. Thus the reduction of
$B=\qCBd(\sfp)P_\tau$ is invertible over $\mathbb F_2(\sfp)$.
Applying this identity to $B=\qCBd(\sfp)P_\tau$ and the diagonal entry
indexed by $\tau(i)$, clearing
denominators, and using $f(\sfp)^2=f(\sfp^2)$ in characteristic $2$
gives \eqref{eq:bK0-characteristic-congruence}.
By symmetry, the terms of \eqref{eq:bK0-source-order-parity}
indexed by distinct pairs $(i,a)\neq(j,b)$ cancel modulo $2$.
Together with \eqref{eq:bK0-characteristic-congruence}, this gives
\[
\operatorname{ord}_{t=t_0}r^\pm(t)
\equiv
\sum_{i,r}M_{ii,r}
\sum_{\substack{a\in\bbC^\times\\\zeta a^2=q^rt_0}}
(n_{i,a}^2-n_{i,a})
\equiv0\pmod2.
\]
Thus every scalar eigenvalue of $R^\pm(t)$ has zeroes and poles of even
order, and the square-root criterion above proves that
$\exp(\mathcal S_{\bfQ,\bfE}^\pm(t))$ is rational.\\

Thus both functions in \eqref{eq:bK0-meromorphic-sources} are
meromorphic, holomorphic at the corresponding expansion point, and
equal to $1$ there. The reduction in \ref{sss:pf-mero-1} gives the
meromorphicity of
\[
\exp(\Upsilon_{\bfb,\bfG}^\pm(s))
\qquad\text{and}\qquad
\exp\left(
\Upsilon_{\bfa}^\pm(s)+
\Upsilon_{\bfb,\bfE}^\pm(s)
\right).
\]
All the operators involved commute, and the product of these two
functions is $\exp(\Upsilon^\pm(s))$. This completes the proof of
Proposition~\ref{prop:bK0-meromorphic}.\hfill\qed

\subsection{Algebraic ansatz}\label{ss:bK0-ansatz}
The expression of $\Upsilon^\pm(s)$ is essentially dictated by the required properties of the abelian K-matrix.

\begin{lemma}\label{lem:bK0-algebraic-ansatz}
Let $\fkB^\pm(s)\in s^{\pm1}\fext{\bfU^{0}}{s^{\pm1}}$ be such that the
operator 
\[\bKK^{0,\pm}_\fkB(s)=q^{\hc^\pm}\bK{\pm\theta}{}\exp(\fkB^\pm(s))\] 
satisfies
\eqref{eq:bK0-shift} and \eqref{eq:bK0-coproduct}. Then, 
$\fkB^\pm(s)=\Upsilon^\pm_\bfa(s)+\Upsilon^\pm_\bfb(s)$, where 
$\Upsilon^\pm_\bfb(s)$ is defined as in \eqref{eq:bK0-Upsilon-b}
for some $\bfb_{i,k}^\pm\in\bbC$. Moreover, $\bKK^{0,\pm}_\fkB(s)$ satisfies the twisted unitarity property \eqref{eq:bK0-unitary} if and only if 
\begin{equation}\label{eq:bK0-b-unitarity}
\bfb_{i,k}^{+}=\zeta^{-k}\bfb_{\tau(i),k}^{-}.
\end{equation}
\end{lemma}

\begin{pf}
First, we observe that the properties \eqref{eq:bK0-shift} and \eqref{eq:bK0-coproduct} force the condition 
$\fkB^\pm(s)\in\fext{\bfU^{0,\pm}}{s^{\pm1}}$.
The factor $q^{\hc^\pm}$ is group-like and therefore does not
affect the coproduct equation \eqref{eq:bK0-coproduct}. Hence, by
\eqref{eq:coproduct-operator-k}, equation \eqref{eq:bK0-coproduct} reduces
to
\begin{equation}\label{eq:bK0-coboundary-Ups}
	\DD{s_2/s_1}(\fkB^{\pm}(s_1))
	=
	\fkB^{\pm}(s_1)\ten 1
	+1\ten\fkB^{\pm}(s_2)
	+(\bpsi\ten\id)(\Lambda^{\pm}(s_1s_2)).
\end{equation}
A direct calculation, using \eqref{eq:log-R}, shows that
$\Upsilon_{\bfa}^{\pm}(s)$ is a particular solution of \eqref{eq:bK0-coboundary-Ups}. 
It follows that $\fkB^\pm(s)-\Upsilon_{\bfa}^\pm(s)$ must be primitive, and therefore, 
by \eqref{eq:bK0-shift}, of the form
\[
\sum_{k>0}s^{\pm k}
\sum_{i\in\bfI}\bfb_{i,k}^{\pm}h_{i,\pm k}.
\]
for some $\bfb_{i,k}^{\pm}\in\bbC$, \ie $\fkB^\pm(s)=\Upsilon_{\bfa}^\pm(s)+\Upsilon_{\bfb}^\pm(s)$ for some $\bfb_{i,k}^{\pm}\in\bbC$.
Finally, one readily verifies that $-\Upsilon_{\bfa}^+(s)=\bpsi(\Upsilon_{\bfa}^-(s^{-1}))$. Hence, $\bKK^{0,\pm}_\fkB(s)$ satisfies \eqref{eq:bK0-unitary} if and only if \eqref{eq:bK0-b-unitarity} holds.
\end{pf}

\begin{corollary}\label{cor:bK0-algebraic-ansatz}
The operator $\bKK^{0,\pm}(s)$ satisfies the twisted unitarity property \eqref{eq:bK0-unitary}.	
\end{corollary}
\begin{pf}
It is enough to observe that, by Lemma~\ref{lem:bK0-coeff-b}, the coefficients $\bfb_{i,k}^\pm\in\bbC$ defined  as in \eqref{eq:bK0-coefficients} satisfy the condition  \eqref{eq:bK0-b-unitarity}, \ie $\bfb_{i,k}^{+}=\zeta^{-k}\bfb_{\tau(i),k}^{-}$. 
Hence, the result follows by Lemma~\ref{lem:bK0-algebraic-ansatz}.
\end{pf}

\subsection{Decoupling the twisted intertwining equation}\label{ss:bK0-intertwining}
We shall now prove that \eqref{eq:bK0-intertwining} uniquely determines the coefficients $\bfb_{i,k}^\pm$ and consequently $\bKK^{0,\pm}(s)$. We start by studying the decoupling of \eqref{eq:bK0-intertwining}, its consistency, and its compatibility with \eqref{eq:bK0-unitary}.

For every $\ell\in\bfI$ and $r\in\bbZ$, set
\begin{equation}
\begin{aligned}
	\fkX_{\ell,r}(s)&=\sum_{a,b\ges 0}
	\pg_{\ell, -a}x^+_{\tau(\ell),\,r+a+b}\phi^+_{\ell,b}
	\,\zeta^{-b}s^{r+a+2b},\\
	\fkY_{\ell,r}(s)&=\sum_{a\ges 0}\ps_{\ell,r+a}\phi^+_{\ell,a}\zeta^{-a}s^a,\\
	\fkX^{\bpsi}_{\ell,r}(s)&=-\sum_{a,b\ges 0}
	\pg_{\ell,-a}x^-_{\ell,\,-r-a-b}\phi^-_{\tau(\ell),-b}
	\,\zeta^{r+a+b}s^{-r-a-2b},\\
	\fkY^\bpsi_{\ell,r}(s)&=\sum_{a\ges 0}\ps_{\ell,r+a} \phi^-_{\tau(\ell),-a}s^{-a}.
\end{aligned}
\end{equation}
Note that $\fkX_{\ell,r}^{\bpsi}(s)=\bpsi(\fkX_{\ell,r}(s^{-1}))$ and $\fkY_{\ell,r}^{\bpsi}(s)=\bpsi(\fkY_{\ell,r}(s^{-1}))$.
By \eqref{eq:B-generators}, we have
\begin{equation}
\begin{aligned}
	\SH_{s}(y_{\ell,r})&=
	x^-_{\ell,-r}\zeta^r s^{-r}
	+\fkX_{\ell,r}(s)
	+\fkY_{\ell,r}(s),\\
	\bpsi(\SH_{s^{-1}}(y_{\ell,r}))
	&=-x^+_{\tau(\ell),r}s^{r}
	+\fkX_{\ell,r}^{\bpsi}(s)
	+\fkY_{\ell,r}^{\bpsi}(s).
\end{aligned}
\end{equation}
By \eqref{eq:bK0-form}, conjugation by $\bKK^{0,\pm}(s)$ is weight preserving, since $\Upsilon^{\pm}(s)\in\fext{\bfU^{0,\pm}}{s^{\pm1}}$.
Hence, \eqref{eq:bK0-intertwining} is equivalent to the following intertwining equations:
\begin{align}
\label{eq:int-AdbK-1}
\Ad(\bKK^{0,\pm}(s))(x^-_{\ell,-r}\zeta^r s^{-r})&=\fkX^{\bpsi}_{\ell,r}(s),\\
\label{eq:int-AdbK-2}
\Ad(\bKK^{0,\pm}(s))(\fkX_{\ell,r}(s))&=-x^+_{\tau(\ell),r}s^{r},\\
\label{eq:int-AdbK-3}
\Ad(\bKK^{0,\pm}(s))(\fkY_{\ell,r}(s))&=\fkY^{\bpsi}_{\ell,r}(s).
\end{align}	

\begin{remark}
We readily observe that the intertwining equation \eqref{eq:bK0-intertwining}
makes sense only at the level of category $\sO$ representations. For instance, \eqref{eq:int-AdbK-3} reduces to the identity $\fkY_{\ell,r}(s)=\fkY^{\bpsi}_{\ell,r}(s)$, where $\fkY_{\ell,r}(s)\in\fext{\bfU^{0}}{s}$ while $\fkY^{\bpsi}_{\ell,r}(s)\in\fext{\bfU^{0}}{s^{-1}}$. However, we shall prove that
the operator $\bKK^{0,\pm}(s)$ does satisfy the identities
\eqref{eq:int-AdbK-1}--\eqref{eq:int-AdbK-2}--\eqref{eq:int-AdbK-3}
when evaluated on a category $\sO$ representation.
\hfill\rmkend
\end{remark}

\begin{lemma}\label{lem:bK0-intertwining}
\hfill
\begin{enumerate}[leftmargin=2em]\itemsep0.25cm
	\item For every $V\in\sO_\zeta$ and $\lambda\in P(V)$, one has
	\begin{equation}
	\fkY_{\ell,r}(s)\mid_{V_\lambda}
	=\fkY^{\bpsi}_{\ell,r}(s)\mid_{V_\lambda}
	\end{equation} 
	for every $\ell\in\bfI$ and $r\in\bbZ$.
	\item 
	Assume $\bKK^{0,+}(s)$ satisfies \eqref{eq:bK0-unitary}, \ie $\bKK^{0,+}(s)^{-1}=\bpsi(\bKK^{0,-}(s^{-1}))$. Then \eqref{eq:int-AdbK-1} is equivalent to \eqref{eq:int-AdbK-2}.
\end{enumerate}
\end{lemma}

\begin{pf}
(1) First, we observe that the identity is trivial whenever $\tau(\ell)\neq\ell$, since
$\ps_{\ell, r+a}=0$ and therefore $\fkY_{\ell,r}(s)=0=\fkY^{\bpsi}_{\ell,r}(s)$. Assume now $\tau(\ell)=\ell$. Then,
\begin{align}
	\fkY_{\ell,r}(s)&=\sum_{a\ges0}\ps_{\ell,r+a}\phi^+_{\ell,a}\zeta^{-a}s^a\\
	&=\oint_{\La_{\ell,\lambda}^\SF}\SF_\ell(w)\left(\sum_{a\ges0}\phi_{\ell,a}^+\zeta^{-a}s^aw^{-a}\right)\frac{dw}{w^{r+1}}\\
	&=\oint_{\La_{\ell,\lambda}^\SF}\SF_\ell(w)\Phif_\ell(\zeta^{-1}sw^{-1})\frac{dw}{w^{r+1}}
\end{align}
and
\begin{align}
	\fkY^{\bpsi}_{\ell,r}(s)&=\sum_{a\ges0}\ps_{\ell,r+a}\phi^-_{\ell,-a}s^{-a}\\
	&=\oint_{\La_{\ell,\lambda}^\SF}\SF_\ell(w)\left(\sum_{a\ges0}\phi_{\ell,-a}^-s^{-a}w^{-a}\right)\frac{dw}{w^{r+1}}\\
	&=\oint_{\La_{\ell,\lambda}^\SF}\SF_\ell(w)\Phif_\ell(sw)\frac{dw}{w^{r+1}}.
\end{align}
We claim that the last two contour integrals coincide. It is enough to prove
this for generic $s$, since both sides are meromorphic functions of $s$.
For such an $s$, we may choose $\La_{\ell,\lambda}^\SF$ to enclose the
possible poles of $\SF_\ell(w)$ and none of the poles of
$\Phif_\ell(\zeta^{-1}sw^{-1})$ or $\Phif_\ell(sw)$. Every pole $p$ of $\SF_\ell(w)$
satisfies $p^2=\zeta^{-1}$, \ie $p=\zeta^{-1}p^{-1}$. Hence, the operator
\[\Phif_\ell(s\zeta^{-1}w^{-1})-\Phif_\ell(sw)\]
vanishes at every such $p$. More explicitly, it can be written locally as
\[
\Phif_\ell(s\zeta^{-1}w^{-1})-\Phif_\ell(sw)
=
(\zeta^{-1}w^{-1}-w)F_s(w)
=
\frac{1-\zeta w^2}{\zeta w}F_s(w),
\]
where $F_s(w)$ is holomorphic at $p$. Therefore, by \eqref{eq:SF},
\[
\SF_\ell(w)
\left(\Phif_\ell(\zeta^{-1}sw^{-1})-\Phif_\ell(sw)\right)
\frac{1}{w^{r+1}}
\]
is holomorphic inside $\La_{\ell,\lambda}^\SF$, and its contour integral
vanishes. Hence, for generic
$s$, $\fkY_{\ell,r}(s)=\fkY^{\bpsi}_{\ell,r}(s)$. The result follows by meromorphic continuation.\\

\noindent
(2) We prove that \eqref{eq:int-AdbK-1} implies \eqref{eq:int-AdbK-2}. One has
\begin{align}
	\Ad(\bKK^{0,\pm}(s))(\fkX_{\ell,r}(s))&=\bpsi(\Ad(\bKK^{0,\mp}(s^{-1}))^{-1}(\fkX_{\ell,r}^{\bpsi}(s^{-1})))\\
	&=\bpsi(x_{\ell,-r}^-\zeta^rs^{r})\\
	&=-x^+_{\tau(\ell),r}s^r,
\end{align}
where the first equality follows from \eqref{eq:bK0-unitary}, the identity $\fkX_{\ell,r}^{\bpsi}(s^{-1})=\bpsi(\fkX_{\ell,r}(s))$, and $\bpsi^2=\id$, while the second one follows from \eqref{eq:int-AdbK-1}. The converse is similar.
\end{pf}

\subsection{The twisted intertwining equation}
\label{ss:bK0-intertwining-main}
By Lemma~\ref{lem:bK0-intertwining}, it remains to prove \eqref{eq:int-AdbK-1}. 
Together with Proposition~\ref{prop:bK0-meromorphic},
and Corollary~\ref{cor:bK0-algebraic-ansatz}, this will complete the proof of
Theorem~\ref{thm:basic-K}.

\begin{proposition}
Let $\bfb_{i,k}^\pm\in\bbC$ be defined as in \eqref{eq:bK0-coefficients}.
For every $V\in\sO_\zeta$, the operators $\bKK^{0,\pm}_V(s)$ satisfy \eqref{eq:int-AdbK-1}.
\end{proposition}

The proof is carried out in \ref{sss:bK0-int-1}--\ref{sss:bK0-int-4}.
Fix $\lambda\in P(V)$, and $\ell\in\bfI$.

\subsubsection{Setting}\label{sss:bK0-int-1}
We shall prove the result for generic $s\in\bbC^\times$ and then extend it meromorphically.
For generic $s$, choose a Jordan curve $\C_{\ell,\lambda}(s)$ which
encloses $s\pole_\ell^\zeta(V_\lambda)$ and none of the points in
\begin{equation}
\{0\}\cup\pole_{\ell,\GF}\cup
s^{-1}\pole_{\tau(\ell)}(V_\lambda).
\end{equation}
Thus $\C_{\ell,\lambda}(s)$ encloses the poles of
$\Xf_\ell^-(s\zeta^{-1}w^{-1})$, while $\GF_\ell(w)$ and
$\Phif_{\tau(\ell)}(sw)$ are holomorphic inside it. The contour formulae
described in \ref{ss:QL-contour}, followed by the change of variable
$u=s\zeta^{-1}w^{-1}$, give
\begin{equation}\label{eq:bK0-contour-negative-mode}
x^-_{\ell,-r}\zeta^rs^{-r}
=
-\oint_{\C_{\ell,\lambda}(s)}
\Xf_\ell^-(s\zeta^{-1}w^{-1})\frac{dw}{w^{r+1}}
\end{equation}
on $V_\lambda$.
Indeed, $du=-u\,dw/w$, while the holomorphic inversion preserves the
clockwise orientation of a contour which does not enclose $0$. Thus the
minus sign from the differential remains.
Similarly, relying on the definitions of $\pg_\ell(z)$ and
$\phi^-_{\tau(\ell)}(z)$ as the expansions at infinity of
$\GF_\ell(z)$ and $\Phif_{\tau(\ell)}(z)$, respectively, one has
\begin{equation}\label{eq:bK0-contour-Xpsi}
\fkX_{\ell,r}^{\bpsi}(s)
=
\oint_{\C_{\ell,\lambda}(s)}
\GF_\ell(w)
\Xf_\ell^-(\zeta^{-1}sw^{-1})
\Phif_{\tau(\ell)}(sw)
\frac{dw}{w^{r+1}}
\end{equation}
on $V_\lambda$.
Indeed, the identity \eqref{eq:bK0-contour-negative-mode}, with $r$ replaced by
$r+a+b$, shows that the expansion of the integrand in
\eqref{eq:bK0-contour-Xpsi} yields
\begin{equation}
-\sum_{a,b\ges0}
\pg_{\ell,-a}
x^-_{\ell,-r-a-b}
\phi^-_{\tau(\ell),-b}
\zeta^{r+a+b}s^{-r-a-2b},
\end{equation}
which is precisely $\fkX_{\ell,r}^{\bpsi}(s)$.

\subsubsection{Commutation relation for $h_{i,k}$}\label{sss:bK0-int-2}
For rational operator-valued functions of $w$, we write
$F(w)\equiv_\C G(w)$ if $F(w)-G(w)$ is holomorphic inside
$\C_{\ell,\lambda}(s)$. Since $0$ lies outside $\C_{\ell,\lambda}(s)$, this
implies
\begin{equation}\label{eq:Gamma-equivalence-moments}
\oint_{\C_{\ell,\lambda}(s)} F(w)\frac{dw}{w^{r+1}}
=
\oint_{\C_{\ell,\lambda}(s)} G(w)\frac{dw}{w^{r+1}}
\end{equation}
for every $r\in\bbZ$. Set $u=s\zeta^{-1}w^{-1}$. Then, we observe that the commutation relation 
\eqref{eq:adjoint-h} implies
\begin{equation}\label{eq:h-X-principal-part}
[h_{i,\pm k},\Xf^-_\ell(u)]
\equiv_{\C}
-c_{i\ell,k}u^{\mp k}\Xf^-_\ell(u)
\end{equation}
where
\begin{equation}
c_{i\ell,k}=\frac{[ka_{i\ell}]_{q_i}}{k}
=\frac{1}{k}\frac{q^k-q^{-k}}{q_i-q_i^{-1}}\qCB_{i\ell}(q^k).
\end{equation}
Indeed, by expanding at $u=0$ and
$u=\infty$, respectively, and then applying \eqref{eq:adjoint-h}, one has
\begin{align}
[h_{i,k},\Xf^-_\ell(u)]
&=
-c_{i\ell,k}u^{-k}\Xf^-_\ell(u)
+c_{i\ell,k}\sum_{m=0}^{k-1}x^-_{\ell,m}u^{m-k},
\label{eq:h-X-boundary-plus}\\
[h_{i,-k},\Xf^-_\ell(u)]
&=
-c_{i\ell,k}u^k\Xf^-_\ell(u)
-c_{i\ell,k}\sum_{m=-k}^{-1}x^-_{\ell,m}u^{m+k}.
\label{eq:h-X-boundary-minus}
\end{align}
Since the Laurent polynomial components are holomorphic inside $\C_{\ell,\lambda}(s)$,
this yields \eqref{eq:h-X-principal-part}.

\subsubsection{Commutation relation for $\Upsilon^\pm(s)$}\label{sss:bK0-int-3}

From \eqref{eq:h-X-principal-part}, we readily obtain the commutation relations of $\Upsilon^\pm_{\bfa}(s)$ and $\Upsilon^\pm_{\bfb}(s)$ with $\Xf^-_\ell(s\zeta^{-1}w^{-1})$.

\begin{lemma}\label{lem:comm-Upsilon}
One has
\begin{align}
[\Upsilon_{\bfa}^\pm(s),\Xf^-_\ell(s\zeta^{-1}w^{-1})]
&\equiv_{\C}\Xf^-_\ell(s\zeta^{-1}w^{-1})\cdot
\Big(
\pm(q_\ell-q_\ell^{-1})
\sum_{k>0}h_{\tau(\ell),\pm k}(sw)^{\pm k}\\
&\hspace{2.5cm}\mp\frac{1}{2}
\sum_{k>0}
\frac{q_\ell^{ka_{\ell\tau(\ell)}}
-q_\ell^{-ka_{\ell\tau(\ell)}}}{k}
\zeta^{\pm k}w^{\pm2k}\Big),\\
[\Upsilon_{\bfb}^\pm(s),\Xf^-_\ell(s\zeta^{-1}w^{-1})]
&\equiv_{\C}\Xf^-_\ell(s\zeta^{-1}w^{-1})\cdot\Big(
-\sum_{k>0}\zeta^{\pm k}w^{\pm k}
\sum_{i\in\bfI}c_{i\ell,k}\bfb_{i,k}^\pm
\Big).
\end{align}
\end{lemma}

\begin{pf}
The identity for $\Upsilon_{\bfb}^\pm(s)$ is obvious. In the case of $\Upsilon_{\bfa}^\pm(s)$, 
it is enough to observe that, by \eqref{eq:def-coeff-a} and since $\qCC(q^k)\qCB(q^k)=1$, 
\begin{equation}\label{eq:basic-contraction-a-c}
\sum_{j\in\bfI}\bfa_{ij,k}c_{j\ell,k}
=(q_i-q_i^{-1})\delta_{i\ell},
\end{equation}
hence
\begin{align}
\sum_{i,j\in\bfI}
\bfa_{ij,k}c_{j\ell,k}h_{\tau(i),\pm k}
&=(q_\ell-q_\ell^{-1})h_{\tau(\ell),\pm k},
\label{eq:contraction-a-c-1}\\
\sum_{i,j\in\bfI}
\bfa_{ij,k}c_{\tau(i)\ell,k}h_{j,\pm k}
&=(q_\ell-q_\ell^{-1})h_{\tau(\ell),\pm k},
\label{eq:contraction-a-c-2}\\
\sum_{i,j\in\bfI}
\bfa_{ij,k}c_{\tau(i)\ell,k}c_{j\ell,k}
&=
(q_\ell-q_\ell^{-1})\frac{[ka_{\ell\tau(\ell)}]_{q_\ell}}{k}.
\label{eq:contraction-a-c-3}
\end{align}
Finally, since 
\[\zeta^{\mp k}s^{\pm 2k}u^{\mp k}
=(sw)^{\pm k}\aand\zeta^{\mp k}s^{\pm2k}u^{\mp2k}
=\zeta^{\pm k}w^{\pm 2k},\] 
the result follows.
\end{pf}

\subsubsection{Commutation relation for $\bKK^{0,\pm}(s)$}\label{sss:bK0-int-4}
Set
\[
\GF_\ell^+=\GF_\ell(0)\aand \GF_\ell^-=\GF_\ell(\infty).
\]
Then, we have
\begin{equation}\label{eq:bK0-zero-mode-GF}
\begin{aligned}
\Ad(q^{\hc^\pm}\bK{\pm\theta}{})
(\Xf^-_\ell(u))
&\equiv_{\C}
q^{-\alpha_\ell(\hc^\pm)}
q_\ell^{\mp a_{\ell\tau(\ell)}/2}
\Xf^-_\ell(u)\phi^\pm_{\tau(\ell),0}\\
&\equiv_{\C}
-\GF_\ell^\pm\Xf^-_\ell(u)\phi^\pm_{\tau(\ell),0},
\end{aligned}
\end{equation}
where the first equality follows from \eqref{eq:adjoint-operator-k} and the second one from the definition of $\hc^\pm$ in \eqref{eq:bK0-H-minus}--\eqref{eq:bK0-H-plus}.

From Lemma~\ref{lem:comm-Upsilon} and \eqref{eq:bK0-zero-mode-GF}, it follows that
\begin{equation}\label{eq:adjoint-K-partial}
\Ad(\bKK^{0,\pm}(s))(\Xf^-_\ell(u))
\equiv_\C
-\GF_{\ell}^{\pm}
\exp(\bfg_\ell^\pm(w))
\Xf^-_\ell(u)\Phif_{\tau(\ell)}(sw),
\end{equation}
where
\begin{equation}\label{eq:def-G-bK}
\bfg_\ell^\pm(w)
=
\mp\frac{1}{2}
\sum_{k>0}
\frac{q_\ell^{ka_{\ell\tau(\ell)}}
-q_\ell^{-ka_{\ell\tau(\ell)}}}{k}
\zeta^{\pm k}w^{\pm2 k}
-
\sum_{k>0}\zeta^{\pm k}w^{\pm k}
\sum_{i\in\bfI}c_{i\ell,k}\bfb_{i,k}^\pm.
\end{equation}

\begin{lemma}\label{lem:b-gamma}
The series $\GF_{\ell}^{\pm}\exp(\bfg_\ell^\pm(w))$ are, respectively, the expansions at $w=0$ and $w=\infty$ of the rational function $\GF_{\ell}(w)$.
\end{lemma}

\begin{pf}
The proof relies on the explicit definition of the coefficients $\bfb_{i,k}^\pm$
given in \eqref{eq:bK0-coefficients}. Since
$\qCB(q^k)\qCC(q^k)=1$, one has
\begin{equation}\label{eq:linear-system-bK-b}
\sum_{i\in\bfI}c_{i\ell,k}\bfb_{i,k}^\pm
=
-\zeta^{\mp k}d_{\ell,k}^\pm
\end{equation}
for every $\ell\in\bfI$ and $k>0$. Hence, by \eqref{eq:def-d-bK}
\[\bfg^\pm_\ell(w)=\sum_{k>0}g_{\ell,k}^\pm w^{\pm k},\]
where $g^\pm_{\ell, k}\in\bbC$ are defined as in \eqref{eq:log-GF-pm}.
The result follows.
\end{pf}

Finally, we prove \eqref{eq:int-AdbK-1}. Namely, on $V_\lambda$, we have
\begin{align}
\Ad(\bKK_V^{0,\pm}(s))
(	x^-_{\ell,-r}\zeta^rs^{-r})
&=
-\oint_{\C_{\ell,\lambda}(s)}
\Ad(\bKK^{0,\pm}(s))
\left(\Xf^-_\ell(\zeta^{-1}sw^{-1})\right)
\frac{dw}{w^{r+1}}\\
&=
\oint_{\C_{\ell,\lambda}(s)}
\GF_\ell(w)\Xf^-_\ell(\zeta^{-1}sw^{-1})
\Phif_{\tau(\ell)}(sw)\frac{dw}{w^{r+1}}\\
&=
\fkX^{\bpsi}_{\ell,r}(s),
\end{align}
where the first equality follows from \eqref{eq:bK0-contour-negative-mode}, the second one from \eqref{eq:adjoint-K-partial} together with Lemma~\ref{lem:b-gamma}, and the third one from \eqref{eq:bK0-contour-Xpsi}.

\section{The tensor abelian K-matrix and boundary $q$-characters}
\label{s:tensor-abelian-K-matrix}
In this section, we show that the construction of the abelian K-matrix from the previous section naturally gives rise to a tensor abelian K-matrix twisted-intertwining the action of $\bfU^\imath$ on any tensor product $M\ten V$, with $M\in\sM_\zeta$
\footnote{More precisely, $M$ is assumed to be a subquotient of a restricted representation in $\sO_\zeta$, see \ref{eq:tensor-K-thm}.} and $V\in\sO_\zeta$. 
The resulting theory of boundary $q$-characters recovers and extends, within a common framework, earlier constructions in \cite{appel-vlaar-25c} and \cite{Przez-23, li-przezdziecki-25, LP-eval, li-przezdziecki-26}.

\subsection{The series $\Xi_{\bfa}^\pm$}\label{eq:operator-Xi}
Fix $\zeta\in\bbC^\times$. Let $\ps_i(z),\pg_i(z)$ ($i\in\bfI$) be 
defined as in \eqref{eq:SF+GF}, together with
\eqref{eq:serre-orthogonal-parameters} and \eqref{eq:serre-sigma-parity}. 
We consider the resulting complex parameters and rational functions as in Section~\ref{ss:cat-O-zeta}.
Let $\veps=\veps_{\ps,\pg}$ be the
corresponding character \eqref{eq:character-normalization}, and set
$\mathbf c=(\pga{i})_{i\in\bfI}$. We shall consider the central block
$\sM_{\zeta,\mathbf c}$ and the normalized operator
\begin{equation}\label{eq:eta-normalization}
\weta_i(z)=\veps(\wetad_i(z))^{-1}\wetad_i(z),
\end{equation}
whose constant terms act on $\sM_{\zeta,\mathbf c}$ as
$\weta_{i,0}=q^{d_ih_i-d_{\tau(i)}h_{\tau(i)}}$. Define
\begin{equation}\label{eq:tensor-normalized-logarithms}
\weta_i(z)=\weta_{i,0}
\exp\left((q_i-q_i^{-1})\sum_{k>0}\wLeta_{i,k}z^k\right),
\end{equation}
see Remark~\ref{rmk:eta-normalization}. With the coefficients
$\bfa_{ij,k}$ from \eqref{eq:def-coeff-a}, set
\begin{equation}\label{eq:tK0-Xi-a}
\begin{aligned}
\Xi_{\bfa}^-(s)
&=-\sum_{k>0}s^{-k}\sum_{i,j\in\bfI}
\bfa_{ij,k}\wLeta_{i,k}\ten h_{j,-k},\\
\Xi_{\bfa}^+(s)
&=-\sum_{k>0}s^k\zeta^{-k}\sum_{i,j\in\bfI}
\bfa_{ij,k}\wLeta_{\tau(i),k}\ten h_{j,k}.
\end{aligned}
\end{equation}
These are formal series in $s^{-1}$ and $s$, respectively, with
coefficients in the tensor product of the Cartan subalgebras.

The two formulas reflect the two expansions of the quantum loop
Cartan currents. Indeed, \eqref{eq:xi} gives
\begin{equation}\label{eq:tensor-logarithm-realization}
\xi_\veps(\wLeta_{i,k})=h_{i,k}-\zeta^k h_{\tau(i),-k}.
\end{equation}
For the negative series, its negative is
$-(\id+\bpsi)(h_{i,k})$. For the positive series, replacing $i$ by
$\tau(i)$ and multiplying by $-\zeta^{-k}$ gives
$h_{i,-k}-\zeta^{-k}h_{\tau(i),k}=(\id+\bpsi)(h_{i,-k})$.
Here $\bpsi(h_{i,\pm k})=-\zeta^{\pm k}h_{\tau(i),\mp k}$.
Consequently, by \eqref{eq:log-R}, one has
\begin{equation}\label{eq:Xi-Lambda}
\begin{aligned}
(\xi_\veps\ten\id)(\Xi^\pm_{\bfa}(s))
&=((\id+\bpsi)\ten\id)(\Lambda^\pm(s))\\
&=\Lambda^\pm(s)-(\id\ten\bpsi)(\Lambda^\mp(s^{-1})).
\end{aligned}
\end{equation}
The second equality follows from
$(\bpsi\ten\bpsi)(\Lambda^\pm(s))=-\Lambda^\mp(s^{-1})$
and $\bpsi^2=\id$ on the Cartan subalgebra.

\subsection{The abelian tensor K-matrix}\label{eq:tensor-K-thm}
Let $\sM_{\zeta,\mathbf c}^{\mathrm{res}}$ be the full subcategory
of $\sM_{\zeta,\mathbf c}$ whose objects are isomorphic to subquotients
of modules obtained from $\sO_\zeta$ by restriction.
It follows from Theorem~\ref{thm:basic-K} and the passage from basic
to tensor K-matrices in \cite{kolb-20,appel-vlaar-25b,LBG-26}.
We shall consider the operator
\begin{equation}
\tkt^\pm=\tk^{\pm(\id+\theta)}\cdot \id\ten q^{\hc^\pm}\bK{\pm\theta}{},
\end{equation}
see \ref{ss:operator-C} and \ref{ss:Casimir}. Namely, for every $V,W\in\sO$, $\lambda\in P(V)$, $\mu\in P(W)$, $\tkt^\pm$ acts on $V_\lambda\ten W_\mu$ as
\[
\tkt^\pm_{V_\lambda W_\mu}=q^{\pm(\lambda,(\id+\theta)(\mu))+\mu(\hc^\pm)\pm(\theta(\mu),\mu)/2}.
\]
Note that, since every object in $\sM_{\zeta,\mathbf c}$ is semisimple over $\h^\theta$, $\tkt^\pm$
acts more generally on any tensor product of the form $M\ten V$ with $M\in\sM_{\zeta,\mathbf c}$, $V\in\sO_\zeta$.

\begin{theorem}\label{thm:tensor-K}
Assume Conditions~\ref{cond:qCartan-non-vanishing} and
\ref{cond:qCartan-dominant}, as in Theorem~\ref{thm:appel-gautam-26}.
Set $\Xi^\pm(s)=\Xi_{\bfa}^\pm(s)+1\ten\Upsilon^\pm(s)$ where 
$\Xi^\pm_\bfa(s)$ and $\Upsilon^\pm(s)=\Upsilon_{\bfa}^\pm(s)+\Upsilon_{\bfb}^\pm(s)$ are defined as in \eqref{eq:tK0-Xi-a} and  \eqref{eq:bK0-Upsilon-a}--\eqref{eq:bK0-Upsilon-b}, respectively. 
The operators
\begin{equation}\label{eq:tensor-K}
	\KK^{0,\pm}(s)={\tkt}^\pm\exp(\Xi^{\pm}(s)),
\end{equation}
induce two meromorphic braidings on the module category $\sM_{\zeta,\mathbf c}^{\mathrm{res}}$ related by twisted unitarity.
Specifically, for every $M\in\sM_{\zeta,\mathbf c}^{\mathrm{res}},V\in\sO_\zeta$, $\KK_{MV}^{0,\pm}(s)$
define two meromorphic $\End(M\otimes V)$--valued functions
with the following properties:
\begin{enumerate}[font=\upshape]\itemsep0.25cm
	\item\label{K0:int} For every $M\in\sM_{\zeta,\mathbf c}^{\mathrm{res}}, V\in\sO_\zeta$,
	\begin{equation}\label{eq:K0-intertwining}
	\KK^{0,\pm}_{M,V}(s) : M\DT{s} V \to
	M\DT{s^{-1}}V^{\bpsi}
	\end{equation}
	is a $\bfU^{\imath}$--intertwiner.
	\item\label{K0:shift} For every $M\in\sM_{\zeta,\mathbf c}^{\mathrm{res}},V\in\sO_\zeta$
	and $a\in\bbC^\times$, one has
	\begin{equation}\label{eq:K0-shift}
	\KK^{0,\pm}_{M,V(a)}(s) = \KK^{0,\pm}_{M,V}(sa).
	\end{equation}
	\item\label{K0:cabling} For every $M\in\sM_{\zeta,\mathbf c}^{\mathrm{res}}$, $V_1,V_2\in\sO_\zeta$, one has
	\begin{align}
		\label{eq:cabling-K-left}
		\KK^{0,\pm}_{M\DT{s_1}V_1,V_2}(s_2) &= \RR^{0,\mp}_{V_1,V_2^{\bpsi}}(1/(s_1s_2))^{-1}\KK^{0,\pm}_{M,V_2}(s_2)
		\RR^{0,\pm}_{V_1,V_2}(s_2/s_1), \\
		\label{eq:cabling-K-right}
		\KK^{0,\pm}_{M,V_1\DT{s_2/s_1}V_2}(s_1) &= \KK^{0,\pm}_{M,V_2}(s_2)
		\RR^{0,\pm}_{V_1^\bpsi,V_2}(s_1s_2)\KK^{0,\pm}_{M,V_1}(s_1).
	\end{align}
	\item\label{K0:unitary} For every $M\in\sM_{\zeta,\mathbf c}^{\mathrm{res}}$, $V\in\sO_\zeta$, one has
	\begin{equation}\label{eq:K0-unitarity}
	\KK^{0,+}_{M,V}(s)^{-1} =
	\KK^{0,-}_{M,V^{\bpsi}}(s^{-1}).
	\end{equation}
	\item\label{K0:rat} For all highest weight representations $V_1,V_2\in\sO^\bbQ_{\zeta}$, there are meromorphic functions $g^{\pm}_{V_1,V_2}(s)$ such that the normalized operator
	\begin{equation}\label{eq:K0-rational}
	\KK^{0,\norm}_{V_1,V_2}(s)=g^{\pm}_{V_1,V_2}(s)^{-1}\KK^{0,\pm}_{V_1,V_2}(s)
	\end{equation}
	is independent of $\pm$, rational and unitary in $s$.
\end{enumerate}
\end{theorem}

\begin{pf}
We first take $M$ to be the restriction of $V_1\in\sO_\zeta$.
The logarithmic identity \eqref{eq:Xi-Lambda} and the definition of
$\tkt^\pm$ give
\begin{equation}\label{eq:tensor-basic-factorization}
\KK^{0,\pm}_{V_1,V_2}(s)
=\RR^{0,\mp}_{V_1,V_2^{\bpsi}}(s^{-1})^{-1}
(1\ten\bKK^{0,\pm}_{V_2}(s))
\RR^{0,\pm}_{V_1,V_2}(s).
\end{equation}
Indeed, for the Cartan component, the two R-matrix factors contribute
$\tk^{\pm(\id+\theta)}$, while the basic K-matrix contributes
$1\ten q^{\hc^\pm}\bK{\pm\theta}{}$.

The right-hand side is the tensor K-matrix associated to the basic
K-matrix, see \cite{kolb-20,appel-vlaar-25b}. The R-matrix
intertwining relation and \eqref{eq:bK0-intertwining} therefore give
\eqref{eq:K0-intertwining}. The shift identities for R and the basic
K-matrix give \eqref{eq:K0-shift}. Their coproduct identities give
\eqref{eq:cabling-K-left}--\eqref{eq:cabling-K-right}, and their
unitarity identities give \eqref{eq:K0-unitarity}. Meromorphicity follows
from the factorization and Theorems~\ref{thm:appel-gautam-26}
and~\ref{thm:basic-K}.

Let $M=S/T$, with $T\subseteq S\subseteq V_1$ restricted $\bfU^\imath$-submodules. Since each component of the formal operator \eqref{eq:tensor-K} preserves $S\ten V_2$ and $T\ten V_2$, then so does the corresponding meromorphic operator.\footnote{
In particular, we rely on the condition \ref{cond:M4}, which is naturally satisfied by restricted modules.} 
Hence the operators descend to $M\ten V_2$, depend only on its $\bfU^\imath$-module structure, and satisfy all previous identities.

Finally, for highest weight modules $V_1,V_2\in\sO_\zeta^{\bbQ}$, normalize
the two R-matrix factors and the basic K-matrix in
\eqref{eq:tensor-basic-factorization}.\footnote{In this case, by normalization we means division by the scalar action on $v_1\ten v_2$, where $v_i$ is a highest weight vector of the original
module $V_i$. For the exponential factors, we use
\eqref{eq:Xi-Lambda}.} Their scalar normalizations
combine to give $g^\pm_{V_1,V_2}(s)$. The result follows.
\end{pf}

\begin{remarks}
\hfill
\begin{enumerate}[leftmargin=2em]\itemsep0.25cm
\item 
The formulae \eqref{eq:tensor-normalized-logarithms}--\eqref{eq:tensor-K}
are formal and therefore make sense on the full central block 
$\sM_{\zeta,\mathbf c}$.
However, in order to extend Theorem~\ref{thm:tensor-K} to $\sM_{\zeta,\mathbf c}$, one has to directly prove the twisted intertwining equation \eqref{eq:K0-intertwining} and the meromorphicity of the operator $\KK^{0,\pm}_{M,V}$ for $M\in\sM_{\zeta,\mathbf c}$ and $V\in\sO_\zeta$.
The proof relies on the same strategy used in Sections~\ref{s:abelian-R-matrix}--\ref{s:basic-abelian-K-matrix}.
\item
Let $M\in\sM_\zeta$ be a cyclic module generated by $m\in M$ and $V\in\sO_\zeta^\bbQ$ a highest weight module with highest weight vector $v\in V$.
Following the proof of Theorem~\ref{thm:generic-submodules-cyclicity-restriction}, we can readily show that
$M\DT{s}V$ is generated over $\bfU^\imath$ by $M\ten v$, although this argument does not establish generation by $m\ten v$. Hence, the normalization property \eqref{eq:K0-rational} can only be stated for restricted modules.
\hfill\rmkend
\end{enumerate}
\end{remarks}

\subsection{The reduced tensor K-matrix}
\label{ss:reduced-tensor-K0}
Following \cite[Sec.~6]{appel-vlaar-25c}, we consider the
{reduced abelian tensor K-matrix}
\begin{equation}\label{eq:reduced-tensor-K}
\dKK^{0,\pm}(s)
=
\bigl(1\ten\bKK^{0,\pm}(s)\bigr)^{-1}
\KK^{0,\pm}(s)
=
\tk^{\pm(\id+\theta)}
\exp\bigl(\Xi_{\bfa}^{\pm}(s)\bigr).
\end{equation}
Note that $\dKK^{0,\pm}(s)$ still depends on a choice of parameters, since
it relies on the normalized generators $\weta_i(z)$ and $\wLeta_i(z)$, see \eqref{eq:eta-normalization}. Moreover, for every $V_1,V_2\in\sO_\zeta$,
\begin{equation}
{\dKK}^{0,\pm}_{V_1(s_1),V_2}(s_2)
=
\RR^{0,\mp}_{V_1,V_2^\bpsi}(1/(s_1s_2))^{-1}
\RR^{0,\pm}_{V_1,V_2}(s_2/s_1),
\end{equation}
where $V_1(s_1)$ is regarded as a $\bfU^\imath$--module by restriction through
$\xi_{\ps,\pg}$.

\begin{proposition}\label{prop:reduced-tensor-K-cabling}
Let $M\in\sM_{\zeta,\mathbf c}^{\mathrm{res}}$ and
$V_1,V_2\in\sO_\zeta$.
One has
\begin{equation}\label{eq:reduced-tensor-K-R-factorization}
	{\dKK}^{0,\pm}_{V_1,V_2}(s)
	=
	\RR^{0,\mp}_{V_1,V_2^{\bpsi}}(s^{-1})^{-1}
	\RR^{0,\pm}_{V_1,V_2}(s),
\end{equation}
and
\begin{align}
	\label{eq:reduced-K-cabling-right}
	{\dKK}^{0,\pm}_{M,V_1\DT{s_2/s_1}V_2}(s_1)
	&=
	{\dKK}^{0,\pm}_{M,V_1}(s_1)
	{\dKK}^{0,\pm}_{M,V_2}(s_2),
	\\
	\label{eq:reduced-K-cabling-left}
	{\dKK}^{0,\pm}_{M\DT{s_1}V_1,V_2}(s_2)
	&=
	{\dKK}^{0,\pm}_{M,V_2}(s_2)
	{\dKK}^{0,\pm}_{V_1(s_1),V_2}(s_2),
\end{align}
where, on the left hand side of the first identity and the right hand side of the third identity, $V_1$ and $V_1(s_1)$ are
regarded as a $\bfU^\imath$--module by restriction through
$\xi_{\ps,\pg}$.
\end{proposition}

\begin{pf}
The first identity follows directly from
\eqref{eq:tensor-basic-factorization} and \eqref{eq:reduced-tensor-K}. 	
The second identity follows from \eqref{eq:cabling-K-right} and 
\eqref{eq:bK0-coproduct}.
Similarly, the third identity follows from
\eqref{eq:cabling-K-left} and
\[
{\dKK}^{0,\pm}_{V_1(s_1),V_2}(s_2)
=
\RR^{0,\mp}_{V_1,V_2^\bpsi}(1/(s_1s_2))^{-1}
\RR^{0,\pm}_{V_1,V_2}(s_2/s_1),
\]
see \eqref{eq:reduced-tensor-K-R-factorization}.
\end{pf}

In Sections \ref{ss:q-characters}--\ref{ss:boundary-left-trace}, we shall explain how the coproduct formulae \eqref{eq:reduced-K-cabling-left}--\eqref{eq:reduced-K-cabling-right}
reflect the multiplicative properties of the boundary $q$-characters associated to $\bfU^\imath$.

\subsection{q-characters and the abelian R-matrix \cite{frenkel-reshetikhin-98,hernandez-05}}
\label{ss:fusion}
\label{ss:q-characters}
Recall that an $\ell$-weight is a tuple $\omega=(\lambda,\Psi)$, where
$\lambda\in\h^*$ and $\Psi=(\Psi_i(z))_{i\in\bfI}$ consists of
rational functions regular at $0$ and $\infty$, with
$\Psi_i(0)=q_i^{\lambda(h_i)}$ and
$\Psi_i(\infty)=q_i^{-\lambda(h_i)}$.
Write $\Psi_i^\pm(z)=\sum_{r\geq0}\Psi_{i,\pm r}^\pm z^{\pm r}$
for their expansions at $0$ and $\infty$. 

For $V\in\sO$, let
$V_\omega\subseteq V_\lambda$ be the simultaneous generalized
eigenspace on which every $\phi_{i,\pm r}^\pm$ has eigenvalue
$\Psi_{i,\pm r}^\pm$. Note that, for every $u\in\bbC^\times$, the $\ell$-weights 
of $V(u)$ are given by $\omega(u)=(\lambda,(\Psi_i(uz))_{i\in\bfI})$.
The q-character of $V$ is
\begin{equation}\label{eq:q-character-definition}
\chi_q(V)=\sum_\omega\dim V_\omega\,[\omega],
\end{equation}
where $[\omega]$ is a formal symbol.\footnote{
The symbols multiply by $[\omega][\omega']=[\omega\omega']$, where
$\omega\omega'=(\lambda+\lambda',\Psi\Psi')$.
We allow sums with finitely many labels at each ordinary weight
and support in finitely many cones $\lambda_a-\sfQ_+$.
}\\

We regard $\omega$ also as a character of
$\bfU^0$, by $\omega(q^h)=q^{\lambda(h)}$ and
$\omega(\phi_{i,\pm r}^\pm)=\Psi_{i,\pm r}^\pm$.
Assume Condition~\ref{cond:qCartan-non-vanishing}, and set
\begin{equation}\label{eq:bulk-character-substitution}
\RR_\omega(s)=(\id\ten\omega)\RR^{0,-}(s)
=q^{-\tinv{\lambda}}
\exp\left(-\sum_{k>0}s^{-k}\sum_{i,j\in\bfI}
\bfa_{ij,k}\omega(h_{j,-k})h_{i,k}\right),
\end{equation}
where $\tinv{\lambda}=\nu^{-1}(\lambda)\in\h$ is defined by
$\iip{\tinv{\lambda}}{h}=\lambda(h)$, as in
Section~\ref{ss:sym-km-summary}. Note that the values $\omega(h_{j,-k})$
are determined by $\Psi_j^-(z)$.\\

Let $\cH^+$ be the commutative algebra generated by $\{q^h\}_{h\in\h}$ and $\{h_{i,k}\}_{i\in\bfI,k>0}$.
Let $\cE(\cH^+)$ denote the weight series
$\sum_\lambda a_\lambda e^\lambda$ ($a_\lambda\in\cH^+$) supported on finitely many cones, and consider the weighted trace
\begin{equation}\label{eq:bulk-weighted-trace}
\mathcal T_R(V;s)=\sum_\lambda e^\lambda
(\id\ten\Tr_{V_\lambda})\RR^{0,-}(s)
\in\fext{\mathcal E(\cH^+)}{s^{-1}}.
\end{equation}

\begin{proposition}\label{prop:bulk-trace-q-character}
For every $V\in\sO$, one has
\begin{equation}\label{eq:bulk-trace-character-identification}
	\mathcal T_R(V;s)=\chi_q(V)\vert_{[\omega]\mapsto e^\lambda\RR_\omega(s)}.
\end{equation}
Moreover, for every $V,W\in\sO$ and $u\in\bbC^\times$ such that 	
$V\DT{u}W$ is defined, one has
\begin{align}
\chi_q(V\DT{u}W)&=\chi_q(V)\chi_q(W(u)),
\label{eq:q-character-Drinfeld-product}\\
\mathcal T_R(V\DT{u}W;s)
&=\mathcal T_R(V;s)\mathcal T_R(W;su).
\label{eq:bulk-trace-Drinfeld-product}
\end{align}
\end{proposition}

\begin{pf}
The identification \eqref{eq:bulk-trace-character-identification}
follows from simultaneous triangularization on $V_\lambda$,
since the trace on $V_\omega$ is precisely $\dim V_\omega\,\RR_\omega(s)$.
The formula  \eqref{eq:q-character-Drinfeld-product} follows from 
\eqref{eq:Drinfeld-fields}, while 
\eqref{eq:bulk-trace-Drinfeld-product} follows from the right coproduct identity
satisfied by $\RR^{0,-}$, see Theorem~\ref{thm:appel-gautam-26} \eqref{R0:cabling}.
\end{pf}

\begin{remark}
By \cite[Sec.~5.1]{hernandez-07}, the deformed Drinfeld tensor product of 
integrable category $\O$ modules can be suitably specialized at $u=1$. This yields
a ring structure on the Grothendieck group for which the q-character $\chi_q$ is properly multiplicative.
\hfill\rmkend
\end{remark}

\subsection{Boundary q-characters of \texorpdfstring{$\bfU$}{U}-modules}
\label{ss:boundary-q-characters}
\label{ss:boundary-right-trace}

Fix the central block $\bfU^\imath_{\zeta,\mathbf c}$ and the
normalization \eqref{eq:eta-normalization}.
For every $V\in\sO$, its \emph{boundary $q$-character} is
\begin{equation}\label{eq:right-boundary-q-character}
\chi^\imath_{q,\bfU}(V;s)
=\sum_\lambda e^\lambda(\id\ten\Tr_{V_\lambda})\dKK^{0,-}(s)
\in\fext{\cE(\bfU^\imath_{\zeta,\mathbf c})}{s^{-1}},
\end{equation}
where $\cE(\bfU^\imath_{\zeta,\mathbf c})$ is defined in the same way as $\cE(\cH^+)$.

\begin{proposition}\label{prop:boundary-right-multiplicativity}
For every $V,W\in\sO$ and $u\in\bbC^\times$ such that 
$V\DT{u}W$ is defined,
\begin{equation}\label{eq:right-boundary-q-character-multiplicativity}
\chi^\imath_{q,\bfU}(V\DT{u}W;s)
=\chi^\imath_{q,\bfU}(V;s)
\chi^\imath_{q,\bfU}(W;su).
\end{equation}
Moreover, 
\begin{equation}
\chi^\imath_{q,\bfU}(V;s)
=\varphi(\chi_q(V)\vert_{[\omega]\mapsto e^\lambda\RR_\omega(s)})
\end{equation} 
where $ \varphi\colon\cE(\cH^+)\longrightarrow\cE(\bfU^\imath_{\zeta,\mathbf c})$ is defined by
\begin{equation}\label{eq:boundary-Cartan-substitution}
\varphi(q^h)= q^{h+\theta(h)},
\qquad \varphi(h_{i,k})=\wLeta_{i,k},
\qquad \varphi(e^\lambda)=e^{\lambda},
\end{equation}
and its extension to $\fext{\cE(\cH^+)}{s^{-1}}$ acts coefficientwise.
\end{proposition}

\begin{pf}
The factorization formula \eqref{eq:right-boundary-q-character-multiplicativity}
follows from the right coproduct identity
\eqref{eq:reduced-K-cabling-right}. Note indeed that \eqref{eq:reduced-K-cabling-right}
is formally satisfied by the operator \eqref{eq:reduced-tensor-K}. 
The map $\varphi$ is an algebra homomorphism, and, by \eqref{eq:log-R} and
\eqref{eq:reduced-tensor-K}, $\dKK^{0,-}=(\varphi\ten\id)\RR^{0,-}$.
The result follows.
\end{pf}

Hence the boundary $q$-character of a $\bfU$-module is just the usual $q$-character 
with the symbols $[\omega]$ replaced by
$e^\lambda(\id\ten\omega)\dKK^{0,-}(s)=e^\lambda\varphi(\RR_\omega(s))$.
This multiplicative substitution however need not be injective.

\begin{remark}
For the $q$-Onsager algebra embedded in $U_q(L\mathfrak{sl}_2)$, 
Baseilhac and Gainutdinov proposed in \cite[Sec.~7]{baseilhac-gainutdinov-26}
a refined version of the  Frenkel--Reshetikhin $q$-characters via a 
suitable Harish--Chandra projection
of the normalized boundary transfer matrix defined in \cite{appel-vlaar-25a}.
Comparing their proposal with \eqref{eq:right-boundary-q-character} would require,
in particular, proving that the standard and Drinfeld coactions on the
$q$-Onsager algebra are twist equivalent as we explain below.
\hfill\rmkend
\end{remark}

\subsection{Boundary q-characters of \texorpdfstring{$\bfU^\imath$}{U-iota}-modules}
\label{ss:boundary-left-trace}

Let $M\in\sM_{\zeta,\mathbf c}$ be a finite-dimensional $\bfU^\imath$-module. Then, its
\emph{boundary $q$-character} is
\begin{equation}\label{eq:left-boundary-q-character}
\chi^\imath_{q,\bfU^\imath}(M;s)=(\Tr_M\ten\id)\dKK^{0,-}(s)
\in\fext{\bfU^0}{s^{-1}}.
\end{equation}
As before, a boundary $\ell$-weight is a tuple $\varpi=(\mu,\Theta)$, where
$\mu\in(\h^\theta)^*$ and $\Theta=(\Theta_i(z))_{i\in\bfI}$ consists of
rational functions regular at $0$, with
$\Theta_i(0)=q_i^{\mu(h_i+\theta(h_i))}$.
Write $\Theta_i(z)=\sum_{r\geq0}\Theta_{i,r} z^{r}$
for their expansions at $0$.

For $M\in\sM_{\zeta,\mathbf c}$, let
$M_\varpi\subseteq M_\mu$ be the simultaneous generalized
eigenspace on which every $\weta_{i,r}\in\bfU^{\imath,0}$ has eigenvalue
$\Theta_{i,r}$. We regard $\varpi$ as a character of $\bfU^{\imath,0}$, by
$\varpi(q^{h})=q^{\mu(h)}$ ($h\in\h^\theta$) and $\varpi(\weta_{i,r})=\Theta_{i,r}$.
Then, we set
\begin{equation}\label{eq:boundary-ell-character}
\chi^\imath_q(M)=\sum_\varpi\dim M_\varpi[\varpi]
\aand K_\varpi(s)=(\varpi\ten\id)\dKK^{0,-}(s).
\end{equation}
Note that $(\varpi\ten\id)\tk^{-(\id+\theta)}=q^{-\tinv{(\id+\theta)\bar\mu}}$, see \eqref{eq:reduced-tensor-K}, hence it does not depend on the extension $\bar\mu\in\h^*$ of $\mu\in(\h^\theta)^*$.
Note also that $K_{\varpi\varpi'}=K_\varpi K_{\varpi'}$.

\begin{remark}
Up to the normalization \eqref{eq:eta-normalization}, $\chi^\imath_q(M)$
coincides with the boundary q-character introduced in \cite{li-przezdziecki-25, li-przezdziecki-26}.\hfill\rmkend 
\end{remark}

For an ordinary $\ell$-weight $\omega=(\lambda,\Psi)$, define its folding
by
\begin{equation}\label{eq:boundary-ell-weight-folding}
\mathfrak f_u(\omega)=
\left(\lambda|_{\h^\theta},
\bigl(\Psi_i^+(uz)
\Psi_{\tau(i)}^-(\zeta^{-1}uz^{-1})\bigr)_{i\in\bfI}\right).
\end{equation}
The product of series on the right is expanded at $z=0$.
Both substitutions
\begin{equation}\label{eq:boundary-folding-map}
[\omega]\mapsto[\mathfrak f_u(\omega)]
\aand
[\omega]\mapsto K_{\mathfrak f_u(\omega)}(s)
\end{equation}
are multiplicative on finite sums of $\ell$-weight symbols.

\begin{proposition}\label{prop:boundary-left-module-property}
For every finite-dimensional $M\in\sM_{\zeta,\mathbf c}$, one has
\begin{equation}\label{eq:left-boundary-character-weights}
\chi^\imath_{q,\bfU^\imath}(M;s)
=\chi^\imath_q(M)\vert_{[\varpi]\mapsto K_\varpi(s)}.
\end{equation}
Moreover, if $M\in\sM_{\zeta,\mathbf c}$ and $V\in\sO$ are finite-dimensional,
then, for every $u\in\bbC^\times$ such that $M\DT{u}V$ is defined,
\begin{align}
\chi^\imath_q(M\DT{u}V)
&=\chi^\imath_q(M)\, (\chi_q(V)\vert_{[\omega]\mapsto[\mathfrak f_u(\omega)]}),
\label{eq:boundary-character-Drinfeld-action}\\
\chi^\imath_{q,\bfU^\imath}(M\DT{u}V;s)
&=\chi^\imath_{q,\bfU^\imath}(M;s)\,
(\chi_q(V)\vert_{ [\omega]\mapsto K_{\mathfrak f_u(\omega)}(s)}).
\label{eq:left-boundary-q-character-module-property}
\end{align}
\end{proposition}

\begin{pf}
The formula \eqref{eq:left-boundary-character-weights} follows by
simultaneous triangularization on $M$. To prove \eqref{eq:boundary-character-Drinfeld-action},
it is enough to observe that, by \eqref{eq:DDi-weta-field}, one has
\[
\DDi{u}(\weta_i(z))=\weta_i(z)\ten
\phi_i^+(uz)\phi_{\tau(i)}^-(\zeta^{-1}uz^{-1}),
\]
whose generalized eigenvalues are of the form $\varpi\mathfrak f_u(\omega)$ with
their multiplicities  $\dim M_\varpi\dim V_\omega$.
Applying $[\varpi]\mapsto K_\varpi(s)$ immediately proves
\eqref{eq:left-boundary-q-character-module-property}, which is
precisely the trace of the left coproduct identity
\eqref{eq:reduced-K-cabling-left}.
\end{pf}

\begin{remarks}
\hfill
\begin{enumerate}[leftmargin=2em]\itemsep0.25cm
\item 
For infinite-dimensional modules, the same formulas hold
\emph{whenever they are well-defined}. This requires finite-dimensional
boundary generalized weight spaces, finite folded multiplicities, and locally
finite products of q-characters. One simply considers a generalized weighted trace, replacing
each $K_\varpi(s)$ by $[\varpi]K_\varpi(s)$.

\item
In order to interpret 
\eqref{eq:boundary-character-Drinfeld-action}--\eqref{eq:left-boundary-q-character-module-property} 
in terms of Grothendieck groups, one should introduce a \emph{fusion action} on $\bfU^\imath$-modules, compatible with Hernandez's fusion product. It is natural to expect that such an action should arise 
by a suitable specialization of the Drinfeld coaction at $u=1$.

\item \label{rmk:boundary-finite-type}
When $\g$ is finite-dimensional, and therefore $\bfU$ is the corresponding quantum affine algebra, 
the Drinfeld coproduct is rationally twist-equivalent to the usual Kac--Moody coproduct
\cite[Sec.~3]{enriquez-khoroshkin-pakuliak-07}.
On finite-dimensional modules, conjugation by the twist preserves
partial traces. Hence, the preceding identities readily imply
\[
\begin{aligned}
\chi_q(V\DT{\scsop{KM}}W)&=\chi_q(V)\chi_q(W),\\
\chi^\imath_{q,\bfU}(V\DT{\scsop{KM}}W;s)
&=\chi^\imath_{q,\bfU}(V;s)\chi^\imath_{q,\bfU}(W;s).
\end{aligned}
\]
The corresponding argument for boundary $q$-characters of $\bfU^\imath$-modules would require proving that
the Drinfeld coaction is rationally twist-equivalent with the standard coaction associated to the
quantum affine symmetric pair.
\hfill\rmkend
\end{enumerate}
\end{remarks}

\section{Shifted generalizations and GKLO representations}\label{s:shifted-generalizations}
\newcommand{\DDs}{\underset{s}{\Delta}{}_{\mu_1,\mu_2}}
\newcommand{\diff}{\mathcal{A}}
\newcommand{\sw}{\mathsf{w}}
\newcommand{\sz}{\mathsf{z}}
\newcommand{\arxiv}[1]{\href{http://arxiv.org/abs/#1}{\tt arXiv:\nolinkurl{#1}}}

In this section, we extend the Drinfeld comodule structure to \emph{shifted} 
$\imath$quantum  loop algebras,
yielding a surprisingly simple new construction of $\imath$GKLO representations.\\ 

We shall need some additional notations complementing those from
Section~\ref{ss:sym-km-summary}. 
Let $\Lambda$ and $\Lambda^\vee$ denote the weight and coweight lattices, respectively, with fundamental weights $\omega_i$ and coweights $\omega_i^\vee$. We endow $\Lambda^\vee$ with the standard partial order and denote the set of dominant coweights by $\Lambda^\vee_+$. 

\subsection{Shifted quantum loop algebras} 

Let $\mu \in \Lambda^\vee$ be \emph{any} coweight and set $\mu_i = \alpha_i(\mu)$. 
Consider the following generating series 
\[
x_i^{\pm}(z)=\sum_{k\in\bbZ}x_{i,k}^{\pm}z^k,
\qquad
\phi_i^{+}(z)=\sum_{k\ges 0}\phi_{i,k}^{+}z^k,
\qquad
\phi_i^{-}(z)=\sum_{k\les \mu_i}
\phi_{i,k}^{-}z^{k}. 
\] 

Following \cite{FT19}, the $\mu$-shifted quantum loop algebra $\bfU_\mu=U_q(L\g)_\mu$ is the unital associative $\bbC$-algebra generated by the elements
\[
\{\phi^+_{i, k}\}_{i\in\bfI,\ k\in\bbZ_{\geqslant0}},
\qquad
\{\phi^-_{i,k}\}_{i\in\bfI,\ k \leqslant \mu_i},
\qquad
\{x^\pm_{i,k}\}_{i\in\bfI,\ k\in\bbZ},
\qquad
\{ (\phi_{i, 0}^+)^{-1}, (\phi_{i, \mu_i}^-)^{-1}   \}_{i\in\bfI}
\]
subject to relations \eqref{eq:commutative-rel-field}--\eqref{eq:serre-relations} and, additionally, 
\[
(\phi_{i,0}^+)^{-1}\phi_{i,0}^+ = \phi_{i,0}^+(\phi_{i,0}^+)^{-1}=1, \qquad 
(\phi_{i,\mu_i}^-)^{-1}\phi_{i,\mu_i}^- = \phi_{i,\mu_i}^-(\phi_{i,\mu_i}^-)^{-1}=1. 
\]

\begin{remark}
	At the level of currents, the relations are literally the same as in the unshifted case, but 
	the shift in the range of $\phi_i^{-}(z)$ radically affects the structure of the algebra. 
	The definition above corresponds to the algebra which is denoted in \cite[Sec.~5]{FT19} 
	as $\mathcal{U}_{0,\mu}^{sc}$. 
	While the original definition in \emph{loc.\ cit.} contains two shifts $\mu_-, \mu_+$, it is 
	well known that $\mathcal{U}_{\mu_+,\mu_-}^{sc} \cong \mathcal{U}_{0,\mu_+ + 
		\mu_-}^{sc}$. Our series conventions are opposite to those in \cite{FT19}, \ie they 
	correspond under $z\mapsto z^{-1}$, but match the ones in \cite{hernandez-05}. 
	Note also that, in contrast with Section~\ref{s:QL}, we do not include the operators $q^h$ ($h\in\h$) among the generators of $\bfU_\mu$.
	\hfill\rmkend
\end{remark}

It was observed in \cite[Sec.~10(i)]{FT19} that the construction of the Drinfeld coproduct straightforwardly generalizes to this setting.\footnote{Since the Drinfeld--Serre relations involve at most two vertices, the argument in \cite[Sec.~10(i)]{FT19} applies whenever $a_{ij}a_{ji}\les3$ for $i\neq j$. Using \cite[Thm.~6.16]{damiani-24}, the same construction extends to arbitrary symmetrizable generalized Cartan matrices.}
More precisely, given any $\mu_1, \mu_2 \in \Lambda^\vee$, there exists an algebra homomorphism
\[
\DDs \colon \bfU_{\mu_1 + \mu_2} \to \Lext{\bfU_{\mu_1} \otimes \bfU_{\mu_2}}{s^{-1};s}, 
\]
given by the formulae \eqref{eq:Drinfeld-fields} as in the unshifted case. 

\subsection{GKLO representations} 

Assume that $\sfA$ is a symmetrizable generalized Cartan matrix of finite or affine type other than $\sfA_1^{(1)}$ and $\sfA_2^{(2)}$.  
Let $\Gamma$ be the graph associated to the Cartan matrix 
$\sfA$ with the edge multiplicities forgotten, \ie $\Gamma$ has vertex set $\bfI$ and an edge between $i$ and $j$ if and only if $a_{ij}<0$.
Fix an orientation of $\Gamma$ and write $i \sim j$ whenever $i$ and $j$ are connected 
by an edge, and $i \to j$ or $j \to i$ if there is an arrow in the given direction.

Fix a \emph{dominant} coweight $\lambda \in \Lambda^\vee_+$ with $\lambda \ges \mu$ and write 
\[
\lambda - \mu = \sum_{i\in \mathbf{I}} m_i  h_{i}, 
\]
where $m_i = \omega_i(\lambda - \mu)$.
In particular, for $\lambda_i=\alpha_i(\lambda)$, one has
\[
\mu_i=\lambda_i-\sum_{j\in\bfI} a_{ji}m_j.
\] 
The \emph{trigonometric algebra of difference operators} $\diff^\lambda_\mu$ is the $\mathbb{C}[\sz_{i,r}^{\pm 1}]_{i \in \mathbf{I}}^{1 \les r \les \lambda_i}$-algebra generated by
$\{D_{i,r}^{\pm 1}, \sw_{i,r}^{\pm \frac12}\}_{i\in \mathbf{I}}^{1\les r\les m_i}$ 
subject to the relations:
\begin{align*}
	[D_{i,r},D_{j,s}]=[\sw^{1/2}_{i,r},\sw^{1/2}_{j,s}]=0,& \qquad 
	D_{i,r}^{\pm 1}D_{i,r}^{\mp 1}=\sw_{i,r}^{\pm 1/2}\sw_{i,r}^{\mp 1/2}=1,\\ 
	D_{i,r}\sw^{1/2}_{j,s} &= q_i^{\delta_{ij}\delta_{rs}}\sw^{1/2}_{j,s}D_{i,r},
\end{align*}
for all $i,j\in \mathbf{I},\ 1\les r\les m_i,\ 1\les s\les m_j$.
Let $\widetilde{\diff}^\lambda_\mu$ be the localization of $\diff^\lambda_\mu$ at the 
multiplicative set generated by 
$$\{\sw_{i,r}-q_i^k\sw_{i,s}\}_{i\in \mathbf{I}, k\in \bbZ}^{1\les r\ne s\les m_i}.$$ 
Define
\begin{align*}
	Z_i(z) &=\prod_{s=1}^{\lambda_i} \left(1-\frac{q_i \sz_{i,s}}{z}\right), \\ 
	W_i(z) &=\prod_{r=1}^{m_i} \sw_{i,r}^{-\frac12} \left(1-\frac{\sw_{i,r}}{z}\right),\\ 
	W_{i,r}(z) &= \sw_{i,r}^{-\frac12} \prod_{1\les s\les m_i}^{s\ne r} \sw_{i,s}^{-\frac12} 
	\left(1-\frac{\sw_{i,s}}{z}\right).
\end{align*}
By \cite[Thm.~7.1]{FT19}, stated originally in \cite{GKLOq} in the unshifted case, there 
exists a unique $\mathbb{C}[\sz_{i,r}^{\pm 1}]_{i \in \mathbf{I}}^{1 \les r \les 
	\lambda_i}$-algebra homomorphism\footnote{The proof in \cite{FT19} concerns quantum loop algebras associated to Cartan matrices of finite type. However, since all the relations involve at most two vertices, the proof generalizes to affine types with the exception of $\sfA_1^{(1)}$ and $\sfA_2^{(2)}$.}
\begin{equation} \label{eq: GKLO homo}
	{\Phi}^{\lambda}_\mu \colon
	\mathbf{U}_\mu[\sz_{i,r}^{\pm 1}]_{i \in \mathbf{I}}^{1 \les r \les \lambda_i} \longrightarrow \widetilde{\diff}^\lambda_\mu, 
\end{equation}
called the \emph{GKLO homomorphism}, such that
\begin{align*}
	x_i^+(z) &\mapsto \frac{-q_i}{1-q_i^2} 
	\sum_{r=1}^{m_i} \delta\left(\sw_{i,r}z\right)\frac{Z_i(\sw_{i,r})}{W_{i,r}(\sw_{i,r})} \prod_{j\to i} \prod_{p=1}^{-a_{ji}} W_j(q_j^{-a_{ji}-2p}z^{-1})D_{i,r}^{-1}, \\ 
	x_i^-(z) &\mapsto \frac{q_i}{1-q_i^2} \cdot
	\sum_{r=1}^{m_i} \delta\left(q_i^2\sw_{i,r}z\right)\frac{1}{W_{i,r}(\sw_{i,r})}   \prod_{j\leftarrow i}\prod_{p=1}^{-a_{ji}} W_j(q_j^{-a_{ji}-2p}z^{-1})D_{i,r}, \\ 
	\phi^\pm_i(z) &\mapsto \left(\frac{Z_i(z^{-1})}{W_i(z^{-1})W_i(q_i^{-2}z^{-1})}
	\prod_{j \sim i} \prod_{p=1}^{-a_{ji}} W_j(q_j^{-a_{ji}-2p}z^{-1})\right)^\pm. 
\end{align*}
As before, $f(z)^\pm$ denote the expansion of a rational function $f(z)$ in $z^{\pm 
	1}$, respectively. 
Pulling back the natural representation of $\widetilde{\diff}^\lambda_\mu$ on a localized polynomial algebra via the GKLO homomorphism ${\Phi}^{\lambda}_\mu$ yields a \emph{GKLO representation}.  

\subsection{Shifted quantum loop symmetric pairs} \label{ss: sh iloop thm}

Consider the generating series 
\begin{equation}
	\wetad_{i}(z)=\sum_{k\ges - \mu_{\tau(i)}}\wetad_{i,k}z^k \aand \wy_i(z) = 
	\sum_{r\in\bbZ}\wy_{i,r}z^r,
\end{equation}
and let $\bfU^\imath_\mu$ be the algebra generated over $\bbC[\FC^{\pm1}]$ by 
$\{\wetad_{i,k}\}_{i\in\bfI,k\geqslant - \mu_{\tau(i)}}$, 
$\{\wetad_{i,-\mu_{\tau(i)}}^{-1}\}_{i\in\bfI}$, and $\{\wy_{i,r}\}_{i\in\bfI,r\in\bbZ}$, 
subject to relations 
\eqref{eq:weight-rel-eta-y-abs}--\eqref{eq:drinfeld-serre-relations-abs} and 
\[
\wetad_{i,- \mu_{\tau(i)}}\wetad_{i,- \mu_{\tau(i)}}^{-1}=\wetad_{i,- 
	\mu_{\tau(i)}}^{-1}\wetad_{i,- \mu_{\tau(i)}}=1, \qquad
[\wetad_{i,k},\wetad_{j,l}]=0.
\]
This algebra was introduced in the simply-laced split case in \cite{LP-gklo}, and in the 
more general quasi-split $\mathsf{ADE}$ case in \cite{lu-wang-weekes-26}. 

Many properties of $\bfU^\imath$ outlined in Section~\ref{s:comodule-realization} 
generalize straightforwardly to $\bfU^\imath_\mu$ as they only depend on the relations 
in the series form. In particular, we have the following analogue of Theorem 
\ref{ss:comodule-structure}. 

\begin{theorem}\label{thm:shifted-coaction}
	For any $\mu_1, \mu_2 \in \Lambda^\vee$, there is a morphism of algebras
	\[
	\DDi{s}{}_{\mu_1,\mu_2} \colon\bfU^\imath_{\mu_1+\mu_2} \longrightarrow
	\Lext{(\bfU^\imath_{\mu_1}\ten\bfU_{\mu_2})[\FC^{\pm1}]}{s^{-1};s}
	\]
	given by \eqref{eq:DDi-weta-field}--\eqref{eq:DDi-wy-field}. 
	Moreover, for every $\mu_3\in\Lambda^\vee$,
	\begin{equation}
		(\DDi{s_1}{}_{\mu_1,\mu_2}\ten\id)\circ\DDi{s_1s_2}{}_{\mu_1+\mu_2,\mu_3}
		=(\id\ten\DD{s_2}{}_{\mu_2,\mu_3})\circ\DDi{s_1}{}_{\mu_1,\mu_2+\mu_3}.
		\label{app:coassociativity}
	\end{equation}
\end{theorem} 

The proof is carried out in Appendix \ref{app:shifted-coaction}.\\

Following the coideal approach of Section~\ref{s:coideal-QLSP}, we define $\bfB_\mu$ 
to be the subalgebra of $\Lext{\bfU_\mu}{\FC;\FC^{-1}}$ generated over 
$\bbC[\FC^{\pm1}]$ by the elements
$\{\eta_{i,r}\}_{i\in\bfI,r\ges - \mu_{\tau(i)}}$, $\{\eta_{i,-\mu_{\tau(i)}}^{-1}\}_{i\in\bfI}$ 
and $\{y_{i,k}\}_{i\in\bfI, k\in\bbZ}$ encoded by the generating series 
\begin{equation}
	\eta_{i}(z)=\sum_{r\ges-\mu_{\tau(i)}}\eta_{i,r}z^r\aand 
	y_i(z)=\sum_{k\in\bbZ}y_{i,k}z^k
\end{equation}
defined as in \eqref{eq: eta defi}--\eqref{eq: y defi}. As a straightforward generalization of Proposition and  Remark~\ref{ss:central-specializations}, there is a surjective algebra homomorphism $\mathbf{U}^\imath_\mu \twoheadrightarrow \mathbf{B}_\mu$ given by
\[
\wetad_{i}(z) \mapsto \deta_i(z)=\pg_i(\FC^{-1}z^{-1})\eta_i(z), \quad 
\wy_i(z) \mapsto y_i(z). 
\]
More explicitly, we obtain a homomorphism
\[
\Psi_\mu \colon \mathbf{U}^\imath_\mu \to \Lext{\bfU_\mu}{\FC;\FC^{-1}}
\]
sending 
\begin{align}\label{eq: igklo i emb}
	\begin{split} 
		\wetad_{i}(z) &\mapsto \pg_i(\FC^{-1}z^{-1}) \phi_{i}^{+}(z) \phi_{\tau(i)}^{-}(\FC^{-1}{z}^{-1}),\\ 
		\wy_i(z) &\mapsto x_{i}^{-}(\FC^{-1}{z}^{-1}) + \pg_{i}(z) x_{\tau(i)}^{+}(z) \phi_{i}^{+}(\FC^{-1} z^{-1}) + \ps_{i}(z) \phi_{i}^{+}(\FC^{-1}z^{-1}).
	\end{split}
\end{align}

\subsection{New construction of $\imath$GKLO representations} \label{ss: iGKLO}

By taking completions on both sides, we can extend the GKLO homomorphism 
\eqref{eq: GKLO homo} to the completed algebra $\Lext{\bfU_\mu}{\FC;\FC^{-1}}$. 
Let ${}^\imath\widetilde{\diff}^\lambda_\mu$ be the localization of 
$\widetilde{\diff}^\lambda_\mu[\FC^{\pm1}]$ at the multiplicative subset generated by 
\[
1-q_i^n\FC^{-1}\mathsf w_{i,r}\mathsf w_{\tau(i),s}
\aand
1-q_i^n\FC^{-1}\mathsf w_{i,r}^{\,2}
\]
for all valid indices and $n\in\bbZ$.\\ 

Choose an orientation on $\Gamma$ satisfying, for each $i \neq \tau(i)$, the property: if $i \to j$ then $\tau(j)\to\tau(i)$, and if  $j \to i$ then $\tau(i)\to\tau(j)$. 
Set 
\begin{align} 
	\mathbf W_i(u)
	=W_i(q_i^{-1}u)\,
	W_{\tau(i)}(q_i^{-1}\FC u^{-1}),
	\aand 
	\mathbf Z_i(u)
	=Z_i(u)\,Z_{\tau(i)}(\FC u^{-1}).
\end{align}
These satisfy
\begin{align} \label{eq: igklo1}
	\mathbf W_i(u)=\mathbf W_{\tau(i)}(\FC u^{-1}),
	\qquad
	\mathbf Z_i(u)=\mathbf Z_{\tau(i)}(\FC u^{-1}), 
\end{align}
as well as 
\begin{align} \label{eq: igklo2}
	\prod_{p=1}^{-a_{ji}}
	W_j(q_j^{-a_{ji}-2p}u)\,
	W_{\tau(j)}(q_j^{-a_{ji}-2p}\FC u^{-1}) = 
	\prod_{p=1}^{-a_{ji}}
	\mathbf W_j(q_j^{-a_{ji}-2p+1}u).
\end{align}

\begin{corollary}
	Assume that $\sfA$ is a symmetrizable generalized Cartan matrix of finite or affine type other than $\sfA_1^{(1)}$ and $\sfA_2^{(2)}$.
	Then there exists a unique homomorphism 
	\[
	{}^\imath{\Phi}^{\lambda}_\mu \colon
	\mathbf{U}_\mu^\imath[\sz_{i,r}^{\pm 1}]_{i \in \mathbf{I}}^{1 \les r \les \lambda_i} \longrightarrow {}^\imath\widetilde{\diff}^\lambda_\mu, 
	\]
	given, for $i\in\bfI$ and $k=\tau(i)$, by
	\begin{align*}
		\wetad_i(z)
		&\mapsto
		\pg_i(\FC^{-1}z^{-1})
		\Biggl(
		\frac{\mathbf{Z}_i(z^{-1})}
		{\mathbf{W}_i(q_i z^{-1})\mathbf{W}_i(q_i^{-1}z^{-1})}
		\prod_{j\sim i}\prod_{p=1}^{-a_{ji}}
		\mathbf{W}_j(q_j^{-a_{ji}-2p+1}z^{-1})
		\Biggr)^+,
		\\[2mm]
		\wy_i(z)
		&\mapsto
		\frac{q_i}{1-q_i^2}
		\sum_{r=1}^{m_i}
		\delta\left(\frac{q_i^2\sw_{i,r}}{\FC z}\right)
		\frac{1}{W_{i,r}(\sw_{i,r})}
		\prod_{j\leftarrow i}\prod_{p=1}^{-a_{ji}}
		W_j(q_j^{-a_{ji}-2p}\FC z)D_{i,r}
		\\[1mm]
		&\quad-
		\frac{q_i}{1-q_i^2}\pg_i(z) \sum_{r=1}^{m_k}
		\frac{\delta(\sw_{k,r}z)}
		{W_{k,r}(\sw_{k,r})}
		\Biggl(
		\frac{
			(q_i^{-a_{ik}}-\FC^{-1}z^{-2})\mathbf{Z}_i(\FC z)}
		{
			(1-q_i^{-a_{ik}}\FC^{-1}z^{-2})
			W_i(\FC z)W_i(q_i^{-2}\FC z)}
		\\*
		&\qquad\times	
\prod_{\substack{j\sim i\\ \tau(j)\not\to k}}
\prod_{p=1}^{-a_{ji}}
W_j\bigl(q_j^{-a_{ji}-2p}\FC z\bigr)
\prod_{\substack{j\sim i\\ \tau(j)\to k}}
\prod_{p=1}^{-a_{ji}}
\mathbf W_j\bigl(q_j^{-a_{ji}-2p+1}\FC z\bigr)
		\Biggr)^- D_{k,r}^{-1}
		\\[1mm]
		&\quad+
		\ps_i(z)
		\Biggl(
		\frac{Z_i(\FC z)}
		{W_i(\FC z)W_i(q_i^{-2}\FC z)}
		\prod_{j\sim i}\prod_{p=1}^{-a_{ji}}
		W_j(q_j^{-a_{ji}-2p}\FC z)
		\Biggr)^-. 
	\end{align*} 
\end{corollary} 

\begin{proof}
	It is enough to observe that ${}^\imath{\Phi}^{\lambda}_\mu = {\Phi}^{\lambda}_\mu \circ 
	\Psi_\mu$. The result then follows by a direct computation, using 
	\eqref{eq: igklo1}--\eqref{eq: igklo2} and basic commutation relations. One readily 
	checks that the image of ${}^\imath{\Phi}^{\lambda}_\mu$ lies in the localization 
	${}^\imath\widetilde{\diff}^\lambda_\mu$, relying on the identities
	\begin{align}
		\delta(\mathsf w_{\tau(i),r}z)
		\left(\frac{1}{1-\mathsf w_{i,s}/(\FC z)}\right)^-
		&= 
		\frac{\delta(\mathsf w_{\tau(i),r}z)}
		{1-\FC^{-1}\mathsf w_{i,s}\mathsf w_{\tau(i),r}}, \\ 
		\delta(\FC z^2)
		\left(\frac{1}{1-a/(\FC z)}\right)^-
		&=
		\frac{1+az}{1-\FC^{-1} a^2}\,\delta(\FC z^2), 
	\end{align}
	for $a=\mathsf w_{i,s}, \ q_i^2\mathsf w_{i,s}$.
\end{proof} 

\begin{remarks}
	\hfill
	
	\begin{enumerate}[leftmargin=2em]\itemsep0.25cm
		\item While the $\imath$GKLO homomorphism constructed above does not match the 
		earlier constructions in \cite{LP-gklo, lu-wang-weekes-26} exactly, it has a remarkably 
		similar form, especially when $\tau = \id$. In particular, in that case the different nature 
		of the first two summands of \eqref{eq: igklo i emb} explains the intrinsic asymmetry 
		between the coefficients on $D_{i,r}$ and $D_{i,r}^{-1}$ observed in \cite[Remark 
		3.9]{LP-gklo}. Moreover, the third summand of \eqref{eq: igklo i emb} explains the 
		appearance of a corresponding third constant (no difference operator) term in the 
		$\imath$GKLO formulae of \cite{LP-gklo, lu-wang-weekes-26}. Remarkably, the 
		condition on $\sigma_i(z)$ established in Proposition \ref{ss:SF+GF}(2) above is 
		precisely the restriction on the appearance of the constant term observed in 
		\cite[(3.2)]{LP-gklo} and \cite[(3.2)]{lu-wang-weekes-26}.  
		\item Our construction provides a new conceptual explanation for the existence of 
		$\imath$GKLO homomorphisms. While the idea of restricting the usual GKLO 
		homomorphism for a quantum affine algebra to its coideal subalgebra may appear as an 
		obvious ansatz, this strategy is of little use when applied to the ordinary non-topological 
		subalgebra structure. Indeed, with a view to geometric applications, it is natural to 
		desire to formulate an $\imath$GKLO homomorphism in terms of the Drinfeld 
		presentation. 
		However, the relationship between the Drinfeld presentations of a coideal subalgebra 
		and a quantum affine algebra is complicated and can only be described approximately 
		\cite{Przez-23, li-przezdziecki-25, li-przezdziecki-26}. In particular, a GKLO 
		representation restricted with respect to the usual subalgebra structure would not admit 
		a direct description in terms of the Drinfeld presentation of $\mathbf{U}^\imath$. For 
		this reason, $\imath$GKLO homomorphisms were constructed in \cite{LP-gklo, 
			lu-wang-weekes-26} \emph{directly} without any reference to the GKLO 
		homomorphisms for $\mathbf{U}$. Corollary \ref{ss: iGKLO}  above gives a new 
		construction of $\imath$GKLO homomorphisms, which, surprisingly, connects them 
		to ordinary GKLO representations. 
		\item Corollary \ref{ss: iGKLO} gives a uniform construction which works for any symmetrizable generalized Cartan matrix of finite or affine type other than $\sfA_1^{(1)}$ and $\sfA_2^{(2)}$,
		extending, in particular, the 
		results of \cite{LP-gklo, lu-wang-weekes-26} to non-simply-laced and toroidal settings. 
		\hfill\rmkend
	\end{enumerate}
\end{remarks}

\appendix
\renewcommand{\endpf}{\qed\par\medskip}
\section{The Drinfeld coaction in the shifted setting}
\label{app:shifted-coaction}
We give a direct proof of Theorem~\ref{thm:shifted-coaction} for an arbitrary 
symmetrizable Cartan matrix. We first check that the coaction formulae are well defined 
in the shifted setting, then verify the defining relations and coassociativity.

\subsection{The Drinfeld coaction}
Let $\sfA=(a_{ij})_{i,j\in\bfI}$ be the symmetrizable Cartan matrix and $\tau$ 
the diagrammatic involution fixed in the main text. Recall that $q_i=q^{d_i}$ and
$d_{\tau(i)}=d_i$, and $a_{\tau(i)\tau(j)}=a_{ij}$.
Fix $\mu_1,\mu_2\in\Lambda^\vee$, set $\mu=\mu_1+\mu_2$, and write 
$\mu_{a,i}=\alpha_i(\mu_a)$ for $a=1,2$. We suppress the shifts from the notation 
$\DDi{s}{}_{\mu_1,\mu_2}$.
For the reader's convenience, we recall the main definitions.
The generating series of $\bfU^\imath_\mu$ and $\bfU_{\mu_2}$ are
\[
\wetad_i(z)=\sum_{r\ges-\mu_{\tau(i)}}\wetad_{i,r}z^r,\quad
\wy_i(z)=\sum_{r\in\bbZ}\wy_{i,r}z^r,\quad
x_i^\pm(z)=\sum_{r\in\bbZ}x_{i,r}^\pm z^r,
\]
\[
\phi_i^+(z)=\sum_{r\ges0}\phi_{i,r}^+z^r,\qquad
\phi_i^-(z)=\sum_{r\les\mu_{2,i}}\phi_{i,r}^-z^r.
\]
Recall also that
\[
\RF_{ij}(z,w)=\frac{q_i^{a_{ij}}z-w}{z-q_i^{a_{ij}}w},
\qquad
\fdel(z)=\sum_{r\in\bbZ}z^r,
\]
and, since $d_i a_{ij}=d_j a_{ji}$, one has
\[
\RF_{ij}(z,w)=\RF_{ji}(z,w)
\aand
\RF_{ij}(z,w)^{-1}=\RF_{ij}(w,z).
\]
As before, we write $\RF_{ij}(z/w)=\RF_{ij}(z,w)$, with the prescribed formal expansion.
We shall use the following relations in $\bfU_{\mu_2}$:
\begin{align}
	\phi_i^\pm(u)x_j^+(v)&=\RF_{ij}(v,u)x_j^+(v)\phi_i^\pm(u),\label{app:phi-x-plus}\\
	\phi_i^\pm(u)x_j^-(v)&=\RF_{ij}(u,v)x_j^-(v)\phi_i^\pm(u),\label{app:phi-x-minus}\\
	[x_i^+(u),x_j^-(v)]&=\frac{\delta_{ij}}{q_i-q_i^{-1}}
	\bigl(\fdel(v/u)\phi_i^+(v)-\fdel(u/v)\phi_i^-(u)\bigr).
	\label{app:mixed-current}
\end{align}
The weight relations are expanded at $u=0$ for $\phi_i^+(u)$ and at
$u=\infty$ for $\phi_i^-(u)$. The quadratic root-current relations are
used in polynomial form:
\begin{align}
	(q_i^{a_{ij}}u-v)x_i^+(u)x_j^+(v)
	&=(u-q_i^{a_{ij}}v)x_j^+(v)x_i^+(u),\notag\\
	(u-q_i^{a_{ij}}v)x_i^-(u)x_j^-(v)
	&=(q_i^{a_{ij}}u-v)x_j^-(v)x_i^-(u).
	\label{app:ordinary-quadratic}
\end{align}
Moreover, for $i\ne j$ and $N=1-a_{ij}$, the Drinfeld--Serre relations read
\begin{equation}
	\Sym_{\ul z}\sum_{p=0}^{N}(-1)^p
	\begin{bmatrix}N\\p\end{bmatrix}_{q_i}
	x_i^\pm(z_1)\cdots x_i^\pm(z_p)x_j^\pm(w)
	x_i^\pm(z_{p+1})\cdots x_i^\pm(z_N)=0.
	\label{app:ordinary-serre}
\end{equation}
In $\bfU^\imath_\mu$, the weight and quadratic relations read
\begin{equation}
	\wetad_i(z)\wy_j(w)=\RF_{ij}(\FC zw)\RF_{\tau(i)j}(w/z)\wy_j(w)\wetad_i(z),
	\label{app:eta-y}
\end{equation}
where both rational factors are expanded at $z=0$, and
\begin{align}
	&(q_i^{a_{ij}}z-w)\wy_i(z)\wy_j(w)
	+(q_i^{a_{ij}}w-z)\wy_j(w)\wy_i(z)\notag\\
	&\qquad=
	\delta_{i,\tau(j)}\frac{\fdel(\FC zw)}{q_i-q_i^{-1}}
	\Bigl((q_i^{a_{i\tau(i)}}z-w)\wetad_{\tau(i)}(z)
	+(q_i^{a_{i\tau(i)}}w-z)\wetad_i(w)\Bigr).
	\label{app:y-y}
\end{align}
As in \eqref{eq:abstract-serre-cartan}, set
\begin{equation}
	\wetadd_i(z)=\frac{q_i^{a_{i\tau(i)}}\FC z^2-1}{q_i^{a_{i\tau(i)}}-\FC z^2}\,\wetad_{\tau(i)}(z),
	\label{app:ddot-eta}
\end{equation}
where the scalar factor is expanded at $z=0$.
We recall the recursive form of the $\imath$Drinfeld--Serre relations from
\ref{ss:i-quantum-group}. For a tuple $\ul v=((i_1,z_1),\ldots,(i_n,z_n))$, let $\ul v^k$ 
denote the tuple obtained by removing its $k$th entry. Multiple superscripts indicate 
simultaneous deletions, preserving the order of the remaining entries. Set 
$\ww(\varnothing)=1$ and
\begin{equation}
	\small
	\ww(\ul v)=\wy_{i_1}(z_1)\ww(\ul v^1)-\sum_{k=2}^n
	\delta_{i_1,\tau(i_k)}\frac{\fdel(\FC z_1z_k)}{q_{i_1}-q_{i_1}^{-1}}
	\left(\prod_{r=2}^{k-1}\RF_{i_k i_r}(z_k/z_r)\right)
	\wetadd_{i_k}(z_k)\ww(\ul v^{1,k}).
	\label{app:W-recursion}
\end{equation}
The rational functions in the product are expanded at $z_k/z_r=0$.  
Note that the recursion is well-defined even in the shifted setting.
Then, for $i\ne j$, 
$N=1-a_{ij}$, and
\[
\ul v_p=((i,z_1),\ldots,(i,z_p),(j,w),
(i,z_{p+1}),\ldots,(i,z_N)),
\]
the defining relation is
\begin{equation}
	\Sym_{\ul z}\sum_{p=0}^{N}(-1)^p
	\begin{bmatrix}N\\p\end{bmatrix}_{q_i}\ww(\ul v_p)=0.
	\label{app:iota-serre}
\end{equation}
Finally, the assignments in Theorem~\ref{thm:shifted-coaction} defining the Drinfeld coaction on $\bfU^\imath_\mu$ are
\begin{align}
	\DDi{s}(\FC)&=\FC\ten1,\label{app:coaction-zeta}\\
	\DDi{s}(\wetad_i(z))
	&=\wetad_i(z)\ten
	\phi_i^+(zs)\phi_{\tau(i)}^-({\FC^{-1}z^{-1}} s),\label{app:coaction-eta}\\
	\DDi{s}(\wy_i(z))
	&=1\ten x_i^-({\FC^{-1}z^{-1}} s)+\wy_i(z)\ten\phi_i^+({\FC^{-1}z^{-1}} s)\notag\\
	&\quad+\wetad_i({\FC^{-1}z^{-1}})\ten x_{\tau(i)}^+(zs)\phi_i^+({\FC^{-1}z^{-1}} s).
	\label{app:coaction-y}
\end{align}

\begin{lemma}\label{app:completion-lemma}
	The coefficients of the right-hand sides of \eqref{app:coaction-eta} and 
	\eqref{app:coaction-y} belong to
	$\Lext{(\bfU^\imath_{\mu_1}\ten\bfU_{\mu_2})[\FC^{\pm1}]}{s^{-1};s}$.
	Moreover, the lowest $z$-degree in \eqref{app:coaction-eta} is $-\mu_{\tau(i)}$, and its 
	coefficient is invertible.
\end{lemma}

\begin{pf}
	Write
	\[
	\phi_i^+(zs)=\sum_{b\ges0}\phi_{i,b}^+z^bs^b
	\aand
	\phi_{\tau(i)}^-({\FC^{-1}z^{-1}} s)=\sum_{c\les\mu_{2,\tau(i)}}
	\FC^{-c}\phi_{\tau(i),c}^-z^{-c}s^c.
	\]
	A term in \eqref{app:coaction-eta} has $z$-degree $a+b-c\ges-\mu_{\tau(i)}$, where 
	$a\ges-\mu_{1,\tau(i)}$, $b\ges0$, and $c\les\mu_{2,\tau(i)}$. For 
	$n\ges-\mu_{\tau(i)}$, the coefficient of $z^n$ is therefore
	\begin{equation}
		\DDi{s}(\wetad_{i,n})
		=\sum_{\substack{b\ges0,\ c\les\mu_{2,\tau(i)}\\
				n-b+c\ges-\mu_{1,\tau(i)}}}
		\FC^{-c}s^{b+c}\wetad_{i,n-b+c}\ten
		\phi_{i,b}^+\phi_{\tau(i),c}^-.
		\label{app:eta-mode}
	\end{equation}
	Since $0\les b\les n+\mu_{\tau(i)}$ and $b-n-\mu_{1,\tau(i)}\les c\les\mu_{2,\tau(i)}$, 
	this sum is finite. For $n=-\mu_{\tau(i)}$, the only term has $a=-\mu_{1,\tau(i)}$, $b=0$, 
	and $c=\mu_{2,\tau(i)}$. Hence
	\begin{equation}
		\DDi{s}(\wetad_{i,-\mu_{\tau(i)}})
		=\FC^{-\mu_{2,\tau(i)}}s^{\mu_{2,\tau(i)}}
		\wetad_{i,-\mu_{1,\tau(i)}}\ten
		\phi_{i,0}^+\phi_{\tau(i),\mu_{2,\tau(i)}}^-.
		\label{app:lowest-eta-mode}
	\end{equation}
	All factors on the right are invertible. Thus the lowest degree is $-\mu_{\tau(i)}$, and the assignment extends to the inverse of the lowest Cartan coefficient. This holds for arbitrary coweights $\mu_1,\mu_2$.
	
	For \eqref{app:coaction-y}, the coefficient of $z^r$ is
	\begin{align}
		\DDi{s}(\wy_{i,r})={}&
		\FC^rs^{-r}(1\ten x_{i,-r}^-)
		+\sum_{c\ges0}\FC^{-c}s^c\wy_{i,r+c}\ten\phi_{i,c}^+\notag\\
		&+\sum_{\substack{a\ges-\mu_{1,\tau(i)}\\c\ges0}}
		\FC^{-a-c}s^{r+a+2c}\wetad_{i,a}\ten
		x_{\tau(i),r+a+c}^+\phi_{i,c}^+.
		\label{app:y-mode}
	\end{align}
	The three terms have $s$-degrees bounded below by $-r$, $0$, and 
	$r-\mu_{1,\tau(i)}$, respectively. In the last sum, the coefficient of $s^m$ receives 
	contributions only from pairs satisfying $a+2c=m-r$, $a\ges-\mu_{1,\tau(i)}$, and 
	$c\ges0$. There are finitely many such pairs, so the expression belongs to the stated 
	completion.
\end{pf}

\begin{remark}
	The scalar factor in \eqref{app:ddot-eta} is a power series in $z^2$ with constant term 
	$-q_i^{-a_{i\tau(i)}}$. Hence
	\[
	\wetadd_{i,-\mu_i}=-q_i^{-a_{i\tau(i)}}\wetad_{\tau(i),-\mu_i},
	\]
	where $\wetadd_i(z)=\sum_{r\ges-\mu_i}\wetadd_{i,r}z^r$.
	In particular,
	\begin{equation}
		\DDi{s}(\wetadd_i(z))
		=\wetadd_i(z)\ten
		\phi_{\tau(i)}^+(zs)\phi_i^-({\FC^{-1}z^{-1}} s).
		\label{app:ddot-image}
	\end{equation}
	Indeed, this follows from \eqref{app:coaction-eta} and \eqref{app:ddot-eta}, since $\DDi{s}$ fixes the scalar factor and $\tau^2=\id$. In particular,
	\[
	\DDi{s}(\wetadd_{i,-\mu_i})=\FC^{-\mu_{2,i}}s^{\mu_{2,i}}
	\wetadd_{i,-\mu_{1,i}}\ten\phi_{\tau(i),0}^+\phi_{i,\mu_{2,i}}^-,
	\]
	so the lowest $z$-degree in \eqref{app:ddot-image} is $-\mu_i$.
	\hfill\rmkend
\end{remark}

\subsection{Weight and quadratic relations}
\label{app:quadratic-prop}
Until compatibility with the defining relations is established, $\DDi{s}$ denotes the 
multiplicative extension of \eqref{app:coaction-zeta}--\eqref{app:coaction-y} to the free 
algebra on the current coefficients.
Set
\begin{equation}
	\begin{aligned}
		A_i(z)&=1\ten x_i^-({\FC^{-1}z^{-1}} s),\\
		B_i(z)&=\wy_i(z)\ten\phi_i^+({\FC^{-1}z^{-1}} s),\\
		C_i(z)&=\wetad_i({\FC^{-1}z^{-1}})\ten x_{\tau(i)}^+(zs)\phi_i^+({\FC^{-1}z^{-1}} s).
	\end{aligned}
	\label{app:ABC}
\end{equation}
Thus $\DDi{s}(\wy_i(z))=A_i(z)+B_i(z)+C_i(z)$. In products of these images, we label 
each entry by the chosen summand $A$, $B$, or $C$.
Set also
\begin{equation}
	\kappa_{ij}(z,w)=
	\delta_{i,\tau(j)}
	\frac{\fdel(\FC zw)}{q_i-q_i^{-1}},
	\label{app:kappa}
\end{equation}
and define
\[
\Theta_i^\pm(z)=\wetadd_i(z)\ten\phi_{\tau(i)}^+(zs)\phi_i^\pm({\FC^{-1}z^{-1}} s).
\]
By \eqref{app:ddot-image}, $\Theta_i^-(z)=\DDi{s}(\wetadd_i(z))$.

\begin{lemma}
	For every $i,j\in\bfI$, one has the following identities, interpreted coefficientwise in $z,w$ with values in
		$\Lext{(\bfU^\imath_{\mu_1}\ten\bfU_{\mu_2})[\FC^{\pm1}]}{s^{-1};s}$,
		where every factor $\RF_{ij}(w/z)$ is expanded at $w/z=0$:
	\begin{align}
		B_i(z)A_j(w)&=\RF_{ij}(w/z)A_j(w)B_i(z),\label{app:BA}\\
		C_i(z)B_j(w)&=\RF_{ij}(w/z)B_j(w)C_i(z),\label{app:CB}
	\end{align}
	and
	\begin{equation}
		C_i(z)A_j(w)=\RF_{ij}(w/z)A_j(w)C_i(z)
		-\kappa_{ij}(z,w)\bigl(\Theta_j^+(w)-\Theta_j^-(w)\bigr).
		\label{app:CA}
	\end{equation}
	Consequently,
	\begin{equation}
		\begin{aligned}
			&B_i(z)B_j(w)+C_i(z)A_j(w)-\kappa_{ij}(z,w)\Theta_j^-(w)=\\[-1mm]
			&\quad=\bigl(\wy_i(z)\wy_j(w)-\kappa_{ij}(z,w)\wetadd_j(w)\bigr)\ten \phi_i^+({\FC^{-1}z^{-1}} s)\phi_j^+({\FC^{-1}w^{-1}} s)\\
			&\hspace{5cm}
			+\RF_{ij}(w/z)A_j(w)C_i(z).
		\end{aligned}
		\label{app:local-cancellation}
	\end{equation}
\end{lemma}

\begin{pf}
	Identity \eqref{app:BA} follows from the weight relation for $\phi_i^+$ and $x_j^-$ in the second tensor factor. For \eqref{app:CB}, apply \eqref{app:eta-y} in the first tensor factor and the weight relation for $\phi_j^+$ and $x_{\tau(i)}^+$ in the second. The factors involving $\RF_{\tau(i)j}(\FC zw)$ cancel, leaving $\RF_{ij}(w/z)$.

	If $i\neq\tau(j)$, the root currents in $C_i(z)A_j(w)$ commute, and \eqref{app:CA} 
	follows in the same way. Suppose now that $i=\tau(j)$. Then $j=\tau(i)$, $q_j=q_i$, and
	\[
	q_j^{a_{j\tau(j)}}=q_j^{a_{ji}}=q_i^{a_{ij}}.
	\]
	On the support of $\fdel(\FC zw)$, substitute
	$z={\FC^{-1}w^{-1}}$, so that ${\FC^{-1}z^{-1}}=w$. Then
	\[
	\RF_{ij}(w/z)\wetad_i({\FC^{-1}z^{-1}})
	=\frac{q_i^{a_{ij}}\FC w^2-1}{\FC w^2-q_i^{a_{ij}}}\,\wetad_i(w)
	=-\wetadd_j(w).
	\]
	Thus,
	\begin{equation}
		\fdel(\FC zw)\RF_{ij}(w/z)\wetad_i({\FC^{-1}z^{-1}})
		=-\fdel(\FC zw)\wetadd_j(w).
		\label{app:eta-ddot-support}
	\end{equation}
		Moving $\phi_i^+({\FC^{-1}z^{-1}} s)$ past
		$x_j^-({\FC^{-1}w^{-1}} s)$ and applying \eqref{app:mixed-current} gives
		\[
		\begin{aligned}
			&C_i(z)A_j(w)-\RF_{ij}(w/z)A_j(w)C_i(z)=\\
			&\quad=\kappa_{ij}(z,w)\RF_{ij}(w/z)\wetad_i({\FC^{-1}z^{-1}})
			\ten\bigl(\phi_j^+({\FC^{-1}w^{-1}} s)-\phi_j^-(zs)\bigr)
			\phi_i^+({\FC^{-1}z^{-1}} s).
		\end{aligned}
		\]
		On the delta support, $\phi_j^-(zs)=\phi_j^-({\FC^{-1}w^{-1}} s)$ and
		$\phi_i^+({\FC^{-1}z^{-1}} s)=\phi_{\tau(j)}^+(ws)$. Thus
		\eqref{app:eta-ddot-support} gives
		$-\kappa_{ij}(z,w)(\Theta_j^+(w)-\Theta_j^-(w))$.
	This proves \eqref{app:CA}. Finally, the substitution $w={\FC^{-1}z^{-1}}$ gives
	\[
	\fdel(\FC zw)\phi_i^+({\FC^{-1}z^{-1}} s)\phi_j^+({\FC^{-1}w^{-1}} s)=\fdel(\FC zw)\phi_{\tau(j)}^+(ws)\phi_j^+({\FC^{-1}w^{-1}} s),
	\]
	and therefore
	\[
	\begin{aligned}
		&B_i(z)B_j(w)-\kappa_{ij}(z,w)\Theta_j^+(w)\\
		&\qquad=\bigl(\wy_i(z)\wy_j(w)
		-\kappa_{ij}(z,w)\wetadd_j(w)\bigr)
		\ten \phi_i^+({\FC^{-1}z^{-1}} s)\phi_j^+({\FC^{-1}w^{-1}} s).
	\end{aligned}
	\]
	Adding this identity to \eqref{app:CA} yields
	\eqref{app:local-cancellation}.
\end{pf}

\begin{proposition}
	The assignments \eqref{app:coaction-zeta}--\eqref{app:coaction-y} preserve the commutativity of the Cartan currents and the weight and quadratic relations \eqref{app:eta-y}--\eqref{app:y-y}.
\end{proposition}

\begin{pf}
	The images of the Cartan currents commute by \eqref{app:coaction-eta}. Relation \eqref{app:eta-y} follows by applying the weight relations in the two tensor factors to each of $A_j(w)$, $B_j(w)$, and $C_j(w)$.
	
	To prove compatibility with \eqref{app:y-y}, set $\alpha=q_i^{a_{ij}}z-w$ and $\beta=q_i^{a_{ij}}w-z$.
	The quadratic relations for the root currents give
	\[
	\alpha A_i(z)A_j(w)+\beta A_j(w)A_i(z)=0,
	\qquad
	\alpha C_i(z)C_j(w)+\beta C_j(w)C_i(z)=0.
	\]
	For the mixed $A,B$ and $B,C$ terms, combine \eqref{app:BA} and
	\eqref{app:CB} with the identity $\alpha\RF_{ij}(w/z)+\beta=0$.
	Thus these contributions to the left-hand side of \eqref{app:y-y} vanish.
	
	The contribution involving two $B$-terms is
	\[
	\kappa_{ij}(z,w)
	\bigl(\alpha\wetad_{\tau(i)}(z)+\beta\wetad_i(w)\bigr)
	\ten \phi_i^+({\FC^{-1}z^{-1}} s)\phi_j^+({\FC^{-1}w^{-1}} s),
	\]
	and \eqref{app:CA}, together with the identity obtained by interchanging $(i,z)$ and $(j,w)$, gives the mixed $A,C$ contribution
	\[
	-\kappa_{ij}(z,w)
	\bigl(\alpha(\Theta_j^+(w)-\Theta_j^-(w))
	+\beta(\Theta_i^+(z)-\Theta_i^-(z))\bigr).
	\]
	Both expressions vanish if $i\ne\tau(j)$. If $i=\tau(j)$, then $w={\FC^{-1}z^{-1}}$ on the support of $\fdel(\FC zw)$, and \eqref{app:ddot-eta} gives
	\[
	\beta\wetadd_i(z)=\alpha\wetad_j(z),\qquad
	\alpha\wetadd_j(w)=\beta\wetad_i(w).
	\]
	Moreover, $\phi_i^+({\FC^{-1}z^{-1}} s)\phi_j^+({\FC^{-1}w^{-1}} s)=\phi_i^+(ws)\phi_j^+(zs)$ on the same support. Thus the terms involving $\Theta^+$ cancel the contribution involving two $B$-terms. The remaining terms satisfy on the same delta support
	\[
	\alpha\Theta_j^-(w)+\beta\Theta_i^-(z)
	=\alpha\DDi{s}(\wetad_j(z))+\beta\DDi{s}(\wetad_i(w)).
	\]
	Multiplying by $\kappa_{ij}(z,w)$ gives the required right-hand side.
\end{pf}

\subsection{The $\imath$Drinfeld--Serre relations}
For a tuple $\ul v=((i_1,z_1),\ldots,(i_n,z_n))$, set $\cW_s(\ul v)=\DDi{s}(\ww(\ul v))$, where $\ww(\ul v)$ is regarded as an expression in the free algebra. Writing $X_a=A_{i_a}(z_a)+B_{i_a}(z_a)+C_{i_a}(z_a)$, the recursion \eqref{app:W-recursion} gives
\begin{equation}
	\cW_s(\ul v)=X_1\cW_s(\ul v^1)-\sum_{k=2}^n
	\kappa_{i_1i_k}(z_1,z_k)\left(\prod_{r=2}^{k-1}\RF_{i_k i_r}(z_k/z_r)\right)
	\Theta^-_{i_k}(z_k)\cW_s(\ul v^{1,k}).
	\label{app:Ws-recursion}
\end{equation}

\subsubsection{Commutation with $A_l(v)$ and $C_l(v)$}
\label{app:Bblock-A}
For $\ul u=((j_1,u_1),\ldots,(j_m,u_m))$, set
\[
\bbW(\ul u)=\ww(\ul u)\ten \phi_{j_1}^+({\FC^{-1}u_1^{-1}} s)\cdots \phi_{j_m}^+({\FC^{-1}u_m^{-1}} s)
\aand
\bbW(\varnothing)=1.
\]
The dependence on $s$, through the currents $\phi_j^+({\FC^{-1}u^{-1}} s)$, 
is omitted from the notation.
We shall use the following commutation relations.
\begin{lemma}
	Let $\ul u=((j_1,u_1),\ldots,(j_m,u_m))$.  Then
	\begin{equation}
		\bbW(\ul u)A_l(v)=\left(\prod_{a=1}^m\RF_{j_a l}(v/u_a)\right)A_l(v)
		\bbW(\ul u).
		\label{app:Bblock-A-eq}
	\end{equation}
	Moreover,
	\begin{equation}
		C_l(v)\bbW(\ul u)=\left(\prod_{a=1}^m\RF_{l j_a}(u_a/v)\right)
		\bbW(\ul u)C_l(v).
		\label{app:C-Bblock-eq}
	\end{equation}
\end{lemma}

The rational functions in the two identities above are expanded at $v/u_a=0$ and
$u_a/v=0$, respectively. Both identities hold coefficientwise with values in
	$\Lext{(\bfU^\imath_{\mu_1}\ten\bfU_{\mu_2})[\FC^{\pm1}]}{s^{-1};s}$.

\begin{pf}
Since $A_l(v)=1\ten x_l^-({\FC^{-1}v^{-1}} s)$, it is enough to commute $x_l^-({\FC^{-1}v^{-1}} s)$ past the currents $\phi_{j_a}^+({\FC^{-1}u_a^{-1}} s)$ in the second tensor factor. The weight relations give the factors in \eqref{app:Bblock-A-eq}.

For \eqref{app:C-Bblock-eq}, apply the weight relation to
$\wetad_l({\FC^{-1}v^{-1}})\ww(\ul u)$. Induction on the recursion gives the
product of the weight factors for all entries: on the support of the delta functions
appearing in a lower degree term, the factors for its two entries multiply to $1$.
Commuting $x_{\tau(l)}^+(vs)$ past the Cartan currents in the second
tensor factor cancels the factors $\RF_{\tau(l)j_a}(\FC vu_a)$.
The remaining factors are those in \eqref{app:C-Bblock-eq}.
\end{pf}

\subsubsection{Homogeneous quadratic relations}
\label{app:internal-exchange}
We next consider the exchange of adjacent entries with the same label.
\begin{lemma}
	For $i,j\in\bfI$, one has
	\begin{equation}
		\begin{aligned}
			(q_i^{a_{ij}}z-w)A_i(z)A_j(w)+(q_i^{a_{ij}}w-z)A_j(w)A_i(z)&=0,\\
			(q_i^{a_{ij}}z-w)C_i(z)C_j(w)+(q_i^{a_{ij}}w-z)C_j(w)C_i(z)&=0.
		\end{aligned}
		\label{app:AC-internal-exchange}
	\end{equation}
	Moreover, for any ordered tuples $\ul u,\ul u'$, one has
	\begin{equation}
		\begin{aligned}
			&(q_i^{a_{ij}}z-w)\bbW\bigl(\ul u,(i,z),(j,w),\ul u'\bigr)\\
			&\qquad+(q_i^{a_{ij}}w-z)\bbW\bigl(\ul u,(j,w),(i,z),\ul u'\bigr)=0.
		\end{aligned}
		\label{app:B-internal-exchange}
	\end{equation}
\end{lemma}

The proof of these identities uses only the weight and quadratic relations.

\begin{pf}
	The first identity in \eqref{app:AC-internal-exchange} is the negative
	quadratic relation evaluated at ${\FC^{-1}z^{-1}} s,{\FC^{-1}w^{-1}} s$.
	For the second, move the Cartan currents to the right and apply the
	positive quadratic relation. The additional factor is symmetric in
	$(i,z),(j,w)$, since $q_i^{a_{i,\tau(j)}}=q_j^{a_{j,\tau(i)}}$.
	
	For \eqref{app:B-internal-exchange}, it is enough to prove the same
	relation for $\ww$, since the Cartan factors commute. Recall the
	polynomials $\ZF_{ij}(z,w)=q_i^{a_{ij}}z-w$ from \ref{ss:QL}.
	We induct on the length of the tuple. For length two, the assertion
	follows from \eqref{app:y-y} and \eqref{app:ddot-eta}.
	
	The weight relation gives
	\begin{equation}
		\begin{aligned}
			&\kappa_{ik}(z,u)\wetadd_k(u)\wy_j(w)\\
			&\qquad=\kappa_{ik}(z,u)\RF_{ij}(w/z)\RF_{jk}(u/w)^{-1}
			\wy_j(w)\wetadd_k(u),
		\end{aligned}
		\label{app:d-y-exchange}
	\end{equation}
	with the expansions inherited from \eqref{app:eta-y}. Indeed, on the
	support of $\kappa_{ik}(z,u)$, one has $i=\tau(k)$ and $z={\FC^{-1}u^{-1}}$.
	
	Suppose the exchanged entries occupy positions $\ell,\ell+1$ with $\ell>1$.
		Expand \eqref{app:W-recursion} for both orders. The leading terms
		$\wy_{i_1}(z_1)\ww(\ul v^1)$ and the summands indexed by
		$k\notin\{\ell,\ell+1\}$ cancel by induction: the exchanged entries remain
		adjacent in the shorter tuple, and their scalar prefactors agree.
		For $k=\ell,\ell+1$, group the terms from the two expansions according
		to which entry is removed. Each pair cancels by
	\[
	\ZF_{ij}(z,w)+\ZF_{ji}(w,z)\RF_{ij}(z/w)=0
	\]
	and the identity obtained by interchanging $(i,z)$ and $(j,w)$.
	
	It remains to exchange the first two entries of $\ul v$. Expand the
	recursion twice in
	\[
	\begin{aligned}
		&\ZF_{i_1i_2}(z_1,z_2)\ww(\ul v)\\
		&\qquad+\ZF_{i_2i_1}(z_2,z_1)
		\ww\bigl((i_2,z_2),(i_1,z_1),\ul v^{1,2}\bigr).
	\end{aligned}
	\]
	The terms with no delta functions or with the pair $(1,2)$ cancel by the
	length-two case. Those with a single delta function cancel by
	\eqref{app:d-y-exchange}.
	
	The remaining terms contain two disjoint delta functions, pairing $1$ with
	$k$ and $2$ with $l$, where $k,l\ges3$. Pair the corresponding terms in
	the two expansions by interchanging $k,l$. The Cartan factors commute,
	and on the common delta support one has
	\[
	i_1=\tau(i_k),\quad i_2=\tau(i_l),\quad
	z_1={\FC^{-1}z_k^{-1}},\quad z_2={\FC^{-1}z_l^{-1}},
	\]
	so that
	\[
	\RF_{i_2i_k}(z_k/z_2)=\RF_{i_1i_l}(z_l/z_1).
	\]
	The products of intermediate kernels in \eqref{app:W-recursion} agree
	except for the factor associated to $k,l$. Thus the paired terms have
	common factors multiplying, respectively,
	\begin{equation}
		\begin{cases}
			\ZF_{i_1i_2}(z_1,z_2)
			+\ZF_{i_2i_1}(z_2,z_1)\RF_{i_ki_l}(z_l/z_k),&k<l,\\[1mm]
			\ZF_{i_1i_2}(z_1,z_2)
			+\ZF_{i_2i_1}(z_2,z_1)\RF_{i_li_k}(z_k/z_l)^{-1},&k>l.
		\end{cases}
		\label{app:d-d-factors}
	\end{equation}
	Both expressions vanish after the delta substitutions, by
	symmetrizability and $\tau$-invariance. The exchange polynomial cancels
	the scalar denominator before expansion, so this holds with either of
	the prescribed expansions. This completes the induction.
\end{pf}

\subsubsection{Expansion of the recursive expression}
\label{app:reduction-lemma}
Order the labels by $A<B<C$. For 
$\epsilon=(\epsilon_1,\ldots,\epsilon_n)\in\{A,B,C\}^n$, let $\ul v_B^\epsilon$ be the 
subtuple of entries with $\epsilon_a=B$, in their original order. Set
\[
A_\epsilon=\prod_{\epsilon_a=A}^{\longrightarrow}A_{i_a}(z_a),\qquad
C_\epsilon=\prod_{\epsilon_a=C}^{\longrightarrow}C_{i_a}(z_a),\qquad
\chi_\epsilon(\ul v)=
\prod_{\substack{a<b\\ \epsilon_a>\epsilon_b}}
\RF_{i_ai_b}(z_b/z_a),
\]
where the products defining $A_\epsilon$ and $C_\epsilon$ preserve the original order. 
The factors in $\chi_\epsilon$ have the expansions prescribed by the current relations.

\begin{lemma}
	For every tuple $\ul v=((i_1,z_1),\ldots,(i_n,z_n))$,
	\begin{equation}
		\cW_s(\ul v)=
		\sum_{\epsilon\in\{A,B,C\}^n}
		\chi_\epsilon(\ul v)\,
		A_\epsilon\,\bbW(\ul v_B^\epsilon)\,C_\epsilon .
		\label{app:induction-formula}
	\end{equation}
	Each term is defined coefficientwise with values in
		$\Lext{(\bfU^\imath_{\mu_1}\ten\bfU_{\mu_2})[\FC^{\pm1}]}{s^{-1};s}$.
\end{lemma}

\begin{pf}
	We proceed by induction on $n$, writing $A_a=A_{i_a}(z_a)$ and similarly for 
	$B_a,C_a$. The cases $n=0,1$ follow from $\cW_s(\varnothing)=1$ and 
	$X_1=A_1+B_1+C_1$. For $n=2$, reorder $B_1A_2$ and $C_1B_2$ using
	\eqref{app:BA}--\eqref{app:CB}. The only contributions still to be treated are
	\[
	B_1B_2+C_1A_2-\kappa_{i_1i_2}(z_1,z_2)\Theta^-_{i_2}(z_2),
	\]
	which give the terms indexed by $(B,B)$ and $(C,A)$ in \eqref{app:induction-formula}, 
	by \eqref{app:local-cancellation}.
	
	Suppose that \eqref{app:induction-formula} holds for tuples of length less than $n$, and 
	apply it to $\ul v^1$ in \eqref{app:Ws-recursion}. The contribution of $A_1$ is already in 
	the required order. Moving $B_1$ past the $A$-terms by \eqref{app:BA} gives the 
	factors in $\chi_\epsilon(\ul v)$ with $\epsilon_1=B$. Likewise, moving $C_1$ past the 
	$A$-terms and the $\bbW$-factor gives the required homogeneous factors by 
	\eqref{app:CA} and \eqref{app:C-Bblock-eq}.
	
	For the nonhomogeneous terms, fix $k>1$ with $i_1=\tau(i_k)$ and proceed by 
	induction on $k-1$. The case $k=2$ is \eqref{app:local-cancellation}. For $1<r<k$, the 
	weight relations for $\phi^+$ and $\phi^-$ and the substitution through the delta 
	functions give
	\begin{equation}
		\begin{aligned}
			&\fdel(\FC z_1z_k)\RF_{i_ki_r}(z_k/z_r)
			\Theta_{i_k}^-(z_k)A_{i_r}(z_r)\\
			&\qquad=\fdel(\FC z_1z_k)\RF_{i_1i_r}(z_r/z_1)
			A_{i_r}(z_r)\Theta_{i_k}^-(z_k).
		\end{aligned}
		\label{app:theta-A-ordering}
	\end{equation}
	Here $i_1=\tau(i_k)$ and $z_1={\FC^{-1}z_k^{-1}}$ on the delta support.
	The factor contributed by $\phi_{i_k}^-(\FC^{-1}z_k^{-1}s)$ is the inverse of
	$\RF_{i_ki_r}(z_k/z_r)$ with the same expansion at $z_k/z_r=0$.
	Thus the two scalar factors cancel coefficientwise.
	
	Move the $A$-term at position $k$ to the left, using Lemma~\ref{app:Bblock-A} for 
	intermediate $B$-terms. Each additional nonhomogeneous term arising from 
	\eqref{app:CA} involves a pair at smaller distance. The factor associated to $(1,k)$ is
	\[
	P_{1k}=\prod_{r=2}^{k-1}\RF_{i_ki_r}(z_k/z_r).
	\]
	Applying \eqref{app:local-cancellation} to $(1,k)$ gives the homogeneous $(C,A)$ term 
	and the recursive correction for the induced $B$-subtuple.
	The $\Theta^-$ contribution cancels the corresponding term of
	\eqref{app:Ws-recursion}, using \eqref{app:theta-A-ordering}.
	On the delta support, one has
	\[
	\Theta^+_{i_k}(z_k)=\wetadd_{i_k}(z_k)\ten
	\phi_{i_1}^+({\FC^{-1}z_1^{-1}} s)\phi_{i_k}^+({\FC^{-1}z_k^{-1}} s),
	\]
	so its term in \eqref{app:CA} supplies the lower degree term in the induced
	$\bbW$-factor.
	For an intermediate $B$-entry, the corresponding factor in $P_{1k}$ is the one in the recursion for this subtuple. For an intermediate $A$-entry, \eqref{app:theta-A-ordering} gives the required factor in $\chi_\epsilon$. For an intermediate $C$-entry, the factor already has this form. Thus the total coefficient is $\chi_\epsilon(\ul v)$. Since each additional nonhomogeneous term involves a shorter pair and the recursion removes the paired entries, the induction closes.
\end{pf}

\subsubsection[Cancellation of the mixed terms]{Cancellation of the mixed terms}
\label{app:mixed-cancellation}
Fix $i\ne j$ and $N=1-a_{ij}$. To apply \eqref{app:induction-formula} to the Serre sum, set
\[
\scrE=\{1,\ldots,N,\star\},\qquad
(i_r,z_r)=(i,z_r)\ (1\les r\les N),\qquad
(i_\star,z_\star)=(j,w).
\]
For $0\les p\les N$, let
\[
\ul v_p=(1,\ldots,p,\star,p+1,\ldots,N),
\]
where $r$ denotes $(i,z_r)$ and $\star$ denotes $(j,w)$. Fix a map
\[
\epsilon:\scrE\longrightarrow\{A,B,C\},\qquad A<B<C,
\]
and let $\ul v_{p,A}^\epsilon$, $\ul v_{p,B}^\epsilon$, and $\ul v_{p,C}^\epsilon$ be the subtuples with the respective labels, in the order induced by $\ul v_p$. For example, suppose $a_{ij}=-1$, so $N=2$, and take
\[
\epsilon(1)=A,\qquad \epsilon(2)=B,\qquad
\epsilon(\star)=B.
\]
For $p=1$, the ordered tuple is
$\ul v_1=((i,z_1),(j,w),(i,z_2))$, and hence
\[
\ul v_{1,A}^\epsilon=((i,z_1)),\qquad
\ul v_{1,B}^\epsilon=((j,w),(i,z_2)),\qquad
\ul v_{1,C}^\epsilon=\varnothing.
\]
For $p=2$, the ordered tuple is
$\ul v_2=((i,z_1),(i,z_2),(j,w))$, and hence
\[
\ul v_{2,A}^\epsilon=((i,z_1)),\qquad
\ul v_{2,B}^\epsilon=((i,z_2),(j,w)),\qquad
\ul v_{2,C}^\epsilon=\varnothing.
\]
Set
\begin{equation}
	\chi_{\epsilon,p}=
	\prod_{\substack{a\text{ precedes }b\text{ in }\ul v_p\\
			\epsilon(a)>\epsilon(b)}}
	\RF_{i_ai_b}(z_b/z_a).
	\label{app:chi-ep}
\end{equation}
Thus $\chi_{\epsilon,p}$ is the factor obtained by ordering the labels as $A<B<C$ while preserving the order within each subtuple.

For an ordered tuple $\ul u=((j_1,u_1),\ldots,(j_m,u_m))$, write
\[
A(\ul u)=A_{j_1}(u_1)\cdots A_{j_m}(u_m),\qquad
C(\ul u)=C_{j_1}(u_1)\cdots C_{j_m}(u_m),
\]
with value $1$ for the empty tuple.  Formula \eqref{app:induction-formula}, applied to $\ul v_p$, reads
\begin{equation}
	\cW_s(\ul v_p)=
	\sum_{\epsilon:\scrE\to\{A,B,C\}}
	\chi_{\epsilon,p}\,
	A(\ul v_{p,A}^\epsilon)\,
	\bbW(\ul v_{p,B}^\epsilon)\,
	C(\ul v_{p,C}^\epsilon).
	\label{app:induction-serre-form}
\end{equation}
By Lemma~\ref{app:internal-exchange}, each of the three factors satisfies the same homogeneous quadratic relation.
\begin{lemma}
Let $\scrE_{\rm mix}$ be the set of nonconstant maps $\epsilon:\scrE\to\{A,B,C\}$. Then
\begin{equation}
\begin{aligned}
0={}&\Sym_{\ul z}\sum_{p=0}^{N}(-1)^p
\begin{bmatrix}N\\p\end{bmatrix}_{q_i}
\sum_{\epsilon\in\scrE_{\rm mix}}\chi_{\epsilon,p}\,
A(\ul v_{p,A}^\epsilon)\,
\bbW(\ul v_{p,B}^\epsilon)\,
C(\ul v_{p,C}^\epsilon).
\end{aligned}
\label{app:mixed-operator-cancellation}
\end{equation}
\end{lemma}

\begin{pf}
Let $\widetilde{\scrQ}$ be the algebra generated by the coefficients of
$x_i^-(z)=\sum_{r\in\bbZ}x_{i,r}^-z^r$ and
$\phi_i^+(z)=\sum_{r\ges0}\phi_{i,r}^+z^r$, together with
$(\phi_{i,0}^+)^{-1}$. The Cartan coefficients commute, and we impose
only the ordinary weight and quadratic relations, without the Serre relations:
\begin{equation}\label{app:quadratic-word-algebra}
\begin{split}
(q_i^{a_{ij}}w-z)x_i^-(z)x_j^-(w)
 &+(q_i^{a_{ij}}z-w)x_j^-(w)x_i^-(z)=0,\\
\phi_i^+(z)x_j^-(w)
 &=\RF_{ij}(z,w)x_j^-(w)\phi_i^+(z),
\end{split}
\end{equation}
where the weight relation is expanded at $z/w=0$.
We consider the completion
\[
\widehat{\widetilde{\scrQ}^{\otimes3}}
=(\widetilde{\scrQ}^{\otimes3})[\FC^{\pm1}]
[\negthinspace[t_1,t_2]\negthinspace][t_1^{-1},t_2^{-1}],
\]
with $\FC$ central and the powers of $t_1,t_2$ bounded below separately.

Let $\scrQ\subset\widetilde{\scrQ}$ be the subalgebra generated by the coefficients of 
$x^-_i(z)$, and set
\[
e_{i,r}=x^-_{i,-r}
\aand
e_i(z)=\sum_{r\in\bbZ}e_{i,r}z^r=x_i^-(z^{-1}).
\]
Then,
\[
(q_i^{a_{ij}}z-w)e_i(z)e_j(w)
+(q_i^{a_{ij}}w-z)e_j(w)e_i(z)=0.
\]
Endow $\scrQ$ with the $\mathsf Q_+$-grading given by
$\deg e_{i,r}=\alpha_i$, so that
$\scrQ=\bigoplus_{\beta\in\mathsf Q_+}\scrQ_\beta$.
We use three copies $\scrQ_L$, for $L\in\{A,B,C\}$, and
set 
\[
\scrQ_{\operatorname{mix}} = \bigoplus_{\substack{\beta_A + \beta_B + \beta_C = N\alpha_i + \alpha_j\\ |\{L \mid \beta_L \neq 0\}| \geq 2}} \scrQ_{A, \beta_A} \otimes \scrQ_{B, \beta_B} \otimes \scrQ_{C, \beta_C}.  
\]
We also equip $\scrQ$ with the loop grading $\deg e_{i,r}=r$,
equivalently $\deg x^-_{i,r}=-r$.
For an elementary tensor $v\in\scrQ_{\operatorname{mix}}$ of multidegree $(r_A,r_B,r_C)$, set
\[
d=r_A+r_B+r_C,\qquad h_1=-r_A-r_B,\qquad h_2=-r_A.
\]
Write $(\scrQ_{\operatorname{mix}})_{r_A,r_B,r_C}$ for the component of multidegree $(r_A,r_B,r_C)$. 
We consider $\widehat{\scrQ_{\operatorname{mix}}}=\bigoplus_{d\in\bbZ}[\widehat{\scrQ_{\operatorname{mix}}}]_d$,
where $[\widehat{\scrQ_{\operatorname{mix}}}]_d$
allows formal sums
\[
\sum_{\substack{h_1\ges L_1\\h_2\ges L_2}}v_{h_1,h_2}
\quad\mbox{where}\quad
v_{h_1,h_2}\in
(\scrQ_{\operatorname{mix}})_{-h_2,\,h_2-h_1,\,d+h_1},
\]
for some integers $L_1,L_2$.\footnote{The indices $h_1,h_2$
record the powers of $t_1,t_2$ after the Cartan projection described below.}
For an ordered tuple
$\ul u=((i_1,z_1),\ldots,(i_m,z_m))$, write
\[
[\ul u]_L=e_{i_1}(z_1)\cdots e_{i_m}(z_m).
\]
We consider the \emph{linear} maps
\[
\widehat{\widetilde{\scrQ}^{\otimes 3}} \ \overset{\Phi}{\longleftarrow} \ \widehat{\scrQ_{\operatorname{mix}}} \ \overset{\Upsilon}{\longrightarrow} \ \Lext{(\bfU^\imath_{\mu_1}\ten\bfU_{\mu_2})[\FC^{\pm1}]}{s^{-1};s}
\]
given by 
\begin{align}
	\Upsilon\colon [\ul u_A]_A\ten[\ul u_B]_B\ten[\ul u_C]_C
	&\mapsto A(\ul u_A)\bbW(\ul u_B)C(\ul u_C),\\
	\Phi\colon [\ul u_A]_A\ten[\ul u_B]_B\ten[\ul u_C]_C
	&\mapsto F_1(\ul u_C)\ten F_{t_1}(\ul u_B)H_{t_1}(\ul u_C)\notag\\
	&\qquad\ten F_{t_1t_2}(\ul u_A)
	H_{t_1t_2}(\ul u_B)H_{t_1t_2}(\ul u_C),
\end{align}
where 
\[
F_t({\ul u}) =\prod_{a=1}^m x_{i_a}^-({\FC^{-1}z_a^{-1}} t), \qquad 
H_t({\ul u}) =\prod_{a=1}^m\phi_{i_a}^+({\FC^{-1}z_a^{-1}} t), 
\]
for an ordered tuple ${\ul u} = ((i_1, z_1), \cdots, (i_m, z_m))$.
On the algebraic tensor product, $\Phi$ is well-defined by the ordinary
weight and quadratic relations, and $\Upsilon$ by
Lemma~\ref{app:internal-exchange}. Both maps are readily verified to be continuous, hence they extend to the completions.
We now claim that: (i) the map $\Phi$ is injective, and (ii) the identity 
\begin{equation} \label{eqapp: Sym 2 id}
\Sym_{\ul z}\sum_{p=0}^{N}(-1)^p
\begin{bmatrix}N\\p\end{bmatrix}_{q_i}
\sum_{\epsilon\in\scrE_{\rm mix}}\chi_{\epsilon,p}\,
\Phi \big( [\ul{v}_{p,A}^\epsilon]_A \otimes 
[\ul{v}_{p,B}^\epsilon]_B \otimes 
[\ul{v}_{p,C}^\epsilon]_C \big) = 0
\end{equation} 
holds coefficientwise in $\widehat{\widetilde{\scrQ}^{\otimes3}}$.
These two claims, together with the definition of $\Upsilon$, immediately
imply \eqref{app:mixed-operator-cancellation}.

Let us prove the first claim. Let $\mathbf{U}^-$ and $\widetilde{\mathbf{U}}^-$ be the quotients of $\scrQ$ and $\widetilde{\scrQ}$, respectively, by the two-sided ideal generated by the Drinfeld--Serre relations \eqref{eq:serre-relations}. Consider the canonical projection $\pi \colon \scrQ \to \mathbf{U}^-$. We first observe that $\pi$ is injective on the graded components $\scrQ_{m\alpha_i + \alpha_j}$ ($0 \les m < N=1-a_{ij}$) and $\scrQ_{m\alpha_i}$ ($m \ges 0$). Indeed, \eqref{eq:serre-relations} for the ordered pair $(i,j)$ has degree $N\alpha_i + \alpha_j$, so its degree is too large to contribute to the kernel of $\pi$ in degree $m \alpha_i + \alpha_j$. For the reversed pair $(j,i)$, \eqref{eq:serre-relations} has degree $(1-a_{ji})\alpha_j + \alpha_i$. If $a_{ji} < 0$, then $\alpha_j$ has multiplicity at least $2$; if $a_{ji} = 0$ then $a_{ij} = 0$, implying $m = 0$. So this relation cannot contribute to the kernel either. Since the Drinfeld--Serre ideal is homogeneous, it has zero intersection with the stated homogeneous components. Next, we use the standard triangular decomposition of $\widetilde{\mathbf{U}}^-$, and define the projection 
\[
P \colon \widetilde{\mathbf{U}}^- \to \mathbf{U}^-, \qquad P(xh) = x \varepsilon(h)
\]   
for $x \in \mathbf{U}^-$, element $h$ in the Cartan subalgebra, and counit morphism 
$\varepsilon$. Let $\operatorname{rev}$ denote the reversal of tensor factors. 
For $v\in\scrQ_{\operatorname{mix}}$ of multidegree $(r_A,r_B,r_C)$,
one has
\[
P^{\otimes3}\widetilde\pi^{\otimes3}\Phi(v)
=\FC^d t_1^{h_1}t_2^{h_2}
\operatorname{rev}\bigl(\pi^{\otimes3}(v)\bigr),
\]
where $\widetilde\pi\colon\widetilde{\scrQ}\to\widetilde{\bfU}^-$
is the canonical projection.
By the previous argument, the RHS above is injective when restricted to 
$\wh{\scrQ}_{\operatorname{mix}}$. Hence, $\Phi$ is injective.

Let us prove the second claim. The main idea is to use the iterated coproduct $(\DD{t_1}\ten\id)\DD{t_1t_2}$.
Observe that 
\[(\DD{t_1}\ten\id)\DD{t_1t_2}(x_k^-(z^{-1}))
=\sfA_k(z)+\sfB_k(z)+\sfC_k(z),\]
where
\[
\begin{gathered}
\begin{aligned}
\sfA_k(z)&=1\ten1\ten x_k^-(z^{-1}t_1t_2),\\
\sfB_k(z)&=1\ten x_k^-(z^{-1}t_1)\ten\phi_k^+(z^{-1}t_1t_2),\\
\sfC_k(z)&=x_k^-(z^{-1})\ten\phi_k^+(z^{-1}t_1)
                         \ten\phi_k^+(z^{-1}t_1t_2).
\end{aligned}
\end{gathered}
\]
Apply $(\DD{t_1}\ten\id)\DD{t_1t_2}$ to the LHS of \eqref{app:ordinary-serre}, evaluated at $\FC^{-1}z_1^{-1},\ldots,\FC^{-1}z_N^{-1},\FC^{-1}w^{-1}$ and
considered as an element of $\widetilde{\scrQ}$. The resulting expression can be 
written as a polynomial in $\sfA_k, \sfB_k, \sfC_k$ for $k \in \{i,j\}$. We can group the 
monomials by the position $p$ of $x_j(w)$ and the number of factors $\sfA_k, \sfB_k$ 
and $\sfC_k$ occurring in them. 
By
\cite[Lem.~4.4]{lu-pan-26}, based on
\cite[Thm.~6.16]{damiani-24}, the Drinfeld--Serre currents are
skew-primitive before imposing Serre relations. Applying this
identity twice shows that the sum of the mixed monomials
vanishes. Note that skew-primitivity for finite rank-two types was
established by Enriquez \cite[Cor.~4.1]{enriquez-00} and
Gross\'{e} \cite[Sec.~1.4]{grosse-07}, and conjectured in general
by Hernandez \cite[Sec.~6.3.1]{hernandez-07}.
Moreover, for each choice of $p$ and each map $\epsilon\in\scrE_{\rm mix}$, if we move the Cartan currents (in each tensor factor) in the corresponding mixed monomial to the right, we obtain the coefficient $\chi_{\epsilon,p}$. This identifies the sum of the mixed monomials with the LHS of \eqref{eqapp: Sym 2 id}, and proves the vanishing of the latter. 
\end{pf}

\subsubsection[Compatibility with the Serre relations]{Compatibility with the Serre relations}
\label{app:serre-prop}
\begin{proposition}
	The assignments \eqref{app:coaction-zeta}--\eqref{app:coaction-y} preserve the $\imath$Drinfeld--Serre relations \eqref{app:iota-serre}.
\end{proposition}

\begin{pf}
	Fix $i\ne j$ and set $N=1-a_{ij}$. By the definition of $\cW_s$, it is enough to prove
	\begin{equation}
		\Sym_{\ul z}\sum_{p=0}^{N}(-1)^p
		\begin{bmatrix}N\\p\end{bmatrix}_{q_i}
		\cW_s(\ul v_p)=0.
		\label{app:target-serre}
	\end{equation}
	Substituting \eqref{app:induction-serre-form}, the terms indexed by nonconstant 
	$\epsilon$ cancel by Lemma~\ref{app:mixed-cancellation}.
	
	It remains to consider the three constant maps. For $\epsilon\equiv A$, the contribution 
	vanishes by the negative Drinfeld--Serre relation \eqref{app:ordinary-serre}. For 
	$\epsilon\equiv B$, one has $\chi_{\epsilon,p}=1$ and $\bbW(\ul v_p)=
	\ww(\ul v_p)\ten \phi_j^+({\FC^{-1}w^{-1}} s)\prod_{r=1}^N \phi_i^+({\FC^{-1}z_r^{-1}} s)$. Since the 
	Cartan currents commute, the second tensor factor is independent of $p$ and 
	symmetric in $\ul{z}$. The contribution therefore vanishes by \eqref{app:iota-serre} in 
	the first tensor factor.

	Finally, suppose that $\epsilon\equiv C$. Moving $\phi_a^+({\FC^{-1}u^{-1}} s)$ past $x_{\tau(b)}^+(vs)$ gives the factor $\RF_{a,\tau(b)}(\FC uv)$, expanded at $(\FC uv)^{-1}=0$. By symmetrizability and $\tau$-invariance, $d_a a_{a,\tau(b)}=d_b a_{b,\tau(a)}$, so this factor is symmetric in the two entries, with the same prescribed expansion. Set
	\[
	\begin{gathered}
		E=\wetad_j({\FC^{-1}w^{-1}})\prod_{r=1}^N\wetad_i({\FC^{-1}z_r^{-1}}),\qquad
		\Phi=\phi_j^+({\FC^{-1}w^{-1}} s)\prod_{r=1}^N \phi_i^+({\FC^{-1}z_r^{-1}} s),\\
		G=\prod_{\{a,b\}\subset\scrE}\RF_{i_a,\tau(i_b)}(\FC z_a z_b),
	\end{gathered}
	\]
	where either order of each pair gives the same factor in $G$. The expressions $E$, $\Phi$, and $G$ are independent of $p$ and symmetric in $\ul z$. Let $X_p^+$ be the ordered product obtained from $\ul v_p$ by replacing each entry $(a,u)$ with $x_{\tau(a)}^+(us)$. The contribution is
	\[
	GE\ten\left(\Sym_{\ul z}\sum_{p=0}^N(-1)^p\begin{bmatrix}N\\p\end{bmatrix}_{q_i}X_p^+\right)\Phi=0,
	\]
	by the positive Drinfeld--Serre relation for $(\tau(i),\tau(j))$, evaluated at $z_rs$ and $ws$. Here we use $q_{\tau(i)}=q_i$ and $a_{\tau(i)\tau(j)}=a_{ij}$.
\end{pf}

\subsection{Coassociativity}

\leavevmode We can now complete the proof of Theorem~\ref{thm:shifted-coaction}.

\begin{pf}
	By Lemma~\ref{app:completion-lemma} and Propositions~\ref{app:quadratic-prop} 
	and~\ref{app:serre-prop}, the assignments 
	\eqref{app:coaction-zeta}--\eqref{app:coaction-y} take values in the required 
	completion and preserve all defining relations. The inverse Cartan generators are 
	accounted for by \eqref{app:lowest-eta-mode}. Hence the assignments extend uniquely 
	to a morphism of algebras.

	It remains to verify the coassociativity \eqref{app:coassociativity}. By the Drinfeld 
	coproduct formulae \eqref{eq:Drinfeld-fields}, the identity is immediate on $\FC$ and 
	$\wetad_i(z)$. For $\wy_i(z)$, set $u={\FC^{-1}z^{-1}} s_1$ and $v=zs_1$. The terms arising 
	from the first two summands of \eqref{app:coaction-y} agree directly. For the last 
	summand, one has
	\[
	\DD{s_2}\bigl(x_{\tau(i)}^+(v)\phi_i^+(u)\bigr)
	=x_{\tau(i)}^+(v)\phi_i^+(u)\ten\phi_i^+(us_2)
	+\phi_{\tau(i)}^-(v)\phi_i^+(u)\ten x_{\tau(i)}^+(vs_2)\phi_i^+(us_2).
	\]
	The first term agrees with the contribution of $C_i(z)$ obtained by applying $\DDi{s_1}$ 
	to the middle summand on the left-hand side of \eqref{app:coassociativity}. The second 
	agrees with the last summand, since 
	$\DDi{s_1}(\wetad_i({\FC^{-1}z^{-1}}))=
	\wetad_i({\FC^{-1}z^{-1}})\ten\phi_i^+(u)\phi_{\tau(i)}^-(v)$, by 
	${\FC^{-1}(\FC^{-1}z^{-1})^{-1}}=z$ and the commutativity of the Cartan currents.
	
	For each coefficient of $z$, both expressions belong to the common
	completion
	\[
	\bigl(\bfU^\imath_{\mu_1}\ten\bfU_{\mu_2}\ten\bfU_{\mu_3}\bigr)
	[\FC^{\pm1}][\negthinspace[s_1,s_2]\negthinspace][s_1^{-1},s_2^{-1}].
	\]
	Indeed, the only term involving a negative Cartan current together
	with a root current is
	\[
	\wetad_i({\FC^{-1}z^{-1}})\ten
	\phi_i^+({\FC^{-1}z^{-1}} s_1)\phi_{\tau(i)}^-(zs_1)\ten
	x_{\tau(i)}^+(zs_1s_2)\phi_i^+({\FC^{-1}z^{-1}} s_1s_2).
	\]
	Let $a\ges-\mu_{1,\tau(i)}$, $b,e\ges0$, and
	$c\les\mu_{2,\tau(i)}$ be the indices of its four Cartan currents,
	in their displayed order. For the coefficient of $z^r$, the
	$s_1,s_2$-degrees are
	\[
	M_1=r+a+2b+2e,\qquad M_2=r+a+b+2e-c.
	\]
	Both are bounded below, and fixing them leaves finitely many
	choices of $a,b,c,e$. The other terms satisfy the same conclusion
	directly from the bounds in Lemma~\ref{app:completion-lemma}.
	This proves \eqref{app:coassociativity}.
\end{pf}


\nc{\etalchar}[1]{$^{#1}$}

\end{document}